\documentclass[reqno,10pt]{amsart}
\usepackage[foot]{amsaddr}
\usepackage{amssymb,amsmath,amsthm,amsfonts}
\usepackage{mathrsfs,dsfont,comment,mathscinet,mathtools}

\usepackage[utf8]{inputenc}
\usepackage[normalem]{ulem}
\usepackage[left=2.5cm,right=2.5cm,top=2.5cm,bottom=2.5cm]{geometry}

\usepackage{graphicx,tikz,color}
\usetikzlibrary{quotes,angles}

\usepackage[format=hang,labelfont=bf]{caption}

\usepackage{enumerate,esint}

\usepackage{ wasysym }

\usepackage{epstopdf}
\usepackage[colorlinks=true, pdfstartview=FitV, linkcolor=blue, 
            citecolor=blue, urlcolor=blue]{hyperref}

\allowdisplaybreaks
\numberwithin{equation}{section}
\theoremstyle{plain}
\newtheorem{theorem}{Theorem}[section]
\newtheorem{proposition}[theorem]{Proposition}
\newtheorem{lemma}[theorem]{Lemma}
\newtheorem{example}[theorem]{Example}
\newtheorem{corollary}[theorem]{Corollary}
\theoremstyle{definition}
\newtheorem{definition}[theorem]{Definition}
\newtheorem{remark}[theorem]{Remark}

\newtheorem*{theorem*}{Theorem}

\makeatletter
\def\th@plain{%
  \thm@notefont{}
  \itshape 
}
\def\th@definition{%
  \thm@notefont{}
  \normalfont 
}
\makeatother

\definecolor{mblue}{HTML}{13439b}

\newcommand\R{\mathbb R}

\newcommand\N{\mathbb N}
\newcommand\Z{\mathbb Z}

\renewcommand{\d}{\mathrm{d}}

\newcommand{\restr}[1]{|_{#1}}

\makeatletter
\newcommand{\subsetsim}{\mathrel{\mathpalette\subset@sim\relax}}
\newcommand{\subset@sim}[2]{%
	\vbox{\offinterlineskip\m@th
		\ialign{\hfil$#1##$\hfil\cr
			\sim\cr\subset\cr
		}%
	}%
}

\def\Xint#1{\mathchoice
{\XXint\displaystyle\textstyle{#1}}%
{\XXint\textstyle\scriptstyle{#1}}%
{\XXint\scriptstyle\scriptscriptstyle{#1}}%
{\XXint\scriptscriptstyle\scriptscriptstyle{#1}}%
\!\int}
\def\XXint#1#2#3{{\setbox0=\hbox{$#1{#2#3}{\int}$ }
\vcenter{\hbox{$#2#3$ }}\kern-.6\wd0}}

\def\dashint{\Xint-}

\makeatletter
\newcommand{\subalign}[1]{%
	\vcenter{%
		\Let@ \restore@math@cr \default@tag
		\baselineskip\fontdimen10 \scriptfont\tw@
		\advance\baselineskip\fontdimen12 \scriptfont\tw@
		\lineskip\thr@@\fontdimen8 \scriptfont\thr@@
		\lineskiplimit\lineskip
		\ialign{\hfil$\m@th\scriptstyle##$&$\m@th\scriptstyle{}##$\hfil\crcr
			#1\crcr
		}%
	}%
}
\makeatother

\definecolor{bblue}{HTML}{3C3C9F}

\newcommand{\loc}{\mathrm{loc}}

\newcommand{\tr}{\mathrm{tr}}
\newcommand{\cof}{\mathrm{cof}}
\newcommand{\adj}{\mathrm{adj}}

\renewcommand{\deg}{\mathrm{deg}}
\newcommand{\Det}{\mathrm{Det}}
\newcommand{\per}{\mathrm{Per}}
\newcommand{\domg}{\mathrm{dom}_{\rm G}}
\newcommand{\imt}{\mathrm{im}_{\rm T}}
\newcommand{\timt}{\widetilde{\mathrm{im}}_{\rm T}}
\newcommand{\img}{\mathrm{im}_{\rm G}}
\newcommand{\imc}{\mathrm{im}_{\rm C}}
\newcommand{\imm}{\mathrm{im}_{\rm M}}
\renewcommand{\div}{\mathrm{div}}
\newcommand{\dist}{\mathrm{dist}}
\newcommand\wk{\rightharpoonup}
\newcommand{\rnn}{\R^{N\times N}}
\newcommand{\wks}{\overset{\ast}{\rightharpoonup}}

\newcommand*\closure[1]{\overline{#1}}

\newcommand{\leb}{\mathscr{L}^N}
\newcommand{\haus}{\mathscr{H}^{N-1}}
\newcommand{\sgn}{\mathrm{sgn}}

\newcommand{\mres}{\mathbin{\vrule height 1.6ex depth 0pt width
		0.13ex\vrule height 0.13ex depth 0pt width 1.3ex}}
		
\DeclareMathOperator*{\aplim}{\mathrm{ap\,lim}}

\def\Xint#1{\mathchoice
{\XXint\displaystyle\textstyle{#1}}%
{\XXint\textstyle\scriptstyle{#1}}%
{\XXint\scriptstyle\scriptscriptstyle{#1}}%
{\XXint\scriptscriptstyle\scriptscriptstyle{#1}}%
\!\int}
\def\XXint#1#2#3{{\setbox0=\hbox{$#1{#2#3}{\int}$ }
\vcenter{\hbox{$#2#3$ }}\kern-.6\wd0}}

\def\dashint{\Xint-}

\newcommand{\widebar}[1]{\overline{\rule{0pt}{.5em}{#1}}}

\title[Eulerian-Lagrangian variational models allowing for material failure]{Variational models with Eulerian-Lagrangian formulation allowing for material failure}
\author[M. Bresciani]{Marco Bresciani${}^{*}$}
\address{* Department Mathematik, Friedrich-Alexander-Universit\"{a}t Erlangen-N\"{u}rnberg, Cauerstrasse 11, 91058 Erlangen (DE)}
\email{marco.bresciani@fau.de}
\author[M. Friedrich]{Manuel Friedrich${}^{*}$}
\email{manuel.friedrich@fau.de}
\author[C. Mora-Corral]{Carlos Mora-Corral${}^\dagger$}
\address{$\dagger$ Departamento de Matem\'{a}ticas, Universidad Aut\'{o}noma de Madrid, Calle Francisco Tom\'{a}s y Valiente 7, 28049 Madrid (ES) and Instituto de Ciencias Matem\'{a}ticas, CSIC-UAM-UC3M-UCM, Calle Nicol\'{a}s Cabrera 13--15, 28049 Madrid (ES)}
\email{carlos.mora@uam.es}

\date{\today}

\begin{document}

\setlength\parindent{0pt}

\vskip .2truecm
\begin{abstract}
We investigate the existence of minimizers of variational models with  Eulerian-Lagrangian formulations.  We consider energy functionals  depending  on the deformation of a body,  defined  on its reference configuration, and  an Eulerian map defined  on the unknown deformed configuration in the actual space. Our existence theory moves beyond the purely elastic setting and accounts for material failure by addressing   free-discontinuity problems   where  both  deformations and Eulerian fields are allowed to jump.     To do so, we build upon the work of Henao and Mora-Corral regarding the variational modeling of cavitation and fracture in nonlinear elasticity. Two main settings are considered by modeling deformations  as Sobolev  and $SBV$-maps, respectively.  The regularity of Eulerian maps is specified in each of these two settings according to the geometric and topological properties of the deformed configuration. 
We present some applications  to specific models of liquid crystals, phase transitions, and ferromagnetic elastomers.
Effectiveness and limitations of the theory are illustrated by means of explicit examples.

\end{abstract}
\maketitle


\section{Introduction}

\subsection{Motivation and overview}
\label{subsec:intro-motivation}
Variational  models featuring  Eulerian-Lagrangian formulations arise naturally in many multiphysics problems, where finite elasticity is coupled with other effects. In such models, equilibrium states correspond to minimizers of energy  functionals   comprising terms in both Eulerian and Lagrangian coordinates. The energy depends on at least two variables: the deformation of the body, classically defined on the reference configuration, and an Eulerian map defined on the deformed configuration in the actual space   which is      often subject to     nonlinear constraints. Concrete examples concern the modeling of nematic and ferromagnetic elastomers, where the Eulerian map represents the nematic director and the magnetization field, respectively, see, e.g.,    \cite{barchiesi.desimone,barchiesi.henao.moracorral,henao.stroffolini,warner.tarentjev} and \cite{barchiesi.henao.moracorral,bresciani,bresciani.davoli.kruzik,brown,kruzik.stefanelli.zeman,rybka.luskin}.  Other instances  regard the fields of plasticity   \cite{kruzik.melching.stefanelli,stefanelli}  and piezoelectricity \cite{rogers},   where the plastic deformation and the polarization vector are modeled as Eulerian fields.

The rigorous analysis of such variational models faces many difficulties which are mainly due to the fact that the domain of Eulerian maps, that is, the deformed configuration,  constitutes  by itself one of the unknowns. Additionally, the energy functional often depends on the composition of Eulerian fields and deformations, which is generally hard to handle by variational methods. 
In principle,  these  issues could be circumvented by exploiting the invertibility of deformations and rewriting the whole energy functional in Lagrangian coordinates. However, this approach  has  the  drawback of requiring the postulation of factitious constitutive assumptions in order to  make the analysis amenable. Moreover, specific  models might account for nonlocal effects,  for instance  long-range  and self-interactions  determined by the stray field in magnetic materials  \cite{brown,rogers}, which cannot be effectively described by means of Lagrangian coordinates only.    Summarizing, in various fields of applicative relevance  it is unavoidable to  cope with a mixed   Eulerian-Lagrangian structure.

 The question of proving  existence of minimizers for  Eulerian-Lagrangian energies has attained   increasing attention in recent years with various works focusing on specific physical models.    The most recent contributions provide  satisfactory answers within the setting of pure elasticity.  We provide  a brief review of the most relevant literature  in    Subsection \ref{subsec:literature} below.     Our aim is to move beyond the purely elastic setting and to establish  existence theories accounting for failure phenomena such as   cavitation and fracture \cite{henao.moracorral.invertibility}.   To our knowledge, the present work constitutes the first contribution in this direction.  Apart from being a challenging mathematical problem, this topic is also relevant from the mechanical point of view. Indeed, the widespread application of active materials in engineering has raised the question of  their  reliability  under different loading conditions.

In the present work, we combine the energetic approach for Eulerian-Lagrangian problems with the variational modeling of cavitation and fracture in nonlinear elasticity proposed  in \cite{henao.moracorral.invertibility}. 
  Specifically, we adopt the setting in \cite{henao.moracorral.lusin,mueller.spector} for the modeling of cavitation and   the one for   fracture in \cite{henao.moracorral.invertibility}.   
This leads to the formulation of free-discontinuity problems with energies comprising bulk and surface terms in both Eulerian and Lagrangian coordinates. 
More precisely, letting   $\Omega \subset \R^N$ be the reference configuration of an elastic body subject to a deformation $\boldsymbol{y}\colon \Omega \to \R^N$,  we consider an Eulerian map $\boldsymbol{v}\colon \boldsymbol{y}(\Omega) \to \R^M$ defined on the deformed configuration $\boldsymbol{y}(\Omega)$. Neglecting possible lower-order terms, we deal with energies of the form
\begin{equation}
	\label{eqn:intro-EL}
	(\boldsymbol{y},\boldsymbol{v})\mapsto \int_\Omega W(D\boldsymbol{y},\boldsymbol{v}\circ \boldsymbol{y})\,\d\boldsymbol{x}+\mathcal{S}(\boldsymbol{y})+\haus(J_{\boldsymbol{y}})+\int_{\boldsymbol{y}(\Omega)}|D \boldsymbol{v}|^2\,\d\boldsymbol{\xi}+\haus(J_{\boldsymbol{v}}). 
\end{equation} 
The first term in \eqref{eqn:intro-EL}   represents  the elastic energy and  exhibits  a coupling between the two variables  $\boldsymbol{y}$ and $\boldsymbol{v}$.  The second  term   accounts for possible failure phenomena such as cavitation and fracture. The functional  $\mathcal{S}$ has been   introduced in  \cite{henao.moracorral.invertibility} and  has an intricate definition, but morally stands for 
\begin{equation}\label{eqn:intro-surface}
	\mathcal{S}(\boldsymbol{y})  =  \haus(\partial \boldsymbol{y}(\Omega))-\haus(\boldsymbol{y}(\partial \Omega)),
\end{equation}
measuring  the new surface created by the deformation $\boldsymbol{y}$. More concretely, $\mathcal{S}(\boldsymbol{y})$ gives the sum of the perimeter of the cavities opened by $\boldsymbol{y}$ and the area of the fractured surface in $\boldsymbol{y}(\Omega)$.  The third   term    also models  the formation of  cracks,   expressed  in terms of the jump set   $J_{\boldsymbol{y}}$ of $\boldsymbol{y}$, and measures the breaking of atomic bonds. The fourth term corresponds to the Dirichlet energy of $\boldsymbol{v}$, which is given by an integral over the unknown deformed configuration $\boldsymbol{y}(\Omega)$. Eventually, the last term  penalizes  possible jumps of $\boldsymbol{v}$.

Concerning the energy terms in  \eqref{eqn:intro-EL}, the first and the third one  unambiguously qualify as Lagrangian, while the last two terms are distinctly Eulerian. In contrast, the nature of the  term $\mathcal{S}(\boldsymbol{y})$ is subject to interpretation: on the one hand, its expression involves a supremum of integrals over the reference configuration (see Definition~\ref{def:surface-energy}) which makes it a Lagrangian term; on the other hand, $\mathcal{S}(\boldsymbol{y})$ measures the area of the new surface created by $\boldsymbol{y}$ in the actual space according to \eqref{eqn:intro-surface} and, hence, it features an Eulerian character.

As we will discuss in Section~\ref{sec:appl}, the energy in \eqref{eqn:intro-EL}  accounts for many different physical models. It agrees with the Oseen-Frank energy of a liquid-crystal elastomer where the Dirichlet integral on the deformed set corresponds to the nematic term in the one-parameter approximation \cite{barchiesi.desimone,barchiesi.henao.moracorral,desimone.teresi,henao.stroffolini,warner.tarentjev}. 
 When  $\boldsymbol{v}$ is interpreted as a phase indicator, \eqref{eqn:intro-EL} describes phase transitions with Eulerian interfaces \cite{silhavy.proc,silhavy}. 
Eventually, neglecting the effect of the stray-field, the functional in  \eqref{eqn:intro-EL} expresses the magnetoelastic energy of a deformable ferromagnet with the fourth term  standing for    its exchange energy \cite{bresciani,bresciani.davoli.kruzik,brown,kruzik.stefanelli.zeman,rybka.luskin}.   

The aim of this paper is to prove the existence of minimizers 
  for the energy in \eqref{eqn:intro-EL} under Dirichlet boundary conditions  on the deformations and physical constraints on the Eulerian fields.    We refer to  Subsection \ref{subsec:intro-main} below for a more detailed account of our results.

  For the moment,  let us describe some of the most delicate mathematical issues. 
 Starting from the modeling,  the set $\boldsymbol{y}(\Omega)$ in \eqref{eqn:intro-EL} needs to be suitably interpreted according to the regularity assumptions  on  $\boldsymbol{y}$. Indeed, deformations belonging to Sobolev or $SBV$-spaces, as customary in nonlinear elasticity, may be only defined  almost everywhere, thus making the definition of the set $\boldsymbol{y}(\Omega)$ ambiguous. In particular,  the images of $\Omega$ under two  different  representatives of the same deformation  may  disagree on a set of positive measure  if  Lusin's condition (N) is violated, see \cite{ponomarev}.  Also, the specification of the regularity of $\boldsymbol{v}$ requires special care given that distributional derivatives can clearly be defined only on open sets. Furthermore, the composition $\boldsymbol{v}\circ \boldsymbol{y}$ in  \eqref{eqn:intro-EL} may be undefined on a set of positive measure whenever   $\boldsymbol{v}$ is  only defined almost everywhere and    $\boldsymbol{y}$ does not satisfy Lusin's condition (N${}^{-1}$).

  Two central questions in the analysis are  the compactness  of  (1)  the deformed configurations and  of   (2) the compositions of Eulerian fields and deformations with respect to the relevant  topologies.  Our results pertaining these two questions  essentially   encompass all the ones previously obtained in the literature  and,  to a  large extent,   hold in the  general framework of approximately differentiable maps.  Fine  properties and  invertibility  of deformations  play a crucial role within  the  proofs which require an extensive use of  Federer's change-of-variable formula.
Other challenges come from the fact that we are dealing with free-discontinuity problems as both the possible cavitation points and cracks are not prescribed, but rather constitute unknowns of the problem itself.  In addition, we observe that fictitious jumps  of  Eulerian fields may arise in the case of deformations exhibiting self-contact at the boundary. 

In the setting of Sobolev deformations, the topological degree represents a fundamental tool and the interplay between geometric and topological image lies at the core of our arguments. Fractures are clearly excluded, but deformations are allowed to create cavities. The location of the cavitation points and the volume of the cavities is described by the distributional determinant, which is a Radon measure on the reference configuration   \cite{henao.moracorral.lusin,mueller.spector}.  This feature and the weak continuity of distributional determinants are exploited  in our techniques  for proving the compactness of Eulerian maps which  extend the ones conceived in  \cite{barchiesi.henao.moracorral} for   the case of deformations that do not create cavities.   In our arguments,  we make the most of the knowledge on the relationship between geometric image, topological image, and surface energy achieved in  \cite{henao.moracorral.lusin}   by means of the results recalled in Theorem~\ref{thm:INV-top-im} below.  In this regard, the treatment of the last term in \eqref{eqn:intro-EL} requires special attention for both the proof of coercivity and lower semicontinuity.  

 In the setting of $SBV$-deformations, we profit by the generality of the results in Subsection~\ref{subsec:conv-approx}. As the deformations under consideration  exhibit discontinuities,  the topological degree is not available, so that most of the techniques devised for Sobolev maps cannot be adapted.  
Instead, the analysis relies on the convergence properties of the images of the deformations, their inverses, and the compositions with them.

Several examples are included in our paper to illustrate the effectiveness and limitations of our existence theory.

\subsection{Literature on Eulerian-Lagrangian energies}
\label{subsec:literature}
Without  claim  of completeness, we briefly review the most relevant literature concerning the existence of minimizers for Eulerian-Lagrangian energies.

A first contribution  in the framework of nonsimple materials was given in \cite{rybka.luskin}   for  a model of magnetoelasticty. The existence problem has been addressed for homeomorphic deformations   in $W^{1,p}(\Omega;\R^N)$  for   $p>N$   with integrable distortion in several papers concerning models of plasticity \cite{kruzik.melching.stefanelli,stefanelli},   viscoelasticity \cite{chiesa.kruzik.stefanelli},   electroelasticity \cite{davoli.molchanova.stefanelli}, and phase transitions \cite{brazda.lruzik.rupp.stefanelli,grandi.etal,grandi.etal2}. Continuous, but not necessarily homeomorphic  deformations   in $W^{1,p}(\Omega;\R^N)$   have been considered in \cite{barchiesi.desimone}   with $p=N$   and \cite{bresciani.davoli.kruzik,kruzik.stefanelli.zeman}  with $p>N$ for nematic  and  magnetic elastomers, respectively. The papers \cite{bresciani.davoli.kruzik,chiesa.kruzik.stefanelli,kruzik.melching.stefanelli,kruzik.stefanelli.zeman} also address the corresponding quasistatic evolution.

In \cite{barchiesi.henao.moracorral}, a class of possibly discontinuous deformations in $W^{1,p}(\Omega;\R^N)$  for  $p>N-1$ excluding the formation of  cavities has been introduced.  For this class of  deformations,  the existence of minimizers for Eulerian-Lagrangian models of liquid crystals and ferromagnets has been established. The analysis in \cite{barchiesi.henao.moracorral} has been extended in \cite{henao.stroffolini} by enlarging the class of admissible deformations   to the scale of Orlicz-Sobolev spaces.  An extension  of the magnetoelasticity model in  \cite{barchiesi.henao.moracorral} to the quasistatic setting  has been performed  in \cite{bresciani} within   the same class of deformations.  

As mentioned in Subsection~\ref{subsec:intro-motivation}, so far, the existence of minimizers of Eulerian-Lagrangian energies has been proved only in the setting of pure elasticity, i.e., by excluding both cavitation and fracture.

\subsection{Main results}
\label{subsec:intro-main}
In order to  state our main results, we need to specify the setting. Given a bounded Lipschitz domain $\Omega \subset \R^N$, we consider deformations $\boldsymbol{y}\colon \Omega \to \R^N$ that are almost everywhere approximately differentiable and almost everywhere injective. For such maps, we have the notion of \emph{geometric image}     introduced in \cite{henao.moracorral.invertibility}   which is  a subset  $\img(\boldsymbol{y},\Omega) \subset \boldsymbol{y}(\Omega)$  with  full measure.   We consider Eulerian maps $\boldsymbol{v}\in L^2(\img(\boldsymbol{y},\Omega);Z)$ for some given measurable set $Z\subset \R^M$ that are almost everywhere approximately differentiable,  with their approximate gradient denoted by  $\nabla \boldsymbol{v}$.  The set $Z$ embodies possible constraints that Eulerian fields have to comply with.

In this setting, it is possible to define the functional in \eqref{eqn:intro-EL} rigorously so that its value does not depend on representatives of $\boldsymbol{y}$ and $\boldsymbol{v}$.   Since compactness results for approximately differentiable maps are not available, we  specialize the setting  by modeling deformations as  Sobolev or $SBV$-maps. Accordingly, we  incorporate small modifications in \eqref{eqn:intro-EL} and we require Eulerian fields to satisfy suitable regularity conditions.    

\textbf{Sobolev deformations.} For a given exponent $p>N-1$, we consider  the  class of admissible deformations
\begin{equation*}
	\mathcal{Y}^{\rm cav}_p(\Omega)\coloneqq \left\{ \boldsymbol{y}\in W^{1,p}(\Omega;\R^N):\: \det D \boldsymbol{y}\in L^1_+(\Omega), \: \text{$\boldsymbol{y}$ satisfies condition (INV)}  \right\}.
\end{equation*}
The restriction on the exponent $p$ is crucial and allows us to employ the topological degree. The invertibility condition (INV) has been introduced in \cite{mueller.spector} and entails almost everywhere injectivity. Roughly speaking, this condition excludes the possibility  that  a cavity created at one point is filled by material coming from elsewhere.
Apart from excluding pathological behaviors, condition (INV)   enables   us to resort to the results in \cite{henao.moracorral.lusin,mueller.spector}. Following these works, for $\boldsymbol{y}\in \mathcal{Y}^{\rm cav}_p(\Omega)$, we define the \emph{topological image} of $\Omega$ under $\boldsymbol{y}$, denoted by $\imt(\boldsymbol{y},\Omega)\subset \R^N$, which is an open set  independent of  representatives. Also, we define the set $C_{\boldsymbol{y}}$ of cavitation points     associated to  $\boldsymbol{y}$    and, for $\boldsymbol{a}\in C_{\boldsymbol{y}}$, the corresponding cavity $\imt(\boldsymbol{y},\boldsymbol{a})$, which is a compact set. As  shown in \cite{henao.moracorral.lusin},  it holds that 
\begin{equation}
	\label{eqn:intro-img-imt}
	\imt(\boldsymbol{y},\Omega) \cong \img(\boldsymbol{y},\Omega) \cup \bigcup_{\boldsymbol{a}\in C_{\boldsymbol{y}}} \imt(\boldsymbol{y},\boldsymbol{a}).
\end{equation} 
where $\cong$ denotes equality almost everywhere.
In our first main result, we consider Eulerian maps $\boldsymbol{v}\in L^2(\img(\boldsymbol{y},\Omega);Z)$ as above such that their extension to $\imt(\boldsymbol{y},\Omega)$ by zero   enjoys suitable regularity.    The precise statement is given in Theorem~\ref{thm:cav}.

\begin{theorem}[Sobolev deformations]
	\label{thm:intro-Sobolev}
	Let $p>N-1$. Under  standard  continuity,  coercivity,   
	and polyconvexity assumptions on $W$, Dirichlet boundary conditions on the deformations,   and physical constraints on the Eulerian fields,  the functional
	\begin{equation*}
		(\boldsymbol{y},\boldsymbol{v})\mapsto \int_\Omega W(D\boldsymbol{y},\boldsymbol{v}\circ \boldsymbol{y})\,\d\boldsymbol{x}+\mathcal{S}(\boldsymbol{y})+\int_{\img(\boldsymbol{y},\Omega)} |\nabla \boldsymbol{v}|^2\,\d \boldsymbol{\xi} + \haus(J_{\boldsymbol{v}} \cap \imt(\boldsymbol{y},\Omega))
	\end{equation*}
	 admits minimizers in   the class   of admissible states $(\boldsymbol{y},\boldsymbol{v})$ with $\boldsymbol{y}\in\mathcal{Y}_p^{\rm cav}(\Omega)$ and $\boldsymbol{v}\in L^2(\img(\boldsymbol{y},\Omega);Z)$ almost everywhere approximately differentiable such that 
	 its extension  to $\imt(\boldsymbol{y},\Omega)$ by zero  enjoys suitable regularity.  
\end{theorem}

 Note that here  we do not account for jumps of $\boldsymbol{v}$ arising by self-contact at the boundary. Indeed, in contrast  to  \eqref{eqn:intro-EL}, the energy controls  $\haus(J_{\boldsymbol{v}} \cap \imt(\boldsymbol{y},\Omega))$ rather than $\haus(J_{\boldsymbol{v}})$,  see  Example \ref{ex:jump-not-contained-imt}  below  for an illustration  of the    difference between these two sets.  This  represents  an appreciable feature of the  setting,   for jumps of $\boldsymbol{v}$ induced by self-contact  are    fictitious.

In Theorem~\ref{thm:intro-Sobolev},  jumps of $\boldsymbol{v}$ can be excluded by restricting to the class of admissible states $(\boldsymbol{y},\boldsymbol{v})$   such that the extension of $\boldsymbol{v}$ to $\imt(\boldsymbol{y},\Omega)$ by zero belongs to $ W^{1,2}(\imt(\boldsymbol{y},\Omega);\R^M)$.   In view of \eqref{eqn:intro-img-imt}, this assumption implicitly enforces homogeneous Dirichlet boundary conditions for $\boldsymbol{v}$ on the boundary of the cavities. However, such boundary conditions can be inadequate for specific physical models, e.g.,  in the modeling of nematic elastomers, where the Eulerian field $\boldsymbol{v}$  representing the nematic director is constrained to have unit length.  

Therefore, we propose an alternative formulation of our variational model by restricting ourselves to the case of deformations creating    a finite number of cavities. Suppose that $\boldsymbol{y}\in \mathcal{Y}_p^{\rm cav}(\Omega)$ satisfies $\mathscr{H}^0(C_{\boldsymbol{y}})<+\infty$. We define the \emph{material image} of $\boldsymbol{y}$ by setting
\begin{equation*}
	\imm(\boldsymbol{y},\Omega)\coloneqq \imt(\boldsymbol{y},\Omega) \setminus \bigcup_{\boldsymbol{a}\in C_{\boldsymbol{y}}} \imt(\boldsymbol{y},\boldsymbol{a}).
\end{equation*}
 Bearing in mind \eqref{eqn:intro-img-imt}, we have $\imm(\boldsymbol{y},\Omega)\cong \img(\boldsymbol{y},\Omega)$, but, in contrast to $\img(\boldsymbol{y},\Omega)$, we can guarantee that $\imm(\boldsymbol{y},\Omega)$ is open  and   independent  of the representative of $\boldsymbol{y}$. Thus,  we can consider Eulerian maps satisfying $\boldsymbol{v}\in W^{1,2}(  \imm(\boldsymbol{y},\Omega);  Z)$.  Unfortunately,   it turns out that the resulting class of states is not  closed with respect to the relevant topology.  Yet,  we are able to solve this issue up to imposing a further constraint in the form of a lower bound on the cavity volumes.  
For a given $\kappa>0$, we consider  the  class of  deformations
\begin{equation*}
	\mathcal{Y}^{\rm cav}_{p,\kappa}(\Omega)\coloneqq \left\{ \boldsymbol{y}\in\mathcal{Y}_p^{\rm cav}(\Omega): \: \inf_{\boldsymbol{a}\in C_{\boldsymbol{y}}} \leb(\imt(\boldsymbol{y},\boldsymbol{a}))\geq \kappa  \right\}.
\end{equation*}
  As a consequence of the representation formula for $\mathcal{S}$ in \cite{henao.moracorral.lusin}, deformations in this class create  at most   a finite number of cavities.  
The second main result is given by Theorem~\ref{thm:cavk}. A simplified version of its statement reads as follows.

\begin{theorem}[Sobolev deformations with lower bound of the cavity volumes]
	\label{thm:intro-sobolev-k}
	Let $p>N-1$ and $\kappa>0$. Under  standard  continuity,  coercivity,   
	and polyconvexity assumptions on $W$, Dirichlet boundary conditions on the deformations,  and physical constraints on the Eulerian fields,    the functional
	\begin{equation*}
		(\boldsymbol{y},\boldsymbol{v})\mapsto \int_\Omega W(D\boldsymbol{y},\boldsymbol{v}\circ \boldsymbol{y})\,\d\boldsymbol{x}+\mathcal{S}(\boldsymbol{y})+\int_{\img(\boldsymbol{y},\Omega)} |D \boldsymbol{v}|^2\,\d \boldsymbol{\xi}
	\end{equation*}
	 admits minimizers in   the class of admissible states $(\boldsymbol{y},\boldsymbol{v})$ with $\boldsymbol{y}\in \mathcal{Y}^{\rm cav}_{p,\kappa}(\Omega)$ and $\boldsymbol{v}\in W^{1,2}(\imm(\boldsymbol{y},\Omega);Z)$.
\end{theorem}

As already mentioned, the class of admissible states  $(\boldsymbol{y},\boldsymbol{v})$  without   a lower bound on the cavity volumes  is not closed with respect to the relevant topology. This  is due to the fact that cavities can close   along converging sequences of deformations for which sequences of Eulerian maps with Sobolev regularity may exhibit jumps in the limit,  see Example \ref{ex:cavity-segment-relaxation} below.

\textbf{Deformations with bounded variation.} For the modeling of brittle materials, we fix $p>N-1$ and $b>0$, and we consider the class of deformations
\begin{equation*}
	\mathcal{Y}_{p,b}^{\rm frac}(\Omega)\coloneqq \left\{ \boldsymbol{y}\in SBV^p(\Omega;\R^N): \: \|\boldsymbol{y}\|_{L^\infty(\Omega;\R^N)}\leq b, \quad  \det \nabla \boldsymbol{y}\in L^1_+(\Omega), \quad  \text{$\boldsymbol{y}$ a.e.~injective} \right\}.
\end{equation*}
The confinement condition given by $b$ is a standard assumption in fracture   models  that allows us to work with  deformations   in $SBV^p(\Omega;\R^N)$. As the topological degree is not available in this setting, we can only resort to the geometric image. For $\boldsymbol{y}\in \mathcal{Y}_{p,b}^{\rm frac}(\Omega)$, we consider Eulerian maps $\boldsymbol{v}\in L^2(\img(\boldsymbol{y},\Omega);Z)$ as above such that their extension  to the whole space by zero   enjoys suitable regularity.   This  modeling assumption  requires the control of the boundary of the geometric image. Thus, we include the term $\per(\img(\boldsymbol{y},\Omega))$ into the energy,  whereby the stretching of the outer boundary is  also   penalized.   Indeed, recalling \eqref{eqn:intro-surface}, we  roughly have 
\begin{equation}
	\label{eqn:intro-perimeter}
	\per (\img(\boldsymbol{y},\Omega)) = \haus(\partial \boldsymbol{y}(\Omega)) = \mathcal{S}(\boldsymbol{y})+\haus(\boldsymbol{y}(\partial \Omega)).
\end{equation}

 The third main result can be   formulated   as follows. The precise statement is given in Theorem~\ref{thm:frac}.

\begin{theorem}[$\boldsymbol{SBV}$-deformations]
	\label{thm:intro-SBV}
	Let $p>N-1$ and $b>0$. Under  standard  continuity,  coercivity,   
	and polyconvexity assumptions on $W$,  Dirichlet boundary conditions on the deformations,   and physical constraints on the Eulerian fields,  the functional
	\begin{equation*}
		(\boldsymbol{y},\boldsymbol{v})\mapsto \int_\Omega W(\nabla \boldsymbol{y},\boldsymbol{v}\circ \boldsymbol{y})\,\d\boldsymbol{x}+\mathcal{S}(\boldsymbol{y})+\per \left( \img(\boldsymbol{y},\Omega) \right)+\haus(J_{\boldsymbol{y}})+\int_{\img(\boldsymbol{y},\Omega)} |\nabla \boldsymbol{v}|^2\,\d \boldsymbol{\xi} + \haus(J_{\boldsymbol{v}})
	\end{equation*}
	 admits minimizers in  the class of admissible states $(\boldsymbol{y},\boldsymbol{v})$ with $\boldsymbol{y}\in \mathcal{Y}_{p,b}^{\rm frac}(\Omega)$ and $\boldsymbol{v}\in L^2(\img(\boldsymbol{y},\Omega);Z)$ almost everywhere approximately differentiable such that the extension of $\boldsymbol{v}$ to the whole space by zero enjoys suitable regularity.
\end{theorem}

  This result is  less satisfactory compared  to the ones  for Sobolev deformations.  In fact,  in contrast  to  the energy in Theorem~\ref{thm:intro-Sobolev}, the functional considered here also penalizes  jumps of Eulerian fields determined by self-contact at the boundary. This feature seems to be unavoidable as no set that plays the role of the topological image in Theorem~\ref{thm:intro-Sobolev} is available for $SBV$-maps.   Additionally,   in contrast  to  Theorem~\ref{thm:intro-sobolev-k}, there is no variant of  Theorem~\ref{thm:intro-SBV} which excludes   jumps of $\boldsymbol{v}$.   Indeed,    simple counterexamples  show that  fractured pieces of the body  can be in contact in the limit leading to discontinuities of $\boldsymbol{v}$,   see Example~\ref{ex:jump-frac} below.

\subsection{Structure of the paper}
In Section \ref{sec:preliminaries}, we revise some known results concerning approximately differentiable maps, Sobolev maps, and maps with (generalized) bounded variation.  In particular, for Sobolev deformations,  we prove in Proposition~\ref{prop:top-im-inv} that the  topological images of nested domains are nested.  
In Section \ref{sec:conv}, we first   recast   the convergence properties of  deformed configurations  and then we establish new results in the case of Sobolev deformations creating cavities.  Section \ref{sec:EL} is devoted to the study of  Eulerian-Lagrangian energies and contains our main results. In Section \ref{sec:appl}, as an application of  our  main   results,   we establish the existence of minimizers for specific Eulerian-Lagrangian models. Precisely, we discuss the extensions of known models for nematic and magnetic elastomers and phase transitions to the setting of material failure. 
Eventually, results on radial deformations are contained in  the  Appendix.

The reader interested in the main results of the paper can move directly to Section \ref{sec:EL}, where all the necessary definitions are recalled. The paper contains  numerous  remarks that are relevant for our discussion but, in a first read,  can be skipped.

\section{Preliminaries} \label{sec:preliminaries}

This section collects some results mostly available in the literature. Sometimes these results are reformulated according to our specific needs. We  consider general approximately differentiable maps, Sobolev maps, and maps with bounded variation. Eventually, we recall a classical lower semicontinuity result.

\subsection{Notation} 

The integers $N,M\in \N$ satisfy $N \geq 2$ and $M \geq 1$.    The exponent $p>N-1$ is fixed and $\Omega \subset \R^N$ denotes a bounded Lipschitz domain. By domains we mean open and connected sets.

The set of $M\times N$ real matrices is denoted by $\R^{M \times N}$ and $\rnn_+$ stands for the set of matrices $\boldsymbol{F}\in \rnn$ with $\det \boldsymbol{F}>0$.  The unit matrix and the null matrix  in $\rnn$ are denoted by $\boldsymbol{I}$ and $\boldsymbol{O}$, respectively. The symbol $\boldsymbol{id}$ stands for the identity map in $\R^N$.
For every $\boldsymbol{F}\in \rnn$, we define $\adj \boldsymbol{F}$ as the unique matrix in $\rnn$ satisfying $\boldsymbol{F}(\adj \boldsymbol{F})=(\det \boldsymbol{F})\boldsymbol{I}$, and we set $\cof \boldsymbol{F}\coloneqq (\adj\boldsymbol{F})^\top$. 
We adopt the notation in \cite[Section~5.4]{dacorogna}:  for  every $1\leq r\leq N$, we define the  matrix $\adj_r\boldsymbol{F}$ in $ \R^{\binom{N}{r}\times \binom{N}{r}}$ given by the minors of $\boldsymbol{F}$ of order $r$ with suitable signs. With this notation, we have $\boldsymbol{F}=\adj_1 \boldsymbol{F}$, $\adj \boldsymbol{F}=\adj_{N-1} \boldsymbol{F}$, and $\det \boldsymbol{F}=\adj_N \boldsymbol{F}$.  

For $1\leq \nu \leq\infty$, we consider the standard $\nu$-norm defined as
\begin{equation}\label{eqn:q-norm}
	|\boldsymbol{x}|_\nu\coloneqq \begin{cases}
		\left (\sum_{i=1}^{N}|x_i|^\nu \right )^{1/\nu } & \text{if $1\leq \nu <\infty$,}\\
		\max \{|x_i|:\:i=1,\dots,N\} & \text{if $\nu =\infty$.}
	\end{cases} \quad \text{for all $\boldsymbol{x}\in\R^N$.}
\end{equation}
Accordingly, given $\boldsymbol{x}_0\in \R^N$ and $0<r<R$, we  define balls and annuli by  
\begin{equation}\label{eqn:q-ball-annulus}
	B_\nu(\boldsymbol{x}_0,r)\coloneqq \{\boldsymbol{x}\in \R^N:\hspace{2pt}|\boldsymbol{x}-\boldsymbol{x}_0|_\nu<r  \}, \quad A_\nu(  \boldsymbol{x}_0,  r,R)\coloneqq \{ \boldsymbol{x}\in \R^N: \hspace{2pt} r<|\boldsymbol{x}  - \boldsymbol{x}_0|_\nu   <R  \}.
\end{equation} 
We denote the closure of $B_\nu (\boldsymbol{x}_0,r)$ by $\closure{B}_\nu (\boldsymbol{x}_0,r)$, while we write
\begin{equation}\label{eqn:q-sphere}
	S_\nu (\boldsymbol{x}_0,r)\coloneqq \{\boldsymbol{x}\in\R^N: \:|\boldsymbol{x}-\boldsymbol{x}_0|_\nu =r\}
\end{equation}
for its boundary. For $\boldsymbol{x}_0=\boldsymbol{0}$, we omit the center by simply writing $B_\nu(r)  :=  B_\nu(\boldsymbol{0},r)$, $\closure{B}_\nu(r)\coloneqq \closure{B}_\nu(\boldsymbol{0},r)$,  $A_\nu(r,R)\coloneqq A_\nu(\boldsymbol{0},r,R)$, and $S_\nu(r)\coloneqq S_\nu(\boldsymbol{0},r)$. Additionally,
for $r=1$, we set $B_\nu \coloneqq B_\nu(1)$  and $S_\nu\coloneqq S_\nu(1)$. In all these  notations,  we omit the subscript for $\nu=2$. In particular, we use the notation
\begin{equation}\label{eqn:B-S}
	B\coloneqq B_2(1), \quad   \quad S\coloneqq S_2(1).
\end{equation}

Given a set $E\subset \R^N$, we denote its  interior, closure, and  boundary as $E^\circ$, $\closure{E}$, and  $\partial E$, respectively.   Given $F\subset \R^N$, we write $E\subset \subset F$ whenever $\closure{E}\subset F$. The characteristic function of a measurable set $A\subset \R^N$ is denoted by $\chi_A$.  When $\boldsymbol{u}\colon A \to \R^M$ is measurable, we use the symbol $\chi_A \boldsymbol{u}$ for the extension of this map to the whole space by zero.  
For $k\in \N$, we denote by $\mathscr{L}^k$ the Lebesgue measure on $\R^k$ and, for $\alpha>0$, we use the  notation $\mathscr{H}^\alpha $  for the $\alpha$-dimensional Hausdorff measure on $\R^N$.   We use  the dashed integral to denote the integral average, i.e., the value of the integral divided by the measure of the set of integration.   

 Given $E,F\subset \R^N$,   we write $E \cong F$  whenever $\leb(E \triangle F)=0$, and  $E \simeq F$ whenever $\haus(E \triangle F)=0$. Here, $E \triangle F\coloneqq (E \setminus F) \cup (F \setminus E)$ is the symmetric difference. Given two measurable functions $\boldsymbol{u},\boldsymbol{v}\colon A \subset \R^N \to \R^M$, we write $\boldsymbol{u}\cong \boldsymbol{v}$ and $\boldsymbol{u}\simeq \boldsymbol{v}$ if $\leb(\{\boldsymbol{x}\in A:\; \boldsymbol{u}(\boldsymbol{x})\neq \boldsymbol{v}(\boldsymbol{x}) \})=0$ and $\haus(\{\boldsymbol{x}\in A:\; \boldsymbol{u}(\boldsymbol{x})\neq \boldsymbol{v}(\boldsymbol{x}) \})=0$, respectively.  

We use standard notation for sets of finite perimeter (see, e.g., \cite[Chapter~4]{ambrosio.fusco.pallara} and  \cite[Chapters~12--16]{maggi}).   The reduced boundary of a measurable set $E\subset \R^N$ is denoted by $\partial^* E$, so that its perimeter $\per(E)$ satisfies $\per(E)= \haus(\partial^* E)$. Given an open set $O \subset \R^N$,  the relative perimeter is indicated as $\per(E;O)$, so that $\per(E;O)= \haus(\partial^* E \cap O)$. The Lebesgue density of $E$ at $\boldsymbol{x}_0\in\R^N$ is denoted by $\Theta^N(E,\boldsymbol{x}_0)$ and we set $E^{(t)}\coloneqq \{ \boldsymbol{x}\in \R^N:\:\Theta^N(E,\boldsymbol{x}_0)=t \}$ for $0 \leq t \leq 1$. Moreover, the essential boundary $\partial^- E$ of $E$ is defined as the complement of $E^{(0)} \cup E^{(1)}$.  

We use standard notation for Lebesgue and Sobolev spaces, for spaces of continuous and differentiable functions, and for the space of (generalized) maps of bounded variation as well, as for their local counterparts. Domain and codomain are separated by semicolon and the codomain is omitted when it is given by $\R$. In the case of continuous functions, the subscripts `b' and `c' are used for bounded and compactly supported maps,  respectively.    We denote by $L^1_+(\Omega)$ the set of maps $u\in L^1(\Omega)$ satisfying $u>0$ almost everywhere. Given two measurable sets $A\subset \R^N$ and $Z\subset \R^M$, we denote by $L^2(A;Z)$ the set of maps $\boldsymbol{u}\in L^2(A;\R^M)$ satisfying $\boldsymbol{u}(\boldsymbol{x})\in Z$ for almost every $\boldsymbol{x}\in A$. We use analogous notations also for the Sobolev space $W^{1,2}(A;Z)$ when $A$ is an open set. We denote the distributional gradient of maps by the symbol $D$. 
In Subsection \ref{subsec:mag}, we will consider the homogeneous Sobolev space defined as
\begin{equation}\label{beppo}
	V^{1,2}(\R^N)\coloneqq \left\{ u \in L^2_{\rm loc}(\R^N): \hspace{4pt} D u \in L^2(\R^N;\R^N) \right\}.
\end{equation}
The space of bounded Radon measures on $\Omega$ is denoted by $M_{\rm b}(\Omega)$.

We remark that in this paper we adopt the definition of equi-integrability in \cite{fonseca.leoni}:  a sequence $(\boldsymbol{u}_n)_n$ of measurable maps $\boldsymbol{u}_n\colon A \subset \R^N \to \R^M$ is termed equi-integrable whenever for every $\varepsilon>0$ there exists $\delta=\delta(\varepsilon)>0$ such that
\begin{equation*}
	\sup_{n \in \N} \int_E |\boldsymbol{u}_n(\boldsymbol{x})|\,\d\boldsymbol{x}<\varepsilon \quad \text{for all measurable sets $E\subset A$ with $\leb(E)<\delta$.}
\end{equation*}    
For functions on sets with finite measure, this definition is equivalent to other ones available in the literature, see, e.g., \cite[Definition~1.26]{ambrosio.fusco.pallara}, but otherwise this is not the case. In particular, for  results related to equi-integrability such as Vitali's convergence theorem \cite[Theorem~2.24]{fonseca.leoni}, the Dunford-Pettis theorem \cite[Theorem~2.54]{fonseca.leoni}, and the De la Vall\'{e}e Poussin criterion \cite[Theorem~2.29]{fonseca.leoni}, we resort to the formulations given in \cite{fonseca.leoni}.

\subsection{Approximately differentiable maps}
We adopt the definitions of approximate limit, continuity, and  differentiability given in \cite[Definition 1]{henao.moracorral.fracture} and \cite[Definition 2]{henao.moracorral.xu}. In particular, the approximate limit of a measurable map  $\boldsymbol{u}\colon A \subset \R^N \to \R^M$ as $\boldsymbol{x}\to \boldsymbol{x}_0$, denoted by $\aplim_{\boldsymbol{x}\to \boldsymbol{x}_0} \boldsymbol{u}(\boldsymbol{x})$,   is defined only at the points   $\boldsymbol{x}_0\in A^{(1)}$.   When  $\boldsymbol{u}$   is approximately differentiable at $\boldsymbol{x}_0 \in A^{(1)}$, we denote its approximate gradient by $\nabla \boldsymbol{u}(\boldsymbol{x}_0)$.  For more information on these notions, we refer to \cite[Chapter 3, Section 1.4]{cartesian.currents}.

We  register the following simple fact about the approximate differentiability of extensions. For a similar result, see \cite[Lemma 1]{henao.moracorral.xu}. Its proof is an immediate consequence of Lebesgue's density theorem.

\begin{lemma}
	\label{lem:approx-diff-extension}
Let  $\boldsymbol{u}\colon A \subset \R^N\to \R^M$ be almost everywhere approximately differentiable. Then, $\widebar{\boldsymbol{u}}\coloneqq \chi_A \boldsymbol{u}$ is almost everywhere approximately differentiable  with $\nabla \widebar{\boldsymbol{u}}\cong\chi_A \nabla \boldsymbol{u}$ in $\R^N$. 
\end{lemma}

Recall that a measurable map  defined on a subset of $\R^N$   is said to satisfy  Lusin's condition (N) when  it  maps $\leb$-negligible sets to  $\leb$-negligible  sets, while it is said to satisfy Lusin's condition (N${}^{-1}$) when the preimage of each $\leb$-negligible set is $\leb$-negligible.  For maps defined on hypersurfaces, analogous definitions apply with $\haus$ in place of $\leb$.  

We state Federer's  area and  change-of-variable   formulas  \cite[Theorem~1, p.~220]{cartesian.currents}.

\begin{proposition}[Area and change-of-variable formulas]
	\label{prop:change-of-variable}
	Let  $\boldsymbol{u}\colon A \subset \R^N \to \R^N$ be almost everywhere approximately differentiable. Denote by $D$ the set of approximate differentiability points of $\boldsymbol{u}$. Then, for every measurable set $E\subset A$, the function $\boldsymbol{\xi}\mapsto  \mathscr{H}^0(\{\boldsymbol{x}\in E \cap D:\hspace{2pt}\boldsymbol{u}(\boldsymbol{x})=\boldsymbol{\xi}  \})$ is measurable and the  area formula
	\begin{equation}
		\label{eqn:area-formula}
		\int_E |\det \nabla \boldsymbol{u}(\boldsymbol{x})|\,\d\boldsymbol{x}=\int_{\boldsymbol{u}(E \cap D)} \mathscr{H}^0(\{\boldsymbol{x}\in E \cap D:\hspace{2pt}\boldsymbol{u}(\boldsymbol{x})=\boldsymbol{\xi}  \})\,\d\boldsymbol{\xi}
	\end{equation}
	 holds. 	Moreover, for every measurable map $\boldsymbol{\zeta}\colon \boldsymbol{u}(D)\to \R^M$, the function $\boldsymbol{x}\mapsto \boldsymbol{\zeta}(\boldsymbol{u}(\boldsymbol{x}))|\det \nabla \boldsymbol{u}(\boldsymbol{x})|$ is measurable and the  change-of-variable formula 
	\begin{equation}
		\label{eqn:change-of-variable-formula}
		\int_E \boldsymbol{\zeta}( \boldsymbol{u}(\boldsymbol{x}))|\det \nabla \boldsymbol{u}(\boldsymbol{x})|\,\d\boldsymbol{x}=\int_{\boldsymbol{u}(E \cap D)}  \boldsymbol{\zeta}(\boldsymbol{\xi}) \mathscr{H}^0(\{\boldsymbol{x}\in E \cap D:\hspace{2pt}\boldsymbol{u}(\boldsymbol{x})=\boldsymbol{\xi}  \})\,\d\boldsymbol{\xi}
	\end{equation}
	 holds, 
	whenever one of the two integrals exists. 
\end{proposition}

\begin{remark} \label{rem:federer}
\begin{enumerate}[(a)]
	\item  By \cite[Chapter 3, Subsection 1.5, Proposition 1(i)--(ii)]{cartesian.currents}, the map $\boldsymbol{u}\restr{D}$ satisfies Lusin's condition (N) and the set $\boldsymbol{u}(E \cap D)$ is measurable.  The first claim follows immediately from \eqref{eqn:area-formula}.
\item  If   $\det \nabla \boldsymbol{u}\neq 0$ almost everywhere, then $\boldsymbol{u}$ satisfies  Lusin's condition (N${}^{-1}$). In that case,  $\boldsymbol{\zeta}\circ \boldsymbol{u}$ is well defined and measurable.  To check the Lusin condition, let $Y\subset \R^N$  with $\leb(Y)=0$ and set $X\coloneqq \boldsymbol{u}^{-1}(Y)$. At this point, we do not know whether the set $X$ is measurable, so that we cannot apply formula \eqref{eqn:change-of-variable-formula} for this set. By  regularity, there exists a Borel set $F\subset \R^N$ such that $Y\subset F$ and $\leb(F)=\leb(Y)=0$. The set $E\coloneqq \boldsymbol{u}^{-1}(F)$ is measurable with $X\subset E$ and $\leb(\boldsymbol{u}(E))=0$. Thus, by \eqref{eqn:area-formula}, we have
	\begin{equation*}
	\int_E |\det \nabla \boldsymbol{u}|\,\d\boldsymbol{x}=\int_{\boldsymbol{u}(E\cap D)} \mathscr{H}^0(\{\boldsymbol{x}\in  E  \cap D: \hspace{2pt}\boldsymbol{u}(\boldsymbol{x})=\boldsymbol{\xi}  \})\,\d\boldsymbol{\xi}=0,
	\end{equation*}
	where the integral on the right-hand side equals zero since $\leb(\boldsymbol{u}(E \cap D))=0$.
	Since, by assumption, the integrand on the left-hand side of the previous equation is almost everywhere positive, we conclude that $\leb(E)=0$ and, in turn, $\leb(X)=0$.
	\item  In \eqref{eqn:change-of-variable-formula}, the composition $\boldsymbol{\zeta}\circ \boldsymbol{u}$ may not be defined on a set of positive measure. However, if we set $\boldsymbol{\zeta}\circ \boldsymbol{u}|\det \nabla \boldsymbol{u}|$  to be zero on  $\{\det \nabla \boldsymbol{u}=0\}$ even if $\boldsymbol{\zeta}\circ \boldsymbol{u}$ is not defined, then the resulting function is measurable.  Its  equivalence class is uniquely determined by those of $\boldsymbol{\zeta}$ and $\boldsymbol{u}$. 
\end{enumerate}
\end{remark}

We give the definition of geometric domain and geometric image.  We refer to  \cite[Definition 3]{henao.moracorral.fracture} for further explanations and motivations.   The terminology  of geometric domain has been proposed in \cite{bresciani},  but the notion was already introduced in \cite{mueller.spector}.   Henceforth, $\Omega \subset \R^N$ is assumed to be a bounded Lipschitz domain.

\begin{definition}[Geometric domain and geometric image]
	\label{def:geometric-image}
	Let $\boldsymbol{y}\colon \Omega \to \R^N$ be almost everywhere approximately differentiable with $\det \nabla \boldsymbol{y}\neq 0$ almost everywhere. We define the geometric domain of $\boldsymbol{y}$  as the set $\domg(\boldsymbol{y},\Omega)$ of points $\boldsymbol{x}_0\in \Omega$ such that $\boldsymbol{y}$ is approximately differentiable at $\boldsymbol{x}_0$ with $\det\nabla\boldsymbol{y}(\boldsymbol{x}_0)\neq 0$, and there exist a compact set $K\subset \Omega$ with  $\Theta^N(K,\boldsymbol{x}_0)=1$ and a map $ \boldsymbol{w}\in C^1(\R^N;\R^N)$ satisfying $ \boldsymbol{w}\restr{K}=\boldsymbol{y}\restr{K}$ and $D \boldsymbol{w}\restr{K}=\nabla\boldsymbol{y}\restr{K}$. For every $E\subset \Omega$ measurable, we set
	\begin{equation*}
		\domg(\boldsymbol{y},E)\coloneqq E \cap \domg(\boldsymbol{y},\Omega)
	\end{equation*}
	and we define
	the geometric image of $E$ under $\boldsymbol{y}$  as
	\begin{equation*}
		\img(\boldsymbol{y},E)\coloneqq \boldsymbol{y}(E \cap \domg(\boldsymbol{y},\Omega)).
	\end{equation*}
\end{definition}
\begin{remark}\label{rem:geom-dom-im}
	\begin{enumerate}[(a)]
		\item The set $\domg(\boldsymbol{y},E)$ is measurable and $\domg(\boldsymbol{y},E)\cong E$,  see \cite[pp. 582--583]{henao.moracorral.fracture}.
		\item The map $\boldsymbol{y}\restr{\domg(\boldsymbol{y},E)}$ has the Lusin property (N)  and the set $\img(\boldsymbol{y},E)$ is measurable because of Remark~\ref{rem:federer}(a).
		\item The geometric image   depends  on the representative of the deformation. However, if $\boldsymbol{y}_1 \cong \boldsymbol{y}_2$, then $\img(\boldsymbol{y}_1,\Omega)\cong \img(\boldsymbol{y}_2,\Omega)$  thanks to the Lusin's condition (N).
	\end{enumerate}
\end{remark}

We recall that a  map defined on a subset of $\R^N$  is called almost everywhere injective if its restriction to the complement of  an   $\leb$-negligible set   is an injective map.  For maps defined on hypersurfaces, an analogous definition applies by replacing the measure $\leb$ with $\haus$.

We will work with the following class of admissible deformations: 
\begin{equation}
	\label{eqn:Y}
	\mathcal{Y}(\Omega)\coloneqq \left \{ \boldsymbol{y}\colon \Omega \to \R^N: \hspace*{3pt} \text{$\boldsymbol{y}$ a.e.\ approximately differentiable, \hspace{2pt} $\boldsymbol{y} $ a.e.\ injective, \hspace{2pt} $\det  \nabla \boldsymbol{y}\in L^1_+(\Omega)$}    \right\}. 
\end{equation}
  The following result has been originally proved in \cite[Lemma~2.5]{mueller.spector}. We state its reformulation  from \cite{henao.moracorral.fracture}.

\begin{lemma}[{\cite[Lemma~1]{henao.moracorral.fracture}}]
\label{lem:density}
Let $\boldsymbol{y} \in \mathcal{Y}(\Omega)$. Then, for every measurable set $E\subset \Omega$ and $\boldsymbol{x}_0\in \Omega$ with $\Theta^N(E,\boldsymbol{x}_0)=1$, we have $\Theta^N(\img(\boldsymbol{y},E),\boldsymbol{y}(\boldsymbol{x}_0))=1$,  i.e.,   $\boldsymbol{y}(E^{(1)}\cap \Omega) \subset \img(\boldsymbol{y},E)^{(1)}$.
\end{lemma}

Deformations in the  class $\mathcal{Y}(\Omega)$  are actually injective when restricted to the geometric domain, see \cite{henao.moracorral.fracture}.  We register this fact  in the next lemma together with the differentiability properties of the inverse. 

\begin{lemma}[Differentiability of inverse deformations]
	\label{lem:inverse-differentiable}
	Let  $\boldsymbol{y}\in \mathcal{Y}(\Omega)$.  Then, the following  holds: 
	\begin{enumerate}[(i)]
		\item The map $\boldsymbol{y}\restr{\domg(\boldsymbol{y},\Omega)}$ is injective.
		\item The map $(\boldsymbol{y}\restr{\domg(\boldsymbol{y},\Omega)})^{-1}$ is approximately differentiable in $\img(\boldsymbol{y},\Omega)$ with
		\begin{equation*}
			\nabla (\boldsymbol{y}\restr{\domg(\boldsymbol{y},\Omega)})^{-1}=(\nabla \boldsymbol{y})^{-1}\circ (\boldsymbol{y}\restr{\domg(\boldsymbol{y},\Omega)})^{-1}.
		\end{equation*}
		\item The map $\widebar{\boldsymbol{y}}^{-1}\coloneqq \chi_{\img(\boldsymbol{y},\Omega)}\boldsymbol{y}^{-1}$ is almost everywhere approximately differentiable in $\R^N$ with
		\begin{equation*}
			\text{$\nabla \widebar{\boldsymbol{y}}^{-1}=(\nabla \boldsymbol{y})^{-1} \circ (\boldsymbol{y}\restr{\domg(\boldsymbol{y},\Omega)})^{-1}$ in $\img(\boldsymbol{y},\Omega)$}
		\end{equation*}
		and $\nabla \widebar{\boldsymbol{y}}^{-1}\cong \boldsymbol{O}$ in $\R^N \setminus \img(\boldsymbol{y},\Omega)$. 
	\end{enumerate}
\end{lemma}

\begin{remark}\label{rem:nonzero-jacobian}
More generally, the result  holds for  a class of functions where  the condition $\det  \nabla \boldsymbol{y}\in L^1_+(\Omega)$ in \eqref{eqn:Y} is replaced by   $\det \nabla \boldsymbol{y}\neq 0$ almost   everywhere. 
\end{remark}
\begin{proof}
(i) This was proved in \cite[Lemma 3]{henao.moracorral.fracture}.

(ii)  We refer to \cite[Lemma 3]{henao.moracorral.xu}, where the symbol $\Omega_0$ corresponds to the set $\domg(\boldsymbol{y},\Omega)$ in our notation.

(iii) The claim follows from  (ii) and Lemma \ref{lem:approx-diff-extension}. 
\end{proof}

For the ease of reference, we specialize Proposition \ref{prop:change-of-variable} for  injective deformations. The result is established by applying  the proposition  to the map $\boldsymbol{y}\restr{\domg(\boldsymbol{y},\Omega)}$ and its inverse $(\boldsymbol{y}\restr{\domg(\boldsymbol{y},\Omega)})^{-1}$.  Henceforth, without further mention, the latter map will be simply denoted   by $\boldsymbol{y}^{-1}$.  

\begin{corollary}[Area and change-of-variable formulas for   admissible  deformations]
\label{cor:change-of-variable}
Let  $\boldsymbol{y}\in \mathcal{Y}(\Omega)$.  Then, the following  holds:  
\begin{enumerate}[(i)]
	\item For every measurable set $E\subset \Omega$, we have
	\begin{equation*}
		\leb(\img(\boldsymbol{y},E))  = \int_{\domg(\boldsymbol{y},E)}\det \nabla \boldsymbol{y}(\boldsymbol{x})\,\d\boldsymbol{x}  =\int_E \det \nabla \boldsymbol{y}(\boldsymbol{x})\,\d\boldsymbol{x}.
	\end{equation*}
	Moreover, for every measurable map $\boldsymbol{\psi}\colon \img(\boldsymbol{y},E)\to \R^M$, we have
	\begin{equation*}
		\int_E \boldsymbol{\psi}(\boldsymbol{y}(\boldsymbol{x}))\det \nabla \boldsymbol{y}(\boldsymbol{x})\,\d\boldsymbol{x}=\int_{\img(\boldsymbol{y},E)}\boldsymbol{\psi}(\boldsymbol{\xi})\,\d\boldsymbol{\xi},
	\end{equation*} 
	whenever one of the two integrals exists.
	\item For every measurable set $F\subset \R^N$, we have
	\begin{equation*}
		\leb(\boldsymbol{y}^{-1}\big(F \cap \img(\boldsymbol{y},\Omega))\big)=\int_{F \cap \img(\boldsymbol{y},\Omega)} \det \nabla \boldsymbol{y}^{-1}(\boldsymbol{\xi})\,\d\boldsymbol{\xi}.
	\end{equation*} 	
	Moreover, for every measurable map $\boldsymbol{\varphi}\colon  \boldsymbol{y}^{-1}(F\cap \img(\boldsymbol{y},\Omega))   \to \R^M$, there holds
	\begin{equation*}
		\int_{F \cap \img(\boldsymbol{y},\Omega)} \boldsymbol{\varphi} ( \boldsymbol{y}^{-1}(\boldsymbol{\xi}))\det \nabla \boldsymbol{y}^{-1}(\boldsymbol{\xi})\,\d\boldsymbol{\xi}=\int_{\boldsymbol{y}^{-1}(F \cap \img(\boldsymbol{y},\Omega))} \boldsymbol{\varphi}(\boldsymbol{x})\,\d\boldsymbol{x},
	\end{equation*}
	whenever one of the two integrals exists. 
\end{enumerate}
\end{corollary}

\begin{remark}
	\label{rem:composition-measurable}
\begin{enumerate}[(a)]
	\item Since $\det \nabla \boldsymbol{y}>0$ almost everywhere, the map $\boldsymbol{y}$ satisfies Lusin's condition (N${}^{-1}$) by  Remark \ref{rem:federer}(b).  This entails that the composition $\boldsymbol{\psi}\circ \boldsymbol{y}$ is measurable and its equivalence class  is uniquely determined by the ones of $\boldsymbol{\psi}$ and $\boldsymbol{y}$.    
	\item  By Lemma \ref{lem:inverse-differentiable}(ii), $\nabla \boldsymbol{y}^{-1}=(\nabla \boldsymbol{y})^{-1}\circ \boldsymbol{y}^{-1}$ and, in turn,  $\det \nabla \boldsymbol{y}^{-1}=(\det \nabla \boldsymbol{y})^{-1} \circ \boldsymbol{y}^{-1}>0$  almost everywhere. Hence, $\boldsymbol{y}^{-1}$ satisfies Lusin's condition (N${}^{-1}$) by  Remark \ref{rem:federer}(b),  the composition $\boldsymbol{\varphi}\circ \boldsymbol{y}^{-1}$ is measurable,  and its equivalence class depends only on those of $\boldsymbol{\varphi}$ and $\boldsymbol{y}$. 
	\item In view of Remark \ref{rem:nonzero-jacobian}, the result remains true   when the condition $\det  \nabla \boldsymbol{y}\in L^1_+(\Omega)$ is replaced by   $\det \nabla \boldsymbol{y}\neq 0$ almost  everywhere.    In this case, the terms $\det \nabla \boldsymbol{y}$ and $\det \nabla \boldsymbol{y}^{-1}$ need to be replaced by $|\det \nabla \boldsymbol{y}|$ and $|\det \nabla \boldsymbol{y}^{-1}|$, respectively,  within the formulas. 
\end{enumerate}
\end{remark}

 We present the area and change-of-variable formulas for surface integrals which are also due to Federer. We state a formulation which is adequate to our purposes.  For the definition of approximate tangential differentiability, we refer to \cite[Definition~4]{henao.moracorral.fracture}. 

\begin{proposition}[Area and change-of-variable formula for surface integrals]
	\label{prop:change-of-variable-surface}
Let $\boldsymbol{y}\in \mathcal{Y}(\Omega)$ and $U \subset \subset \Omega$ be a domain of class $C^1$ with   outer unit normal $\boldsymbol{\nu}_U$ satisfying  $\haus(\partial U \setminus \domg(\boldsymbol{y},\partial U))=0$. Then, the following holds:
\begin{enumerate}[(i)]
	\item We have $\img(\boldsymbol{y},\partial U)\subset \partial^* \img(\boldsymbol{y},U)$ and the outer unit normal $\boldsymbol{\nu}_{\img(\boldsymbol{y},U)}$ to  $\img(\boldsymbol{y},U)$  satisfies 
	\begin{equation*}
		\boldsymbol{\nu}_{\img(\boldsymbol{y},U)}(\boldsymbol{y}(\boldsymbol{x}_0))=\frac{(\cof \nabla \boldsymbol{y}(\boldsymbol{x}_0))\boldsymbol{\nu}_U(\boldsymbol{x}_0)}{|(\cof \nabla \boldsymbol{y}(\boldsymbol{x}_0))\boldsymbol{\nu}_U(\boldsymbol{x}_0)|} \quad \text{for all $\boldsymbol{x}_0\in \domg(\boldsymbol{y},\partial U)$.}
	\end{equation*}
	\item  Suppose that the map $\boldsymbol{y}\restr{\partial U}$ is $\haus$-almost everywhere approximately tangentially differentiable with $\nabla^{\partial U} \boldsymbol{y} \simeq (\nabla \boldsymbol{y})\restr{\partial U}(\boldsymbol{I}-\boldsymbol{\nu}_U \otimes \boldsymbol{\nu}_U)$ on $\partial U$. Then, for every $\haus$-measurable subset $E\subset \partial U$, we have
	\begin{equation*} 
		\haus(\img(\boldsymbol{y},E))=\int_E |(\cof \nabla \boldsymbol{y})\boldsymbol{\nu}_U|\,\d\haus.
	\end{equation*}
	Moreover, for any $\haus$-measurable map $\boldsymbol{\psi}\colon \img(\boldsymbol{y},E)\to \R^N$, we have
	\begin{equation*}
		\int_{\img(\boldsymbol{y},E)} \boldsymbol{\psi}(\boldsymbol{\xi}) \cdot \boldsymbol{\nu}_{\img(\boldsymbol{y},U)}(\boldsymbol{\xi})\,\d \haus(\boldsymbol{\xi})=\int_E \boldsymbol{\psi}(\boldsymbol{y}(\boldsymbol{x}))\cdot \left( (\cof \nabla \boldsymbol{y}(\boldsymbol{x}))\boldsymbol{\nu}_U(\boldsymbol{x})  \right)\,\d\haus(\boldsymbol{x}).
	\end{equation*}
\end{enumerate}
\end{proposition}
\begin{proof}
Claim (i) has been proved in \cite[Proposition~6, Claims (i) and (v)]{henao.moracorral.fracture}. For the formula in (ii), we refer to \cite[Proposition~2]{henao.moracorral.fracture}.  
\end{proof}

We will consider the surface energy  introduced in \cite{henao.moracorral.invertibility} and subsequently studied in \cite{henao.moracorral.lusin,henao.moracorral.fracture,henao.moracorral.regularity}.   Loosely speaking,  as already  mentioned  in \eqref{eqn:intro-surface}, this functional measures the surface created by deformations. 

\begin{definition}[Surface energy functional] \label{def:surface-energy}
	Let $\boldsymbol{y}\in \mathcal{Y}(\Omega)$  with $\cof\,\nabla \boldsymbol{y}\in L^1(\Omega;\rnn)$. 	We define
		 \begin{equation*}
		 	\mathcal{S}(\boldsymbol{y})\coloneqq \sup \left\{ \mathcal{S}_{\boldsymbol{y}}(\boldsymbol{\eta}): \hspace{3pt} \boldsymbol{\eta}\in C^\infty_{\rm c}(\Omega \times \R^N;\R^N), \hspace{3pt} \|\boldsymbol{\eta}\|_{C^0(\Omega \times \R^N;\R^N)}\leq 1   \right\},
		 \end{equation*}
		 where we set
		 \begin{equation*}
		 	\mathcal{S}_{\boldsymbol{y}}(\boldsymbol{\eta})\coloneqq \int_\Omega  \big(  \cof\,\nabla \boldsymbol{y}(\boldsymbol{x})\cdot D_{\boldsymbol{x}}\boldsymbol{\eta}(\boldsymbol{x},\boldsymbol{y}(\boldsymbol{x}))+\div_{\boldsymbol{\xi}}\boldsymbol{\eta}(\boldsymbol{x},\boldsymbol{y}(\boldsymbol{x}))\,\det \nabla \boldsymbol{y}(\boldsymbol{x}) \big)\,\d\boldsymbol{x}
		 \end{equation*}
		 for every  function $\boldsymbol{\eta}\in C^\infty_{\rm c}(\Omega \times \R^N;\R^N)$  with variables $(\boldsymbol{x},\boldsymbol{\xi})$. 
\end{definition}

\begin{remark}\label{rem:surface-energy}
\begin{enumerate}[(a)]
	\item In \cite[Definition 9]{henao.moracorral.fracture}, for every $\boldsymbol{y}\in \mathcal{Y}(\Omega)$, the authors defined  the set $\Gamma_{\boldsymbol{y}}$ corresponding to the surface created by $\boldsymbol{y}$  which   is given by the union of two sets  $\Gamma_{\boldsymbol{y}}^{\rm v} \subset \img(\boldsymbol{y},\Omega)^{(1/2)}$ and $\Gamma_{\boldsymbol{y}}^{\rm i}\subset J_{\boldsymbol{y}^{-1}}$  representing the visible and  invisible surface created by $\boldsymbol{y}$, respectively, where $ J_{\boldsymbol{y}^{-1}}$ denotes the set of jump points of the map  
	 $\boldsymbol{y}^{-1}$.  See \cite[p.~577]{henao.moracorral.fracture} for the motivations behind the definition of these sets and \cite[p.~598]{henao.moracorral.fracture} for several examples.
	Moreover, in \cite[Theorem 3]{henao.moracorral.fracture}, the authors proved that $\mathcal{S}(\boldsymbol{y})=\haus(\Gamma_{\boldsymbol{y}}^{\rm v})+2\haus(\Gamma_{\boldsymbol{y}}^{\rm i})$.  This rigorously identifies   $\mathcal{S}(\boldsymbol{y})$ as the measure of the surface created by $\boldsymbol{y}$. 
	\item The identity in \eqref{eqn:intro-surface} might suggest that $\mathcal{S}(\boldsymbol{y})\leq \per \left(\img(\boldsymbol{y},\Omega) \right)$ for all $\boldsymbol{y}\in \mathcal{Y}(\Omega)$.  However, this inequality does not generally hold and fails for deformations creating invisible surface as in (a). Indeed, this type of surface  is originated when two pieces of surface created at different points in the reference configuration are put in contact and, hence, it cannot be detected by the perimeter of the geometric image. In this regard, an instructive example  has been provided in \cite[Section 11, pp.~51--54]{mueller.spector}, where   a sequence of deformations   has   been exhibited creating more and more invisible surface while keeping the perimeter of geometric images uniformly bounded. Explicit computations   have been carried out in \cite[pp.~629--631]{henao.moracorral.invertibility} and \cite[Example~(d), p.~599]{henao.moracorral.fracture}.
\end{enumerate}	
\end{remark}

The surface energy functional has been introduced in \cite{henao.moracorral.invertibility} in connection with the problem of the weak continuity of the Jacobian determinant. The corresponding result reads as follows.

\begin{theorem}[{\cite[Theorem 1]{henao.moracorral.invertibility}}]
	\label{thm:determinant-surface-energy}
	 Let $(\boldsymbol{y}{}_n)_{n}\subset \mathcal{Y}(\Omega)$ be such that $(\cof\,\nabla \boldsymbol{y}_n)_n \subset  L^1(\Omega;\rnn)$. Suppose that there exist an almost everywhere approximately differentiable map $\boldsymbol{y}\colon \Omega \to \R^N$  with $\cof\,\nabla \boldsymbol{y} \in L^1(\Omega;\rnn)$ and  a function $h\in L^1_+(\Omega)$ satisfying 
	\begin{align*}
		\text{$\boldsymbol{y}_n \to \boldsymbol{y}$ a.e.\ in $\Omega$}, \qquad \text{$\cof\, \nabla \boldsymbol{y}_n \wk \cof\, \nabla \boldsymbol{y}$ in $L^1(\Omega;\rnn)$}, \qquad \text{$\det \nabla \boldsymbol{y}_n \wk h$ in $L^1(\Omega)$.}
	\end{align*}
	Additionally, assume that 
	\begin{equation*}
		\sup_{n \in \N} \mathcal{S}(\boldsymbol{y}_n)<+\infty.
	\end{equation*}
	Then, $h\cong\det \nabla \boldsymbol{y}$ and $\boldsymbol{y}$ is almost everywhere injective.  In particular, $\boldsymbol{y}\in\mathcal{Y}(\Omega)$.  Moreover,
	\begin{equation*}
		\mathcal{S}(\boldsymbol{y})\leq \liminf_{n \to \infty} \mathcal{S}(\boldsymbol{y}_n).
	\end{equation*} 
\end{theorem} 
\begin{remark} \label{rem:S-lsc}
\begin{enumerate}[(a)]
	\item We note that Theorem \ref{thm:determinant-surface-energy} requires the weak convergence  of  Jacobian cofactors. In Section \ref{sec:EL}, we will apply this theorem to maps in $W^{1,p}(\Omega;\R^N)$ and $SBV^p(\Omega;\R^N)$ for some $p>N-1$. In these settings, the weak continuity of the Jacobian cofactor is known, see \cite[Theorem 8.20, Part 4]{dacorogna} and \cite[Corollary 5.31]{ambrosio.fusco.pallara}, respectively.
	\item Whenever
	$\boldsymbol{y}_n\to \boldsymbol{y}$ almost everywhere in $\Omega$, $\cof \nabla \boldsymbol{y}_n \wk \cof \nabla \boldsymbol{y}$ in $L^1(\Omega;\rnn)$, and $\det \nabla \boldsymbol{y}_n \wk \det \nabla \boldsymbol{y}$ in $L^1(\Omega)$, the lower semicontinuity of $\mathcal{S}$  can be easily established  by  passing to the limit in  $\mathcal{S}_{\boldsymbol{y}_n}$.
	\item In general, the uniform boundedness of the surface energy cannot be replaced by the one of the perimeter of geometric images. This is in agreement with Remark~\ref{rem:surface-energy}(b). A  counterexample is provided by the sequence of deformations in \cite[Section 11, pp.~51--54]{mueller.spector}. Denoting this sequence by $(\boldsymbol{y}{}_n)_n$, we have $(\boldsymbol{y}{}_n)_n \subset   \mathcal{Y}(\Omega)$ for some bounded set $\Omega \subset \R^2$. By direct computation, one can also check that  $(\det \nabla \boldsymbol{y}_n)_n$ is equi-integrable. The sequence $(\boldsymbol{y}{}_n)_n$ converges almost everywhere towards an almost everywhere approximately deformation $\boldsymbol{y}\colon \Omega \to \R^2$ that violates injectivity on a set of positive measure and satisfies $\det \nabla \boldsymbol{y}>0$ almost everywhere. Therefore, we cannot have $\det \nabla \boldsymbol{y}_n \wk \det \nabla \boldsymbol{y}$ in $L^1(\Omega)$. Indeed, if that is the case, one can prove  by using the Ciarlet-Ne\v{c}as  condition \cite{giacomini.ponsiglione}  that $\boldsymbol{y}$ is almost everywhere injective, thus providing a contradiction.  We refer to \cite[Section 1]{henao.moracorral.invertibility} for further comments.
\end{enumerate}	
\end{remark}

\subsection{Sobolev maps}
\label{subsec:sobolev}
In this  subsection,  we  revise   some results about Sobolev deformations with subcritical integrability. We repeat that, in this paper, the exponent $p>N-1$ is fixed. 

We start by recalling the notion of topological degree for Sobolev maps. For a comprehensive account on the topological degree of continuous mappings we refer to \cite{fonseca.gangbo}. The application of this tool  to  Sobolev maps dates back to \cite{reshetnyak.67} (see also \cite{reshetnyak}) and the use of degree theory for maps  in $W^{1,p}(\Omega;\R^N)$ was initiated in \cite{sverak}.

Given $\boldsymbol{y}\in W^{1,p}(\Omega;\R^N)$, its precise representative $\boldsymbol{y}^*\colon \Omega \to \R^N$ is defined as
\begin{equation*}
	\boldsymbol{y}^*(\boldsymbol{x}_0)\coloneqq \limsup_{r \to 0^+} \dashint_{B(\boldsymbol{x}_0,r)} \boldsymbol{y}(\boldsymbol{x})\,\d\boldsymbol{x}.
\end{equation*}
The set $L_{\boldsymbol{y}}$ of Lebesgue points of $\boldsymbol{y}$ is  given  by the points  $\boldsymbol{x}_0\in \Omega$ for which we have
\begin{equation*}
	\lim_{r \to 0^+} \dashint_{B(\boldsymbol{x}_0,r)} |\boldsymbol{y}(\boldsymbol{x})-\boldsymbol{y}^*(\boldsymbol{x}_0)|\,\d \boldsymbol{x}=0.
\end{equation*}
In that case, the superior limit in  the definition of $\boldsymbol{y}^*(\boldsymbol{x}_0)$ can be replaced by a limit. Moreover, it is known that $\mathscr{H}^1(\Omega \setminus L_{\boldsymbol{y}})=0$, see \cite[Theorem~2, p.~186]{cartesian.currents}.  

\begin{definition}[Topological degree  and topological image of domains] \label{def:top-deg}
	Let $\boldsymbol{y}\in W^{1,p}(\Omega;\R^N)$ and   $U \subset \subset \Omega$ be a domain  with $\partial U\subset L_{\boldsymbol{y}}$   such that $\boldsymbol{y}^*\restr{\partial U}\in C^0(\partial U;\R^N)$. 	We define the topological degree of $\boldsymbol{y}$ on $U$ to be the topological degree of any continuous extension of  $\boldsymbol{y}^*\restr{\partial U}$ to $\closure{U}$. This yields a  map $\deg(\boldsymbol{y},U,\cdot)\colon \R^N \setminus \boldsymbol{y}^*(\partial U)\to \Z$ and we define the topological image of $U$ under $\boldsymbol{y}$ as
	\begin{equation*}
		\imt(\boldsymbol{y},U)\coloneqq \left\{ \boldsymbol{\xi}\in \R^N \setminus \boldsymbol{y}^*(\partial U): \hspace{3pt} \deg(\boldsymbol{y},U,\boldsymbol{\xi})\neq 0 \right\}.
	\end{equation*} 
\end{definition}

\begin{remark} \label{rem:top-deg}
	\begin{enumerate}[(a)]
		\item    The choice of considering  domains $U$ satisfying $\partial U\subset L_{\boldsymbol{y}}$ and $\boldsymbol{y}^*\restr{\partial U}\in C^0(\partial U;\R^N)$  is motivated just by notational simplicity and   allows us to use the precise representative whenever we refer to pointwise properties of deformations or to their restriction on lower-dimensional sets. The same choice is made also in \cite{henao,mueller.spector,sivaloganathan.spector,sivaloganathan.spector.tilakraj}.  The abundance of domains satisfying these conditions   will be demonstrated in Lemma~\ref{lem:abundance} below.      
		\item  If $\boldsymbol{y}^*\restr{\partial U}$ is continuous, then it admits a continuous extension to $\closure{U}$ by Tiezte's theorem  \cite[Theorem 1.15]{fonseca.gangbo}. Since the topological degree depends only on the boundary values  \cite[Theorem 2.4]{fonseca.gangbo}, the map $\deg(\boldsymbol{y},U,\cdot)$ is well defined.   
		\item The map $\deg(\boldsymbol{y},U,\cdot)$ is  locally constant \cite[Theorem 2.3(3)]{fonseca.gangbo} and hence   continuous. As consequences,   $\imt(\boldsymbol{y},U)$ is open and $\partial\, \imt(\boldsymbol{y},U)\subset \boldsymbol{y}^* (\partial U)$.
		\item  We have $\deg(\boldsymbol{y},U,\cdot)=0$ in the unique unbounded component of $\R^N \setminus \boldsymbol{y}^*(\partial U)$. Thus,   $\imt(\boldsymbol{y},U)$ is bounded. The claim can be proved as follows. Let $\boldsymbol{\varphi}\in C^0(\closure{U};\R^N)$ be an extension of $\boldsymbol{y}^*\restr{\partial U}$ and $V_0$ be the  unbounded component of $\R^N \setminus \boldsymbol{y}^*(\partial U)$. By contradiction, suppose that  $\deg(\boldsymbol{y},U,\boldsymbol{\xi})=\deg(\boldsymbol{\varphi},U,\boldsymbol{\xi})\neq 0$ at some point $\boldsymbol{\xi}\in V_0$. Then, the same holds at any other point in $V_0$ by (c), so that $V_0\subset \boldsymbol{\varphi}(U)$ by  \cite[Theorem~2.1]{fonseca.gangbo}. As $V_0$ is unbounded and $\boldsymbol{\varphi}(U)$ is bounded, this provides a contradiction.   
	\end{enumerate}
\end{remark}

As a consequence of the continuity of the topological degree, we have the following characterization of the asymptotic  behavior  of the topological images for weakly converging  sequences   of deformations. Proofs are given in \cite[Lemma~3.6]{barchiesi.henao.moracorral} and \cite[Lemma~2.10]{bresciani}.

\begin{lemma}\label{lem:top-deg-conv}
	Let $(\boldsymbol{y}{}_n)_{n}\subset W^{1,p}(\Omega;\R^N)$ and $\boldsymbol{y}\in W^{1,p}(\Omega;\R^N)$. Also, let $U \subset \subset \Omega$ be a domain  with $\partial U \subset \bigcap_{n\in \N} L_{\boldsymbol{y}_n}\cap L_{\boldsymbol{y}}$. Suppose that 
	\begin{equation*}
		 \text{$(\boldsymbol{y}^*_n\restr{\partial U})_n\subset C^0(\partial U;\R^N)$,    \qquad $\boldsymbol{y}_n^* \to \boldsymbol{y}^*$ uniformly on $\partial U$.}
	\end{equation*}
	Then, the following properties hold:
	\begin{enumerate}[(i)]
		\item For every compact set $ H  \subset  \imt(\boldsymbol{y},U)$, we have $ H  \subset \imt(\boldsymbol{y}_n,U)$ for $n \gg 1$ depending on $ H  $.
		\item For every compact set $ H  \subset \R^N \setminus \left( \imt(\boldsymbol{y},U) \cup \boldsymbol{y}^*(\partial U) \right)$, we have $ H   \subset \R^N \setminus \left( \imt(\boldsymbol{y}_n,U) \cup \boldsymbol{y}_n^*(\partial U) \right)$ for $n \gg 1$ depending on $ H  $.
		\item We have $\chi_{\imt(\boldsymbol{y}_n,U)}\to \chi_{\imt(\boldsymbol{y},U)}$ in $L^1(\R^N)$.
	\end{enumerate}
\end{lemma} 

We introduce some notation.   Given a domain $U \subset \subset \Omega$ of class $C^2$, the signed distance function $d_U\colon \R^N \to \R$ is defined by setting
\begin{equation*}
	d_U(\boldsymbol{x})\coloneqq \begin{cases}
		\dist(\boldsymbol{x};\partial U) & \text{if $\boldsymbol{x}\in U,$} \\
		0 & \text{if $\boldsymbol{x}\in \partial U,$} \\
		-\dist(\boldsymbol{x};\partial U) & \text{if $\boldsymbol{x}\in \R^N \setminus \closure{U}.$}
	\end{cases}
\end{equation*}
Given the regularity of $U$, there exists $0<\delta<\dist(\partial U;\partial \Omega)$ such that $d_U$ is of class $C^2$ in the tubular neighborhood $\{ \boldsymbol{x}\in \Omega: \hspace{3pt} -\delta<d_U(\boldsymbol{x})<\delta\}$. In particular,  the  set $U_s\coloneqq \{\boldsymbol{x}\in \R^N:\:d_U(\boldsymbol{x})>s\}$ is a domain of class $C^2$ for every $s \in (-\delta,\delta)$. We refer to \cite[Section~4]{ambrosio} for the proofs of these claims.   The  outer unit normal of  $U$  will be denoted by  $\boldsymbol{\nu}_U$ and, analogously,  by $\boldsymbol{\nu}_{U_s}$  for $U_s$.  

The next result has been originally established in \cite[Lemma 2.9]{mueller.spector}. We state its formulation from \cite{barchiesi.henao.moracorral} and we combine it  with claim (iii)  of Lemma \ref{lem:top-deg-conv}. 

\begin{lemma}[{\cite[Lemma 2.24]{barchiesi.henao.moracorral}}]
	\label{lem:conv-boundaries}
	Let $(\boldsymbol{y}{}_n)_n\subset W^{1,p}(\Omega;\R^N)$ and $\boldsymbol{y}\in W^{1,p}(\Omega;\R^N)$ be such that
	\begin{equation*}
		\text{$\boldsymbol{y}_n \wk \boldsymbol{y}$ in $W^{1,p}(\Omega;\R^N)$.}
	\end{equation*}
	Let $U \subset \subset \Omega$ be a domain of class $C^2$. Then, there exists $\delta>0$ such that,  for almost every $s\in (-\delta,\delta)$, we have  $(\boldsymbol{y}_n^*\restr{\partial U})_n\subset C^0(\partial U_s;\R^N)$ and $\boldsymbol{y}^*\restr{\partial U} \in C^0(\partial U_s;\R^N)$.   Additionally, for almost every $s$ as above, there exists a not relabeled subsequence possibly depending on $s$, for which 
	\begin{equation*}
		\text{$\boldsymbol{y}_{n}^* \to \boldsymbol{y}^*$ uniformly on $\partial U_s$}, \qquad \text{$\chi_{\imt(\boldsymbol{y}_n,U_s)}\to \chi_{\imt(\boldsymbol{y},U_s)}$ in $L^1(\R^N)$.}
	\end{equation*}
\end{lemma}

We recall the invertibility condition (INV) introduced in \cite{mueller.spector}.   Morally, this condition excludes the possibility that a cavity created at one point can be filled by material coming from elsewhere.   

\begin{definition}[Invertibility condition]\label{def:INV}
	Let $\boldsymbol{y}\in W^{1,p}(\Omega;\R^N)$. Then, $\boldsymbol{y}$  satisfies condition (INV) whenever, for every $\boldsymbol{a}\in \Omega$ and  almost all $r \in (0,\dist(\boldsymbol{a};\partial \Omega))$   the following  holds: 
	\begin{enumerate}[(i)]
		\item $S(\boldsymbol{a},r)\subset L_{\boldsymbol{y}}$ and  $\boldsymbol{y}^*\in C^0(S(\boldsymbol{a},r);\R^N)$;  
		\item $\boldsymbol{y}(\boldsymbol{x})\in \imt(\boldsymbol{y},B(\boldsymbol{a},r))$ for almost every $\boldsymbol{x}\in B(\boldsymbol{a},r)$;
		\item $\boldsymbol{y}(\boldsymbol{x})\notin \imt(\boldsymbol{y},B(\boldsymbol{a},r))$ for almost every $\boldsymbol{x}\in \Omega \setminus B(\boldsymbol{a},r)$.
	\end{enumerate}
\end{definition} 
\begin{remark}\label{rem:INV}
From  (ii) and Remark \ref{rem:top-deg}(d),  we deduce that $\boldsymbol{y}\in L^\infty_{\loc}(\Omega;\R^N)$. As $\boldsymbol{y}$ is approximately continuous at any point of $\domg(\boldsymbol{y},\Omega)$, we deduce that $\domg(\boldsymbol{y},\Omega)\subset L_{\boldsymbol{y}}$ and $\boldsymbol{y}\restr{\domg(\boldsymbol{y},\Omega)}=\boldsymbol{y}^*\restr{\domg(\boldsymbol{y},\Omega)}$.  This follows from \cite[Proposition~3.65]{ambrosio.fusco.pallara}.  (Note that  the reference  employs  a different terminology for approximate limits).
\end{remark}

We observe that condition (INV) is sufficient for almost everywhere injectivity but not necessary,  see  \cite[Remark 2 after Definition 3.2]{mueller.spector}.  

\begin{lemma}[{\cite[Lemma 3.4]{mueller.spector}}]
	\label{lem:inv-injectivity}
	Let $\boldsymbol{y}\in W^{1,p}(\Omega;\R^N)$ with $\det D \boldsymbol{y} \neq 0$ almost everywhere satisfy condition {\rm (INV)}. Then, $\boldsymbol{y}$ is almost everywhere injective.	
\end{lemma}

The next result ensures the stability of condition (INV) for weakly converging sequences of deformations.

\begin{lemma}[{\cite[Lemma 3.3]{mueller.spector}}]
	\label{lem:INV-stability}
	Let $(\boldsymbol{y}{}_n)_n \subset W^{1,p}(\Omega;\R^N)$ and $\boldsymbol{y}\in W^{1,p}(\Omega;\R^N)$ be such that
	\begin{equation*}
		\text{$\boldsymbol{y}_n \wk \boldsymbol{y}$ in $W^{1,p}(\Omega;\R^N)$.}
	\end{equation*}
	Suppose that each $\boldsymbol{y}_n$ satisfies condition {\rm (INV)}.  Then, also $\boldsymbol{y}$ satisfies condition {\rm (INV)}.
\end{lemma}

We include condition (INV) within our definition of admissible deformations in order to exclude pathological behaviors. From a technical point of view, this choice is motivated by the relationship between geometric image, topological image, and surface energy illustrated in Theorem \ref{thm:INV-top-im} below, as well as the representation formula for the distributional determinant provided by   Theorem \ref{thm:Det} below.  Both   results are  proved in \cite{henao.moracorral.lusin}, require condition (INV), and  will be essential for our analysis.

We introduce the following class of maps  
\begin{equation}
	\label{eqn:Y-cav}
	\mathcal{Y}_p^{\rm cav}(\Omega)\coloneqq \left\{ \boldsymbol{y}\in W^{1,p}(\Omega;\R^N):\;\det D \boldsymbol{y}\in L^1_+(\Omega),\;\text{$\boldsymbol{y}$ satisfies (INV)} \right\},
\end{equation}
where the superscript ``cav'' stands for cavitation.    From Lemma \ref{lem:inv-injectivity}, we  see that $\mathcal{Y}_p^{\rm cav}(\Omega)\subset \mathcal{Y}(\Omega)$.

For technical reasons,  we introduce the class of  good domains.   The next definition is a modification of  \cite[Definition 2.17]{barchiesi.henao.moracorral}.    In addition to the properties required in \cite{barchiesi.henao.moracorral},   we include property (ii)  which is motivated by Remark \ref{rem:top-deg}(a) and property (iv) which corresponds to requiring that an analog of   (INV) is satisfied on each good domain.      

\begin{definition}[Good domains]
	\label{def:regular-subdomains}
	Let $\boldsymbol{y}\in W^{1,p}(\Omega;\R^N)$. We define  $\mathcal{U}_{\boldsymbol{y}}$ as the  family  of domains $U\subset \subset \Omega$ of  class  $C^2$ satisfying the following conditions:
	\begin{enumerate}[(i)]
		\item $\boldsymbol{y}\restr{\partial U}\in W^{1,p}(\partial U;\R^N)$ and $(\cof\,\nabla \boldsymbol{y})\restr{\partial U}\in L^1(\partial U;\rnn)$.
		\item $\partial U \subset L_{\boldsymbol{y}}$ and $\boldsymbol{y}^*\restr{\partial U}\in C^0(\partial U;\R^N)$.
		\item $\haus(\partial U \setminus \domg(\boldsymbol{y},\Omega))=0$ and $\nabla^{\partial U}{\boldsymbol{y}}\restr{\partial U}\simeq\nabla \boldsymbol{y}(\boldsymbol{I}-\boldsymbol{\nu}_U\otimes \boldsymbol{\nu}_U)$  on $\partial U$.
		\item   $\boldsymbol{y}(\boldsymbol{x})\in \imt(\boldsymbol{y},U)$ for almost every $\boldsymbol{x}\in U$ and $\boldsymbol{y}(\boldsymbol{x})\notin \imt(\boldsymbol{y},U)$ for almost every $\boldsymbol{x}\in \Omega \setminus U$.  
		\item There holds
		\begin{equation*}
			\displaystyle \lim_{\delta \to 0^+} \dashint_{-\delta}^\delta \left |\int_{\partial U_s} |\cof\hspace{1pt}\nabla \boldsymbol{y}|\,\d\haus- \int_{\partial U} |\cof\hspace{1pt}\nabla \boldsymbol{y}|\,\d\haus\right |\,\d s=0.
		\end{equation*}
		\item For every $\boldsymbol{\psi}\in C^1_{\rm c}(\R^N;\R^N)$, there holds
		\begin{equation*}
			\begin{split}
				&\lim_{\delta \to 0^+} \dashint_{-\delta}^\delta \bigg | \int_{\partial U_s} \boldsymbol{\psi}\circ\boldsymbol{y}\cdot (\cof \hspace{1pt}\nabla\boldsymbol{y})\boldsymbol{\nu}_{U_s}\,\d\haus- \int_{\partial U} \boldsymbol{\psi}\circ\boldsymbol{y}\cdot (\cof \hspace{1pt}\nabla\boldsymbol{y})\boldsymbol{\nu}_{U}\,\d\haus          \bigg |\,\d s=0.
			\end{split}
		\end{equation*}
	\end{enumerate}
\end{definition}

\begin{remark}
	\label{rem:regular-subdomains}
	\begin{enumerate}[(a)]
		\item  By (i)--(ii),  $\boldsymbol{y}^*\restr{\partial U}$ is the continuous representative of  $\boldsymbol{y}\restr{\partial U}$. Thus,  $\boldsymbol{y}^*\restr{\partial U}$ satisfies Lusin's condition (N) by \cite[Corollary~1]{marcus.mizel}. Also, as $\boldsymbol{y}\restr{\domg(\boldsymbol{y},\partial U)}=\boldsymbol{y}^*\restr{\domg(\boldsymbol{y},\partial U)}$ by Remark~\ref{rem:INV},  the map $\boldsymbol{y}^*\restr{\domg(\boldsymbol{y},\partial U)}$ is injective   by    Lemma~\ref{lem:inverse-differentiable}(i). Therefore, in view of (iii), $\boldsymbol{y}^*\restr{\partial U}$ is almost everywhere injective.
		\item  As a consequence of (i), the map $\boldsymbol{y}\restr{\partial U}$ is almost {everywhere} approximately tangentially differentiable. Then, by (iii) and Proposition~\ref{prop:change-of-variable-surface}(ii), we have 
		\begin{equation*}
			\haus(\img(\boldsymbol{y},\partial U))=\int_{\partial U} |(\cof \nabla \boldsymbol{y})\boldsymbol{\nu}_U|\,\d\haus \leq \|\cof \nabla \boldsymbol{y}\|_{L^1(\partial U;\rnn)},
		\end{equation*}
		where the right-hand side is finite thanks to (i). As $\boldsymbol{y}\restr{\domg(\boldsymbol{y},\partial U)}=\boldsymbol{y}^*\restr{\domg(\boldsymbol{y},\partial U)}$ by Remark~\ref{rem:INV},  we see that $\img(\boldsymbol{y},\partial U)\subset \boldsymbol{y}^*(\partial U)$. Then, thanks to (iii) and the fact that $\boldsymbol{y}^*\restr{\partial U}$ satisfies Lusin's condition (N)
		by item (a), we deduce $\boldsymbol{y}^*(\partial U)\simeq\img(\boldsymbol{y},\partial U)$.  In particular, $\haus(\boldsymbol{y}^*(\partial U))<+\infty$. 
		\item We observe that the class $\mathcal{U}_{\boldsymbol{y}}$  depends on the specific representative of $\boldsymbol{y}$.
	\end{enumerate}
\end{remark}

The following lemma states that, roughly speaking, almost every subdomain of $\Omega$ belongs to the class $\mathcal{U}_{\boldsymbol{y}}$.

\begin{lemma}
	\label{lem:abundance}
	Let $\boldsymbol{y}\in \mathcal{Y}^{\rm cav}_p(\Omega)$ and $U \subset \subset \Omega$ be a domain of class $C^2$. Then, there exists $\delta>0$  such that,  for almost every $s\in (-\delta,\delta)$, we have $U_s\in\mathcal{U}_{\boldsymbol{y}}$. 
\end{lemma}
\begin{proof}
 Let $\boldsymbol{y}\in \mathcal{Y}^{\rm cav}_p(\Omega)$ and $U \subset \subset \Omega$ be a domain of class $C^2$. Consider $\delta>0$ for which the tubular neighborhood $\{\boldsymbol{x}\in \Omega: \: -s<d_U(\boldsymbol{x})<s \}$ is well defined. Then, for almost every $s\in (-\delta,\delta)$, the set  $U_s$ satisfies properties (i)--(iii) and (v)--(vi)  of Definition \ref{def:regular-subdomains} by  \cite[Lemma~2.20]{barchiesi.henao.moracorral} and  property (iv) by  \cite[Theorem~9.1 and p.~48]{mueller.spector} (see also \cite[Theorem~A.1]{sivaloganathan.spector.tilakraj}).   
\end{proof}

In the next lemma, we begin to examine the relationship between geometric and topological image. The claims in \eqref{eqn:img-imt} have already been proved in \cite{barchiesi.henao.moracorral,conti.delellis,henao.moracorral.lusin} for different classes of admissible deformations. The observation in \eqref{eqn:img-imt2} was made in \cite[Lemma~7.2]{mueller.spector}.

\begin{lemma}
	\label{lem:img-imt}
	Let $\boldsymbol{y}\in\mathcal{Y}(\Omega)$ and $U\in\mathcal{U}_{\boldsymbol{y}}$. Then, the following holds:
	\begin{equation}
		\label{eqn:img-imt}
		\img(\boldsymbol{y},U)\subset \imt(\boldsymbol{y},U)^{(1)}, \quad \img(\boldsymbol{y},\Omega \setminus U)\subset \R^N \setminus \imt(\boldsymbol{y},U)^{(1)}, \quad \partial\,\imt(\boldsymbol{y},U)=\boldsymbol{y}^*(\partial U).
	\end{equation}
	Moreover,
	\begin{equation}
		\label{eqn:img-imt2}
		\boldsymbol{y}^*(U \cap L_{\boldsymbol{y}}) \subset \closure{\imt(\boldsymbol{y},U)}, \qquad \boldsymbol{y}^*((\Omega \setminus U) \cap L_{\boldsymbol{y}}) \subset \R^N \setminus \imt(\boldsymbol{y},U).
	\end{equation}
\end{lemma} 
\begin{proof}
To prove  the  first inclusion in \eqref{eqn:img-imt}, we argue as in \cite[Lemma~3.8]{conti.delellis}. Denote by $E_1$ the set of points $\boldsymbol{x}\in \domg(\boldsymbol{y},U)$ for which $\boldsymbol{y}(\boldsymbol{x})\in \imt(\boldsymbol{y},U)$. By Remark \ref{rem:geom-dom-im}(a) and 
Definition \ref{def:regular-subdomains}(iv), we have  $\domg(\boldsymbol{y},U) \cong E_1$. Take any $\boldsymbol{x}\in \domg(\boldsymbol{y},U)$. Then, $\Theta^N(E_1,\boldsymbol{x})=1$, so that  $\Theta^N(\boldsymbol{y}(E_1),\boldsymbol{y}(\boldsymbol{x}))=1$ by Lemma \ref{lem:density}. As $\boldsymbol{y}(E_1)\subset \imt(\boldsymbol{y},U)$, we deduce $\Theta^N(\imt(\boldsymbol{y},U),\boldsymbol{y}(\boldsymbol{x}))=1$. 

The second inclusion in \eqref{eqn:img-imt} is proved similarly:  denote  by $E_2$  the set of points $\boldsymbol{x}\in \domg(\boldsymbol{y},\Omega \setminus U)$ such that $\boldsymbol{y}(\boldsymbol{x})\notin \imt(\boldsymbol{y},U)^{(1)}$. Thus, $\domg(\boldsymbol{y},\Omega \setminus {U})\cong E_2$ by Remark \ref{rem:geom-dom-im}(a) and Definition \ref{def:regular-subdomains}(iv). Here, we also use that  $\boldsymbol{y}$ satisfies Lusin's condition (N${}^{-1}$) by Remark \ref{prop:change-of-variable}(b).  Consider any $\boldsymbol{x}\in \domg(\boldsymbol{y},\Omega \setminus {U})$. Then, $\Theta^N(E_2,\boldsymbol{x})=1$ and $\Theta^N(\boldsymbol{y}(E_2),\boldsymbol{y}(\boldsymbol{x}))=1$ by Lemma~\ref{lem:density}. As $\boldsymbol{y}(E_2) \subset \R^N \setminus \imt(\boldsymbol{y},U)^{(1)}$, we deduce $\Theta^N(\R^N \setminus \imt(\boldsymbol{y},U)^{(1)},\boldsymbol{y}(\boldsymbol{x}))=1$. This yields  $\Theta^N(\imt(\boldsymbol{y},U),\boldsymbol{y}(\boldsymbol{x}))=\Theta^N(\imt(\boldsymbol{y},U)^{(1)},\boldsymbol{y}(\boldsymbol{x}))=0$ and hence $\boldsymbol{y}(\boldsymbol{x})\notin \imt(\boldsymbol{y},U)^{(1)}$. 

We look at the  third property   in \eqref{eqn:img-imt}. In view of   Remark~\ref{rem:top-deg}(c),  we only have to show the inclusion $ \boldsymbol{y}^*(\partial U)\subset \partial \,\imt(\boldsymbol{y},U)$. We argue similarly to \cite[Lemma~5.4(b)]{barchiesi.henao.moracorral}. Let $\boldsymbol{x}_0\in \partial U$. By Remark \ref{rem:geom-dom-im}(a) and the fact that $\boldsymbol{x}_0\in L_{\boldsymbol{y}}$, there exists a sequence $(\boldsymbol{x}_n)_n\subset \domg(\boldsymbol{y},U)$ with $\boldsymbol{x}_n \to \boldsymbol{x}_0$ such that $\boldsymbol{y}(\boldsymbol{x}_n)\to \boldsymbol{y}^*(\boldsymbol{x}_0)$. Since $(\boldsymbol{y}(\boldsymbol{x}_n))_n \subset  \closure{\imt(\boldsymbol{y}, U  )}$ by \eqref{eqn:img-imt}, we conclude that $\boldsymbol{y}^*(\boldsymbol{x}_0)\in \closure{\imt(\boldsymbol{y},U)}$. Being $\imt(\boldsymbol{y}, U )$ an open subset of $\R^N \setminus \boldsymbol{y}^*(\partial U)$ by  Remark \ref{rem:top-deg}(c), the claim follows. 

Finally, to prove \eqref{eqn:img-imt2}, we observe that \eqref{eqn:img-imt} entails
\begin{equation*}
	\img(\boldsymbol{y},U)\subset \closure{\imt(\boldsymbol{y},U)}, \qquad \img(\boldsymbol{y},\Omega\setminus U) \subset \R^N \setminus \imt(\boldsymbol{y},U).
\end{equation*}
Then, \eqref{eqn:img-imt2} is deduced from these two {inclusions} thanks to Remark~\ref{rem:geom-dom-im}(a) and Remark~\ref{rem:INV} by exploiting the approximate continuity of $\boldsymbol{y}^*$ on $L_{\boldsymbol{y}}$.
\end{proof}

In \cite[Lemma~3.5]{mueller.spector} it has been shown that, assuming condition (INV) and nonzero Jacobian determinant almost everywhere,  the positivity of the Jacobian determinant is equivalent to that of the topological degree. 
In the next proposition, we present  one of the two implications and we highlight some simple consequences. See \cite[Proposition~2.2 and Lemma~2.4]{mueller.spector.tang} for similar results.

\begin{lemma}[Degree of Sobolev maps]
	\label{lem:orientation}
Let $\boldsymbol{y}\in \mathcal{Y}_p^{\rm cav}(\Omega)$ and $U\in \mathcal{U}_{\boldsymbol{y}}$. Then,
\begin{equation}\label{eqn:deg-nonnegative}
	\deg(\boldsymbol{y},U,\boldsymbol{\xi}) = \chi_{\imt(\boldsymbol{y},U)}(\boldsymbol{\xi}) \quad \text{for all $\boldsymbol{\xi}\in \R^N \setminus \boldsymbol{y}^*(\partial U)$.}
\end{equation}
Moreover,
\begin{equation*}
	\per(\imt(\boldsymbol{y},U))<+\infty, \qquad \partial^* \,\imt(\boldsymbol{y},U) \simeq \boldsymbol{y}^*(\partial U)\simeq \img(\boldsymbol{y},\partial U).
\end{equation*}
\end{lemma}
\begin{proof}
The proof works as the one of \cite[Lemma~3.5]{mueller.spector} by considering the set $U$ in place of a ball, see also \cite[Lemma~3.10]{conti.delellis}. For the convenience of the reader, we briefly recall the  strategy.

Set $d(\boldsymbol{\xi})\coloneqq \deg(\boldsymbol{y},U,\boldsymbol{\xi})$ for all $\boldsymbol{\xi}\in \R^N \setminus \boldsymbol{y}^*(\partial U)$.  By extending $d$ arbitrarily to zero, we consider it as  an integer-valued map defined on the whole $\R^N$.  Let $\boldsymbol{\psi}\in C^1_{\rm c}(\R^N;\R^N)$. By  \cite[Proposition~2.1]{mueller.spector},  Proposition~\ref{prop:change-of-variable-surface}, and items (a) and (b) in Remark~\ref{rem:regular-subdomains}, we have
\begin{equation*}
	\int_{\R^N} d(\boldsymbol{\xi})\,\div \boldsymbol{\psi}(\boldsymbol{\xi})\,\d\boldsymbol{\xi}=\int_{\partial U} \boldsymbol{\psi} (\boldsymbol{y}^*(\boldsymbol{x}))\cdot ((\cof \nabla \boldsymbol{y}^*(\boldsymbol{x}))\boldsymbol{\nu}_U(\boldsymbol{x}))\,\d\haus(\boldsymbol{x})=\int_{\boldsymbol{y}^*(\partial U)} \boldsymbol{\psi}(\boldsymbol{\xi})\cdot \boldsymbol{\nu}_{\img(\boldsymbol{y},U)}(\boldsymbol{\xi})\,\d \boldsymbol{\xi}.
\end{equation*}
This shows that $d\in BV(\R^N;\Z)$ with distributional gradient 
\begin{equation*}
	Dd=-\boldsymbol{\nu}_{\img(\boldsymbol{y},U)} \haus\mres \boldsymbol{y}^*(\partial U).
\end{equation*}
By fine properties of functions with bounded variation (see \cite[Lemma~3.10]{conti.delellis} or \cite[Lemma~3.5]{mueller.spector} for details), we find  that $d$ is a characteristic function. This along with the definition of $d$ shows   $d=\chi_{\imt(\boldsymbol{y},U)}$ and thus   \eqref{eqn:deg-nonnegative}.  Observe that $\imt(\boldsymbol{y},U)=\{d>0\}$ has finite perimeter thanks to the coarea formula \cite[Theorem~3.40]{ambrosio.fusco.pallara}. Hence, 
\begin{equation*}
	Dd=-\boldsymbol{\nu}_{\imt(\boldsymbol{y},U)}\haus \mres \partial^*\,\imt(\boldsymbol{y},U).
\end{equation*}
This proves that $\partial^*\,\imt(\boldsymbol{y},U)\simeq \boldsymbol{y}^*(\partial U)$. As $\img(\boldsymbol{y},\partial U)\simeq \boldsymbol{y}^*(\partial U)$ by Remark~\ref{rem:regular-subdomains}(b), this  concludes the proof. 
\end{proof}

Some of the claims in the next proposition have  already  been  proved in \cite{mueller.spector}.  The fact that topological image is monotone in the  sense  that it preserves inclusions is stated without proof in \cite[Proposition~2.14(iii)]{henao.moracorral.regularity}. Therein, the authors refer to the results in \cite{conti.delellis,henao.moracorral.lusin}, where a different definition of topological image is employed (see Remark~\ref{rem:top-im-inv}(a) below). In \cite{barchiesi.henao.moracorral}, this issue is circumvented by means of a technical expedient, see \cite[Lemma~5.18(a)]{barchiesi.henao.moracorral}. Here, we provide a simple proof of the monotonicity of topological images by showing that the topological image is a regular open set, i.e., it coincides with the interior of its closure. However, our result requires condition (INV), so that it does not hold in the setting of \cite{barchiesi.henao.moracorral}, where deformations are only locally injective.

\begin{proposition}
	\label{prop:top-im-inv}	
	Let $\boldsymbol{y}\in\mathcal{Y}_p^{\rm cav}(\Omega)$.  Then, the following holds: 
	\begin{enumerate}[(i)]
		\item For every $U\in \mathcal{U}_{\boldsymbol{y}}$,  we have 
		\begin{equation}
			\label{eqn:top-prop-imt}
			\haus \left ( \imt(\boldsymbol{y},U)^{(1)} \cap \boldsymbol{y}^*(\partial U) \right )=0, \qquad \imt(\boldsymbol{y},U)= \left( \closure{\imt(\boldsymbol{y},U)} \right)^\circ.
		\end{equation}
		\item For every $U_1,U_2\in \mathcal{U}_{\boldsymbol{y}}$  with $\haus(\partial U_1 \cap \partial U_2)=0$,   we have
		\begin{align}
			\label{eqn:imt-inclusion}
			\text{$\imt(\boldsymbol{y},U_1)\subset \imt(\boldsymbol{y},U_2)$ \:  and \:   $\closure{\imt(\boldsymbol{y},U_1)}\subset \closure{\imt(\boldsymbol{y},U_2)}$} \qquad &\text{whenever \: $U_1 \subset U_2$,}\\
			\label{eqn:imt-disjoint}
			\text{$\imt(\boldsymbol{y},U_1)\cap \imt(\boldsymbol{y},U_2)=\emptyset$} \qquad &\text{whenever \:  $U_1 \cap U_2 =\emptyset$.}	
		\end{align} 
	\end{enumerate}
\end{proposition}
\begin{remark}\label{rem:top-im-inv}
	\begin{enumerate}[(a)]
		\item In \cite{conti.delellis,henao.moracorral.lusin}, the authors define the topological image of  $\boldsymbol{y}$ under $U$ as
		\begin{equation}\label{eq: tildeT}
			\timt(\boldsymbol{y},U)\coloneqq \imt(\boldsymbol{y},U)^{(1)}.
		\end{equation}
		In general, the set $\timt(\boldsymbol{y},U)$ is not open.  We have $\imt(\boldsymbol{y},U)\subset \timt(\boldsymbol{y},U)$, where the inclusion can be strict, see Example~\ref{ex:dolby} below.  
		However,
		\begin{equation*}
			\timt(\boldsymbol{y},U)\setminus \imt(\boldsymbol{y},U)=\timt(\boldsymbol{y},U) \cap \partial\,\imt(\boldsymbol{y},U)=\timt(\boldsymbol{y},U) \cap \boldsymbol{y}^*(\partial U)
		\end{equation*}
		 by  \eqref{eqn:img-imt}, so that \eqref{eqn:top-prop-imt} yields
		\begin{equation*}
			\haus\left (\timt(\boldsymbol{y},U)\setminus \imt(\boldsymbol{y},U)\right )=0.  
		\end{equation*}   
		\item By \eqref{eqn:img-imt}, $\img(\boldsymbol{y},  U   )\subset \timt(\boldsymbol{y},U)$ and $\img(\boldsymbol{y},\Omega \setminus U)\subset \R^N \setminus \timt(\boldsymbol{y},U)$. In particular,
		\begin{equation*}
			\haus(\img(\boldsymbol{y}, U  ) \setminus \imt(\boldsymbol{y}, U  ))=0
		\end{equation*}
		 by item (a).  
	\end{enumerate}
\end{remark}

\begin{proof}
 (i) For convenience, set  $V\coloneqq \imt(\boldsymbol{y},U)$. By  Lemma \ref{lem:img-imt} and  Lemma \ref{lem:orientation}, the set $V$ has finite perimeter and $\haus(\partial V \setminus \partial^*V)=0$. By the Federer-Volpert theorem \cite[Theorem~16.2]{maggi}, this yields $\haus(\partial V \setminus \partial^-V)=0$ and, in particular, $\haus( V^{(1)} \cap \partial V)=0$.   Using $\partial\, V  =\boldsymbol{y}^*(\partial U)$ by \eqref{eqn:img-imt}, this proves the first identity in \eqref{eqn:top-prop-imt}.

For the second identity, set $\widetilde{V}\coloneqq \left (\closure{V}\right )^\circ$. Clearly, $V \subset \widetilde{V}$ and we claim that these two sets actually coincide. We observe that $\widetilde{V}\setminus V \subset V^{(1)} \cap  \partial  V$, so that $\haus (\widetilde{V}\setminus V )=0$. 
Now, by contradiction, we  suppose that there exists $\boldsymbol{\xi}_0 \in \widetilde{V}\setminus V$. Being $\widetilde{V}$ open, we have $B(\boldsymbol{\xi}_0,\varepsilon) \subset \widetilde{V}$ for some $\varepsilon>0$.  Using \eqref{eqn:img-imt} we let   $\boldsymbol{x}_0\in \partial U$ be such that $\boldsymbol{\xi}_0=\boldsymbol{y}^*(\boldsymbol{x}_0)$. Given the continuity of $\boldsymbol{y}^*\restr{\partial U}$, we have $\boldsymbol{y}^*(\partial U \cap B(\boldsymbol{x}_0,\delta))\subset \partial V \cap B(\boldsymbol{\xi}_0,\varepsilon)$ for some $\delta>0$. As $U$ is of class $C^2$, the set $E\coloneqq \partial U \cap B(\boldsymbol{x}_0,\delta)$ satisfies $\haus(E)>0$ and $\boldsymbol{y}^*(E)\subset \widetilde{V} \cap \partial V=\widetilde{V}\setminus V$. Moreover, by Proposition \ref{prop:change-of-variable-surface}, we have
\begin{equation*}
	\haus(\boldsymbol{y}^*(E))=\int_E |(\cof D\boldsymbol{y})\boldsymbol{\nu}_U|\,\d\haus>0,
\end{equation*} 
which provides a contradiction  to $\haus (\widetilde{V}\setminus V )=0$.   Thus, $\widetilde{V}\setminus V=\emptyset$, as desired.

(ii) First of all, we observe that the assumption $\haus(\partial U_1 \cap \partial U_2)=0$ entails
\begin{equation}
	\label{eqn:int}
	\haus(\boldsymbol{y}^*(\partial U_1) \cap \boldsymbol{y}^*(\partial U_2))=0.
\end{equation}
Indeed, since $\boldsymbol{y}^*\restr{\partial U_1 \cup \partial U_2}$ satisfies Lusin's condition (N), $\boldsymbol{y}\restr{\domg(\boldsymbol{y},\partial U_1 \cup \partial U_2)}=\boldsymbol{y}^*\restr{\domg(\boldsymbol{y},\partial U_1 \cup \partial U_2)}$ in view of Remark~\ref{rem:regular-subdomains}(a), and $\boldsymbol{y}\restr{\domg(\boldsymbol{y},\partial U_1 \cup \partial U_2)}$ is injective by Lemma~\ref{lem:inverse-differentiable}(i), we have
\begin{equation*}
\boldsymbol{y}^*(\partial U_1) \cap \boldsymbol{y}^*(\partial U_2) \simeq \boldsymbol{y}^*(\domg(\boldsymbol{y},\partial U_1)) \cap \boldsymbol{y}^*(\domg(\boldsymbol{y},\partial U_2)) = \boldsymbol{y}^*(\domg(\boldsymbol{y},\partial U_1 \cap \partial U_2))
\end{equation*}
with $\haus (\domg(\boldsymbol{y},\partial U_1 \cap \partial U_2) )= 0$.  

We look at \eqref{eqn:imt-inclusion} and we argue as in \cite[Lemma~7.3(i)]{mueller.spector}. Let $U_1 \subset U_2$. 
Recalling Definition~\ref{def:regular-subdomains}(ii)   we have 
\begin{equation}
	\label{eqn:inc}
	\boldsymbol{y}^*(\partial U_1) \subset   \closure{\imt(\boldsymbol{y}, U_2 )}, 
	 \qquad \boldsymbol{y}^*(\partial U_2) \subset \R^N \setminus \imt(\boldsymbol{y},U_1),
\end{equation}
 where we applied  \eqref{eqn:img-imt2}, first for $U_2$ and then for $U_1$.   
   From  \eqref{eqn:img-imt} and  \eqref{eqn:int}--\eqref{eqn:inc}, by applying \cite[Lemma~A.1]{mueller.spector} with $A=\imt(\boldsymbol{y},U_1)$ and $D=\R^N \setminus \closure{\imt(\boldsymbol{y},U_2)}$, we obtain  $A \cap D = \emptyset$ and thus  $\imt(\boldsymbol{y},U_1)\subset \closure{\imt(\boldsymbol{y},U_2)}$.   Taking the closure  on  both sides, we get $\closure{\imt(\boldsymbol{y},U_1)}\subset \closure{\imt(\boldsymbol{y},U_2)}$. From this, by considering the interiors, we have $\imt(\boldsymbol{y},U_1)\subset \imt(\boldsymbol{y},U_2)$ thanks to \eqref{eqn:top-prop-imt}.

Eventually, \eqref{eqn:imt-disjoint} is proved by arguing as  in \cite[Lemma~7.3(ii)]{mueller.spector}.   Suppose that $U_1 \cap U_2 = \emptyset$. By  the second item in \eqref{eqn:img-imt2} and Definition \ref{def:regular-subdomains}(ii),  we have
\begin{equation}
	\label{eqn:inc2}
	\boldsymbol{y}^*(\partial U_1)\subset \R^N \setminus \imt(\boldsymbol{y},U_2), \qquad \boldsymbol{y}^*(\partial U_2)\subset \R^N \setminus \imt(\boldsymbol{y},U_1).
\end{equation}
Given \eqref{eqn:img-imt}, \eqref{eqn:int}, and \eqref{eqn:inc2}, by applying \cite[Lemma~A.1]{mueller.spector} with $A=\imt(\boldsymbol{y},U_1)$ and $D=\imt(\boldsymbol{y},U_2)$, we obtain \eqref{eqn:imt-disjoint}. 
\end{proof}

\begin{figure}
		\begin{tikzpicture}[scale=2,baseline, remember picture]
		\draw[xshift=-50pt,fill=black!20] (-1,-1) rectangle (1,1);
		\draw[xshift=-50pt]  (0,1) node[above] {$U$};
		\draw[xshift=-50pt,->, thick] (1.25,0)  to [out=30,in=150] node[above,midway] {$\boldsymbol{y}$} (2.25,0);
		\draw[xshift=50pt,fill=black!20] (-1,1) arc(90:-90:1) -- cycle;
		\draw[xshift=50pt,fill=black!20] (1,1) arc(90:270:1) -- cycle;
		\draw[xshift=50pt]  (-.5,1) node[above] {$V^-$};
		\draw[xshift=50pt]  (.5,1) node[above] {$V^+$};
		\draw[xshift=50pt]  (0,0) node[above right] {$\boldsymbol{0}$};
		\fill[black,xshift=50pt] (0,0) circle (1pt);
	\end{tikzpicture}
	\caption{The deformation $\boldsymbol{y}$ in Example \ref{ex:dolby}. }
	\label{fig:dolby}
\end{figure}

\begin{example}[Density points of topological image]{\label{ex:dolby}}
Let $N=2$ and $U\coloneqq (-1,1)^2$. Consider the deformation $\boldsymbol{y} \in W^{1,\infty}(U;\R^2)$ defined as
\begin{equation*}
	\boldsymbol{y}(\boldsymbol{x})\coloneqq \left (x_1,x_2 \sqrt{1-(|x_1|-1)^2}\right )^\top.
\end{equation*} 
The map $\boldsymbol{y}$ is almost everywhere injective and transforms $U$ into the union of two half balls touching at one point (see Figure \ref{fig:dolby}). 
One can check that $\boldsymbol{y}$ can be extended as a map in $\mathcal{Y}^{\rm cav}_p(\Omega)$ for some bounded Lipschitz domain $\Omega \subset \R^2$ with $U \subset \subset \Omega$. Note that $\det D \boldsymbol{y}(\boldsymbol{x})= \sqrt{1-(|x_1|-1)^2}>0$ for almost all $\boldsymbol{x}\in U$. By Lemma~\ref{lem:degree-supercritical}, we have 
$\imt(\boldsymbol{y},U)=V^- \cup V^+$, where
\begin{align*}
	V^-&\coloneqq \left\{ (\xi_1,\xi_2) \in (-1,0) \times (-1,1):\: (\xi_1+1)^2+\xi_2^2 <1  \right\},\\
	V^+&\coloneqq \left\{ (\xi_1,\xi_2) \in (0,1) \times (-1,1):\: (\xi_1-1)^2+\xi_2^2 <1  \right\}.
\end{align*}	
Thus, $\timt(\boldsymbol{y},U)=V^- \cup V^+ \cup \{ \boldsymbol{0} \}$. In particular, $\imt(\boldsymbol{y},U)\subset \timt(\boldsymbol{y},U)$ with strict inclusion. 
\end{example}

The topological image of a point has been firstly defined in \cite{sverak} to investigate the fine properties of admissible deformations,  see  also \cite{barchiesi.henao.moracorral,mueller.qi.yan,qi}. Subsequently, the same concept has been employed  to describe cavities created by deformations  in  \cite{conti.delellis,henao,henao.moracorral.lusin, mueller.spector,mueller.spector.tang,sivaloganathan.spector,sivaloganathan.spector.tilakraj}.  

\begin{definition}[Topological image of a point and cavitation points]
	\label{def:topim-point}
	Let $\boldsymbol{y}\in \mathcal{Y}_p^{\rm cav}(\Omega)$ and $\boldsymbol{a}\in \Omega$. The topological image of $\boldsymbol{a}$ under $\boldsymbol{y}$ is defined as
	\begin{equation*}
		\imt(\boldsymbol{y},\boldsymbol{a})\coloneqq \bigcap \left\{ \closure{\imt(\boldsymbol{y},B(\boldsymbol{a},r))}: \hspace{3pt} r>0 \text{  such that }\hspace{3pt}B(\boldsymbol{a},r)\in \mathcal{U}_{\boldsymbol{y}} \right\}.
	\end{equation*}
	The set of cavitation points of $\boldsymbol{y}$ is  defined  by $C_{\boldsymbol{y}}\coloneqq \left\{\boldsymbol{a}\in \Omega: \hspace{2pt} \leb(\imt(\boldsymbol{y},\boldsymbol{a}))>0 \right\}$.
\end{definition}

\begin{remark}
	\label{rem:topim-point}
	\begin{enumerate}[(a)]
		\item In view of   Proposition~\ref{prop:top-im-inv}(ii),   the topological image of a point is given by a decreasing intersection of nonempty compact sets and hence    is also nonempty and compact. 
		\item  Notably, the set $\imt(\boldsymbol{y},\boldsymbol{a})$  depends only on the equivalence class of $\boldsymbol{y}$ and not on the choice of a representative.  This can be shown by arguing as in the case without  cavitations, see  \cite[Remark 5.7(c)]{barchiesi.henao.moracorral}. 
		\item Thanks  to   Proposition~\ref{prop:top-im-inv}(ii),   we have the identities
		\begin{equation*}
			\imt(\boldsymbol{y},\boldsymbol{a})=\bigcap \left\{\closure{\imt(\boldsymbol{y},U)}: \hspace{3pt} U \in \mathcal{U}_{\boldsymbol{y}} \text{ with} \hspace{3pt}\boldsymbol{a}\in U \right \}= \bigcap_{l\in \N}  \closure{\imt(\boldsymbol{y},U_l)},
		\end{equation*}
		where $(U_l)_{l}\subset \mathcal{U}_{\boldsymbol{y}}$ is any  decreasing sequence of sets  with $\boldsymbol{a}\in U_l$  and $U_{l+1} \subset \subset U_l$   for every $l \in \N$ satisfying  $\mathrm{diam}\,U_l \to 0$ as $l \to \infty$. 
		\item In \cite{conti.delellis,henao.moracorral.lusin}, the topological image of $\boldsymbol{a}$ under $\boldsymbol{y}$ is defined as
		\begin{equation*}
			\timt(\boldsymbol{y},\boldsymbol{a})\coloneqq \bigcap \left\{  \timt(\boldsymbol{y},B(\boldsymbol{x}_0,r)): \hspace{3pt} r>0 \text{  such that } \hspace{3pt}B(\boldsymbol{a},r)\in \mathcal{U}_{\boldsymbol{y}} \right\},
		\end{equation*} 
		 where we recall the definition in \eqref{eq: tildeT}.  		The set $\timt(\boldsymbol{y},\boldsymbol{a})$ might be empty  and, in general, is neither open nor closed.
		Clearly, ${\timt}(\boldsymbol{y},\boldsymbol{a})\subset \imt(\boldsymbol{y},\boldsymbol{a})$. Also, $\leb(\imt(\boldsymbol{y},\boldsymbol{a})\setminus \widetilde{\imt}(\boldsymbol{y},\boldsymbol{a}))=0$. To see this, let $(r_l)_l$ be a decreasing sequence  of positive radii  such that $B(\boldsymbol{a},r_l)\in \mathcal{U}_{\boldsymbol{y}}$ for every $l \in \N$ and  $r_l \to 0^+$, as $l \to \infty$. With  Proposition  \ref{prop:top-im-inv}(ii)   in mind, we have  
		\begin{equation*}
	\leb \left (	 {\imt}(\boldsymbol{y},\boldsymbol{a})\right )=\lim_{l \to \infty} \leb\big(\closure{\imt(\boldsymbol{y},B(\boldsymbol{a},r_l))}\big), \qquad  			\leb \left (	{\timt}(\boldsymbol{y},\boldsymbol{a})\right )=\lim_{l \to \infty} \leb\big( \timt(\boldsymbol{y},B(\boldsymbol{a},r_l))\big).
		\end{equation*}
		Now, it suffices to   recall that $\imt(\boldsymbol{y},B(\boldsymbol{a},r_l))\cong\timt(\boldsymbol{y},B(\boldsymbol{a},r_l))$  
 	by  the Lebesgue density theorem.   	
	\end{enumerate}
\end{remark}

 We recall a result from \cite{mueller.spector}.

\begin{lemma}
\label{lem:topim-point-disjoint}	
Let $\boldsymbol{y}\in \mathcal{Y}_p^{\rm cav}(\Omega)$. Then, for every pair of distinct points $\boldsymbol{a}_1,\boldsymbol{a}_2\in \Omega$, there holds $\imt(\boldsymbol{y},\boldsymbol{a}_1)\cap \imt(\boldsymbol{y},\boldsymbol{a}_2)\simeq \emptyset$.	 In particular, $\partial^* \, \imt(\boldsymbol{y},\boldsymbol{a}_1)\cap \partial^* \, \imt(\boldsymbol{y},\boldsymbol{a}_2)\simeq \emptyset$.
\end{lemma}

\begin{proof}
The proof of the first property can be found in \cite[Lemma~7.6]{mueller.spector}. Then, the second follows as  $\partial^* \,  \imt(\boldsymbol{y},\boldsymbol{a}) \subset \imt(\boldsymbol{y},\boldsymbol{a})$ for all $\boldsymbol{a}\in \Omega$ by Remark \ref{rem:topim-point}(a). 
\end{proof}

We define the topological image according to \cite[Equation (7)]{henao.moracorral.regularity}.  

\begin{definition}[Topological image]
	\label{def:topological-image}
	Let $\boldsymbol{y}\in W^{1,p}(\Omega;\R^N)$. 
	The topological image of $\Omega$ under $\boldsymbol{y}$ is defined as
	\begin{equation*}
		\imt(\boldsymbol{y},\Omega)\coloneqq \bigcup_{U \in \mathcal{U}_{\boldsymbol{y}}} \imt(\boldsymbol{y},U).
	\end{equation*}
\end{definition}

\begin{remark}\label{rem:top-im}
	\begin{enumerate}[(a)]
		\item Being the union of open sets, the set $\imt(\boldsymbol{y},\Omega)$ is open.
		\item Arguing as in  the case without cavitation  \cite[Lemma~5.18(b)]{barchiesi.henao.moracorral},  one  can show   that this set depends only on the equivalence class of $\boldsymbol{y}$. This fact is not obvious because of  Remark \ref{rem:regular-subdomains}(c).   
		\item As a consequence of  Remark~\ref{rem:top-im-inv}(b), we have  $\haus(\img(\boldsymbol{y},\Omega)\setminus  \imt(\boldsymbol{y},\Omega))=0$ for $\boldsymbol{y}\in \mathcal{Y}_p^{\rm cav}(\Omega)$. To see this, consider $(U_l)_l\subset \mathcal{U}_{\boldsymbol{y}}$ with $\Omega=\bigcup_{l\in \N} U_l $. Then, $\img(\boldsymbol{y},\Omega)=\bigcup_{l\in\N} \img(\boldsymbol{y},U_l)$. Observing that
		\begin{equation*}
			\haus(\img(\boldsymbol{y},U_l)\setminus \imt(\boldsymbol{y},\Omega))\leq \haus(\img(\boldsymbol{y},U_l)\setminus \imt(\boldsymbol{y},U_l))=0
		\end{equation*}
		for all $l\in \N$ by Remark~\ref{rem:top-im-inv}(b), we deduce the claim. 
		\item By   Proposition \ref{prop:top-im-inv}(ii),  if $\boldsymbol{y}\in \mathcal{Y}_p^{\rm cav}(\Omega)$,  then  
		\begin{equation*}
			\imt(\boldsymbol{y},\Omega)=\bigcup_{l \in \N} \imt(\boldsymbol{y},U_l),
		\end{equation*}
		where $(U_l)_{l} \subset \mathcal{U}_{\boldsymbol{y}}$ is any  increasing sequence with  $U_l \subset  \subset U_{l+1}$ for every $l\in \N$ and  $\Omega=\bigcup_{l \in \N} U_l$.  
	\end{enumerate}
\end{remark}

For notational convenience, we give the following definition.

\begin{definition}[Cavitation image]\label{def:cavitation-image}
	Let $\boldsymbol{y}\in W^{1,p}(\Omega;\R^N)$. For every  subset $E\subset \Omega$, we define  the cavitation image by 
	\begin{equation*}
	\imc(\boldsymbol{y},E)\coloneqq \bigcup_{\boldsymbol{a}\in C_{\boldsymbol{y}}\cap E}\imt(\boldsymbol{y},\boldsymbol{a}).
	\end{equation*}
\end{definition}

\begin{remark}\label{a new remark}
	\begin{enumerate}[(a)]
				\item Thanks to Remark \ref{rem:topim-point}(b), the cavitation image depends only on the equivalence class of the deformations.
		\item Without further assumptions, the set of cavitation points does not need to be finite or countable, see Theorem \ref{thm:INV-top-im}(i) below. In particular, the cavitation image is not necessarily closed.
	\end{enumerate}
\end{remark}

The next theorem is a reformulation of the results in \cite{henao.moracorral.lusin}.    Proofs  have been performed    for maps in $W^{1,N-1}(\Omega;\R^N)\cap L^\infty(\Omega;\R^N)$  but also work for deformations in  $W^{1,p}(\Omega;\R^N)$ as  the arguments  just rely on condition (INV) and  the properties of the degree.   Note that, in \cite{henao.moracorral.lusin},   the topological image of domains and the topological image of points are defined as in Remark~\ref{rem:top-im-inv}(a) and  Remark~\ref{rem:topim-point}(d), respectively. In view of the observation made, this difference does not affect the validity of the results  since $\timt(\boldsymbol{y},U)\simeq \imt(\boldsymbol{y},U)$ and 
 $\imt(\boldsymbol{y},\boldsymbol{a})\cong {\timt}(\boldsymbol{y},\boldsymbol{a})$,   where the latter clearly entails $\partial^*\,\imt(\boldsymbol{y},\boldsymbol{a})= \partial^*\,{\timt}(\boldsymbol{y},\boldsymbol{a})$.  Recall the  functional   $\mathcal{S}$ introduced in Definition \ref{def:surface-energy}. 

\begin{theorem}
	\label{thm:INV-top-im}
	Let $\boldsymbol{y}\in \mathcal{Y}_p^{\rm cav}(\Omega)$ with  $\mathcal{S}(\boldsymbol{y})<+\infty$.  Then, the following holds: 
	\begin{enumerate}[(i)]
		\item  The set $C_{\boldsymbol{y}}$ is countable and $\per(\imt(\boldsymbol{y},\boldsymbol{a}))<+\infty$ for every $\boldsymbol{a}\in C_{\boldsymbol{y}}$. Moreover, we have
		\begin{equation*}
			\mathcal{S}(\boldsymbol{y})=\sum_{\boldsymbol{a} \in C_{\boldsymbol{y}}} \per  \big(  \imt(\boldsymbol{y},\boldsymbol{a})\big).
		\end{equation*}
		\item For every $U\in \mathcal{U}_{\boldsymbol{y}}$, we have 
		\begin{equation}
			\label{eqn:splitting-U}
			\imt(\boldsymbol{y},U)\cong\img(\boldsymbol{y},U) \cup  \imc(\boldsymbol{y},U),
		\end{equation}
		where the union is essentially disjoint,  i.e.,  $\img(\boldsymbol{y},U) \cap \imc(\boldsymbol{y},U)\cong \emptyset $.  
		\item We have 
		\begin{equation}
			\label{eqn:splitting}
			\imt(\boldsymbol{y},\Omega)\cong\img(\boldsymbol{y},\Omega) \cup   \imc(\boldsymbol{y},\Omega),
		\end{equation}
		where the union is essentially disjoint,  i.e., $\img(\boldsymbol{y},\Omega) \cap \imc(\boldsymbol{y},U)\cong \emptyset$.  
		\item For every $U\in  \mathcal{U}_{\boldsymbol{y}}$, we have $\per(\img(\boldsymbol{y},\Omega);\imt(\boldsymbol{y},U))<+\infty$  and $\per(\imc(\boldsymbol{y},\Omega))<+\infty$.  Moreover,  it holds that 
		\begin{equation}
			\label{eqn:bdry-geomim-loc}
			\partial^*\,\img(\boldsymbol{y},\Omega) \cap \imt(\boldsymbol{y},U)\simeq \bigcup_{\boldsymbol{a}\in C_{\boldsymbol{y}} \cap U} \partial^*\,\imt(\boldsymbol{y},\boldsymbol{a}) \simeq \partial^*\,\imc(\boldsymbol{y},U).
		\end{equation}
		\item We have $\per(\img(\boldsymbol{y},\Omega);\imt(\boldsymbol{y},\Omega))<+\infty$ and  $\per(\imc(\boldsymbol{y},\Omega))<+\infty$. Moreover, it holds that  
			\begin{equation}
				\label{eqn:bdry-geomim}
			\partial^*\,\img(\boldsymbol{y},\Omega) \cap \imt(\boldsymbol{y},\Omega)\simeq \bigcup_{\boldsymbol{a}\in C_{\boldsymbol{y}}} \partial^*\,\imt(\boldsymbol{y},\boldsymbol{a})   \simeq \partial^*\,\imc(\boldsymbol{y},\Omega). 
		\end{equation}
	\end{enumerate}	
\end{theorem}

\begin{remark}\label{rem:INV-top-im}
	\begin{enumerate}[(a)]
		\item  From (i) and the isoperimetric inequality, we obtain
		\begin{equation*} \label{eqn:isoperimetric}
			 \sum_{\boldsymbol{a}\in C_{\boldsymbol{y}}} \leb(\imt(\boldsymbol{y},\boldsymbol{a}))\leq C \sum_{\boldsymbol{a}\in C_{\boldsymbol{y}}} \per(\imt(\boldsymbol{y},\boldsymbol{a}))^{\frac{N}{N-1}}\leq C \left( \sum_{\boldsymbol{a}\in C_{\boldsymbol{y}}} \per(\imt(\boldsymbol{y},\boldsymbol{a})) \right)^{\frac{N}{N-1}}=C \mathcal{S}(\boldsymbol{y})^{\frac{N}{N-1}}
		\end{equation*}
		for some constant $C=C(N)>0$.  Given Lemma~\ref{lem:topim-point-disjoint} and the assumption $\mathcal{S}(\boldsymbol{y})<+\infty$, the previous estimate yields $\leb(\imc(\boldsymbol{y},\Omega))<+\infty$. Then, as $\leb(\img(\boldsymbol{y},\Omega))<+\infty$ by Proposition~\ref{prop:change-of-variable}(i), from claim (iii) we deduce $\leb(\imt(\boldsymbol{y},\Omega))<+\infty$.  
		\item From   claims (i)   and  (v)   of Theorem~\ref{thm:INV-top-im} and Lemma \ref{lem:topim-point-disjoint}, it is easy to see that $\mathcal{S}(\boldsymbol{y})\leq \per\left(\img(\boldsymbol{y},\Omega)\right)$ for all $\boldsymbol{y}\in \mathcal{Y}_p^{\rm cav}(\Omega)$, see \cite[Proposition~4.11]{henao.moracorral.lusin}.  In general, without condition (INV), this is  false as observed in Remark~\ref{rem:surface-energy}(b). 
	\end{enumerate}
\end{remark}

\begin{proof}
(i) This  is proved in \cite[Theorem 4.6(ii)]{henao.moracorral.lusin}.

(ii) The claim follows   combining \cite[Proposition 2.17(vi), Theorem 3.2(iii)]{henao.moracorral.lusin}, where the latter can be applied due to  \cite[Theorem 4.6]{henao.moracorral.lusin} and $\mathcal{S}(\boldsymbol{y})<+\infty$.  

(iii) The identity \eqref{eqn:splitting} is stated in  \cite[Equation (9)]{henao.moracorral.regularity}.    It  is deduced from \eqref{eqn:splitting-U} by taking the union over all  $U\in \mathcal{U}_{\boldsymbol{y}}$.  

(iv) The identity \eqref{eqn:bdry-geomim-loc} is established in \cite[Theorem 4.8(i)]{henao.moracorral.lusin}.

(v)   From the first equality in \eqref{eqn:bdry-geomim-loc}, by taking the union over all $U \in \mathcal{U}_{\boldsymbol{y}}$, we obtain the first equality in \eqref{eqn:bdry-geomim}. Thus,  $\per(\img(\boldsymbol{y},\Omega);\imt(\boldsymbol{y},\Omega))<+\infty$ in view of claim (i) and the assumption $\mathcal{S}(\boldsymbol{y})<+\infty$.   The second equality in \eqref{eqn:bdry-geomim} follows as a consequence of Lemma \ref{lem:topim-point-disjoint}.   
\end{proof}

We introduce the distributional determinant. 

\begin{definition}[Distributional determinant]\label{def:distributional-determinant}
	Let $\boldsymbol{y}\in \mathcal{Y}_p^{\rm cav}(\Omega)$. The distributional determinant of $\boldsymbol{y}$ is the distribution  $\Det D \boldsymbol{y}$ defined as
	\begin{equation*}\label{distr-det}
		\langle \Det D \boldsymbol{y}, \varphi \rangle\coloneqq - \frac{1}{N} \int_{\Omega} \left( (\adj D \boldsymbol{y})\boldsymbol{y} \right) 	\cdot D \varphi \,\d\boldsymbol{x}
	\end{equation*}
	for every $\varphi \in C^\infty_{\rm c}(\Omega)$.
\end{definition}
\begin{remark}\label{rem:Det}
Note that $\adj D \boldsymbol{y} \in L^{\frac{p}{N-1}}(\Omega;\rnn)$, while   $\boldsymbol{y}\in L^\infty_\loc(\Omega;\R^N)$  by Remark~\ref{rem:INV}.   Therefore,  $(\adj D \boldsymbol{y})\boldsymbol{y}\in L^1_\loc(\Omega;\R^N)$, and the previous definition sets $\Det D \boldsymbol{y}=N^{-1}\div((\adj D \boldsymbol{y})\boldsymbol{y})$,  where the divergence  on the right-hand side  is understood in the distributional sense. Thus, Definition \ref{def:distributional-determinant} is well posed.
\end{remark}

The following representation formula for the distributional determinant has been firstly established in \cite[Theorem 8.4]{mueller.spector}  under the assumption that the geometric image has finite perimeter, see also \cite[Theorem 4.2]{conti.delellis}. Subsequently,
the same representation was proved  in \cite[Theorem 3.2(i) and Theorem 4.6]{henao.moracorral.lusin} by solely assuming finite surface energy. Also here, the result is proved for deformations in $W^{1,N-1}(\Omega;\R^N)\cap L^\infty(\Omega;\R^N)$, but holds true also in our setting for the same reasons as before.  
  
\begin{theorem}[Distributional determinant representation]
	\label{thm:Det}
Let $\boldsymbol{y}\in \mathcal{Y}_p^{\rm cav}(\Omega)$ with $\mathcal{S}(\boldsymbol{y})<+\infty$. Then, $\Det D \boldsymbol{y}$ belongs to $M_{\rm b}(\Omega)$ and takes the form
\begin{equation}
	\label{eqn:DetDy}
	\Det D \boldsymbol{y}=(\det D \boldsymbol{y})\leb+ \sum_{\boldsymbol{a} \in C_{\boldsymbol{y}}} \leb(\imt(\boldsymbol{y},\boldsymbol{a}))\delta_{\boldsymbol{a}}.
\end{equation}
Moreover,
\begin{equation}
	\label{eqn:DDU}
(\Det D \boldsymbol{y})(U)=\leb(\imt(\boldsymbol{y},U)) \quad \quad  \text{for all 	$U \in \mathcal{U}_{\boldsymbol{y}}$.} 
\end{equation}
\end{theorem}
\begin{proof}
From \cite[Theorem 4.6]{henao.moracorral.lusin}, we deduce a representation formula of the form
\begin{equation*}
	\Det D\boldsymbol{y}=(\det D \boldsymbol{y})\leb+ \sum_{\boldsymbol{a} \in C_{\boldsymbol{y}}} \kappa_{\boldsymbol{a}}\delta_{\boldsymbol{a}}
\end{equation*}
for some constants $(\kappa_{\boldsymbol{a}})_{\boldsymbol{a}\in C_{\boldsymbol{y}}} \subset  (0,+\infty)$. Then, \cite[Theorem 3.2(i)]{henao.moracorral.lusin} yields $\kappa_{\boldsymbol{a}}=\leb(\imt(\boldsymbol{y},\boldsymbol{a}))$ for every $\boldsymbol{a}\in C_{\boldsymbol{y}}$, which proves \eqref{eqn:DetDy}. To show \eqref{eqn:DDU}, let $U\in \mathcal{U}_{\boldsymbol{y}}$. By claim (ii) of Theorem \ref{thm:INV-top-im}, we have
\begin{equation*}
	\leb(\imt(\boldsymbol{y},U))=\leb(\img(\boldsymbol{y},U))+\leb(\imc(\boldsymbol{y},U)). 
\end{equation*}
Then, thanks to Corollary \ref{cor:change-of-variable}(i), Definition \ref{def:cavitation-image}, and Lemma \ref{lem:topim-point-disjoint}, we obtain
\begin{equation*}
	\leb(\imt(\boldsymbol{y},U))=\int_U \det  D  \boldsymbol{y}\,\d\boldsymbol{x}+\sum_{\boldsymbol{a}\in C_{\boldsymbol{y}}  \cap U} \leb(\imt(\boldsymbol{y},\boldsymbol{a})).
\end{equation*}
Comparing the last equation with the representation in \eqref{eqn:DetDy}, the identity \eqref{eqn:DDU} follows.  
\end{proof}

\subsection{Maps of bounded variation}
In this subsection, we collect some results on  (generalized) special maps of bounded  variation.

For the convenience of the reader, we recall the definition of jump points, see \cite[Definition 1(d)]{henao.moracorral.fracture} and \cite[Definition 2(a)]{henao.moracorral.xu}.  To this end, recalling \eqref{eqn:B-S},  given $ \boldsymbol{x}_0  \in \R^N$ and $\boldsymbol{\nu}\in S$,  we introduce the notation   
\begin{equation*}
	H^+( \boldsymbol{x}_0,\boldsymbol{\nu})\coloneqq \left\{\boldsymbol{x}\in \R^N: \hspace{3pt} (\boldsymbol{x}-\boldsymbol{x}_0) \cdot \boldsymbol{\nu} \geq 0  \right\}, \qquad 	H^-(\boldsymbol{x}_0,\boldsymbol{\nu})\coloneqq \left\{\boldsymbol{x}\in \R^N: \hspace{3pt} (\boldsymbol{x}-\boldsymbol{x}_0) \cdot \boldsymbol{\nu} \leq 0  \right\}.
\end{equation*}

\begin{definition}[Jump point, lateral traces, and jump] \label{def:jump}
Let  $\boldsymbol{u}\colon A \subset \R^N \to \R^M$ be measurable.  A point $\boldsymbol{x}_0\in A^{(1)}$ is termed a jump point of $\boldsymbol{u}$ if there exist $\boldsymbol{\nu}\in  S$ and $\boldsymbol{\ell}^+,\boldsymbol{\ell}^-\in \R^N$ with $\boldsymbol{\ell}^+\neq \boldsymbol{\ell}^-$  such that
\begin{equation*}
	\boldsymbol{u}^+(\boldsymbol{x})\coloneqq\aplim_{\substack{\boldsymbol{x}\to \boldsymbol{x}_0 \\ \boldsymbol{x}\in H^+(\boldsymbol{x}_0,\boldsymbol{\nu})}} \boldsymbol{u}(\boldsymbol{x})=\boldsymbol{\ell}^+, \qquad \boldsymbol{u}^-(\boldsymbol{x})\coloneqq \aplim_{\substack{\boldsymbol{x}\to \boldsymbol{x}_0 \\ \boldsymbol{x}\in H^-(\boldsymbol{x}_0,\boldsymbol{\nu})}} \boldsymbol{u}(\boldsymbol{x})=\boldsymbol{\ell}^-.
\end{equation*}
The set of jump points of $\boldsymbol{u}$ is denoted by $J_{\boldsymbol{u}}$. The functions $\boldsymbol{u}^+,\boldsymbol{u}^-\colon J_{\boldsymbol{u}}\to \R^M$ are termed the lateral traces of $\boldsymbol{u}$, and the jump of $\boldsymbol{u}$ is defined as $[\boldsymbol{u}]\coloneqq \boldsymbol{u}^+-\boldsymbol{u}^-$.
\end{definition}

\begin{remark}\label{rem:jump}
\begin{enumerate}[(a)]
	\item If the triplet $(\boldsymbol{\nu},\boldsymbol{\ell}^+,\boldsymbol{\ell}^-)$ satisfies the previous definition, then  $(-\boldsymbol{\nu},\boldsymbol{\ell}^-,\boldsymbol{\ell}^+)$ is the only other triplet satisfying the same conditions \cite[p.~358]{cartesian.currents}. Therefore, the definition of lateral traces and of jump depend on the choice between $\boldsymbol{\nu}$ and $-\boldsymbol{\nu}$. 	We will specify our choice whenever using these definitions.   Note that the set $J_{\boldsymbol{u}}$ is Borel.
Moreover, the map $\boldsymbol{\nu} : J_{\boldsymbol{u}} \to S$ can be chosen to be Borel and, for such a choice, the lateral traces $\boldsymbol{u}^+$, $\boldsymbol{u}^-$ are Borel maps.
We will always work with that choice.  

	\item According to the previous definition, the jump points of $\boldsymbol{u}$ do not necessarily belong to $A$,  but to  $A^{(1)}$. When $A$ is a Lipschitz domain, we have $A=A^{(1)}$ and, in turn,  $J_{\boldsymbol{u}}\subset A$.
	\item If $A_1,A_2\subset \R^N$ satisfy $A_1 \cong {A}_2$, and  $\boldsymbol{u}_1\colon A_1 \to \R^M$ and ${\boldsymbol{u}_2}\colon A_2 \to \R^M$  are measurable with $ \boldsymbol{u}_1  \cong {\boldsymbol{u}}_2$ in $A_1\cap {A}_2$, then $J_{\boldsymbol{u}_1}=J_{{\boldsymbol{u}}_2}$. This property follows since $A_1^{(1)}={A}^{(1)}_2$ and the approximate limits are independent  of the representatives. 
\end{enumerate}
\end{remark}

We make the following  observation about the jump set of the extension   maps defined on sets of finite perimeter.  A similar result was given in \cite[Lemma 2]{henao.moracorral.xu}.  

\begin{lemma}[Jump of extension]
\label{lem:jump-extension}
Let  $D \subset \R^N$ be an open set (possibly unbounded) and let $A\subset D$ be measurable with $\per (A;D) <+\infty$. Let $\boldsymbol{u}\colon A \to \R^M$ be measurable and $\boldsymbol{b}\in \R^M$. Define $\widehat{\boldsymbol{u}}^{\boldsymbol{b}}\colon D \to \R^M$ by setting
\begin{equation*}
	\widehat{\boldsymbol{u}}^{\boldsymbol{b}}(\boldsymbol{x})\coloneqq \begin{cases}
		\boldsymbol{u}(\boldsymbol{x}) & \text{if $\boldsymbol{x}\in A$,} \\ \boldsymbol{b} & \text{if $\boldsymbol{x}\in D \setminus A$.}
	\end{cases}
\end{equation*}
 Then, 
\begin{equation*}
	J_{\widehat{\boldsymbol{u}}^{\boldsymbol{b}}} \cap D \simeq (J_{\boldsymbol{u}} \cap D) \cup ( J_{\widehat{\boldsymbol{u}}^{\boldsymbol{b}}} \cap \partial^* A \cap D).
\end{equation*} 
\end{lemma}

Let $D \subset \R^N$ be an open set (possibly unbounded). For the definition and  basic properties of the spaces $SBV(D;\R^M)$ and $GSBV(D;\R^M)$, we refer to  \cite[Section 4]{ambrosio.fusco.pallara}.  Recall that maps in both these spaces are almost everywhere approximately differentiable.  In few instances, we will  consider the space of piecewise-constant maps that, in view of \cite[Theorem 4.23]{ambrosio.fusco.pallara}, we  define as 
\begin{equation}
	\label{eqn:PC}
PC(D;\R^M)\coloneqq \left\{ \boldsymbol{u}\in SBV(D;\R^M): \hspace{3pt} \nabla \boldsymbol{u}\cong \boldsymbol{O} \right\},
\end{equation}
 where as before $\nabla \boldsymbol{u}$ denotes the approximate gradient of $\boldsymbol{u}$. For $1<q<\infty$, we set
\begin{align}\label{eqn:SBVcapLinftyXXX}
\begin{split}
SBV^q(D;\R^M)&\coloneqq \left\{ \boldsymbol{u}\in SBV(D;\R^M): \hspace{3pt} \nabla \boldsymbol{u}\in L^q(D;\R^{M\times N}), \   \haus(J_{\boldsymbol{u}}) < +\infty    \right\},\\
GSBV^q(D;\R^M)&\coloneqq \left\{ \boldsymbol{u}\in GSBV(D;\R^M): \hspace{3pt} \nabla \boldsymbol{u}\in L^q(D;\R^{M \times N}), \   \haus(J_{\boldsymbol{u}}) < +\infty   \right\}.
\end{split}
\end{align}
As observed in \cite[p.~172]{DalMaso-Francfort-Toader:2005}, there holds
\begin{equation}
	\label{eqn:SBVcapLinfty}
GSBV^q(D;\R^M) \cap L^\infty(D;\R^M)= SBV^q(D;\R^M)\cap L^\infty(D;\R^M).	
\end{equation} 

We will employ the following  compactness  and lower semicontinuity  theorem. Without seeking for generality, we present it in a formulation which is sufficient for our purposes.

\begin{theorem}[{\cite[Theorem 2.2]{ambrosio.compactness}}]
	\label{thm:ambrosio-compactness}
Let $(\boldsymbol{u}_n)_n\subset GSBV^q(D;\R^M)$  be such that
\begin{equation*}
	\sup_{n \in \N} \left\{ \|\boldsymbol{u}_n\|_{L^q(D;\R^M)} + \|\nabla \boldsymbol{u}_n\|_{L^q(D;\R^{M\times N})} + \haus(J_{\boldsymbol{u}_n}) \right\} < + \infty. 
\end{equation*}	
Then, there exists a map $\boldsymbol{u}\in GSBV^q(D;\R^M)$ such that, up to subsequences, we have
\begin{equation*}
	\text{$\boldsymbol{u}_n \to \boldsymbol{u}$ a.e.\ in $D$, \qquad $\boldsymbol{u}_n \wk \boldsymbol{u}$ in $L^q(D;\R^M)$, \qquad $\nabla \boldsymbol{u}_n \wk \nabla \boldsymbol{u}$ in $L^q(D;\R^{M \times N})$.}
\end{equation*}
Moreover,
\begin{align*}
	\haus(J_{\boldsymbol{u}})&\leq \liminf_{n \to \infty} \haus(J_{\boldsymbol{u}_n}),\\
	\|\boldsymbol{u}\|_{L^q(D;\R^M)}&\leq \liminf_{n \to \infty} \|\boldsymbol{u}_n\|_{L^q(D;\R^M)},\\
	\|\nabla \boldsymbol{u}\|_{L^q(D;\R^{M\times N})}&\leq \liminf_{n \to \infty} \|\nabla \boldsymbol{u}_n\|_{L^q(D;\R^{M\times N})}.
\end{align*}
\end{theorem}

In Section \ref{sec:EL}, we will work with the  class of admissible deformations  defined by 
\begin{equation}
\label{eqn:Y-frac}
\mathcal{Y}_{p}^{\rm frac}(\Omega)\coloneqq \big\{ \boldsymbol{y}\in SBV^p(\Omega;\R^N): \hspace{2pt} \det \nabla \boldsymbol{y} \in L^1_+(\Omega), \hspace{4pt} \text{$\boldsymbol{y}$ a.e.\ injective}\big\}
\end{equation}
 for   $p>N-1$, 
where the superscript  ``frac'' stands for ``fracture''. Recalling \eqref{eqn:Y}, we immediately observe that $\mathcal{Y}_{p}^{\rm frac}(\Omega)\subset \mathcal{Y}(\Omega)$.  

\subsection{A lower secontinuity result}

We conclude the section by stating a simple variant of a well-known  result  \cite{ball.currie.olver}, which will be  employed in Section \ref{sec:EL} to establish the lower semicontinuity of elastic energies.

In the case $ Z=\R^M$, the theorem below coincides with \cite[Theorem 5.4]{ball.currie.olver}. Here, just for simplicity, we do not account for dependence of $ \Psi$ on $\boldsymbol{x}\in \Omega$.    The proof of the  theorem  works as the one in \cite{ball.currie.olver}, hence we omit it. 

\begin{theorem}\label{thm:lscT}
Let $\Omega \subset \R^N$ be an open set and $ Z \subset \R^M$ be a measurable set. Also, let $\Psi  \colon \R^P\times  Z  \to [0,+\infty]$ be a continuous map, where $P \in \N$. Suppose that $\boldsymbol{p}\mapsto  \Psi (\boldsymbol{p},\boldsymbol{z})$ is convex for every fixed $\boldsymbol{z}\in  Z $. Define the functional $\mathcal{I}\colon L^1(\Omega;\R^P)\times L^1(\Omega; Z )\to [0,+\infty]$ by setting
\begin{equation*}
	\mathcal{I}(\boldsymbol{f},\boldsymbol{g})\coloneqq \int_\Omega  \Psi(\boldsymbol{f}(\boldsymbol{x}),\boldsymbol{g}(\boldsymbol{x}))\,\d\boldsymbol{x}.
\end{equation*}
Let $(\boldsymbol{f}{}_n)_{ n}\subset L^1(\Omega;\R^P)$ and $(\boldsymbol{g}{}_n)_{ n}\subset L^1(\Omega; Z )$ be such that
\begin{equation*}
	\text{$\boldsymbol{f}_n \wk \boldsymbol{f}$ in $L^1(\Omega;\R^P)$, \qquad $\boldsymbol{g}_n \to \boldsymbol{g}$  in $L^1(\Omega;\R^M)$}
\end{equation*}
for some $\boldsymbol{f}\in L^1(\Omega;\R^P)$ and $\boldsymbol{g}\in L^1(\Omega;\R^M)$. Eventually, assume that $\boldsymbol{g}(\boldsymbol{x})\in  Z$ for almost every $\boldsymbol{x}\in \Omega$. Then,
\begin{equation*}
	\mathcal{I}(\boldsymbol{f},\boldsymbol{g})\leq \liminf_{n \to \infty} \mathcal{I}(\boldsymbol{f}_n,\boldsymbol{g}_n).
\end{equation*}
\end{theorem}
\begin{remark}
\begin{enumerate}[(a)] \label{rem:lscT}
	\item If  $\boldsymbol{f}\colon \Omega \to \R^P$ and $\boldsymbol{g}\colon \Omega \to  Z$ are both measurable, then so is $\boldsymbol{x}\mapsto  \Psi(\boldsymbol{f}(\boldsymbol{x}),\boldsymbol{g}(\boldsymbol{x}))$. In particular, the functional $\mathcal{I}$ is well defined. 
	\item As  stated  in \cite[Theorem 5.4]{ball.currie.olver}, for the maps $\boldsymbol{g}_n$ and $\boldsymbol{g}$ it  would be  sufficient to be measurable (not necessarily integrable) and to converge almost everywhere in $\Omega$.
		\item From $\boldsymbol{g}_n\to \boldsymbol{g}$ almost everywhere in $\Omega$, it follows that $\boldsymbol{g}(\boldsymbol{x})\in  \closure{Z} $ for almost every $\boldsymbol{x}\in \Omega$. Thus, the  assumption $\boldsymbol{g}(\boldsymbol{x})\in  Z $ for almost every $\boldsymbol{x}\in \Omega$ is superfluous whenever $ Z $ is a closed set.
\end{enumerate}
\end{remark}

\section{Convergence results} \label{sec:conv}

In this section, we present some convergence results. First, for general approximately differentiable maps, we investigate the convergence of inverse deformations and that of compositions of Eulerian maps with deformations. Subsequently, for Sobolev maps, we study the convergence properties of deformations creating cavities and we review  convergence properties  of the distributional determinant. We also introduce the concept of material image that will be employed later on.

\subsection{Convergence results for approximately differentiable deformations} \label{subsec:conv-approx}

The results presented in this section are mainly known for Sobolev deformations and we present here extensions to maps which are merely almost everywhere approximately differentiable. In particular, we allow for   deformations  in $SBV$.  

In order to simplify the presentation, we state a preliminary lemma which  has been implicitly proved in \cite[Lemma 2.6]{giacomini.ponsiglione} for $SBV$-maps.

\begin{lemma}
	\label{lem:im-u}
	Let $A\subset \R^N$ be measurable with $\leb(A)<+\infty$, and   $(\boldsymbol{u}_n)_n$ be a sequence of almost everywhere approximately differentiable maps $\boldsymbol{u}_n\colon A \to \R^N$  such that $(\det \nabla \boldsymbol{u}_n)_n$ is equi-integrable.  Also,  let $\boldsymbol{u}\colon A \to \R^N$ be almost everywhere approximately differentiable with $\det \nabla \boldsymbol{u}\in L^1(A)$. Denote  by $D_n$ and $D$ the set of approximate differentiability points of $\boldsymbol{u}_n$ and $\boldsymbol{u}$, respectively.
	Suppose that  there exists a measurable set $E\subset A$ such 
	\begin{equation}
		\label{eqn:l-im-u-s1}
		\text{$\boldsymbol{u}_n \to \boldsymbol{u}$ a.e.\ in $ E$, \qquad $\liminf_{n \to \infty}\leb(\boldsymbol{u}_n( E \cap  D_n))\geq  \leb(\boldsymbol{u}( E \cap  D))$.}
	\end{equation}	
	 Then,
	\begin{equation}
		\label{eqn:l-im-u-s2}
		\text{$\chi_{\boldsymbol{u}_n( E \cap  D_n)}\to \chi_{\boldsymbol{u}( E \cap  D)}$ in $L^1(\R^N)$}.
	\end{equation}
\end{lemma}

Recall that $\boldsymbol{u}_n( E \cap D_n)$ and $\boldsymbol{u}( E \cap D)$ are measurable by Remark \ref{rem:federer}(a). We stress that, in this lemma,  maps  are not required to be injective and their Jacobian determinant can be zero on sets of positive measure.

\begin{proof}
 Let $\varepsilon>0$ be arbitrary. By equi-integrability and by Proposition \ref{prop:change-of-variable}, there exists $\delta=\delta(\varepsilon)>0$ such that
	\begin{equation}
		\label{eqn:l-im-u1}
		\sup_{n \in \N} \leb(\boldsymbol{u}_n( X  \cap D_n))\leq \sup_{n \in \N} \int_{X} |\det \nabla \boldsymbol{u}_n|\,\d\boldsymbol{ x}<\varepsilon/2 \quad \text{for all $X\subset A$ measurable with $\leb(X)<\delta$.}
	\end{equation} 
 Define $\widetilde{E} := E \cap D$.  By Lusin and Egorov theorems, there exists a compact set $ K=K(\delta)$ with $ K \subset \widetilde{E}$ and $\leb(\widetilde{E}\setminus  K )<\delta$ such that  $\boldsymbol{u}_n\restr{ K}$ and $\boldsymbol{u}\restr{ K}$ are continuous for all $n$,  and $\boldsymbol{u}_n \to \boldsymbol{u}$ uniformly on $ K$. By continuity, the set $\boldsymbol{u}( K)$ is compact. Let $V= V(\varepsilon,K)\subset \R^N$ be an open subset  satisfying  $\boldsymbol{u}(K )\subset V$ and $\leb(V \setminus \boldsymbol{u}(K))<\varepsilon/2$. By uniform convergence, $\boldsymbol{u}_n( K)\subset V$ for $n \gg1$ depending on $\varepsilon$  and $K$.  Hence,
	\begin{equation}
		\label{eqn:l-im-u2}
		\leb(\boldsymbol{u}_n(K)\setminus \boldsymbol{u}(\widetilde{E}))\leq \leb(V \setminus \boldsymbol{u}( K ))<  \varepsilon/2, 
	\end{equation}
	where we noted that $\boldsymbol{u}( K )\subset \boldsymbol{u}(\widetilde{E})$.  For $n\gg 1$, setting $\widetilde{E}_n\coloneqq E \cap D_n$, we have  
	\begin{equation*}
		\begin{split}
			\leb(\boldsymbol{u}_n (\widetilde{E}_n)\setminus \boldsymbol{u}(\widetilde{E}))&\leq
			\leb( \boldsymbol{u}_n(\widetilde{E}_n)\setminus \boldsymbol{u}_n(K))+\leb(\boldsymbol{u}_n(K)\setminus \boldsymbol{u}(\widetilde{E}))\\
			  &\leq \leb(\boldsymbol{u}_n(\widetilde{E}_n \setminus K)  )  +\leb(\boldsymbol{u}_n(K)\setminus \boldsymbol{u}(\widetilde{E})).
		\end{split}
	\end{equation*}
	 Note that  $\widetilde{E}_n \setminus K\cong \widetilde{E}\setminus K$ since both $\boldsymbol{u}_n$ and $\boldsymbol{u}$ are almost everywhere approximately differentiable.  Thus,   for the first term on the right-hand side, we obtain $\leb(\boldsymbol{u}_n(\widetilde{E}_n \setminus K))<\varepsilon/2$ thanks to \eqref{eqn:l-im-u1}. This, together with \eqref{eqn:l-im-u2}, yields $\leb(\boldsymbol{u}_n (\widetilde{E}_n) \setminus \boldsymbol{u}(\widetilde{E}))<\varepsilon$.
	As $\varepsilon$ was arbitrary, this proves that
	\begin{equation}
		\label{eqn:l-im-u3}
		\leb(\boldsymbol{u}_n(\widetilde{E}_n)\setminus \boldsymbol{u}(\widetilde{E}))\to 0.
	\end{equation}
	Then, observe that
	\begin{equation*}
		\leb(\boldsymbol{u}(\widetilde{E})\setminus  \boldsymbol{u}_n(\widetilde{E}_n))  =\leb(\boldsymbol{u}(\widetilde{E}))-\leb(\boldsymbol{u}_n(\widetilde{E}_n))+\leb(\boldsymbol{u}_n(\widetilde{E}_n)\setminus \boldsymbol{u}(\widetilde{E})).
	\end{equation*}
	  Passing to the superior limit, as $n \to \infty$, in the previous equation  and  using the second assumption in \eqref{eqn:l-im-u-s1},  we deduce that
	\begin{equation}
		\label{eqn:l-im-u4}
		\leb(\boldsymbol{u}(\widetilde{E})\setminus \boldsymbol{u}_n(\widetilde{E}_n)) \to 0.
	\end{equation}
	The combination of \eqref{eqn:l-im-u3}--\eqref{eqn:l-im-u4} yields \eqref{eqn:l-im-u-s2}.    
\end{proof}

Recall  \eqref{eqn:Y}.  For this class of deformations, the convergence of geometric images has been investigated in \cite[Lemma 4.1 and Remark 2]{mueller.spector} and \cite[Theorem 2]{henao.moracorral.invertibility},  see   also  \cite[Theorem 2.7]{giacomini.ponsiglione} for a similar result for $SBV$-maps. For  Sobolev maps that do not create cavities, the  convergence  of inverses   has been shown in 
\cite[Theorem 6.3]{barchiesi.henao.moracorral}. The next proposition reviews the convergence properties of geometric images and proves   a result for  inverses of deformations in the class $\mathcal{Y}(\Omega)$. In particular, we show the convergence of preimages of deformations which will be exploited  when we prove   the convergence of the compositions of  Eulerian maps with deformations  in   Proposition   \ref{prop:conv-comp} below.

\begin{proposition}[Convergence of inverse deformations]
	\label{prop:approx-diff}
	Let $(\boldsymbol{y}{}_n)_n\subset  \mathcal{Y}(\Omega) $ and $\boldsymbol{y}\in  \mathcal{Y}(\Omega)$.  Suppose that 
	\begin{equation}
		\label{eqn:inv-ass}
		\text{$\boldsymbol{y}_n \to \boldsymbol{y}$ a.e.\ in $\Omega$, \qquad $\det \nabla \boldsymbol{y}_n \wk \det \nabla \boldsymbol{y}$ in $L^1(\Omega)$.}
	\end{equation}
	Then:
	\begin{enumerate}[(i)]
		\item For every $E\subset \Omega$ measurable, we have
			\begin{align}
				\label{eqn:img}
				\text{$\chi_{\img(\boldsymbol{y}_n,E)} \to \chi_{\img(\boldsymbol{y},E)}$ }&\text{in $L^1(\R^N)$},\\
				\label{eqn:img-inv}
				\text{$\chi_{\img(\boldsymbol{y}_n,E)}\boldsymbol{y}^{-1}_n \to \chi_{\img(\boldsymbol{y},E)}\boldsymbol{y}^{-1}$ }&\text{in $L^1(\R^N;\R^N)$.}
			\end{align}
		\item Suppose that the sequence $(\chi_{\img(\boldsymbol{y}_n,\Omega)}\det \nabla \boldsymbol{y}_n^{-1})_n$ is equi-integrable. Then,
		\begin{equation}
			\label{eqn:inv-det-l1}
			\text{$\chi_{\img(\boldsymbol{y}_n,\Omega)}\det \nabla \boldsymbol{y}{}_n^{-1} \wk \chi_{\img(\boldsymbol{y},\Omega)}\det \nabla \boldsymbol{y}^{-1}$ }\text{in $L^1(\R^N)$.}
		\end{equation}
Moreover, for every $F\subset \R^N$ measurable with $\leb(F)<+\infty$, we have
		\begin{equation}
			\label{eqn:inv-preim}
			\text{$\chi_{\boldsymbol{y}_n^{-1}(F \cap\, \img(\boldsymbol{y}_n,\Omega))} \to \chi_{\boldsymbol{y}^{-1}(F \cap\, \img(\boldsymbol{y},\Omega))}$ }\text{in $L^1(\Omega)$.}
		\end{equation}
	\end{enumerate}
\end{proposition}

Recall that deformations in $\mathcal{Y}(\Omega)$ are injective when restricted to their geometric domain,  see  Lemma \ref{lem:inverse-differentiable}. Also, for each $\boldsymbol{y}\in \mathcal{Y}(\Omega)$, we keep denoting the map  $(\boldsymbol{y}\restr{\domg(\boldsymbol{y},\Omega)})^{-1}$ simply by  $\boldsymbol{y}^{-1}$, as in Section \ref{sec:preliminaries}.

\begin{remark}
If we assume the convergence of the Jacobian minors of the deformations, then it is possible to deduce also that  of the Jacobian minors of the inverses. Precisely, suppose \eqref{eqn:inv-ass} and the equi-integrability of $(\chi_{\img(\boldsymbol{y}_n,\Omega)}\det \nabla \boldsymbol{y}{}_n^{-1})_n$. If we have
\begin{equation*}
	\text{$\adj_{N-r}\nabla \boldsymbol{y}_n\wk \adj_{N-r}\nabla \boldsymbol{y}$ in $L^1\left (\Omega;\R^{\binom{N}{r}\times \binom{N}{r}}\right)$}
\end{equation*}
for some $r=1,\dots,N-1$, then 
\begin{equation}\label{eqn:adj-inv}
	\text{$\chi_{\img(\boldsymbol{y}_n,\Omega)}\adj_r(\nabla \boldsymbol{y}_n^{-1})\wk \chi_{\img(\boldsymbol{y},\Omega)}\adj_r(\nabla \boldsymbol{y}^{-1})$ in $L^1\left (\R^N;\R^{\binom{N}{r}\times \binom{N}{r}}\right)$.}
\end{equation}
Here, we implicitly used the identity $\binom{N}{N-r}=\binom{N}{r}$. The case $r=N$, corresponds to \eqref{eqn:inv-det-l1}.   The proof of \eqref{eqn:adj-inv} exploits the structure of the minors of inverse matrices and  proceeds similarly to \cite[Theorem 6.3, Claim (c)]{barchiesi.henao.moracorral} starting from \eqref{eqn:inv-det-l1}. We do not present it as we are not going to employ this result in the following. 
\end{remark}

\begin{proof}  
(i) Let $E \subset \Omega$ be measurable. Denote by $D_n$ and $D$ the set of approximate differentiability points of $\boldsymbol{y}_n$ and $\boldsymbol{y}$, respectively. Note that $\domg(\boldsymbol{y}_n,E)\cong E \cap D_n$ and $\domg(\boldsymbol{y},E)\cong E \cap D$ by Remark~\ref{rem:geom-dom-im}(a), so that $\img(\boldsymbol{y}_n,E)\cong \boldsymbol{y}_n(E \cap D_n)$ and $\img(\boldsymbol{y},E)\cong \boldsymbol{y}(E \cap D)$ by Remark \ref{rem:federer}(a).  Given \eqref{eqn:inv-ass}, using Corollary~\ref{cor:change-of-variable}(i), we obtain
	\begin{equation*}
		\leb(\boldsymbol{y}_n(E \cap D_n))= \int_E \det \nabla \boldsymbol{y}_n\,\d\boldsymbol{x} \to \int_E \det \nabla \boldsymbol{y}\,\d\boldsymbol{x}=\leb(\boldsymbol{y}(E \cap D)).
	\end{equation*}
Thus, Lemma \ref{lem:im-u} yields
\begin{equation*}
	\text{$\chi_{\boldsymbol{y}_n(E\cap D_n)}\to \chi_{\boldsymbol{y}(E\cap D)}$ in $L^1(\R^N)$}
\end{equation*}
which proves \eqref{eqn:img}.

We prove \eqref{eqn:img-inv}.  Let $\varepsilon>0$ be arbitrary. Given \eqref{eqn:inv-ass},   by the  Dunford-Pettis theorem \cite[Theorem 2.54]{fonseca.leoni},  there exists $\delta=\delta(\varepsilon)>0$ such that
	\begin{equation}
		\label{eqn:img7}
		\int_{ X } \det \nabla \boldsymbol{y}\,\d\boldsymbol{x}+\sup_{n \in \N} \int_{ X  } \det \nabla \boldsymbol{y}_n\,\d\boldsymbol{x}<\frac{\varepsilon}{2 \|\boldsymbol{id}\|_{L^\infty(\Omega;\R^N)} } \quad  \text{for all $X\subset {E}$ measurable  with $\leb(X)<\delta$.} 
	\end{equation}
Set $\widetilde{E}\coloneqq \domg(\boldsymbol{y},E) \cap \bigcap_{n=1}^\infty \domg(\boldsymbol{y}_n,E)$. By Lusin and Egorov theorems, there exists a compact set $K=K(\delta) \subset \widetilde{E}$ with $\leb(\widetilde{E}\setminus K)<\delta$ such that $\boldsymbol{y}_n\restr{K}$ and $\boldsymbol{y}\restr{K}$ are continuous  for all $n$,  and  $\boldsymbol{y}_n \to \boldsymbol{y}$ uniformly on $K$. 
	We write
	\begin{equation*}
		\label{eqn:img-8}
		\chi_{\img(\boldsymbol{y}_n,E)}\boldsymbol{y}_n^{-1}-\chi_{\img(\boldsymbol{y},E)}\boldsymbol{y}^{-1}=\chi_{\img(\boldsymbol{y}_n,E\setminus K)}\boldsymbol{y}_n^{-1}  -  \chi_{\img(\boldsymbol{y},E\setminus K)}\boldsymbol{y}^{-1}+\chi_{\boldsymbol{y}_n(K)}\boldsymbol{y}_n^{-1}-\chi_{\boldsymbol{y}(K)}\boldsymbol{y}^{-1}. 
	\end{equation*}
	 Applying Corollary \ref{cor:change-of-variable}(i) with $\psi=|\boldsymbol{y}^{-1}|$ and \eqref{eqn:img7},   we have 
	\begin{equation*}
		\int_{\img(\boldsymbol{y},E\setminus K)} |\boldsymbol{y}^{-1}|\,\d\boldsymbol{\xi}=\int_{E \setminus K} |\boldsymbol{x}|\,\det \nabla \boldsymbol{y}\,\d\boldsymbol{x}  \leq \varepsilon/2
	\end{equation*}
	and analogously
	\begin{equation*}
		\int_{\img(\boldsymbol{y}_n,E\setminus K)} |\boldsymbol{y}_n^{-1}|\,\d\boldsymbol{\xi}=\int_{E \setminus K} |\boldsymbol{x}|\,\det \nabla \boldsymbol{y}_n\,\d\boldsymbol{x}  \leq   \varepsilon/2. 
	\end{equation*}
	Here, while applying \eqref{eqn:img7}, we  used that  $ \leb(\widetilde{E}\setminus K)<\delta$ and $\leb(E \setminus \widetilde{E}) = 0$.    
	Therefore, if we show that
	\begin{equation}
		\label{eqn:img-11}
		\text{$\chi_{\boldsymbol{y}_n(K)}\boldsymbol{y}_n^{-1}\to\chi_{\boldsymbol{y}(K)}\boldsymbol{y}^{-1}$ in $L^1(\R^N;  \R^N)$,}
	\end{equation}
	then 
	\begin{equation*}
		\limsup_{n \to \infty} \|\chi_{\img(\boldsymbol{y}_n,E)}\boldsymbol{y}_n^{-1}-\chi_{\img(\boldsymbol{y},E)}\boldsymbol{y}^{-1}\|_{L^1(\R^N;\R^N)}\leq \varepsilon.
	\end{equation*} 
	As $\varepsilon$ is arbitrary, this proves \eqref{eqn:img-inv}. 
	
	We now   prove \eqref{eqn:img-11}. Taking $K$ in place of $E$ in \eqref{eqn:img}, we deduce $\chi_{\boldsymbol{y}_n(K)}\to \chi_{\boldsymbol{y}(K)}$ in $L^1(\R^N)$ and,   up to subsequences, the same convergence {also holds} almost everywhere. In that case, for almost every $\boldsymbol{\xi}\in \boldsymbol{y}(K)$, we have $\boldsymbol{\xi} \in \boldsymbol{y}_{n}(K)$ for $n \gg 1$ depending on $\boldsymbol{\xi}$.  We claim that $\boldsymbol{y}_{n}^{-1}(\boldsymbol{\xi})\to \boldsymbol{y}^{-1}(\boldsymbol{\xi})$. Since $\boldsymbol{y}\restr{K}$ is continuous and injective  (see Lemma \ref{lem:inverse-differentiable}(i)),  while $K$ is compact, this map is {a} homeomorphism onto its image. Precisely, the map $(\boldsymbol{y}\restr{K})^{-1}$ is uniformly continuous, so that   for $\eta>0$ arbitrary  there exists $\vartheta=\vartheta(\eta)>0$ such that, for every $\widetilde{\boldsymbol{\xi}}_1,\widetilde{\boldsymbol{\xi}}_2\in \boldsymbol{y}(K)$ with $|\widetilde{\boldsymbol{\xi}}_1-\widetilde{\boldsymbol{\xi}}_2|<\vartheta$, we have $|\boldsymbol{y}^{-1}(\widetilde{\boldsymbol{\xi}}_1)-\boldsymbol{y}^{-1}(\widetilde{\boldsymbol{\xi}}_2)|<\eta$. Also, by uniform convergence, $\|\boldsymbol{y}_{n}-\boldsymbol{y}\|_{C^0(K;\R^N)}<\vartheta$ for $n \gg 1$. Now, setting $\boldsymbol{x}_{n}\coloneqq \boldsymbol{y}_{n}^{-1}(\boldsymbol{\xi})\in K$ and $\boldsymbol{\xi}_{n}\coloneqq \boldsymbol{y}(\boldsymbol{x}_{n})\in \boldsymbol{y}(K)$, there holds $|\boldsymbol{\xi}-\boldsymbol{\xi}_{n}|=|\boldsymbol{y}_{n}(\boldsymbol{x}_{n})-\boldsymbol{y}(\boldsymbol{x}_{n})|<\vartheta$ and, in turn,
	\begin{equation*}
		|\boldsymbol{y}^{-1}_{n}(\boldsymbol{\xi})-\boldsymbol{y}^{-1}(\boldsymbol{\xi})|=|\boldsymbol{x}_{n}-\boldsymbol{y}^{-1}(\boldsymbol{\xi})|=|\boldsymbol{y}^{-1}(\boldsymbol{\xi}_{n})-\boldsymbol{y}^{-1}(\boldsymbol{\xi})|<\eta	.
	\end{equation*}
	This proves  $\boldsymbol{y}_{n}^{-1}(\boldsymbol{\xi})\to \boldsymbol{y}^{-1}(\boldsymbol{\xi})$. Thus,    up to subsequences, $\chi_{\boldsymbol{y}_{n}(K)}\boldsymbol{y}_{n}^{-1}\to \chi_{\boldsymbol{y}(K)}\boldsymbol{y}^{-1}$ almost everywhere in $\R^N$. Since $\chi_{\boldsymbol{y}_n(K)}|\boldsymbol{y}^{-1}_n|\leq \|\boldsymbol{id}\|_{L^\infty(\Omega;\R^N)} \chi_{\boldsymbol{y}_n(K)}$ and the sequence on the right-hand side converges in $L^1(\R^N)$, we actually have convergence in $L^1(\R^N;\R^N)$ by the dominated convergence theorem.  The convergence actually holds for the whole sequence thanks to the Urysohn property. This concludes the proof of \eqref{eqn:img-11}.

	(ii) First, similarly to \cite[Theorem 6.3(c)]{barchiesi.henao.moracorral},  we show that 
	\begin{equation}
		\label{eqn:inv-det-cb}
		\text{$\chi_{\img(\boldsymbol{y}_n,\Omega)}\det \nabla \boldsymbol{y}_n^{-1} \wk \chi_{\img(\boldsymbol{y},\Omega)}\det \nabla \boldsymbol{y}^{-1}$ }\text{in $C_{\rm b}^0(\R^N)'$.}
	\end{equation} 
 To this end,  let $\psi\in C^0_{\rm b}(\R^N)$. Applying Corollary \ref{cor:change-of-variable}(ii) for ${\varphi} = \psi \circ \boldsymbol{y}$  and ${\varphi} = \psi \circ \boldsymbol{y}_n$,  and the dominated convergence theorem,  we obtain
	\begin{equation*}
			\int_{\img(\boldsymbol{y}_n,\Omega)} \psi\,\det \nabla \boldsymbol{y}_n^{-1}\,\d\boldsymbol{\xi}= \int_\Omega \psi \circ \boldsymbol{y}_n\,\d\boldsymbol{x} \to  \int_\Omega \psi \circ \boldsymbol{y}\,\d\boldsymbol{x}=\int_{\img(\boldsymbol{y},\Omega)} \psi\,\det \nabla \boldsymbol{y}^{-1}\,\d\boldsymbol{\xi}.
	\end{equation*}
	This proves \eqref{eqn:inv-det-cb}.	
	
	By Corollary \ref{cor:change-of-variable}(ii), we have  
	\begin{equation*}
	  \int_{\img(\boldsymbol{y}_n,\Omega)}\det \nabla \boldsymbol{y}_n^{-1}\,\d\boldsymbol{\xi}=\leb(\Omega),
	\end{equation*}
	so that the sequence $(\chi_{\img(\boldsymbol{y}_n,\Omega)}\det \nabla \boldsymbol{y}_n^{-1})_n$ is bounded in $L^1(\R^N)$.  We are assuming  that the same sequence is also equi-integrable.  Hence, for every $\varepsilon>0$ there exists $\delta=\delta(\varepsilon)>0$ such that
	\begin{equation*}
		 \sup_{n \in \N}\int_{\img(\boldsymbol{y}_n,\Omega)\cap  Y } \det \nabla \boldsymbol{y}_n^{-1}\,\d\boldsymbol{\xi}<\varepsilon\quad  \text{for all $Y\subset \R^N$ measurable with $\leb(Y)<\delta$.}  
	\end{equation*}
	By \eqref{eqn:img}, we have $\leb(\img(\boldsymbol{y}_n,\Omega)\setminus \img(\boldsymbol{y},\Omega))<\delta$ for $n \gg 1$. Thus,
	\begin{equation*}
		\sup_{n \gg 1} \int_{\R^N  \setminus \img(\boldsymbol{y},\Omega)}   \chi_{\img(\boldsymbol{y}_n,\Omega)}   \det \nabla \boldsymbol{y}_n^{-1}\,\d\boldsymbol{\xi}<\varepsilon.
	\end{equation*}
Therefore,   the  Dunford-Pettis theorem  ensures the weak convergence  in $L^1(\R^N)$ of  $(\chi_{\img(\boldsymbol{y}_n,\Omega)}\det \nabla \boldsymbol{y}_n^{-1})_n$ up to subsequences. Hence,
	 \eqref{eqn:inv-det-cb}  and the Urysohn property yield  \eqref{eqn:inv-det-l1}.  
	
	Finally, we prove \eqref{eqn:inv-preim}.  Let   $F\subset \R^N$  be  measurable with $\leb(F)<+\infty$.    For convenience,  set $\widebar{\boldsymbol{y}}_n^{-1}\coloneqq \chi_{\img(\boldsymbol{y}_n,\Omega)} \boldsymbol{y}^{-1}_n$ and   $\widebar{\boldsymbol{y}}^{-1}\coloneqq \chi_{\img(\boldsymbol{y},\Omega)} \boldsymbol{y}^{-1}$   on $\R^N$.    Clearly,
	\begin{equation}\label{eqn:preim}
		\widebar{\boldsymbol{y}}_n^{-1}(F\cap \img(\boldsymbol{y}_n,\Omega))= \boldsymbol{y}_n^{-1}(F \cap \img(\boldsymbol{y}_n,\Omega)), \qquad \widebar{\boldsymbol{y}}^{-1}(F\cap \img(\boldsymbol{y},\Omega))= \boldsymbol{y}^{-1}(F \cap \img(\boldsymbol{y},\Omega)).
	\end{equation}
	By Lemma~\ref{lem:inverse-differentiable}(iii), we have $\nabla \widebar{\boldsymbol{y}}_n^{-1}\cong\chi_{\img(\boldsymbol{y}_n,\Omega)}\nabla \boldsymbol{y}_n^{-1}$  and $\nabla \widebar{\boldsymbol{y}}^{-1}\cong\chi_{\img(\boldsymbol{y},\Omega)}\nabla \boldsymbol{y}^{-1}$. In particular, $\widebar{\boldsymbol{y}}_n^{-1}$ and $\widebar{\boldsymbol{y}}^{-1}$ are  approximately differentiable at any point of $\img(\boldsymbol{y}_n,\Omega)$ and $\img(\boldsymbol{y},\Omega)$, respectively.
Using \eqref{eqn:inv-det-l1}, \eqref{eqn:preim}, and  Corollary \ref{cor:change-of-variable}(ii), we  obtain
	\begin{equation*}
		\begin{split}
		 \leb\big(\overline{\boldsymbol{y}}_n^{-1}(F \cap \img(\boldsymbol{y}_n,\Omega) )\big)&=\int_{F \cap \img(\boldsymbol{y}_n,\Omega)} \det  \nabla  {\boldsymbol{y}}_n^{-1}\,\d\boldsymbol{\xi}\\
		 &\to \int_{F \cap \img(\boldsymbol{y},\Omega)} \det  \nabla  {\boldsymbol{y}}^{-1}\,\d\boldsymbol{\xi}=\leb\big(\overline{\boldsymbol{y}}^{-1}(F  \cap \img(\boldsymbol{y},\Omega) )\big).
		\end{split}
	\end{equation*}
	Observe that $\widebar{\boldsymbol{y}}_n^{-1}\to \overline{\boldsymbol{y}}^{-1}$ in $L^1(\R^N;\R^N)$  by \eqref{eqn:img-inv}. 
	Therefore, by applying Lemma \ref{lem:im-u} for a subsequence of $(\overline{\boldsymbol{y}}_n^{-1}\restr{F})_n$ which converges to $\overline{\boldsymbol{y}}^{-1}\restr{F}$ almost everywhere we conclude that
	\begin{equation*}
		\text{$\chi_{\overline{\boldsymbol{y}}_n^{-1}(F\cap \img(\boldsymbol{y}_n,\Omega) )}\to \chi_{\overline{\boldsymbol{y}}^{-1}(F \cap  \img(\boldsymbol{y},\Omega))}$ in $L^1(\R^N)$.}
	\end{equation*} 
	Thanks to  \eqref{eqn:preim} and the Urysohn property, this proves \eqref{eqn:inv-preim}.   	
\end{proof}

We now address the convergence of compositions of Eulerian maps and deformations. 
The following result extends  \cite[Proposition 3.4]{bresciani.davoli.kruzik} to the case of discontinuous deformations, see also \cite[Subsection 2.2.1]{bresciani.thesis} for the case of Eulerian maps  without a priori $L^\infty$-bounds.   For homeomorphic deformations with bounded distortion,  similar results were proved in \cite[Lemma 5.3]{grandi.etal} and \cite[Lemma 5.3]{grandi.etal2}. A different strategy for proving the convergence of compositions  which relies on the Sobolev  regularity   of deformations and   does not require them to be globally injective   has been conceived in \cite[Theorem 3.2]{bresciani}.

\begin{proposition}[Convergence of compositions]
	\label{prop:conv-comp}
	Let $(\boldsymbol{y}{}_n)_n\subset  \mathcal{Y}(\Omega)$ and $\boldsymbol{y}\in  \mathcal{Y}(\Omega) $. Also, let $(\boldsymbol{v}_n)_n$ be a sequence of measurable  maps $\boldsymbol{v}_n \colon \img(\boldsymbol{y}_n,\Omega) \to \R^M$ and $\boldsymbol{v}\colon  \img(\boldsymbol{y},\Omega) \to  \R^M$.  
	Suppose that \eqref{eqn:inv-ass} holds and both sequences $(\chi_{\img(\boldsymbol{y}_n,\Omega)}\det \nabla \boldsymbol{y}{}_n^{-1})_n$ and 
	$(\boldsymbol{v}_n \circ \boldsymbol{y}{}_n)_n$ are equi-integrable. Eventually, suppose that
	\begin{equation}
		\label{eqn:v-conv-ae}
		\text{$\chi_{\img(\boldsymbol{y}_n,\Omega)}\boldsymbol{v}_n \to \chi_{\img(\boldsymbol{y},\Omega)}\boldsymbol{v}$ a.e.\ in $\R^N$.}
	\end{equation}
	Then,   $\boldsymbol{v}\circ \boldsymbol{y}\in L^1(\Omega;\R^M)$  and we have
	\begin{equation*}
		\text{$\boldsymbol{v}_n \circ \boldsymbol{y}_n \to \boldsymbol{v}\circ \boldsymbol{y}$ in $L^1(\Omega;\R^M)$.}
	\end{equation*}
\end{proposition}

 Recall that the compositions $\boldsymbol{v}_n\circ \boldsymbol{y}_n$ and $\boldsymbol{v}\circ \boldsymbol{y}$ are well defined and measurable by Remark \ref{rem:composition-measurable}(a). 

\begin{proof}
	Let $\varepsilon>0$ be arbitrary. By \eqref{eqn:inv-ass} and Egorov's theorem, there exists a measurable set $E=E(\varepsilon)\subset \Omega$ with $\leb(\Omega \setminus E)<\varepsilon / 2$ such that $\boldsymbol{y}_n \to \boldsymbol{y}$ uniformly on $E$. As $\det \nabla \boldsymbol{y}^{-1}\in L^1(\img(\boldsymbol{y},\Omega))$   by Corollary \ref{cor:change-of-variable}(ii),   there exists $\delta=\delta(\varepsilon)>0$ such that 
	\begin{equation}\label{eqn:eiy}
		\int_Y \det \nabla \boldsymbol{y}^{-1}\,\d\boldsymbol{\xi}<\varepsilon/2 \quad \text{for all $Y\subset \img(\boldsymbol{y},\Omega)$ measurable with $\leb(Y)<\delta$.}
	\end{equation}
	From \eqref{eqn:v-conv-ae}, by applying Lusin and Egorov theorems, we find a measurable set $F=F(\delta)\subset \img(\boldsymbol{y},E)$ with $\leb(\img(\boldsymbol{y},E)\setminus F)<\delta$ such that $\boldsymbol{v}\restr{F}$ is continuous and $\chi_{\img(\boldsymbol{y}_n,\Omega)}\boldsymbol{v}_n \to \boldsymbol{v}$ uniformly on $F$.  In particular,
	\begin{equation*}
		\leb(\domg(\boldsymbol{y},E)\setminus \boldsymbol{y}^{-1}(F))=\leb(\boldsymbol{y}^{-1}(\img(\boldsymbol{y},E)\setminus F))=\int_{\img(\boldsymbol{y},E)\setminus F} \det \nabla \boldsymbol{y}^{-1}\,\d\boldsymbol{\xi}<\varepsilon/2 
	\end{equation*}
	thanks to Corollary \ref{cor:change-of-variable}(ii) and \eqref{eqn:eiy}. Thus, recalling Remark \ref{rem:geom-dom-im}(a), we have
	\begin{equation}\label{eqn:preimageF}
		\leb(\Omega \setminus \boldsymbol{y}^{-1}(F))= \leb(\Omega \setminus E)+\leb(\domg(\boldsymbol{y},E)\setminus \boldsymbol{y}^{-1}(F))<\varepsilon.
	\end{equation}
	
	 We show that, up to subsequences, $\boldsymbol{v}_n \circ \boldsymbol{y}_n \to \boldsymbol{v}\circ \boldsymbol{y}$ almost everywhere in $\boldsymbol{y}^{-1}(F)$. 
	 Let 
	$\eta>0$ be arbitrary. Given our choice of the set $F$, there exist $\vartheta=\vartheta(F,\eta)>0$ and $\widetilde{n}=\widetilde{n}(F,\eta)\in \N$ such that
	\begin{equation} \label{eqn:v-cont}
		\text{$|\boldsymbol{v}(\boldsymbol{\xi}_1)-\boldsymbol{v}(\boldsymbol{\xi}_2)|<\eta/2$ \quad for all $\boldsymbol{\xi}_1,\boldsymbol{\xi}_2 \in F$ with $|\boldsymbol{\xi}_1-\boldsymbol{\xi}_2|<\vartheta$}
	\end{equation}
	  and
	\begin{equation}\label{eqn:unif-v}
		\text{ $\sup_{\boldsymbol{\xi}\in F   \cap \img(\boldsymbol{y}_n,\Omega)   }| \boldsymbol{v}_n(\boldsymbol{\xi}) - \boldsymbol{v}(\boldsymbol{\xi})|<\eta/2$ \quad for all $n\geq \widetilde{n}$.}
	\end{equation} 
 	By Proposition \ref{prop:approx-diff}(ii),  for almost every $\boldsymbol{x}\in \boldsymbol{y}^{-1}(F)$ we have $\boldsymbol{x}\in \boldsymbol{y}_n^{-1}(F \cap \img(\boldsymbol{y}_n,\Omega))$ for $n \gg 1$ depending on $\boldsymbol{x}$ along a not relabeled subsequence. We take any such $n$ with $n\geq \widetilde{n}$ and, noting that $\boldsymbol{y}_n(\boldsymbol{x})\in F \cap \img(\boldsymbol{y}_n,\Omega)$, we write
	\begin{equation*}
		|\boldsymbol{v}_n(\boldsymbol{y}_n(\boldsymbol{x}))-\boldsymbol{v}(\boldsymbol{y}(\boldsymbol{x}))|\leq |\boldsymbol{v}_n(\boldsymbol{y}_n(\boldsymbol{x}))-\boldsymbol{v}(\boldsymbol{y}_n(\boldsymbol{x}))| + |\boldsymbol{v}(\boldsymbol{y}_n(\boldsymbol{x}))-\boldsymbol{v}(\boldsymbol{y}(\boldsymbol{x}))|.
	\end{equation*}
	In the previous equation, the first term on the right-hand side is smaller than $\eta/2$ because of \eqref{eqn:unif-v}. The same holds also for the second one owing to \eqref{eqn:v-cont} by choosing $n\gg 1$ large enough so that $|\boldsymbol{y}_n(\boldsymbol{x})-\boldsymbol{y}(\boldsymbol{x})|<\vartheta$ given the uniform convergence of $\boldsymbol{y}_n$ towards $\boldsymbol{y}$ on $E$. This shows that   $\boldsymbol{v}_n(\boldsymbol{y}_n(\boldsymbol{x}))\to \boldsymbol{v}(\boldsymbol{y}(\boldsymbol{x}))$, as $n\to \infty$.   
	
	Given \eqref{eqn:preimageF}  and the convergence $\boldsymbol{v}_n \circ \boldsymbol{y}_n \to \boldsymbol{v}\circ \boldsymbol{y}$ a.e.\ on $\boldsymbol{y}^{-1}(F)$,  by letting $\varepsilon \to 0^+$ along a sequence and applying a standard diagonal argument, we see that $\boldsymbol{v}_n \circ \boldsymbol{y}_n \to \boldsymbol{v}\circ \boldsymbol{y}$ almost everywhere in $\Omega$ for a not relabeled subsequence. As $(\boldsymbol{v}_n \circ \boldsymbol{y}_n)_n$ is assumed to be equi-integrable, by Vitali's convergence theorem \cite[Theorem 2.24]{fonseca.leoni} we conclude that $\boldsymbol{v} \circ \boldsymbol{y} \in L^1(\Omega;\R^M)$ and $\boldsymbol{v}_n \circ \boldsymbol{y}_n \to \boldsymbol{v}\circ \boldsymbol{y}$ in $L^1(\Omega;\R^M)$, where the latter convergence holds for the whole sequence owing to the Urysohn property.
\end{proof}

\subsection{Convergence results for Sobolev deformations} \label{subsec:conv-sobolev}

In  this  subsection, we focus on the convergence properties of Sobolev deformations. 
The following  result  adapts \cite[Theorem 6.3, Claims (a) and (d)]{barchiesi.henao.moracorral}  to  the more involved  setting allowing for cavities. At the same time, the proof   can also be simplified  because in contrast to \cite{barchiesi.henao.moracorral} deformations  are assumed to     satisfy condition (INV).   We  recall the notion of surface energy in Definition~\ref{def:surface-energy}. 

\begin{proposition}[Convergence of topological images]
	\label{prop:conv-topim}
	Let $(\boldsymbol{y}_n)_n\subset \mathcal{Y}_p^{\rm cav}(\Omega)$ and $\boldsymbol{y}\in \mathcal{Y}_p^{\rm cav}(\Omega)$  with $\mathcal{S}(\boldsymbol{y}) < + \infty$  be such that
	\begin{equation}\label{eqn:weak}
		\text{$\boldsymbol{y}_n \wk \boldsymbol{y}$ in $W^{1,p}(\Omega;\R^N)$. }
	\end{equation}
	Then, the following  holds: 
	\begin{enumerate}[(i)]
		\item For every  compact set $ H\subset \imt(\boldsymbol{y},\Omega)$,  up to subsequences,    we have $H \subset \imt(\boldsymbol{y}_n,\Omega)$ for every $n \in \N$.  
		\item   We have 
		\begin{equation*}
			\text{$\chi_{\imt(\boldsymbol{y}_n,\Omega)}\to \chi_{\imt(\boldsymbol{y},\Omega)}$ in $L^1(\imt(\boldsymbol{y},\Omega))$}
		\end{equation*}
		or equivalently
		\begin{equation*}
			\text{$\chi_{\imt(\boldsymbol{y}_n,\Omega) \cap \imt(\boldsymbol{y},\Omega)} \to \chi_{\imt(\boldsymbol{y},\Omega)}$ in $L^1(\R^N)$.}
		\end{equation*}
	\end{enumerate}
\end{proposition}
\begin{remark}\label{rem:imt-conv}
In contrast  to  the case of deformations that do not create cavities \cite[Theorem 6.3(d)]{barchiesi.henao.moracorral}, in general, we do not have that $\chi_{\imt(\boldsymbol{y}_n,\Omega)}\to \chi_{\imt(\boldsymbol{y},\Omega)}$ in $L^1(\R^N)$. The typical situation in which this convergence fails is when cavitation points move along the sequence towards the boundary of  the domain   as in Example \ref{ex:cavitation-boundary} below.	
\end{remark}
\begin{proof}
(i)  Let $(U_l)_l\subset \mathcal{U}_{\boldsymbol{y}}$ be a  nested   sequence with  $\Omega=\bigcup_{l\in \N} U_l$ as in Remark \ref{rem:top-im}(d). 
  Since $({U}_l)_l$ is increasing and $   $ is compact, there exists $k \in \N$ such that $ H  \subset \imt(\boldsymbol{y},{U}_k)$.  Applying Lemma \ref{lem:abundance} and Lemma \ref{lem:conv-boundaries}, we find $U\in\mathcal{U}_{\boldsymbol{y}}\cap \bigcap_{n\in\N}\mathcal{U}_{\boldsymbol{y}_n}$ with $U_k\subset U$ and a not relabeled subsequence for which $\boldsymbol{y}_n^*\to\boldsymbol{y}^*$ uniformly on $\partial U$. By  Proposition \ref{prop:top-im-inv}(ii),  $ H  \subset \imt(\boldsymbol{y},U_k)\subset \imt(\boldsymbol{y},U)$, so that Lemma \ref{lem:top-deg-conv}(i) yields $ H  \subset \imt(\boldsymbol{y}_n,U)\subset \imt(\boldsymbol{y}_n,\Omega)$ for $n\gg 1$,  or even for all $n \in \N$ after extracting a further subsequence. 
	
(ii)   Let $\boldsymbol{\xi}\in \imt(\boldsymbol{y},\Omega)$. By applying claim (i) with $H=\{ \boldsymbol{\xi}\}$, we find a subsequence indexed by $(n_k)$ depending on $\boldsymbol{\xi}$ for which $\boldsymbol{\xi}\in\imt(\boldsymbol{y}_{n_k},\Omega)$ for all $n\gg 1$. Thanks to the Urysohn property, we actually have $\boldsymbol{\xi}\in \imt(\boldsymbol{y}_n,\Omega)$ for $n\gg 1$. This proves that  $\chi_{\imt(\boldsymbol{y}_n,\Omega)}\to 1$ pointwise in $\imt(\boldsymbol{y},\Omega)$. Therefore,  recalling Remark \ref{rem:INV-top-im}(a),   claim (ii)  follows by the dominated convergence theorem.    
\end{proof}

 Recall Definition \ref{def:cavitation-image}.  The next result addresses the asymptotic  behavior  of the cavitation image.  A result similar to claim (i)  was already  derived in the proof of  \cite[Proposition 3.1]{mora.corral}.

\begin{proposition}[Convergence of  cavitation images] \label{prop:cavities}
	Let $(\boldsymbol{y}{}_n)_n\subset \mathcal{Y}_p^{\rm cav}(\Omega)$    with  $\mathcal{S}(\boldsymbol{y}_n)<+\infty$  for every $n\in\N$ and $\boldsymbol{y}\in \mathcal{Y}_p^{\rm cav}(\Omega)$   with  $\mathcal{S}(\boldsymbol{y})<+\infty$. Suppose that \eqref{eqn:weak} holds and also  
	\begin{equation}\label{eqn:weak-det}
		\text{$\det D \boldsymbol{y}_n \wk \det D \boldsymbol{y}$ in $L^1(\Omega)$. }
	\end{equation}
	 Then, the following  holds: 
	\begin{enumerate}[(i)]
		\item For every domain $U\in \mathcal{U}_{\boldsymbol{y}}\cap \bigcap_{n \in \N} \mathcal{U}_{\boldsymbol{y}_n}$ such that $\boldsymbol{y}_n^* \to \boldsymbol{y}^*$ uniformly on $\partial U$, we have
		\begin{equation*} \label{eqn:imc-loc}
			\text{$\chi_{\imc(\boldsymbol{y}_n,U)}\to \chi_{\imc(\boldsymbol{y},U)}$ in $L^1(\R^N)$.}
		\end{equation*}
	\item Up to subsequences, we have 
	\begin{equation*}
		\text{$\chi_{\imc(\boldsymbol{y}_n,\Omega)}\to \chi_{\imc(\boldsymbol{y},\Omega)}$ in $L^1(\imt(\boldsymbol{y},\Omega))$}	
	\end{equation*}
	or equivalently
	\begin{equation*}
	\text{$\chi_{\imc(\boldsymbol{y}_n,\Omega) \cap \imt(\boldsymbol{y},\Omega)} \to \chi_{\imc(\boldsymbol{y},\Omega)}$ in $L^1(\R^N)$.}	
	\end{equation*}
	\end{enumerate}
\end{proposition} 
\begin{remark}
	 In general,  we do not have $\chi_{\imc(\boldsymbol{y}_n,\Omega)}\to \chi_{\imc(\boldsymbol{y},\Omega)}$ in $L^1(\R^N)$. For example, this convergence fails whenever cavitation points move towards the boundary, as in Example \ref{ex:cavitation-boundary} below.
\end{remark}  
\begin{proof}	
	(i) By Theorem \ref{thm:INV-top-im}(ii),  we have	
	\begin{equation}
		\label{eqn:convcav2}
		\chi_{\imt(\boldsymbol{y},U)}\cong\chi_{\img(\boldsymbol{y},U)}+  \chi_{\imc(\boldsymbol{y},U)},
	\end{equation}
	and analogously
	\begin{equation}
		\label{eqn:convcav1}
		\chi_{\imt(\boldsymbol{y}_n,U)}\cong\chi_{\img(\boldsymbol{y}_n,U)}+\chi_{\imc(\boldsymbol{y}_n,U)}
	\end{equation}
	for every $n \in \N$. By   Lemma \ref{lem:top-deg-conv}(iii),   we have
	\begin{equation*}
		\text{$\chi_{\imt(\boldsymbol{y}_n,U)}\to \chi_{\imt(\boldsymbol{y},U)}$ in $L^1(\R^N)$,}
	\end{equation*}
	while Proposition \ref{prop:approx-diff}(i) yields  
	\begin{equation*}
		\text{$\chi_{\img(\boldsymbol{y}_n,U)}\to \chi_{\img(\boldsymbol{y},U)}$ in $L^1(\R^N)$}
	\end{equation*}
	 for a not relabeled subsequence.  
	Therefore, passing to the limit in \eqref{eqn:convcav1}  and comparing  this  with \eqref{eqn:convcav2}, the desired convergence follows.  Thanks to the Urysohn property, the convergence in (i) holds for the whole sequence.  
	
	(ii) By claim (iii) of Theorem \ref{thm:INV-top-im},  we have	
	\begin{equation}
		\label{eqn:imt-limit}
		\chi_{\imt(\boldsymbol{y},\Omega)}\cong\chi_{\img(\boldsymbol{y},\Omega)}+\chi_{\imc(\boldsymbol{y},\Omega)},
	\end{equation}
	and
	\begin{equation*}
		\chi_{\imt(\boldsymbol{y}_n,\Omega)}\cong\chi_{\img(\boldsymbol{y}_n,\Omega)}+\chi_{\imc(\boldsymbol{y}_n,\Omega)}
	\end{equation*}
	for every $n \in \N$. 
	Thus,
	\begin{equation}
		\label{eqn:imim}
		\chi_{\imt(\boldsymbol{y}_n,\Omega) \cap \imt(\boldsymbol{y},\Omega)} \cong\chi_{\img(\boldsymbol{y}_n,\Omega) \cap \imt(\boldsymbol{y},\Omega)}+\chi_{\imc(\boldsymbol{y}_n,\Omega) \cap \imt(\boldsymbol{y},\Omega)}
	\end{equation}
	for every $n \in \N$. By  Proposition \ref{prop:conv-topim}(ii),  we   get  
	\begin{equation*}
		\text{$\chi_{\imt(\boldsymbol{y}_n,\Omega)\cap \imt(\boldsymbol{y},\Omega)}\to \chi_{\imt(\boldsymbol{y},\Omega)}$ in $L^1(\R^N)$.}
	\end{equation*}
	 On the other hand,  by  Proposition \ref{prop:approx-diff}(i), we have $\chi_{\img(\boldsymbol{y}_n,\Omega)}\to \chi_{\img(\boldsymbol{y},\Omega)}$ in $L^1(\R^N)$  for a not relabeled subsequence, which entails
	\begin{equation*}
		\text{$\chi_{\img(\boldsymbol{y}_n,\Omega)\cap \imt(\boldsymbol{y},\Omega)}\to \chi_{\img(\boldsymbol{y},\Omega)}$ in $L^1(\R^N)$}
	\end{equation*} 
	thanks to   Theorem \ref{thm:INV-top-im}(iii).     Therefore, passing to the limit in \eqref{eqn:imim} and comparing  this  with \eqref{eqn:imt-limit}, the desired convergence follows  thanks to the Urysohn property. 
\end{proof}

The next example illustrates the  phenomenon  of cavitation points moving towards the boundary of the domain along a converging sequence of deformations and clarifies how this behavior can preclude the convergence of topological and cavitation images.

\begin{figure}
	\centering
	\begin{tikzpicture}[scale=2,baseline, remember picture]
		\usetikzlibrary{calc}
		\usetikzlibrary{math}
		\def\a{.15}; 
		\tikzmath{\b=(\a+1)*.5; }
		\draw[xshift=-90pt,fill=black!20] (-1,-1) rectangle (1,1);
		\fill[xshift=-90pt,blue] (1-\a,0) circle (1pt);
		\draw[xshift=-90pt,blue]  (1-\a,0) node[above left] {$\boldsymbol{a}_n$};
		\draw[xshift=-90pt]  (0,1) node[above] {$\Omega$};
		\draw[xshift=-90pt,->, thick] (1.15,0)  to [out=30,in=150] node[above,midway] {$\boldsymbol{f}_n$} (2,0);
		\draw[fill=black!5] (-1,-1) rectangle (1,1);
		\draw[fill=black!20] (-1,-1) rectangle  (\a,1);
		\fill[blue] (0,0) circle (1pt);
		\draw[blue]  (0,0) node[above left] {$\boldsymbol{0}$};
		\draw  (0,1) node[above] {$B_\infty$};
		\draw (0,-1) node[above left] {$\boldsymbol{f}_n(\Omega)$};
		\draw[->, thick] (1.15,0)  to [out=30,in=150] node[above,midway] {$\boldsymbol{u}$} (2,0);
		\draw[xshift=90pt,fill=black!5] (-1,-1) rectangle (1,1);
		\draw[xshift=90pt,fill=black!20] (-1,-1) -- (-1,1) -- (\a,1) --  plot [smooth, samples=20, domain=1:\a, variable=\t] ({((\t+1)/(2*\t))*\a},{(\t+1)/2}) -- (\b,\b) -- (\b,-\b) -- plot [smooth, samples=20, domain=\a:1, variable=\t] ({((\t+1)/(2*\t))*\a},{-(\t+1)/2}) -- (\a,-1) -- (-1,-1);
		\draw[xshift=90pt,fill=white] (-.5,-.5) rectangle (.5,.5);
		\draw[xshift=90pt]  (0,1) node[above] {$A_\infty(1/2,1)$};
		\draw[xshift=90pt] (\a,-1) node[above left] {$\img(\boldsymbol{y}_n,\Omega)$};
	\end{tikzpicture}
	\caption{The deformation $\boldsymbol{y}_n$ in Example \ref{ex:cavitation-boundary}.}
	\label{fig:cavitation-boundary1}
	\bigskip 
	
	\begin{tikzpicture}[scale=2]
	\draw[xshift=-90pt,fill=black!20] (-1,-1) rectangle (1,1);
	\fill[xshift=-90pt,blue] (1,0) circle (1pt);
	\draw[xshift=-90pt,blue]  (1,0) node[above left] {$\boldsymbol{a}$};
	\draw[xshift=-90pt]  (0,1) node[above] {$\Omega$};
	\draw[xshift=-90pt,->, thick] (1.15,0)  to [out=30,in=150] node[above,midway] {$\boldsymbol{f}$} (2,0);
	\draw[fill=black!5] (-1,-1) rectangle (1,1);
	\draw[fill=black!20] (-1,-1) rectangle  (0,1);
	\fill[blue] (0,0) circle (1pt);
	\draw[blue]  (0,0) node[above left] {$\boldsymbol{0}$};
	\draw  (0,1) node[above] {$B_\infty$};
	\draw (0,-1) node[above left] {$\boldsymbol{f}(\Omega)$};
	\draw[->, thick] (1.15,0)  to [out=30,in=150] node[above,midway] {$\boldsymbol{u}$} (2,0);
	\draw[xshift=90pt,fill=black!5] (-1,-1) rectangle (1,1);
	\draw[xshift=90pt,fill=black!20] (-1,-1) rectangle (0,1);
	\draw[xshift=90pt,fill=white] (-.5,-.5) rectangle (.5,.5);
	\draw[xshift=90pt]  (0,1) node[above] {$A_\infty(1/2,1)$};
	\draw[xshift=90pt] (0,-1) node[above left] {$\img(\boldsymbol{y},\Omega)$};
	\end{tikzpicture}
	\caption{The deformation $\boldsymbol{y}$ in Example \ref{ex:cavitation-boundary}. }
	\label{fig:cavitation-boundary2}
\end{figure}

\begin{example}[Cavitation  points moving to the boundary]
	\label{ex:cavitation-boundary}
	Let $N=2$. We employ the notation in \eqref{eqn:q-norm}--\eqref{eqn:q-sphere}.
	Define  $\boldsymbol{u}\colon B_\infty \to \R^2$ as
	\begin{equation*}
		\boldsymbol{u}(\boldsymbol{x})\coloneqq 
		\frac{1}{2}\left(|\boldsymbol{x}|_\infty+1 \right) \frac{\boldsymbol{x}}{|\boldsymbol{x}|_\infty} \quad \text{for all $\boldsymbol{x}\in B_\infty\setminus \{ \boldsymbol{0}  \}$.} 
	\end{equation*}
	This map fixes the boundary $S_\infty$ of $B_\infty$, while opening a square cavity of side $1/2$ at the origin (Figure \ref{fig:cavitation-boundary1}). Precisely, $\boldsymbol{u}(B_\infty\setminus \{\boldsymbol{0}\})=A_\infty(1/2,1)$. For every $n \in \N$, let 
	$f_n\colon [-1,1]\to [-1,2^{-n}  ]$ be given by
	\begin{equation*}
		f_n(t)\coloneqq \begin{cases}
			\dfrac{1}{ 2-2^{-n}}t-\dfrac{1-2^{-n}}{2-2^{-n}}  & \text{if $ \ t \in [-1,1-2^{-n}]$,}\\
			t-(1-2^{-n}) & \text{if $ t \in (1-2^{-n},1]$,}
		\end{cases}
	\end{equation*} 
	and define $\boldsymbol{f}_n\colon B_\infty \to B_\infty$ as $\boldsymbol{f}_n(\boldsymbol{x})\coloneqq (f_n(x_1),x_2)^\top$.  This map is a diffeomorphism that transforms $\closure{B}_\infty$ into $[-1,2^{-n}]\times [-1,1]$. In particular,  $\boldsymbol{f}_n(\boldsymbol{a}_n)=\boldsymbol{0}$ where $\boldsymbol{a}_n\coloneqq (1-2^{-n},0)^\top $.  
 
	Let $\Omega \coloneqq B_\infty$ and consider the deformation $\boldsymbol{y}_n \coloneqq \boldsymbol{u}\circ \boldsymbol{f}_n\colon  \Omega \to \R^2$. Then, $\boldsymbol{y}_n\in \mathcal{Y}_p^{\rm cav}(\Omega)$ by Proposition \ref{prop:radial-coref}. This map  opens a cavity at the point $\boldsymbol{a}_n\in \Omega$. Precisely, we have (Figure \ref{fig:cavitation-boundary1}): 
	\begin{equation*}
		\imt(\boldsymbol{y}_n,\Omega)\subset (-1,b_n)\times (-1,1), \quad \img(\boldsymbol{y}_n,\Omega)\subset \big ( (-1,b_n)\times (-1,1) \big ) \setminus \closure{B}_\infty(1/2),
	\end{equation*}
	where $b_n\coloneqq \frac{1}{2}(2^{-n}+1)$. In particular,
	\begin{equation}\label{eqn:ecavn}
		\text{$\closure{B}_\infty(1/2)=\imt(\boldsymbol{y},\boldsymbol{a}_n)\subset \imt(\boldsymbol{y}_n,\Omega)$ \quad for all $n\in \N$.}
	\end{equation}  
	
	Now, as $n\to \infty$, the sequence $(f_n)_n$ converges pointwise to the function $f\colon [-1,1]\to [-1,0  ]$ given by $	f(t)\coloneqq  \frac{1}{2}t-\frac{1}{2}$. Define $\boldsymbol{f}\colon B_\infty\to B_\infty$ as $\boldsymbol{f}(\boldsymbol{x})\coloneqq (f(x_1),x_2)^\top$. Then, $\boldsymbol{f}$ is a diffeomorphism  that transforms $\closure{B}_\infty$ into $[-1,0]\times[-1,1]$. In particular, $\boldsymbol{f}(\boldsymbol{a})=\boldsymbol{0}$, where $\boldsymbol{a}\coloneqq (1,0)^\top$. 
	We consider the deformation $\boldsymbol{y}\coloneqq \boldsymbol{u}\circ \boldsymbol{f}\colon \Omega \to \R^2$. Thus, $\boldsymbol{y}\in \mathcal{Y}_p^{\rm cav}(\Omega)$  by Proposition \ref{prop:radial-coref}. This map  does not open any cavity since $\boldsymbol{a}\in \partial \Omega$. Moreover, as $b_n \to \frac{1}{2}$, we have (Figure \ref{fig:cavitation-boundary2}): 
	\begin{equation*}
		\imt(\boldsymbol{y},\Omega)=\img(\boldsymbol{y},\Omega) = \big( (-1,0)\times (-1,1) \big) \setminus \closure{B}_\infty(1/2).
	\end{equation*}	
	In particular,
	\begin{equation}\label{eqn:ecav}
		\closure{B}_\infty(1/2) \subset \R^2\setminus \imt(\boldsymbol{y},\Omega).
	\end{equation}

	One checks that both \eqref{eqn:weak}--\eqref{eqn:weak-det} hold. From \eqref{eqn:ecavn}--\eqref{eqn:ecav}, we deduce that $(\chi_{\imt(\boldsymbol{y}_n,\Omega)})_n$ does not converge to $\chi_{\imt(\boldsymbol{y},\Omega)}$ in $L^1(\R^2)$,  while one can show that $\chi_{\imt(\boldsymbol{y}_n,\Omega)\cap \imt(\boldsymbol{y},\Omega)}\to \chi_{\imt(\boldsymbol{y},\Omega)}$ in $L^1(\R^2)$ in agreement with Proposition \ref{prop:conv-topim}(ii). Also, $\imc(\boldsymbol{y}_n,\Omega)=\closure{B}_\infty(1/2)$ for all $n\in \N$, while $\imc(\boldsymbol{y},\Omega)=\emptyset$. Thus, $(\chi_{\imc(\boldsymbol{y}_n,\Omega)})_n$ does not converge to $\chi_{\imc(\boldsymbol{y},\Omega)}$ in $L^1(\R^2)$, while  the convergence in Proposition \ref{prop:cavities}(ii) trivially holds.
\end{example}

Next, we address the convergence properties of the distributional determinant.   Related results have  already  been  proved in the literature,  see, e.g., \cite[Lemma 2.6 and Lemma 2.9]{mora.corral} and   \cite[Lemma 3.2]{sivaloganathan.spector.tilakraj}.

\begin{lemma}[Convergence of distributional determinants]
	\label{lem:cavities-Det}	
	Let $(\boldsymbol{y}{}_n)_n\subset \mathcal{Y}_p^{\rm cav}(\Omega)$  and $\boldsymbol{y}\in \mathcal{Y}_p^{\rm cav}(\Omega)$  satisfy \eqref{eqn:weak}--\eqref{eqn:weak-det}.
	Additionally,  suppose that 
	\begin{equation}\label{eqn:S-bdd}
		\sup_{n\in\N} \mathcal{S}(\boldsymbol{y}_n)<+\infty.  
	\end{equation}
	Then, we have the following: 
	\begin{enumerate}[(i)]
		\item $\Det D \boldsymbol{y}_n \wks \Det D \boldsymbol{y}$ in $M_{\rm b}(\Omega)$
		\item $\displaystyle \sum_{\boldsymbol{a} \in C_{\boldsymbol{y}_n}} \leb(\imt(\boldsymbol{y}_n,\boldsymbol{a})) \delta_{\boldsymbol{a}} \wks \sum_{\boldsymbol{a} \in C_{\boldsymbol{y}}} \leb(\imt(\boldsymbol{y},\boldsymbol{a}))\delta_{\boldsymbol{a}} $ in $M_{\rm b}(\Omega)$.
	\end{enumerate}	
\end{lemma}
\begin{proof}
(i)  We have    $\adj D \boldsymbol{y}_n \wk \adj D \boldsymbol{y}$ in $L^{p/(N-1)}(\Omega;\R^{N \times N})$ by \cite[Theorem 8.20, Part 4]{dacorogna}, while $\boldsymbol{y}_n \to \boldsymbol{y}$ in $L^{q}_{\rm loc}(\Omega;\R^N)$ for every $1\leq q<\infty$  by \cite[Lemma 3.3]{sivaloganathan.spector}.     Hence,  
		\begin{equation*}
			\label{eqn:Det-wk-L1}
			\text{$(\adj D \boldsymbol{y}_n)\boldsymbol{y}_n \wk (\adj D \boldsymbol{y})\boldsymbol{y}$ in $L^1_{\rm loc}(\Omega;\R^N)$,}	
		\end{equation*}
		so that 
		\begin{equation}
			\label{eqn:Det-wks-D'}
			\text{$\Det D \boldsymbol{y}_n \wks \Det D \boldsymbol{y}$ as distributions over $\Omega$}	
		\end{equation}
		 in view of Remark \ref{rem:Det}.  
 
By Theorem~\ref{thm:Det}, we have $(\Det D \boldsymbol{y}_n)_n\subset M_{\rm b}(\Omega)$ with
\begin{equation*}
	\label{eqn:DDet}
	(\Det D \boldsymbol{y}_n)(\Omega)=\|\det D \boldsymbol{y}_n\|_{L^1(\Omega)} + \sum_{\boldsymbol{a} \in C_{\boldsymbol{y}_n}} \leb(\imt(\boldsymbol{y}_n,\boldsymbol{a})) \quad \text{for all $n\in \N$.} 
\end{equation*}
The first term on the right-hand side is uniformly bounded due to \eqref{eqn:weak-det} and the same holds for the second one  by  Remark \ref{rem:INV-top-im}(a) and \eqref{eqn:S-bdd}. 
Thus, the sequence $(\Det D \boldsymbol{y}_n)_n$ is  bounded in $M_{\rm b}(\Omega)$. Therefore, the conclusion follows by combining the compactness theorem for Radon measures with  \eqref{eqn:Det-wks-D'}.

(ii) The claim follows from Theorem \ref{thm:Det},  given  claim (i) and  \eqref{eqn:weak-det}.
\end{proof}

In the literature, a priori upper bounds on the number of cavities have often been imposed, see, e.g. \cite{henao,mora.corral,sivaloganathan.spector,sivaloganathan.spector.tilakraj}. 
Later on,  we will  enforce an even more restrictive constraint by imposing an 
a priori lower bound on the  volume of the cavities.   Indeed, upon assuming the uniform boundedness on the surface energy, a lower bound on the  cavity volumes   entails an upper bound on the number of cavities as a consequence of Remark~\ref{rem:INV-top-im}.

For some given $ \kappa>0$,  we consider the  class of maps  
\begin{equation}
\label{eqn:Y-cav-kappa}
\mathcal{Y}^{\rm cav}_{p,\kappa}(\Omega)\coloneqq \left\{ \boldsymbol{y}\in \mathcal{Y}^{\rm cav}_p(\Omega): \hspace{2pt} \inf_{\boldsymbol{a}\in C_{\boldsymbol{y}}} \leb(\imt(\boldsymbol{y},\boldsymbol{a}))\geq \kappa \right\}
\end{equation}
The next result concerns the convergence of cavitation points under such a priori bounds. The proof employs the same arguments  as  \cite[Proposition 3.1]{mora.corral}.  

\begin{lemma}[Convergence of  cavitation points]
	\label{lem:cavities-Det2}	
	Let $(\boldsymbol{y}{}_n)_n\subset \mathcal{Y}_p^{\rm cav}(\Omega)$ and $\boldsymbol{y}\in \mathcal{Y}_p^{\rm cav}(\Omega)$ satisfy \eqref{eqn:weak}--\eqref{eqn:weak-det} and \eqref{eqn:S-bdd}.
	Then,  the following holds:
	\begin{enumerate}[(i)]
		\item  If $\sup_{n \in \N} \mathscr{H}^0( C_{\boldsymbol{y}_n})\leq  P $ for some $ P  \in \N$, then $\mathscr{H}^0(C_{\boldsymbol{y}})\leq P$. In that case, there exist three  disjoint  set of indices $I_\Omega,I_{\partial \Omega},I_0 \subset \{1,\dots,P\}$ and three corresponding  set of points $\{\boldsymbol{a}_i:\:i\in I_\Omega\}$, $\{\boldsymbol{a}_i:\:i\in I_0\}\subset \Omega$ and $\{\boldsymbol{a}_i:\:i\in I_{\partial \Omega}\}\subset \partial \Omega$ such that, up to subsequences, we can write $C_{\boldsymbol{y}_n}=\{\boldsymbol{a}_i^n: \:i \in  I_\Omega \cup I_{\partial \Omega} \cup I_0\}$   with $\mathscr{H}^0(C_{\boldsymbol{y}_n})=\mathscr{H}^0(I_\Omega)+\mathscr{H}^0(I_{\partial \Omega})+\mathscr{H}^0(I_0)$ for every $n\in \N$  and $C_{\boldsymbol{y}}=\{\boldsymbol{a}_i:\:i\in I_\Omega\}$ in such a way that, as $n \to \infty$, we have
		\begin{align}
			\label{eqn:C1}\boldsymbol{a}^n_i\to \boldsymbol{a}_i \quad &\text{ for all $i\in I_\Omega \cup I_{\partial \Omega} \cup I_0$,}\\
			\label{eqn:C2}\leb(\imt(\boldsymbol{y}_n,\boldsymbol{a}^n_i))\to\leb(\imt(\boldsymbol{y},\boldsymbol{a}_i)) \quad &\text{ for all $i\in I_\Omega$,}\\
			\label{eqn:C3}\leb(\imt(\boldsymbol{y}_n,\boldsymbol{a}^n_i))\to0 \quad &\text{ for all $i\in I_0$.}
		\end{align}
		\item  If  $(\boldsymbol{y}{}_n)_n \subset \mathcal{Y}^{\rm cav}_{p,\kappa}(\Omega)$ for some $\kappa>0$, then $\boldsymbol{y}\in  \mathcal{Y}^{\rm cav}_{p,\kappa}(\Omega)$. In that case,  there  exist two finite {disjoint} sets of indices  $I_\Omega,I_{\partial \Omega}\subset \N$ and two corresponding  sets of points $\{\boldsymbol{a}_i: \: i \in I_\Omega\}\subset \Omega$ and $\{\boldsymbol{a}_i:\: i\in I_{\partial \Omega}  \} \subset \partial \Omega$ such that, up to subsequences, we can  write $C_{\boldsymbol{y}_n}=\{\boldsymbol{a}^n_i:\; i \in I_\Omega \cup I_{\partial \Omega}\}$ {with $\mathscr{H}^0(C_{\boldsymbol{y}_n})=\mathscr{H}^0(I_\Omega)+\mathscr{H}^0(I_{\partial \Omega})$} for every $n \in \N$ and $C_{\boldsymbol{y}}=\{\boldsymbol{a}_i:\;i \in I_\Omega\}$ in such a way that, as $n\to \infty$, we have
		\begin{align}
			\label{eqn:C4}\boldsymbol{a}^n_i\to \boldsymbol{a}_i \quad &\text{ for all $i\in I_\Omega \cup I_{\partial \Omega}$,}\\
			\label{eqn:C5}\leb(\imt(\boldsymbol{y}_n,\boldsymbol{a}^n_i))\to\leb(\imt(\boldsymbol{y},\boldsymbol{a}_i)) \quad &\text{ for all $i\in I_\Omega$.}
		\end{align} 
	\end{enumerate}	
\end{lemma}

\begin{remark} \label{rem:Y-cav-kappa}
\begin{enumerate}[(a)]
	\item  We note  that $\mathscr{H}^0(C_{\boldsymbol{y}})\leq \mathscr{H}^0(I_\Omega)$. The reason for not having an equality is that not all $\boldsymbol{a}_i$ are necessarily different. Repetitions may appear when the coalescence of cavities is observed, i.e., when  two different sequences $(\boldsymbol{a}^n_i)_n$ and $(\boldsymbol{a}^n_j)_n$ of cavitation points  converge to the same point $\boldsymbol{a}_i\in C_{\boldsymbol{y}}$.   
	\item Claim (i) shows that, upon assuming \eqref{eqn:S-bdd}, the class $\{\boldsymbol{y}\in \mathcal{Y}_p^{\rm cav}(\Omega):\: \mathscr{H}^0(C_{\boldsymbol{y}})\leq P\}$ is closed with respect to the convergence in \eqref{eqn:weak}--\eqref{eqn:weak-det} for every $P\in \N$.  The indices  $i\in I_0$ correspond to the cavities  closing in the limit, see Example \ref{ex:cavity-point} and Example \ref{ex:cavity-segment} below.
	\item Claim (ii) shows that, upon assuming \eqref{eqn:S-bdd}, the class $\mathcal{Y}^{\rm cav}_{p,\kappa}(\Omega)$ is closed with respect to the convergence in \eqref{eqn:weak}--\eqref{eqn:weak-det} for every $\kappa>0$. Clearly, the lower bound in \eqref{eqn:Y-cav-kappa} excludes the closing of cavities mentioned  in (b). 
\end{enumerate}
\end{remark}

\begin{proof}
First, we observe that  $\mathcal{S}(\boldsymbol{y})<+\infty$. Indeed, by   Theorem \ref{thm:determinant-surface-energy}  we have
\begin{equation*}
	\mathcal{S}(\boldsymbol{y})\leq \liminf_{n\to \infty} \mathcal{S}(\boldsymbol{y}_n)\leq \sup_{n \in \N} \mathcal{S}(\boldsymbol{y}_n)<+\infty.
\end{equation*}

		(i) By assumption, we have $P_n\coloneqq \mathscr{H}^0(C_{\boldsymbol{y}_n})\leq P$ for every $n\in \N$. Up to subsequences, we have $P_n \to \widebar{P}$ for some $\widebar{P}\in \N$ with $\widebar{P}\leq P$ which entails $P_n=\widebar{P}$ for $n\gg 1$. For any such $n$, we write  $C_{\boldsymbol{y}_n}=\{\boldsymbol{a}^n_1,\dots,\boldsymbol{a}^n_{\widebar{P}}\}$ and $\kappa^n_i\coloneqq \leb(\imt(\boldsymbol{y}_n,\boldsymbol{a}^n_i))$ for $1 \leq i \leq \widebar{P}$.  For convenience, set $I\coloneqq \left\{1,\dots,\widebar{P} \right\}$. For each $i \in I$,  consider the sequences $(\boldsymbol{a}^n_i)_n\subset \Omega$ and $(\kappa^n_i)_n\subset (0,+\infty)$.  In view of \eqref{eqn:S-bdd}, the   latter is uniformly bounded with respect to $n\in \N$ since  
		\begin{equation*}
		\kappa^n_i\leq \sum_{j=1}^{\widebar{P}} \kappa^n_j =\sum_{\boldsymbol{a} \in C_{\boldsymbol{y}_n}} \leb(\imt(\boldsymbol{y}_n,\boldsymbol{a}))\leq C\,\mathcal{S}(\boldsymbol{y}_n)^{\frac{N}{N-1}}  
		\end{equation*}  
		for some $C=C(N)>0$ by Remark \ref{rem:INV-top-im}(a).  
		Therefore, there exist $\boldsymbol{a}_1,\dots,\boldsymbol{a}_{\widebar{P}}\in \closure{\Omega}$ and $\kappa_1,\dots,\kappa_{\widebar{P}}\in [0,+\infty)$ such that $\boldsymbol{a}^n_i \to \boldsymbol{a}_i$ and $\kappa^n_i \to \kappa_i$, as $n \to\infty$, for all $i\in I$ upon the extraction of a subsequence. 
		
		We write $I=I_\Omega \cup I_{\partial \Omega} \cup I_0$, where we set
		\begin{equation*}
			I_\Omega\coloneqq \{i\in I: \hspace{2pt} \boldsymbol{a}_i \in \Omega, \hspace{2pt} \kappa_i>0 \}, \quad I_0\coloneqq \{i\in I:\: \boldsymbol{a}_i\in \Omega,\:\kappa_i=0 \},\quad I_{\partial \Omega}\coloneqq \{i\in I:\:\boldsymbol{a}_i\in \partial \Omega\}.
		\end{equation*}
		Thus, \eqref{eqn:C1}--\eqref{eqn:C3} are proved. Moreover, 
		\begin{equation*}
		\text{$\sum_{\boldsymbol{a}\in C_{\boldsymbol{y}_n}} \leb(\imt(\boldsymbol{y}_n,\boldsymbol{a}))\,\delta_{\boldsymbol{a}}=\sum_{i=1}^{\widebar{P}} \kappa^n_i\,\delta_{\boldsymbol{a}^n_i}\wks \sum_{i \in I_\Omega \cup I_0} \kappa_i\,\delta_{\boldsymbol{a}_i}=\sum_{i \in I_\Omega} \kappa_i\,\delta_{\boldsymbol{a}_i}$ in $M_{\rm b}(\Omega)$.}
		\end{equation*}
		In view of Lemma \ref{lem:cavities-Det}(ii), we conclude that $C_{\boldsymbol{y}}=\{\boldsymbol{a}_i:\hspace{2pt}i \in I_\Omega\}$ and $\kappa_i=\leb(\imt(\boldsymbol{y},\boldsymbol{a}_i))$ for every $i \in I_\Omega$. In particular, $\mathscr{H}^0(C_{\boldsymbol{y}})  \le \mathscr{H}^0(I_{\Omega})     \leq \widebar{P}\leq P$.  
		
		(ii) By assumption, $\leb(\imt(\boldsymbol{y}_n,\boldsymbol{a}))\geq \kappa$ for every $\boldsymbol{a}\in C_{\boldsymbol{y}_n}$ and $n \in \N$. From Remark \ref{rem:INV-top-im}(a), we obtain
		\begin{equation*}
		\kappa \,\mathscr{H}^0(C_{\boldsymbol{y}_n})\leq \sum_{\boldsymbol{a} \in C_{\boldsymbol{y}_n}} \leb(\imt(\boldsymbol{y}_n,\boldsymbol{a})) \leq C \,\mathcal{S}(\boldsymbol{y}_n)^{\frac{N}{N-1}}.
		\end{equation*}  
		Thus, by \eqref{eqn:S-bdd}, we deduce  $\sup_{n \in \N} \mathscr{H}^0(C_{\boldsymbol{y}_n})\leq P$ for some $P \in \N$. By claim (i),  we find three {disjoint} sets $I_\Omega,I_{\partial \Omega}, I_0\subset \N$ and  $\{\boldsymbol{a}_i:\: i\in I_\Omega \cup I_{\partial \Omega } \cup I_0\}\subset \closure{\Omega}$ such that, for a not relabeled subsequence, we have   $C_{\boldsymbol{y}_n}=\{\boldsymbol{a}_i^n: \:i \in  I_\Omega \cup I_{\partial \Omega} \cup I_0\}$ for every $n\in \N$ and $C_{\boldsymbol{y}}=\{\boldsymbol{a}_i:\:i\in I_\Omega\}$. Also, \eqref{eqn:C1}--\eqref{eqn:C3} hold.  Suppose by contradiction that there exists  $i\in I_0$. By assumption $\kappa^n_i \geq \kappa$ for every $n\in \N$ and, in turn, $\liminf_{n \to \infty} \kappa^n_i \geq \kappa>0$. However,  $\kappa_i^n\coloneqq \leb(\imt(\boldsymbol{y}_n,\boldsymbol{a}^n_i))\to 0$ by \eqref{eqn:C3}. This contradiction shows that $I_0=\emptyset$, so that \eqref{eqn:C1}--\eqref{eqn:C2} clearly yield \eqref{eqn:C4}--\eqref{eqn:C5}.  
\end{proof}

Recall  Definition \ref{def:topological-image} and Definition \ref{def:cavitation-image}.   We introduce a new concept of image which is particularly convenient in the case of deformations creating a {finite} number of cavities.

\begin{definition}[Material image] \label{def:material-image}
Let $\boldsymbol{y}\in \mathcal{Y}_p^{\rm cav}(\Omega)$  with   $\mathcal{S}(\boldsymbol{y})<+\infty$.  The material image of $\Omega$ under $\boldsymbol{y}$ is defined as
\begin{equation*}
	\imm(\boldsymbol{y},\Omega)\coloneqq \imt(\boldsymbol{y},\Omega) \setminus \imc(\boldsymbol{y},\Omega).
\end{equation*}
\end{definition}

\begin{remark} \label{rem:matim}
	\begin{enumerate}[(a)]
		\item By Theorem \ref{thm:INV-top-im}(iii), we have that $\imm(\boldsymbol{y},\Omega)\cong \img(\boldsymbol{y},\Omega)$. 
		\item If $\mathscr{H}^0(C_{\boldsymbol{y}})<\infty$, then $	\imm(\boldsymbol{y},\Omega)$ is an open set.  Indeed, in that case, $\imc(\boldsymbol{y},\Omega)$ is  compact being a finite union of such sets.   This is the advantage compared to $\img(\boldsymbol{y},\Omega)$ which is in general not open.  
		\item The material image is uniquely determined by the equivalence class of $\boldsymbol{y}$ as the same property holds for the topological image,  
see    Remark \ref{rem:top-im}(b)  and Remark \ref{a new remark}(a). 
	\end{enumerate}
\end{remark}

The following result improves claim (i) of Proposition \ref{prop:conv-topim} for deformations in the class $\mathcal{Y}^{\rm cav}_{p,\kappa}(\Omega)$. 

\begin{proposition}[Convergence of material images]
	\label{prop:conv-matim}
	Let $(\boldsymbol{y}_n)_n\subset \mathcal{Y}^{\rm cav}_{p,\kappa}(\Omega)$ and $\boldsymbol{y}\in \mathcal{Y}^{\rm cav}_{p,\kappa}(\Omega)$  for some $\kappa>0$ be such that \eqref{eqn:weak}--\eqref{eqn:weak-det} and \eqref{eqn:S-bdd} hold. Then, for every compact set $H \subset \imm(\boldsymbol{y},\Omega)$, up to subsequences, we have $H \subset \imm(\boldsymbol{y}_n,\Omega)$ for every $n \in \N$.
\end{proposition}

\begin{remark}\label{nice new remark}
(a) From Proposition \ref{prop:approx-diff}(i), owing to \eqref{eqn:weak}--\eqref{eqn:weak-det} and Remark \ref{rem:matim}(a), we immediately see that $\chi_{\imm(\boldsymbol{y}_n,\Omega)}\to \chi_{\imm(\boldsymbol{y},\Omega)}$ in $L^1(\R^N)$.

  (b) Note that the  validity of Proposition~\ref{prop:conv-matim}  crucially relies on the lower bound  imposed in \eqref{eqn:Y-cav-kappa},   see   Example~\ref{ex:cavity-point} and Example~\ref{ex:cavity-segment} below.  
\end{remark}

\begin{proof}
	 First of all, note that $\mathcal{S}(\boldsymbol{y})<+\infty$ as in the proof of Lemma \ref{lem:cavities-Det2}. 
	By  Proposition \ref{prop:conv-topim}(i), up to subsequences, we have $H \subset \imt(\boldsymbol{y}_n,\Omega)$ for every $n \in \N$. 
	Applying claim (ii) of Lemma \ref{lem:cavities-Det2},  we select a  not relabeled subsequence  for which we can write
	$C_{\boldsymbol{y}_n}= \{\boldsymbol{a}^n_i:\;  i\in I_\Omega \cup I_{\partial \Omega}  \}$ for every $n \in \N$. Here,  $I_\Omega,I_{\partial \Omega}\subset \N$   are two finite disjoint sets of indices and there  exist two corresponding set of points $\{\boldsymbol{a}_i: \; i \in I_{\Omega}\}\subset \Omega$ and $\{ \boldsymbol{a}_i: \hspace{2pt} i\in I_{\partial \Omega} \}\subset \partial \Omega$  such that $\boldsymbol{a}^n_i\to \boldsymbol{a}_i$ for every  $i \in I_\Omega \cup I_{\partial \Omega}$ and also   $C_{\boldsymbol{y}}=\{\boldsymbol{a}_i:\hspace{2pt} i \in I_\Omega  \}$.  
	
	We have \begin{equation*}
		H \subset \imm(\boldsymbol{y},\Omega)\subset \R^N \setminus \left(  \bigcup_{i \in  I_\Omega  } \imt(\boldsymbol{y},\boldsymbol{a}_i) \right)  =\bigcap_{i \in  I_\Omega } \left( \R^N \setminus \imt(\boldsymbol{y},\boldsymbol{a}_i) \right).
	\end{equation*} 
	Let $i \in  I_\Omega $. Applying   Lemma \ref{lem:conv-boundaries} and Lemma \ref{lem:abundance} together with a diagonal argument,    we select a  decreasing  sequence $(U_l)_l\subset \mathcal{U}_{\boldsymbol{y}} \cap \bigcap_{n \in \N} \mathcal{U}_{\boldsymbol{y}_n}$  with $\boldsymbol{a}_i \in U_l$  and $U_{l+1}\subset \subset U_l$   for every $ l  \in \N$ such that $\mathrm{diam}\,U_l \to 0$, as $l \to \infty$,  and $\boldsymbol{y}^*_n \to \boldsymbol{y}^*$  uniformly on $\partial U_l$, as $n\to \infty$, for all $l\in \N$.  In view of Remark \ref{rem:topim-point}(c), we have
	\begin{equation*}
		\imt(\boldsymbol{y},\boldsymbol{a}_i)=\bigcap_{l \in \N} \closure{\imt(\boldsymbol{y},U_l)}.
	\end{equation*}
Now, since $H \subset \R^N \setminus \imt(\boldsymbol{y},\boldsymbol{a}_i)$  and $H$ is compact,   we necessarily have  $H \subset \R^N \setminus \closure{\imt(\boldsymbol{y},U_k)}$ for some $k\in \N$. By claim (ii) of Lemma \ref{lem:top-deg-conv} and   the equality in \eqref{eqn:img-imt},    this yields $ H \subset \R^N \setminus \closure{\imt(\boldsymbol{y}_n,U_k)}$ for $n \gg 1$.   As $\boldsymbol{a}^n_i \to \boldsymbol{a}_i$, we have $\boldsymbol{a}^n_i\in U_k$ for $n \gg 1$. Therefore, $\imt(\boldsymbol{y}_n,\boldsymbol{a}^n_i) \subset \closure{\imt(\boldsymbol{y}_n,U_k)}$, so that
	  $H \subset \R^N \setminus \imt(\boldsymbol{y}_n,\boldsymbol{a}^n_i)$ for $n \gg 1$.
	  
	   We repeat this argument for every $i \in  I_\Omega $ and we find that 
	 \begin{equation}
	 	\label{eqn:HS}
	 	H \subset \bigcap_{i \in  I_\Omega } (\R^N \setminus \imt(\boldsymbol{y}_n,\boldsymbol{a}^n_i))=\R^N \setminus  \left( \bigcup_{i \in  I_\Omega } \imt(\boldsymbol{y}_n,\boldsymbol{a}^n_i) \right)
	 \end{equation}
 	for  a not relabeled subsequence.  
	
	Let $i \in  I_{\partial \Omega} $. Exploiting Lemma \ref{lem:abundance}, Lemma \ref{lem:conv-boundaries}, and Remark \ref{rem:top-im}(d), we find $U \in \mathcal{U}_{\boldsymbol{y}} \cap \bigcap_{n \in \N}\mathcal{U}_{\boldsymbol{y}_n}$   such that $ H \subset \imt(\boldsymbol{y},U)$ and, up to subsequences, $\boldsymbol{y}_n^* \to \boldsymbol{y}^*$ uniformly on $\partial U$. Thus, Lemma \ref{lem:top-deg-conv}(i) yields  $  H \subset \imt(\boldsymbol{y}_n,U)$ for $n \gg 1$.  
	As $\boldsymbol{a}^n_i \to \boldsymbol{a}_i$ and $\boldsymbol{a}_i \in \partial \Omega$, we have $\boldsymbol{a}^n_i \notin \closure{U}$ for $n \gg 1$. For any such $n \in \N$, let $W^n_i \in  \mathcal{U}_{\boldsymbol{y}_n}  $   with $\boldsymbol{a}^n_i \in W^n_i$ be such that $\closure{U} \cap  \overline{W}^{\,n}_i  =\emptyset$. By Proposition \ref{prop:top-im-inv}(ii), $\imt(\boldsymbol{y}_n,U)\cap \imt(\boldsymbol{y}_n,W^n_i)=\emptyset$. Thus, 
	\begin{equation*}
 		H	\subset  {\imt(\boldsymbol{y}_n,U)} \subset  {\R^N \setminus \imt(\boldsymbol{y}_n,W^n_i)}.
	\end{equation*}
Therefore, $ H  \subset \R^N \setminus \imt(\boldsymbol{y}_n,\boldsymbol{a}^n_i)$.  Repeating the argument for every $i \in  I_{\partial \Omega}$, we find that
	\begin{equation}
		\label{eqn:HT}
		 H \subset \bigcap_{i \in  I_{\partial \Omega} } \left( \R^N \setminus \imt(\boldsymbol{y}_n,\boldsymbol{a}^n_i) \right)=\R^N \setminus \left( \bigcup_{i \in  I_{\partial \Omega}} \imt(\boldsymbol{y}_n,\boldsymbol{a}^n_i) \right)
	\end{equation}
	for $n \gg1 $. The combination of \eqref{eqn:HS}--\eqref{eqn:HT} yields $H \subset \R^N \setminus \imc(\boldsymbol{y}_n,\Omega)$ for $n \gg 1$.  This along with the fact that  $H  \subset \imt(\boldsymbol{y}_n,\Omega)$  concludes the proof.
\end{proof}

We exemplify below the phenomena of cavities closing along converging sequences of deformations. In particular, we demonstrate that the claim in Proposition \ref{prop:conv-matim} fails whenever this behavior occurs. We discuss two examples. 
In the first one, a cavity shrinks to a point while, in the second one, it squashes onto a set of codimension one. 

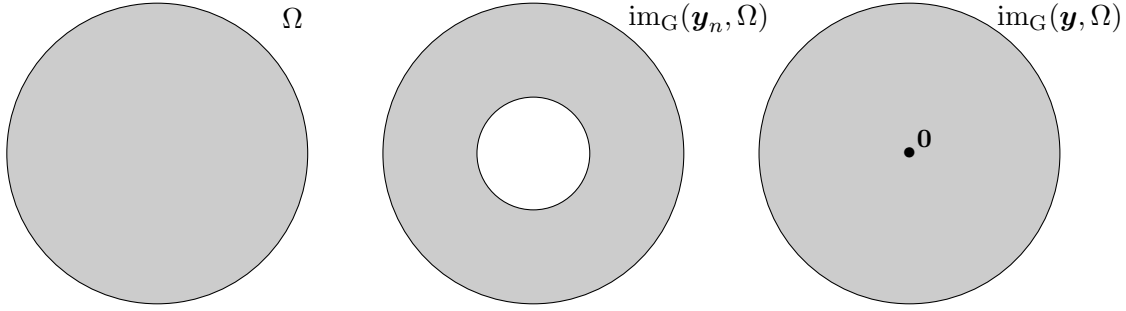
\begin{figure}
	\begin{tikzpicture}[scale=1]
		\draw[fill=black!20] (0,0) circle (2cm);
		\node at (1.8,1.8) {\large $\Omega$};
		\draw[fill=black!20] (5,0) circle (2cm);
		\draw[fill=white] (5,0) circle (.75cm);-n
		\node at (7.2,1.8) {\large $\img(\boldsymbol{y}_n,\Omega)$};
		\draw[fill=black!20] (10,0) circle (2cm);	
		\node at (12,1.8) {\large $\img(\boldsymbol{y},\Omega)$};	
		\node at (10,0) {\textbullet};
		\node at (10.2,.2) {$\boldsymbol{0}$};
	\end{tikzpicture}
	\caption{Deformations in Example \ref{ex:cavity-point}.}
	\label{fig:cavity-point}
\end{figure}

\begin{example}[Cavities shrinking to a point] \label{ex:cavity-point}
Let $N=2$. We employ the notation in \eqref{eqn:q-norm}--\eqref{eqn:B-S}. Set $\Omega\coloneqq B$ and consider the sequence of maps $(\boldsymbol{y}{}_n)_n$ where $\boldsymbol{y}_n\colon \Omega \to \R^2$ is defined as
\begin{equation*}
	\boldsymbol{y}_n(\boldsymbol{x})\coloneqq   \big( (1-2^{-n})|\boldsymbol{x}|+2^{-n} \big ) \frac{\boldsymbol{x}}{ |\boldsymbol{x}|} \quad \text{for all $\boldsymbol{x}\in B\setminus \{ \boldsymbol{0}  \}$.}
\end{equation*} 
Note that $(\boldsymbol{y}_n)_n  \subset \mathcal{Y}_p^{\rm cav}(\Omega)$ by  Lemma~\ref{lem:radial-degree}.  
Each map $\boldsymbol{y}_n$ maps $\Omega$ into itself while opening a spherical cavity of radius $2^{-n}$ at the origin (Figure \ref{fig:cavity-point}). Precisely, we have
\begin{equation*}
		\imm(\boldsymbol{y}_n,\Omega)=\img(\boldsymbol{y}_n,\Omega)=A(2^{-n},1), \quad \imt(\boldsymbol{y}_n,\Omega)=B, \quad \mathcal{S}(\boldsymbol{y}_n)=\mathscr{H}^1(S(2^{-n}))=\pi 2^{1-n}.
\end{equation*}
In particular,  \eqref{eqn:S-bdd} holds. Setting $\boldsymbol{y}\coloneqq \boldsymbol{id}\restr{\Omega}$, one checks \eqref{eqn:weak}--\eqref{eqn:weak-det}. 
Clearly $\imt(\boldsymbol{y},\boldsymbol{0})=\{\boldsymbol{0}\}$, so that $C_{\boldsymbol{y}}=\emptyset$ and  $\imm(\boldsymbol{y},\Omega)=\img(\boldsymbol{y},\Omega)=\imt(\boldsymbol{y},\Omega)=B$. Now,
for all $0<r<1$, we have $ \overline{B}(r)  \subset  \imm(\boldsymbol{y},\Omega)$, but  $ \overline{B}(r)  $ is not contained in $\imm(\boldsymbol{y}_n,\Omega)$ for any $n\in \N$. Therefore, the claim in Proposition \ref{prop:conv-matim} fails for this sequence, while material images converge  in $L^1(\R^N)$,  see Remark~\ref{nice new remark}(a).  
\end{example}  

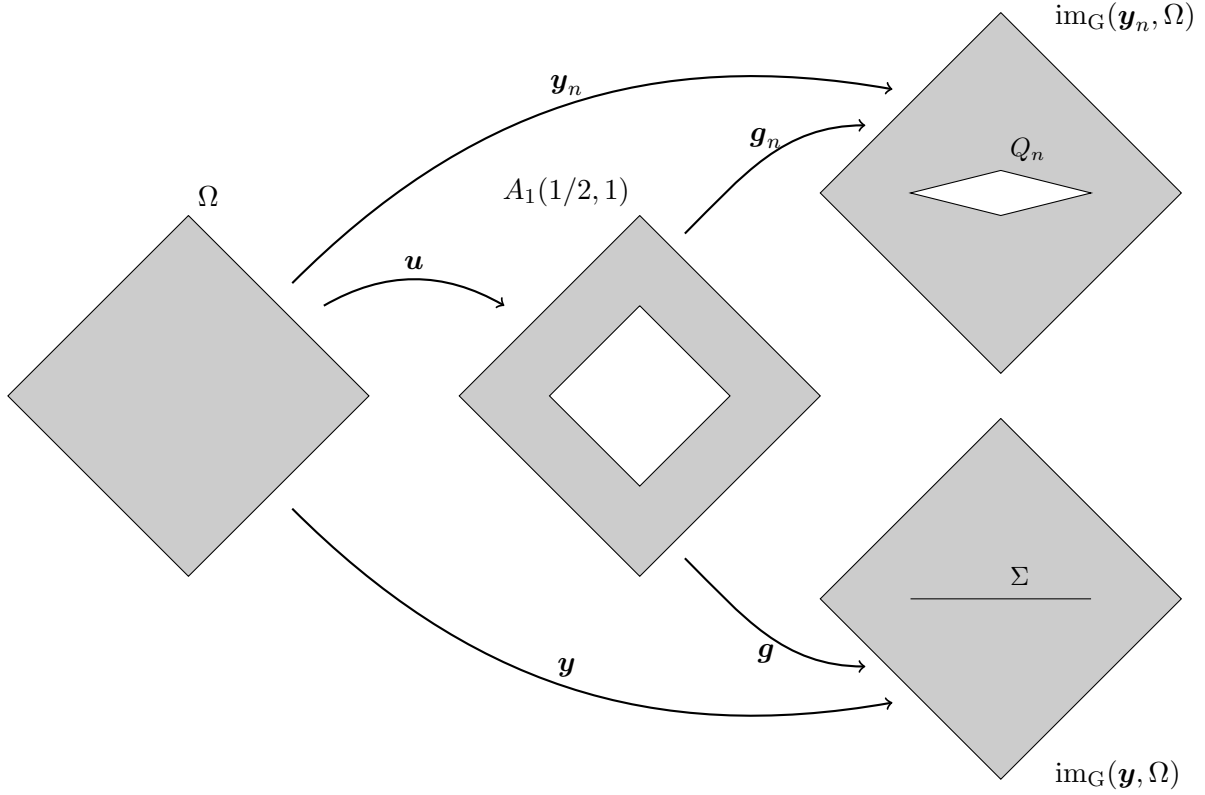
\begin{figure}
	\begin{tikzpicture}[scale=1.2]
		\path[draw=black,fill=black!20] (-5,2) -- (-3,0) -- (-5,-2) -- (-7,0) -- (-5,2) node[above right] {\large $\Omega$};
		\path[draw=black,fill=black!20] (0,2) -- (2,0) -- (0,-2) -- (-2,0) -- (0,2) node[above left] {\large $A_1(1/2,1)$};
		\path[draw=black,fill=white] (0,1) -- (1,0) -- (0,-1) -- (-1,0) -- (0,1);
		\path[draw=black,fill=black!20] (4,4.25) -- (6,2.25) -- (4,.25) -- (2,2.25) -- (4,4.25); 
		\path[draw=black,fill=white] (4,2.5) -- (5,2.25) -- (4,2) -- (3,2.25) -- (4,2.5);
		\path[draw=black,fill=black!20] (4,-.25) -- (6,-2.25) -- (4,-4.25) -- (2,-2.25) -- (4,-.25); 
		\draw (3,-2.25) -- (5,-2.25);
		\node[right] at (4.5,4.2) {\large $\img(\boldsymbol{y}_n,\Omega)$};
		\node[right] at (4.5,-4.2) {\large $\img(\boldsymbol{y},\Omega)$};
		\draw[->, thick] (-3.5,1)  to [out=30,in=150] node[above,midway] {\large $\boldsymbol{u}$} (-1.5,1);
		\draw[->, thick] (.5,1.8)  to [out=45,in=180] node[above,midway] {\large $\boldsymbol{g}_n$} (2.5,3);
		\draw[->, thick] (.5,-1.8)  to [out=315,in=180] node[below,midway] {\large $\boldsymbol{g}$} (2.5,-3);
		\node[above right] at (4,2.5) {$Q_n$};
		\node[above right] at (4,-2.2) {$\Sigma$};
		\draw[->, thick] (-3.85,1.25)  to [out=45,in=170] node[above,midway] {\large $\boldsymbol{y}_n$} (2.8,3.4);
		\draw[->, thick] (-3.85,-1.25)  to [out=-45,in=-170] node[above,midway] {\large $\boldsymbol{y}$} (2.8,-3.4);
	\end{tikzpicture}
	\caption{Deformations in  Example \ref{ex:cavity-segment}.}
	\label{fig:cavity-segment}
\end{figure}

\begin{example}[Cavities squashing onto a segment] \label{ex:cavity-segment}
	Let $N=2$.    We employ the notation in \eqref{eqn:q-norm}--\eqref{eqn:q-ball-annulus}.
	The domain $\Omega \coloneqq B_1$ is a square centered at the origin with side length $\sqrt{2}$  rotated by an angle of $\pi/4$. Consider the map $\boldsymbol{u}\colon \Omega \to \R^2$ defined as
	\begin{equation*}
		\boldsymbol{u}(\boldsymbol{x})\coloneqq  \frac{1}{2}(|\boldsymbol{x}|_1+1)\frac{\boldsymbol{x}}{|\boldsymbol{x}|_1} \quad \text{for all $\boldsymbol{x}\in B_1 \setminus \{ \boldsymbol{0}  \}$.}
	\end{equation*}
	The map $\boldsymbol{u}$ opens a square cavity at the origin with side length $\sqrt{2}/2$ rotated by an angle of $\pi/4$ while keeping the boundary of $\Omega$ fixed (Figure \ref{fig:cavity-segment}). Precisely, we have $\boldsymbol{u}(\Omega \setminus \{\boldsymbol{0}\})=A_1(1/2,1)$. 
	
	 For every $n\in \N$, we consider the map $\boldsymbol{g}_n \colon\R^2 \to \R^2$ defined as 
	\begin{equation*}
		\boldsymbol{g}_n(\boldsymbol{x})\coloneqq \begin{cases}
			(x_1, \sgn(x_2)g_n(|x_1|,|x_2|))^\top & \text{if $\boldsymbol{x}\in \closure{B}_1 \cap \big ([-1/2,1/2]\times \R \big)$}, \\
			\boldsymbol{x} & \text{otherwise},
		\end{cases}
	\end{equation*}
	 where,  setting $D\coloneqq \closure{B}_1 \cap \left( [0,1/2] \times [0,1] \right)$, the function $g_n\colon D \to [0,1]$ in the previous equation is given by
	\begin{equation*}
		g_n(s,t)\coloneqq \begin{cases}
			2^{-n}(1/2-s)^{-1}t & \text{if $t\in [0,1/2-s]$,}\\
			2\left( (1-2^{-n})-s \right)(s+t-1)+1-s & \text{if $t\in (1/2-s,1-s]$.} 
		\end{cases}
	\end{equation*} 
	The map $\boldsymbol{g}_n$ is piecewise smooth and squeezes the cavity $\closure{B}_1(1/2)$ into the (closed) rhombus $Q_n$ of  vertices  $(\pm 1/2,0)^\top$  and $(0,\pm 2^{-n})^\top$ while keeping the outer boundary $S_1$ of $B_1$ fixed.  We consider the deformation $\boldsymbol{y}_n \coloneqq \boldsymbol{g}_n \circ \boldsymbol{u} \colon \Omega \to \R^2$. After verifying all the assumptions, we see that $\boldsymbol{y}_n \in \mathcal{Y}_p^{\rm cav}(\Omega)$ by Proposition~\ref{prop:radial-superposition}. In particular,
	\begin{equation*}
		\imm(\boldsymbol{y}_n,\Omega)=\img(\boldsymbol{y}_n,\Omega)=B_1 \setminus Q_n, \quad \imt(\boldsymbol{y}_n,\Omega)=B_1, \quad \mathcal{S}(\boldsymbol{y}_n)=\mathscr{H}^1(\partial Q_n)=4\sqrt{1/4  +2^{-2n}}.
	\end{equation*}
	Thus, \eqref{eqn:S-bdd} holds.
	As $n\to \infty$,   we have that $g_n \to g$ pointwise in $D$, where $g\colon D \to [0,1]$ is defined as
	\begin{equation*}
		g(s,t)\coloneqq \begin{cases}
			0 & \text{if $t\in [0,1/2-s]$,} \\ 2(1-s)(s+t-1/2) & \text{if $t\in (1/2-s,1-s]$.}
		\end{cases}
	\end{equation*}  
	The map $\boldsymbol{g}\colon \R^2 \to \R^2$  given by
		\begin{equation*}
		\boldsymbol{g}(\boldsymbol{x})\coloneqq \begin{cases}
			(x_1, \sgn(x_2)g(|x_1|,|x_2|))^\top & \text{if $\boldsymbol{x}\in \closure{B}_1 \cap \big( [-1/2,1/2]\times \R \big)$}, \\
			\boldsymbol{x} & \text{otherwise},
		\end{cases}
	\end{equation*}
	squashes the cavity $\closure{B}_1(1/2)$ into the segment  $\Sigma \coloneqq (-1/2,1/2)\times \{0\}$. 
	Now,  define $\boldsymbol{y}\coloneqq \boldsymbol{g} \circ \boldsymbol{u} \colon \Omega \to \R^2$. As before,  $\boldsymbol{y} \in \mathcal{Y}_p^{\rm cav}(\Omega)$ by Proposition \ref{prop:radial-superposition}. For this deformation, we have
	\begin{equation*}
		\img(\boldsymbol{y},\Omega)=B_1 \setminus \Sigma, \quad \imm(\boldsymbol{y},\Omega)=\imt(\boldsymbol{y},\Omega)=B_1.  
	\end{equation*}
	In particular, note that $\imt(\boldsymbol{y},\boldsymbol{0})=\Sigma$, so that $C_{\boldsymbol{y}}=\emptyset$. 
	One checks that \eqref{eqn:weak}--\eqref{eqn:weak-det} hold.  We observe the convergence of material images  in  $L^1(\R^2)$,  but the claim in Proposition \ref{prop:conv-matim} does not hold. Indeed, for all $0<r<1$, we have $ \overline{B}_1(r)  \subset B_1=\imm(\boldsymbol{y},\Omega)$ while $\overline{B}_1(r)$ is not contained in $\imm(\boldsymbol{y}_n,\Omega)$ for any $n\in \N$.   
	
	In passing, we observe that $\mathcal{S}(\boldsymbol{y})=0$ by  Theorem~\ref{thm:INV-top-im}(i).  
\end{example}

\section{Coercivity and lower semicontinuity of Eulerian-Lagrangian energies}
\label{sec:EL}
In this section, we investigate the coercivity and lower semicontinuity of Eulerian-Lagrangian energies with the aim of applying the direct method of the calculus of variations. We analyze different settings corresponding to different regularity assumptions for both deformations and Eulerian maps,  distinguishing again the cases of   Sobolev and  $SBV$-deformations.  

\subsection{Eulerian-Lagrangian energies}
We first describe the general setting. As before, let $\Omega \subset \R^N$ be a bounded Lipschitz domain and  $p>N-1$. Also, let $Z\subset \R^M$ be a measurable set (possibly unbounded). 
Recall  \eqref{eqn:Y}  and Definition~\ref{def:geometric-image}. For each $\boldsymbol{y}\in \mathcal{Y}(\Omega)$, we consider Eulerian  maps $\boldsymbol{v}\in L^2(\img(\boldsymbol{y},\Omega); Z)$ which are almost everywhere approximately differentiable. We will deal with functionals of the form 
\begin{equation}
	\label{eqn:F-EL}
	\begin{split}
		\mathcal{F}(\boldsymbol{y},\boldsymbol{v})\coloneqq& \int_{\Omega} W(\nabla \boldsymbol{y},\boldsymbol{v}\circ \boldsymbol{y})\,\d\boldsymbol{x}+ \lambda_1 \mathcal{S}(\boldsymbol{y})+\lambda_2 \per \left( \img(\boldsymbol{y},\Omega) \right)+\haus(J_{\boldsymbol{y}})\\
		 &+ \int_{\img(\boldsymbol{y},\Omega)} |\nabla \boldsymbol{v}|^2\,\d\boldsymbol{\xi}+\lambda_3 \haus(J_{\boldsymbol{v}}\cap {O}_{\boldsymbol{y}}),
	\end{split}
\end{equation}
where  $W\colon \rnn_+ \times  Z   \to [0,+\infty]$ is an elastic energy density,  $J_{\boldsymbol{y}}$  and $J_{\boldsymbol{v}}$ are  the jump sets  of $\boldsymbol{y}$ and $\boldsymbol{v}$, respectively, according to  Definition~\ref{def:jump},     $\mathcal{S}(\boldsymbol{y})$ is as in Definition \ref{def:surface-energy}, $\lambda_1,\lambda_2,\lambda_3\in [0,+\infty)$, and  ${O}_{\boldsymbol{y}}\subset \R^N$ denotes an open set which is uniquely determined by the equivalence class of $\boldsymbol{y}$.  Specifically, we will consider two settings: 
\begin{align}
	\label{eqn:F-EL2}
		\lambda_1,\lambda_3>0, \hspace{3pt} \lambda_2=0, \qquad O_{\boldsymbol{y}}=\imt(\boldsymbol{y},\Omega) \qquad &\text{for $\boldsymbol{y}\in\mathcal{Y}_p^{\rm cav}(\Omega)$,}\\
	\label{eqn:F-EL3}	
		\lambda_1,\lambda_2,\lambda_3>0, \qquad O_{\boldsymbol{y}}=\R^N \hspace{28pt}\qquad &\text{for $\boldsymbol{y}\in \mathcal{Y}_p^{\rm frac}(\Omega)$.}
\end{align}
Here, we employ the notation in \eqref{eqn:Y-cav} and \eqref{eqn:Y-frac}  as well as Definition \ref{def:topological-image}.

We briefly compare the two settings in \eqref{eqn:F-EL2}--\eqref{eqn:F-EL3}. The  one in \eqref{eqn:F-EL2} is quite natural: the surface term penalizes only the creation of new surface by  deformations and the choice of the set $O_{\boldsymbol{y}}$ ensures that only jump points of $\boldsymbol{v}$ inside of the geometric image are accounted for,   see Remark  \ref{rem:cav}(b) below. Instead, the  setting in \eqref{eqn:F-EL3} is not completely satisfactory. In this case, apart from the creation of new surface by deformations, the energy {also} penalizes the stretching of the boundary. Bearing in mind \eqref{eqn:intro-perimeter}, we can formally write
\begin{equation*}
	\lambda_1 \mathcal{S}(\boldsymbol{y})+\lambda_2 \per \left( \img(\boldsymbol{y},\Omega) \right) =(\lambda_1+\lambda_2 ) \mathcal{S}(\boldsymbol{y})+\lambda_2 \haus(\boldsymbol{y}(\partial \Omega)).
\end{equation*}
In this way, we can think of \eqref{eqn:F-EL3} as a regularized setting where also the measure of the image of the boundary is controlled. 
 The control of both the surface energy and the perimeter of the geometric image is essential for carrying out the analysis.  In particular, the condition $\lambda_1>0$    in \eqref{eqn:F-EL3}  is necessary in order to recover the weak continuity of Jacobian determinants since perimeter bounds do not ensure it,  see Remark~\ref{rem:S-lsc}(c). 
Eventually, the energy in \eqref{eqn:F-EL3} {also} accounts for jumps of $\boldsymbol{v}$ outside of the geometric image, like those arising from self{-}contact at the boundary, see Remark  \ref{rem:cav}(b).

\begin{remark}
The value of the functional in \eqref{eqn:F-EL} does not depend on the representatives of $\boldsymbol{y}$ and $\boldsymbol{v}$, but only on their equivalence classes. More explicitly, let $(\boldsymbol{y}_1,\boldsymbol{v}_1)$ and $(\boldsymbol{y}_2,\boldsymbol{v}_2)$ be two admissible states with $\boldsymbol{y}_1 \cong \boldsymbol{y}_2$ in $\Omega$ and $\boldsymbol{v}_1\cong {\boldsymbol{v}}_2$ in $\img(\boldsymbol{y}_1,\Omega) \cap \img({\boldsymbol{y}}_2,\Omega)$.
Then, $\mathcal{F}(\boldsymbol{y}_1,\boldsymbol{v}_1)=\mathcal{F}({\boldsymbol{y}}_2,{\boldsymbol{v}}_2)$. Indeed, $\nabla \boldsymbol{y}_1\cong \nabla \boldsymbol{y}_2$ and $\boldsymbol{v}_1 \circ \boldsymbol{y}_1\cong \boldsymbol{v}_2 \circ \boldsymbol{y}_2$ thanks to Remark \ref{rem:composition-measurable}(a), while $J_{\boldsymbol{y}_1}\simeq J_{\boldsymbol{y}_2}$ and $J_{\boldsymbol{v}_1}\simeq J_{\boldsymbol{v}_2}$  by    Remark \ref{rem:jump}(c),  and  ${O}_{\boldsymbol{y}_1}={O}_{\boldsymbol{y}_2}$  by Remark \ref{rem:top-im}(b) and  \eqref{eqn:F-EL2}.    The identity $\mathcal{S}(\boldsymbol{y}_1)=\mathcal{S}(\boldsymbol{y}_2)$   follows since   $  \boldsymbol{y}_1\cong   \boldsymbol{y}_2$ and  $\nabla \boldsymbol{y}_1\cong \nabla \boldsymbol{y}_2$, see Definition \ref{def:surface-energy}.    Eventually, the independence of the fourth  and fifth  term in \eqref{eqn:F-EL}  is ensured  as    $\nabla \boldsymbol{v}_1 \cong \nabla \boldsymbol{v}_2$ in $\img(\boldsymbol{y}_1,\Omega) \cap \img({\boldsymbol{y}}_2,\Omega)$  and   $\img(\boldsymbol{y}_1,\Omega)\cong \img(\boldsymbol{y}_2,\Omega)$  by  Remark \ref{rem:geom-dom-im}(c). 
\end{remark}

On the  density $W$  in \eqref{eqn:F-EL}  we make the following  standard  assumptions:
\begin{enumerate}[(W1)]
	\item \textbf{Continuity:} The function $W\colon \rnn_+ \times  Z  \to [0,+\infty]$ is continuous.
	\item \textbf{Coercivity:} There exist a constant $C>0$ and a Borel function $\gamma \colon (0,+\infty) \to [0,+\infty]$  satisfying
	\begin{equation}
		\label{eqn:growth-gamma}
		\lim_{h \to 0^+} \gamma(h)=+\infty, \qquad \lim_{h \to +\infty} \frac{\gamma(h)}{h}=+\infty,
	\end{equation}
	such that
	\begin{equation}
		\label{eqn:growth}
	    W(\boldsymbol{F},\boldsymbol{z})\geq C |\boldsymbol{F}|^p+\gamma(\det \boldsymbol{F}) \quad \quad  \text{for all $\boldsymbol{F}\in \rnn_+$ and  $\boldsymbol{z}\in  Z $}. 
	\end{equation}
	\item \textbf{Polyconvexity:} There exists a continuous function $\widehat{W}\colon \displaystyle \prod_{r=1}^{N-1} \R^{\binom{N}{r}\times \binom{N}{r}} \times (0,+\infty)\times  Z    \to [0,+\infty]$  such that 
	\begin{equation}
		\label{eqn:polyconvex1}
 \text{$(\boldsymbol{F}_1,\dots,\boldsymbol{F}_{N-1},{F}_N)\mapsto \widehat{W}(\boldsymbol{F}_1,\dots,\boldsymbol{F}_{N-1},{F}_N,\boldsymbol{z})$ is convex \quad  for all   $\boldsymbol{z}\in Z $}
	\end{equation}
	and \
	\begin{equation}
		\label{eqn:polyconvex2}
  W(\boldsymbol{F},\boldsymbol{z})=\widehat{W}(\adj_1\boldsymbol{F},\dots,\adj_{N-1}\boldsymbol{F},\adj_N\boldsymbol{F},\boldsymbol{z}) \quad \quad  \text{for all $\boldsymbol{F}\in \rnn_+$ and $\boldsymbol{z}\in  Z $}.	
	\end{equation}
\end{enumerate}

The first  assumption  ensures the measurability of  $\boldsymbol{x}\mapsto W(\nabla \boldsymbol{y}(\boldsymbol{x}),\boldsymbol{v}(\boldsymbol{y}(\boldsymbol{x})))$. Given the first limit in  \eqref{eqn:growth-gamma},  the extension of $W$ to $\rnn \times  Z $ defined by setting $W(\boldsymbol{F},\boldsymbol{z})=+\infty$ whenever $\det \boldsymbol{F}\leq 0$ is continuous. Also, if  we analogously  extend  $\widehat{W}$   by  setting $\widehat{W}(\boldsymbol{F}_1,\dots,\boldsymbol{F}_{N-1},F_N,  \boldsymbol{z}  )=+\infty$ whenever $F_N\leq 0$, then  \eqref{eqn:polyconvex1}--\eqref{eqn:polyconvex2} are still satisfied by the extension.
We will see examples of densities $W$ satisfying these assumptions in Section~\ref{sec:appl}.

The next two subsections are devoted  to  the study of Eulerian-Lagrangian energies in the case of deformations that are Sobolev maps and special maps of bounded variation, respectively,  see \eqref{eqn:F-EL2} and \eqref{eqn:F-EL3}.

\subsection{Eulerian-Lagrangian energies for Sobolev deformations}

In this subsection,  admissible deformations are Sobolev maps. Recall \eqref{eqn:Y-cav}. For $\boldsymbol{y}\in \mathcal{Y}^{\rm cav}_p(\Omega)$ and $\boldsymbol{v}\in L^2(\img(\boldsymbol{y},\Omega); Z )$  almost everywhere approximately differentiable,  the functional in \eqref{eqn:F-EL}--\eqref{eqn:F-EL2} takes the form
\begin{equation}
	\label{eqn:F-EL-cav}
	{\mathcal{F}^{\rm \, cav}_1}(\boldsymbol{y},\boldsymbol{v})\coloneqq \int_{\Omega} W(D \boldsymbol{y},\boldsymbol{v}\circ \boldsymbol{y})\,\d\boldsymbol{x}+ \lambda_1\mathcal{S}(\boldsymbol{y}) + \int_{\img(\boldsymbol{y},\Omega)} |\nabla \boldsymbol{v}|^2\,\d\boldsymbol{\xi}+ \lambda_3  \haus\big(J_{\boldsymbol{v}} \cap \imt(\boldsymbol{y},\Omega)\big). 
\end{equation}

 Our first  result addresses the coercivity and the lower semicontinuity   of this functional  under the assumption that the extension of $\boldsymbol{v}$ to $\imt(\boldsymbol{y},\Omega)$ by zero belongs to the space $GSBV^2(\imt(\boldsymbol{y},\Omega);\R^M)$,  see \eqref{eqn:SBVcapLinftyXXX}. We refer to   Remark~\ref{rem:cav}(a) for a few comments about this assumption.

\begin{theorem}[Coercivity and lower semicontinuity: Sobolev  deformations] 
	\label{thm:cav}
Let  $ \mathcal{F}^{\rm \, cav}_1 $ be defined as in 	\eqref{eqn:F-EL-cav} with  $W$ satisfying  {\rm (W1)}--{\rm (W3)}  and $\lambda_1\geq \lambda_3>0$.   
Let $(\boldsymbol{y}_n)_n\subset \mathcal{Y}_p^{\rm cav}(\Omega)$ and let $(\boldsymbol{v}_n)_n$ be a sequence of maps $\boldsymbol{v}_n \in L^2(\img(\boldsymbol{y}_n,\Omega);  Z)$ with  $\widetilde{\boldsymbol{v}}_n \in GSBV^2(\imt(\boldsymbol{y}_n,\Omega);\R^M)$, where $\widetilde{\boldsymbol{v}}_n$ denotes the extension of $\boldsymbol{v}_n$ to $\imt(\boldsymbol{y}_n,\Omega)$ by zero.  
Suppose that 
\begin{equation}
	\label{eqn:cav-bound}
	\sup_{n \in \N}  \left \{ \|\boldsymbol{y}_n\|_{L^p(\Omega;\R^N)}+ \|\boldsymbol{v}_n\|_{L^2(\img(\boldsymbol{y}_n,\Omega);\R^M)}+ \mathcal{F}^{\rm \, cav}_1 (\boldsymbol{y}_n,\boldsymbol{v}_n) \right \}<+\infty.
\end{equation}
Eventually, assume that  the sequence $(\boldsymbol{v}_n \circ \boldsymbol{y}{}_n)_n$ is equi-integrable.  Then, there   exist $\boldsymbol{y}\in \mathcal{Y}_{p}^{\rm cav}(\Omega)$ and $\boldsymbol{v}\in L^2( \img(\boldsymbol{y},\Omega);  \R^M)$ satisfying $\boldsymbol{v}\circ \boldsymbol{y}\in L^1(\Omega;\R^M)$ and  $\widetilde{\boldsymbol{v}}\in GSBV^2(\imt(\boldsymbol{y},\Omega);\R^M)$, where $\widetilde{\boldsymbol{v}}$ denotes the extension  of $\boldsymbol{v}$ to $\imt(\boldsymbol{y},\Omega)$ by zero,  such that, up to subsequences, we have: 
\begin{align}
\label{eqn:cav-compactness1}
\text{$\boldsymbol{y}_n \wk \boldsymbol{y}$ in $W^{1,p}(\Omega;\R^N)$}, \qquad &\text{$\det D\boldsymbol{y}_n \wk \det D\boldsymbol{y}$ in $L^1(\Omega)$},\\
\label{eqn:cav-compactness2}
\text{$\chi_{\img(\boldsymbol{y}_n,\Omega)}\boldsymbol{v}_n \to \chi_{\img(\boldsymbol{y},\Omega)}\boldsymbol{v}$ }  &\text{in $L^1(\R^N;\R^M)$,}\\ 
\label{eqn:cav-compactness3}
\text{$\chi_{\img(\boldsymbol{y}_n,\Omega)}\boldsymbol{v}_n \wk \chi_{\img(\boldsymbol{y},\Omega)}\boldsymbol{v}$ }  &\text{in $L^2(\R^N;\R^M)$,}  \\
\label{eqn:cav-compactness4}
\text{$\chi_{\img(\boldsymbol{y}_n,\Omega)}\nabla\boldsymbol{v}_n \wk \chi_{\img(\boldsymbol{y},\Omega)}\nabla\boldsymbol{v}$ } &\text{in $L^2(\R^N;\R^{M \times N})$,}\\
\label{eqn:cav-compactness5}
\text{$\boldsymbol{v}_n \circ \boldsymbol{y}_n \to \boldsymbol{v}\circ \boldsymbol{y}$ } &\text{in $L^1(\Omega;\R^M)$.}
\end{align}
 Moreover, if  $\boldsymbol{v}(\boldsymbol{\xi}) \in Z $ for almost every $\boldsymbol{\xi}\in \img(\boldsymbol{y},\Omega)$,  then  we have
\begin{equation}
\label{eqn:cav-lsc}
\mathcal{F}^{\rm \, cav}_1 (\boldsymbol{y},\boldsymbol{v})\leq \liminf_{n \to \infty} \mathcal{F}^{\rm \, cav}_1 (\boldsymbol{y}_n,\boldsymbol{v}_n). 
\end{equation}
\end{theorem}

\begin{remark} \label{rem:cav}
	\begin{enumerate}[(a)]
		\item The assumption $\widetilde{\boldsymbol{v}}\in GSBV^2(\imt(\boldsymbol{y},\Omega);\R^M)$ is an additional regularity requirement, but the value of the functional in \eqref{eqn:F-EL-cav} is completely determined by  $\boldsymbol{v}$.  Loosely speaking, this assumption means that $\haus(J_{\boldsymbol{v}} )<+\infty$  and  $\boldsymbol{v}$ has finite Dirichlet energy on $\img(\boldsymbol{y},\Omega)$.   
		\item Recall that $J_{\boldsymbol{v}}\subset \img(\boldsymbol{y},\Omega)^{(1)}$ according to  Remark \ref{rem:jump}(b).  In general,   $\img(\boldsymbol{y},\Omega)^{(1)}$  is not contained in $\imt(\boldsymbol{y},\Omega)$,  see Example \ref{ex:jump-not-contained-imt} below.  This  shows  that, compared with $\haus(J_{\boldsymbol{v}})$, the term $\haus(J_{\boldsymbol{v}}\cap \imt(\boldsymbol{y},\Omega))$ in \eqref{eqn:F-EL-cav} is more natural as it does not penalize the jump points of $\boldsymbol{v}$ arising from self-contact at the boundary.   	Clearly, controlling the whole measure of $J_{\boldsymbol{v}}$  would  not affect the compactness result. However, the functional $(\boldsymbol{y},\boldsymbol{v})\mapsto \haus(J_{\boldsymbol{v}})$ is not lower semicontinuous with respect to the convergence in \eqref{eqn:cav-compactness1}--\eqref{eqn:cav-compactness5}, see Example \ref{ex:jump-not-contained-imt} below.  Therefore, the term used in \eqref{eqn:F-EL-cav}  is also advantageous from an analytical viewpoint.   
		\item  The assumption $\lambda_1\geq \lambda_3$ is essential. Indeed,  
		in the proof of Theorem \ref{thm:cav}, we show that
		\begin{align}\label{YYYYY}
			\haus\big(J_{\boldsymbol{v}}  \cap \imt(\boldsymbol{y},\Omega) \big) +\mathcal{S}(\boldsymbol{y}) & \leq \liminf_{n \to \infty}  \lbrace \haus\big(J_{\boldsymbol{v}_n}  \cap \imt(\boldsymbol{y}_n,\Omega) \big)  +\mathcal{S}(\boldsymbol{y}_n)\rbrace,      \\  
			\label{YYYYY2} \mathcal{S}(\boldsymbol{y}) & \leq \liminf_{n \to \infty}   \mathcal{S}(\boldsymbol{y}_n).
		\end{align}
		While the lower semicontinuity of $\mathcal{S}$ follows from Theorem \ref{thm:determinant-surface-energy}, that of  $(\boldsymbol{y},\boldsymbol{v})\mapsto \haus(J_{\boldsymbol{v}}\cap \imt(\boldsymbol{y},\Omega))$  does not generally hold,  and we can only establish a result for the sum.  The typical situation when such lower semicontinuity fails is given by a sequence of deformations opening a cavity that squashes onto a lower-dimensional set the limit,    see  Example \ref{ex:jump-not-lsc} below.    By   multiplying the terms in \eqref{YYYYY}--\eqref{YYYYY2} by $\lambda_3>0$ and $\lambda_1 - \lambda_3 \ge 0$, respectively, the  lower semicontinuity of the two surface terms in \eqref{eqn:F-EL-cav} follows.  
						
		\item The assumption $\boldsymbol{v}(\boldsymbol{\xi})\in  Z $ for almost every $\boldsymbol{\xi}\in \img(\boldsymbol{y},\Omega)$  is only needed to apply Theorem \ref{thm:lscT} and  to  deduce the lower semicontinuity of the elastic energy. As noted in  Remark \ref{rem:lscT}(c), this assumption is superfluous whenever $ Z  $ is closed.
		  
		\item 
		If  $\sup_{n \in \N} \|\boldsymbol{v}_n\|_{L^\infty(\img(\boldsymbol{y}_n,\Omega);\R^M)}<+\infty$, then
		$\widetilde{\boldsymbol{v}}_n \in SBV^2(\imt(\boldsymbol{y}_n,\Omega);\R^M)$ for every $n \in \N$ and also $\widetilde{\boldsymbol{v}}\in SBV^2(\imt(\boldsymbol{y},\Omega);\R^M)$. Additionally, if each $\widetilde{\boldsymbol{v}}_n$ is piecewise-constant, then so is $\widetilde{\boldsymbol{v}}$. The $SBV$-regularity of $\widetilde{\boldsymbol{v}}_n$ comes from \eqref{eqn:SBVcapLinfty}. The one of $\widetilde{\boldsymbol{v}}$ follows analogously by taking  a subsequence that converges almost everywhere  from \eqref{eqn:cav-compactness2} and using the lower semicontinuity of the $L^\infty$-norm. The piecewise-constant case follows thanks to \eqref{eqn:PC}  and \eqref{eqn:cav-compactness4}.  
		\item If $\widetilde{\boldsymbol{v}}_n \in W^{1,2}(\imt(\boldsymbol{y}_n,\Omega);\R^M)$ for every $n \in \N$, then $ \widetilde{\boldsymbol{v}} \in W^{1,2}(\imt(\boldsymbol{y},\Omega);\R^M)$. To check this, consider an open set $V \subset \subset \imt(\boldsymbol{y},\Omega)$. By Proposition \ref{prop:conv-topim}(i), we have $V\subset \subset \imt(\boldsymbol{y}_n,\Omega)$ along a not relabeled subsequence. Combining \eqref{eqn:cav-compactness3}--\eqref{eqn:cav-compactness4},   we deduce that $\widetilde{\boldsymbol{v}}\in W^{1,2}_\loc(\imt(\boldsymbol{y},\Omega);\R^M)$ with $D\widetilde{\boldsymbol{v}}$ given by the extension of $\nabla \boldsymbol{v}$ to $\imt(\boldsymbol{y},\Omega)$ by zero. As $\boldsymbol{v}\in L^2(\img(\boldsymbol{y},\Omega);\R^{M})$ and $\nabla \boldsymbol{v}\in L^2(\img(\boldsymbol{y},\Omega);\R^{M\times N})$, we actually have $\widetilde{\boldsymbol{v}}\in W^{1,2}(\imt(\boldsymbol{y},\Omega);\R^M)$. 	 
	\end{enumerate}
\end{remark}

\begin{proof}[Proof of Theorem \ref{thm:cav}]	
We divide the proof into three steps.	
	
\textbf{Step 1 (Compactness of deformations).}	 
By \eqref{eqn:growth} and \eqref{eqn:cav-bound}, the sequences $(\boldsymbol{y}_n)_n \subset W^{1,p}(\Omega;\R^N)$ and $(\gamma(\det D \boldsymbol{y}_n))_n\subset L^1(\Omega)$ are bounded. Thus,  up to subsequences, we have the first convergence in \eqref{eqn:cav-compactness1} for some $\boldsymbol{y}\in W^{1,p}(\Omega;\R^N)$. Thanks to Lemma \ref{lem:INV-stability}, the map $\boldsymbol{y}$ satisfies condition (INV).
In view of the second condition in \eqref{eqn:growth-gamma},  by  the De la Vallée {Poussin} criterion \cite[Theorem 2.29]{fonseca.leoni}, and the Dunford-Pettis theorem, there exists $h \in L^1(\Omega)$ such that, up to subsequences, $\det D \boldsymbol{y}_n \wk h$ in $L^1(\Omega)$. As each $\boldsymbol{y}_n$ satisfies $\det D \boldsymbol{y}_n>0$ almost everywhere, we have $h\geq 0$ almost everywhere.  To show that $h>0$ almost everywhere, we employ a standard contradiction  argument based on the first condition in \eqref{eqn:growth-gamma}, see, e.g.,  \cite[Theorem 5.1]{mueller.spector}. Suppose that $h=0$ on a measurable set $E\subset \Omega$ with $\leb(E)>0$. Since $\det D \boldsymbol{y}_n>0$ almost everywhere for all $n\in\N$, we have $\det D \boldsymbol{y}_n\to 0$ in $L^1(E)$ which, by  \eqref{eqn:growth-gamma}, yields $\gamma(\det D\boldsymbol{y}_n)\to +\infty$ almost everywhere in $E$ up to subsequences. Thus, \eqref{eqn:growth} and Fatou’s lemma give
\begin{equation*}
	+\infty \leq \liminf_{n\to \infty} \int_E \gamma(\det D \boldsymbol{y}_n)\,\d\boldsymbol{x}\leq \liminf_{n\to \infty} \mathcal{F}^{\rm cav}_1(\boldsymbol{y}_n,\boldsymbol{v}_n),
\end{equation*}
which contradicts \eqref{eqn:cav-bound}. Therefore, $h>0$ almost everywhere.

From the first convergence in  \eqref{eqn:cav-compactness1}, we obtain $\boldsymbol{y}_n\to\boldsymbol{y}$ almost everywhere by  the  Rellich-Kondrachov theorem   and $\adj D \boldsymbol{y}_n \wk \adj D \boldsymbol{y}$ in $L^1(\Omega;\rnn)$ by \cite[Theorem 8.20, Part 4]{dacorogna}. Also, $(\mathcal{S}(\boldsymbol{y}_n))_n$ is bounded by \eqref{eqn:cav-bound}.  
Hence,   by Theorem \ref{thm:determinant-surface-energy}  we deduce the second convergence in  \eqref{eqn:cav-compactness1} as well as the almost everywhere injectivity of $\boldsymbol{y}$   and the finiteness of $\mathcal{S}(\boldsymbol{y})$.     In particular, $\boldsymbol{y}\in \mathcal{Y}_{p}^{\rm cav}(\Omega)$. 
Then, from claim (i) of Proposition~\ref{prop:approx-diff}, we  additionally obtain
\begin{equation}
	\label{eqn:imgeom-conv}
	\text{$\chi_{\img(\boldsymbol{y}_n,\Omega)}\to \chi_{\img(\boldsymbol{y},\Omega)}$ in $L^1(\R^N)$. }
\end{equation}

\textbf{Step 2 (Compactness of Eulerian maps).} From  Lemma \ref{lem:approx-diff-extension} and Theorem \ref{thm:INV-top-im}(iii)   we deduce the identities
\begin{equation}
	\label{eqn:imt=img}
	\chi_{\imt(\boldsymbol{y}_n,\Omega)}\widetilde{\boldsymbol{v}}_n   \cong  \chi_{\img(\boldsymbol{y}_n,\Omega)}\boldsymbol{v}_n, \qquad \chi_{\imt(\boldsymbol{y}_n,\Omega)}\nabla \widetilde{\boldsymbol{v}}_n \cong \chi_{\img(\boldsymbol{y}_n,\Omega)}\nabla\boldsymbol{v}_n \quad \text{for all $n\in \N$.}
\end{equation}
From \eqref{eqn:cav-bound}, we have
\begin{equation}
	\label{eqn:cav-bbound}
	\sup_{n \in \N} \left\{ \int_{\imt(\boldsymbol{y}_n,\Omega)}|\widetilde{\boldsymbol{v}}_n|^2\,\d\boldsymbol{\xi} + \int_{\imt(\boldsymbol{y}_n,\Omega)}|\nabla \widetilde{\boldsymbol{v}}_n|^2\,\d\boldsymbol{\xi}  \right\}<+\infty.
\end{equation}
This entails the existence of two maps $ {\boldsymbol{w}} \in L^2(\R^N;\R^M)$ and $ {\boldsymbol{W}}\in L^2(\R^N;\R^{M \times N})$ for which
\begin{equation}
	\label{eqn:cav-w}
	\text{$\chi_{\imt(\boldsymbol{y}_n,\Omega)}\widetilde{\boldsymbol{v}}_n \wk {\boldsymbol{w}} $ in $L^2(\R^N;\R^M)$}, \qquad \text{$\chi_{\imt(\boldsymbol{y}_n,\Omega)}\nabla \widetilde{\boldsymbol{v}}_n \wk {\boldsymbol{W}}$ in $L^2(\R^N;\R^{M \times N})$}
\end{equation}
along a not relabeled subsequence.  

Let $V \subset \subset \imt(\boldsymbol{y},\Omega)$ be open.  By Proposition \ref{prop:conv-topim}(i), up to subsequences, we have $V \subset \subset \imt(\boldsymbol{y}_n,\Omega)$ for every $n \in \N$.   Thanks to Lemma \ref{lem:jump-extension}, we  get the identity
\begin{equation}
	\label{eqn:jump-imtop}
	J_{\widetilde{\boldsymbol{v}}_n} \cap \imt(\boldsymbol{y}_n,\Omega)= \big( J_{\boldsymbol{v}_n}\cap \imt(\boldsymbol{y}_n,\Omega) \big) \cup \big( J_{\widetilde{\boldsymbol{v}}_n} \cap \partial ^* \img(\boldsymbol{y}_n,\Omega) \cap \imt(\boldsymbol{y}_n,\Omega) \big),
\end{equation}
while claims (i) and (v) of Theorem \ref{thm:INV-top-im}  and Lemma \ref{lem:topim-point-disjoint}	  yield  
\begin{equation*}
	\haus\big(\partial ^* \img(\boldsymbol{y}_n,\Omega) \cap \imt(\boldsymbol{y}_n,\Omega)\big)=\sum_{\boldsymbol{a} \in C_{\boldsymbol{y}_n}} \per \big( \imt(\boldsymbol{y}_n,\boldsymbol{a}) \big) =\mathcal{S}(\boldsymbol{y}_n).
\end{equation*}
Therefore, we estimate
\begin{equation}
	\label{eqn:cav-bbbound}
		\sup_{n \in \N} \haus(J_{\widetilde{\boldsymbol{v}}_n} \cap V)\leq   \sup_{n \in \N} \left\{ \haus\big(J_{\boldsymbol{v}_n}\cap \imt(\boldsymbol{y}_n,\Omega)\big) + \mathcal{S}(\boldsymbol{y}_n) \right\}<+\infty.
\end{equation}
Given  \eqref{eqn:cav-bbound} and \eqref{eqn:cav-bbbound}, Theorem \ref{thm:ambrosio-compactness}  ensures the existence of a map  $\widehat{\boldsymbol{v }}\in GSBV^2(V;\R^M)$ and a not relabeled subsequence, both possibly depending on $V$, such that
\begin{equation}
	\label{eqn:ambrosio-loc1}
	\text{$\widetilde{\boldsymbol{v}}_{n} \to \widehat{\boldsymbol{ v }}$ a.e.\ in $V$},  \qquad \text{$  \widetilde{\boldsymbol{v}}_{n}  \wk \widehat{\boldsymbol{ v }}$ in $L^2(V;\R^M)$,} \qquad \text{$\nabla  \widetilde{\boldsymbol{v}}_{n}  \wk \nabla \widehat{\boldsymbol{ v }}$ in $L^2(V;\R^{M \times N})$.}
\end{equation}
From \eqref{eqn:cav-w}, we deduce $\widehat{\boldsymbol{ v }}\cong {\boldsymbol{w}}$ and $\nabla \widehat{\boldsymbol{ v }}\cong {{\boldsymbol{W}}}$ in $V$. In particular, the map $\widehat{\boldsymbol{ v }}$ does not depend on the set $V$.  Setting $\widetilde{\boldsymbol{v}}\coloneqq   {\boldsymbol{w}} \restr{\imt(\boldsymbol{y},\Omega)}  $  and  using again  \eqref{eqn:cav-bbound} and \eqref{eqn:cav-bbbound},   we obtain $\widetilde{\boldsymbol{v}}\in GSBV(\imt(\boldsymbol{y},\Omega);\R^M)$  with  $\nabla \widetilde{\boldsymbol{v}}\cong {\boldsymbol{W}} \restr{\imt(\boldsymbol{y},\Omega)}$. 

By a diagonal argument, we select a subsequence satisfying \eqref{eqn:ambrosio-loc1}  for every open set $V\subset \subset \imt(\boldsymbol{y},\Omega)$.  In particular, we obtain
\begin{equation}\label{eqn:v-imt}
	\text{$\chi_{\imt(\boldsymbol{y}_n,\Omega)}\widetilde{\boldsymbol{v}}_n \to \widetilde{\boldsymbol{v}}$  a.e.\ in $\imt(\boldsymbol{y},\Omega)$.}
\end{equation}
We claim that $ \widetilde{{\boldsymbol{v}}}\cong\boldsymbol{0}$  in $ \imc(\boldsymbol{y},\Omega)$.  By Proposition \ref{prop:cavities}(ii),   we have 
\begin{equation} \label{eqn:imc-conv}
	\text{$\chi_{\imc(\boldsymbol{y}_n,\Omega)\cap \imt(\boldsymbol{y},\Omega)}\to \chi_{\imc(\boldsymbol{y},\Omega)}$ in $L^1(\R^N)$.}
\end{equation}
 As $\widetilde{\boldsymbol{v}}_n\cong \boldsymbol{0}$ in $\imc(\boldsymbol{y}_n,\Omega)$  by definition,  the claim follows by combining \eqref{eqn:v-imt}--\eqref{eqn:imc-conv}.  
Therefore, setting $\boldsymbol{v}\coloneqq \widetilde{{\boldsymbol{v}}}\restr{\img(\boldsymbol{y},\Omega)} \in L^2( \img(\boldsymbol{y},\Omega);  \R^M)  $,  the map $\widetilde{\boldsymbol{v}}$ coincides almost everywhere with the extension  of $\boldsymbol{v} $ to $\imt(\boldsymbol{y},\Omega)$ by zero because of Theorem \ref{thm:INV-top-im}(iii).  Also, we deduce ${\boldsymbol{w}}\cong \boldsymbol{0}$  and ${\boldsymbol{W}}\cong \boldsymbol{O}$ on $\imc(\boldsymbol{y},\Omega)$.  

Next, we show that ${\boldsymbol{w}}\cong \boldsymbol{0}$ and ${\boldsymbol{W}}\cong \boldsymbol{O}$ on $\R^N\setminus \imt(\boldsymbol{y},\Omega)$.   To do so, we argue as in \cite[Proposition 7.1]{barchiesi.henao.moracorral}. Let $ Y \subset \R^N \setminus  \imt(\boldsymbol{y},\Omega)$ be a bounded and measurable set.  By    \eqref{eqn:imt=img} and  \eqref{eqn:cav-w},  we have
\begin{equation*}
	\left | \int_{ Y}  {{\boldsymbol{w}}}\,\d\boldsymbol{\xi}\right |= \liminf_{n \to \infty} \left | \int_{ Y } \chi_{\img(\boldsymbol{y}_n,\Omega)}\boldsymbol{v}_n\,\d\boldsymbol{\xi} \right | \leq \liminf_{n \to \infty} \int_{\img(\boldsymbol{y}_n,\Omega)\setminus \img(\boldsymbol{y},\Omega)} |\boldsymbol{v}_n|\,\d\boldsymbol{\xi}.
\end{equation*}
As the sequence $(\chi_{\img(\boldsymbol{y}_n,\Omega)}\boldsymbol{v}_n)_n$ is equi-integrable  by \eqref{eqn:cav-bound}  and $\leb(\img(\boldsymbol{y}_n,\Omega)\setminus \img(\boldsymbol{y},\Omega))\to 0$ by \eqref{eqn:imgeom-conv}, the right-hand side equals zero. Since $ Y $ is arbitrary, this proves the claim for ${{\boldsymbol{w}}}$. The proof for ${{\boldsymbol{W}}}$ is analogous.   We conclude that 
\begin{equation*}
	{\boldsymbol{w}}\cong \chi_{\imt(\boldsymbol{y},\Omega)}\widetilde{\boldsymbol{v}}\cong \chi_{\img(\boldsymbol{y},\Omega)}\boldsymbol{v}, \qquad {\boldsymbol{W}}\cong \chi_{\imt(\boldsymbol{y},\Omega)}\nabla \widetilde{\boldsymbol{v}}\cong \chi_{\img(\boldsymbol{y},\Omega)}\nabla \boldsymbol{v}.
\end{equation*}
Thus, \eqref{eqn:imt=img} and \eqref{eqn:cav-w} entail \eqref{eqn:cav-compactness3}--\eqref{eqn:cav-compactness4}.   
 From \eqref{eqn:imgeom-conv}--\eqref{eqn:imt=img} and \eqref{eqn:v-imt}, we see that
\begin{equation*}
	\text{$
		\chi_{\img(\boldsymbol{y}_n,\Omega)}\boldsymbol{v}_n\to \chi_{\img(\boldsymbol{y},\Omega)}\boldsymbol{v}$ a.e.\ in $\R^N$.}
\end{equation*}
Exploiting again the equi-integrability of $(\chi_{\img(\boldsymbol{y}_n,\Omega)}\boldsymbol{v}_n)_n$  together with \eqref{eqn:imgeom-conv},   \eqref{eqn:cav-compactness2} follows by Vitali's convergence theorem. 

We proceed with the proof of \eqref{eqn:cav-compactness5}.  Define   $\widetilde{\gamma}\colon [0,+\infty)\to [0,+\infty]$ by setting $\widetilde{\gamma}(t)\coloneqq t \gamma(1/t)$.   From \eqref{eqn:growth-gamma}, we deduce $\displaystyle \lim_{t \to +\infty} \widetilde{\gamma}(t)/t=+\infty$. Applying Corollary \ref{cor:change-of-variable}(ii), we deduce
\begin{equation*}
 \sup_{n \in \N}\int_{\img(\boldsymbol{y}_n,\Omega)} \widetilde{\gamma}(\det \nabla \boldsymbol{y}_n^{-1})\,\d\boldsymbol{\xi}= \sup_{n \in \N} \int_{\Omega} \gamma(\det D \boldsymbol{y}_n)\,\d\boldsymbol{x}<+\infty,
\end{equation*}
from which the equi-integrability of the sequence $(\chi_{\img(\boldsymbol{y}_n,\Omega)}\det \nabla \boldsymbol{y}_n^{-1})_n$ follows by    the De la Vallée {Poussin} criterion. 
Since the sequence $(\boldsymbol{v}_n \circ \boldsymbol{y}_n)_n$  is  assumed to be equi-integrable  and we have \eqref{eqn:cav-compactness2},   we establish  $\boldsymbol{v}\circ \boldsymbol{y}\in L^1(\Omega;\R^M)$ and  \eqref{eqn:cav-compactness5} by  applying Proposition \ref{prop:conv-comp}.

\textbf{Step 3 (Lower semicontinuity).} From \eqref{eqn:cav-compactness4}, it follows
\begin{equation}
	\label{eqn:cav-exchange}
	\int_{\img(\boldsymbol{y},\Omega)} |\nabla \boldsymbol{v}|^2\,\d\boldsymbol{\xi} \leq \liminf_{n \to \infty}\int_{\img(\boldsymbol{y}_n,\Omega)} |\nabla \boldsymbol{v}_n|^2\,\d\boldsymbol{\xi}. 
\end{equation}
Thanks to \eqref{eqn:cav-compactness1},  the condition $p>N-1$,  and the weak continuity of Jacobian minors \cite[Theorem 8.20, Part 4]{dacorogna},  we have
\begin{equation*}
\text{$\adj_r D\boldsymbol{y}_n \wk \adj_r D \boldsymbol{y}$ in $L^1\left (\Omega;\R^{\binom{N}{r}\times \binom{N}{r}}\right)$ \quad   	 for all    $ r=1,\dots,N$}. 
\end{equation*}
 From the assumption $\boldsymbol{v}(\boldsymbol{\xi})\in  Z$ for almost every $\boldsymbol{\xi}\in \img(\boldsymbol{y},\Omega)$, we deduce that $\boldsymbol{v}(\boldsymbol{y}(\boldsymbol{x}))\in  Z$ for almost every $\boldsymbol{x}\in \Omega$. Here, we exploit the fact that $\boldsymbol{y}$ satisfies Lusin's condition (N${}^{-1}$)  in view of  Remark~\ref{rem:federer}(b).  Therefore, recalling  \eqref{eqn:polyconvex1}--\eqref{eqn:polyconvex2} and   \eqref{eqn:cav-compactness5}, and  invoking  Theorem \ref{thm:lscT}, we obtain
\begin{equation}
	\label{eqn:cav-elastic}
	\int_{\Omega} W(D\boldsymbol{y},\boldsymbol{v}\circ \boldsymbol{y})\,\d\boldsymbol{x}\leq \liminf_{n \to \infty} \int_{\Omega} W(D\boldsymbol{y}_n,\boldsymbol{v}_n \circ \boldsymbol{y}_n)\,\d\boldsymbol{x}.
\end{equation} 
To prove the lower semicontinuity of the remaining energy terms, we proceed as follows.  By  repeating \eqref{eqn:jump-imtop} for ${\widetilde{\boldsymbol{v}}} $ in place of ${\widetilde{\boldsymbol{v}}}_n$    and using  claim (v) of Theorem \ref{thm:INV-top-im}  along with Lemma \ref{lem:topim-point-disjoint},  we have
\begin{equation}
	\label{eqn:jv}
	\begin{split}
		J_{\widetilde{\boldsymbol{v}}} \cap \imt(\boldsymbol{y},\Omega) &\simeq \big( J_{\boldsymbol{v}} \cap \imt(\boldsymbol{y},\Omega) \big) \cup \big( J_{\widetilde{\boldsymbol{v}}} \cap \partial^* \img(\boldsymbol{y},\Omega) \cap \imt(\boldsymbol{y},\Omega)  \big)\\
		&\simeq \big( J_{\boldsymbol{v}} \cap \imt(\boldsymbol{y},\Omega) \big) \cup \bigcup_{\boldsymbol{a}\in C_{\boldsymbol{y}}} \big(J_{\widetilde{\boldsymbol{v}}} \cap \partial^* \imt(\boldsymbol{y},\boldsymbol{a})\big).
	\end{split} 
\end{equation}
Let  $\boldsymbol{a}_0 \in C_{\boldsymbol{y}}$ and denote by  $\boldsymbol{\nu}_{\imt(\boldsymbol{y},\boldsymbol{a}_0)}\colon \partial^*\imt(\boldsymbol{y},\boldsymbol{a}_0) \to S$ the outer unit normal,   where we employ the notation in \eqref{eqn:B-S}.  
Recalling  Definition~\ref{def:jump},  the lateral traces of $\widetilde{\boldsymbol{v}}$ at $\boldsymbol{\xi}_0\in \partial^* \imt(\boldsymbol{y},  \boldsymbol{a}_0  )$ are defined as
\begin{equation}\label{XXX}
	\widetilde{\boldsymbol{v}}^\pm(\boldsymbol{\xi}_0)=\aplim_{\substack{\boldsymbol{\xi}\to \boldsymbol{\xi}_0 \\ \boldsymbol{\xi}\in H^\pm(\boldsymbol{\xi}_0,\boldsymbol{\nu}_0)}} \widetilde{\boldsymbol{v}}(\boldsymbol{\xi}),
\end{equation} 
where $\boldsymbol{\nu}_0\coloneqq \boldsymbol{\nu}_{\imt( \boldsymbol{y}_0,  \boldsymbol{a}_0)}( \boldsymbol{\xi}_0 )$. 
We  observe  that $\widetilde{\boldsymbol{v}}^-=\boldsymbol{0}$, so that the jump of $\widetilde{\boldsymbol{v}}$ across $\partial^*\,\imt(\boldsymbol{y},\boldsymbol{a}_0)$ equals $\left[\widetilde{\boldsymbol{v}} \right]=\widetilde{\boldsymbol{v}}^+$. As the function $\widetilde{\boldsymbol{v}}^+\colon \partial^* \imt(\boldsymbol{y},\boldsymbol{a}_0) \to \R^N$ is  Borel, the sets $\{ \boldsymbol{\xi} \in \partial^* \imt(\boldsymbol{y},\boldsymbol{a}_0):\hspace{4pt} \widetilde{\boldsymbol{v}}^+(\boldsymbol{\xi})=  \boldsymbol{b}\}$ are Borel for all $\boldsymbol{b}\in \R^M$. Since $\haus(\partial^* \imt(\boldsymbol{y},\boldsymbol{a}_0))<+\infty$, we have that for all $\boldsymbol{b}\in \R^M$, except  for an at most countable number, 
\begin{equation*}
	\haus(\{\boldsymbol{\xi} \in \partial^* \imt(\boldsymbol{y},\boldsymbol{a}_0):\hspace{4pt} \widetilde{\boldsymbol{v}}^+(\boldsymbol{\xi})=  \boldsymbol{b}  \})=0 . 
\end{equation*}
 Actually, as $C_{\boldsymbol{y}}$ is countable by claim (i) of Theorem \ref{thm:INV-top-im},  for all $\boldsymbol{b}\in \R^M$, except  for  a countable number,  
\begin{equation}
	\label{eqn:choice-b}
	\haus \left ( \left \{\boldsymbol{\xi} \in \bigcup_{\boldsymbol{a}\in C_{\boldsymbol{y}}} \partial^* \imt(\boldsymbol{y},\boldsymbol{a}):\hspace{4pt} \widetilde{\boldsymbol{v}}^+(\boldsymbol{\xi})=  \boldsymbol{b}  \right \}\right )=0.	
\end{equation}
 Let us choose any such $ \boldsymbol{b}$.  Define $\widetilde{\boldsymbol{v}}^{\,\boldsymbol{b}}\colon \imt(\boldsymbol{y},\Omega) \to \R^M$ by setting
\begin{equation*}
	\widetilde{\boldsymbol{v}}^{\,\boldsymbol{b}}(\boldsymbol{\xi})\coloneqq
	\begin{cases}
		\boldsymbol{v}(\boldsymbol{\xi}) & \text{if $\boldsymbol{\xi}\in \img(\boldsymbol{y},\Omega)$,} \\
		\boldsymbol{b} & \text{if $\boldsymbol{\xi}\in \imt(\boldsymbol{y},\Omega)\setminus \img(\boldsymbol{y},\Omega)$.}
	\end{cases}
\end{equation*}
 Then,   the jump of $\widetilde{\boldsymbol{v}}^{\boldsymbol{b}}$ across $\partial^* \imt(\boldsymbol{y},\boldsymbol{a}_0)$ for a given $\boldsymbol{a}_0\in C_{\boldsymbol{y}}$ equals $[\widetilde{\boldsymbol{v}}^{\boldsymbol{b}}]=\widetilde{\boldsymbol{v}}^+-\boldsymbol{b}$. Using   Lemma~\ref{lem:topim-point-disjoint}, we find that the jump of $\widetilde{\boldsymbol{v}}^{\boldsymbol{b}}$ across $\bigcup_{\boldsymbol{a}\in C_{\boldsymbol{y}}} \partial^* \imt(\boldsymbol{y},\boldsymbol{a})$ satisfies $[\widetilde{\boldsymbol{v}}^{\boldsymbol{b}}]\simeq\widetilde{\boldsymbol{v}}^+-\boldsymbol{b}\neq \boldsymbol{0}$ in $\bigcup_{\boldsymbol{a}\in C_{\boldsymbol{y}}} \partial^* \imt(\boldsymbol{y},\boldsymbol{a})$ owing to \eqref{eqn:choice-b}. In particular
\begin{equation*}
	J_{\widetilde{\boldsymbol{v}}^{\boldsymbol{b}}} \cap \bigcup_{\boldsymbol{a}\in C_{\boldsymbol{y}}} \partial^* \imt(\boldsymbol{y},\boldsymbol{a}) \simeq \bigcup_{\boldsymbol{a}\in C_{\boldsymbol{y}}} \partial^* \imt(\boldsymbol{y},\boldsymbol{a}).
\end{equation*}
 Hence, Theorem~\ref{thm:INV-top-im}(i) gives the identity 
\begin{equation*}
	\sum_{\boldsymbol{a}\in C_{\boldsymbol{y}}}\haus\big( J_{\widetilde{\boldsymbol{v}}^{\boldsymbol{b}}} \cap \partial^* \imt(\boldsymbol{y},\boldsymbol{a})\big)=\sum_{\boldsymbol{a}\in C_{\boldsymbol{y}}} \per \big( \imt(\boldsymbol{y},\boldsymbol{a})\big)=\mathcal{S}(\boldsymbol{y}).
\end{equation*}
 Using    \eqref{eqn:jv} for $\widetilde{\boldsymbol{v}}^{\boldsymbol{b}}$ in place of $\widetilde{\boldsymbol{v}}$,  this yields
\begin{equation}
	\label{eqn:j}
	\haus\big(J_{\widetilde{\boldsymbol{v}}^{\boldsymbol{b}}}\cap \imt(\boldsymbol{y},\Omega)\big)=\haus\big(J_{\boldsymbol{v}}\cap \imt(\boldsymbol{y},\Omega)\big)+\mathcal{S}(\boldsymbol{y}).
\end{equation}
Analogously, we can select a sequence $(\boldsymbol{b}_n)_n\subset \R^M$ in such a way that the map
$\widetilde{\boldsymbol{v}}_n^{\,\boldsymbol{b}_n}\colon \imt(\boldsymbol{y}_n,\Omega) \to \R^M$ given by
\begin{equation*}
	\widetilde{\boldsymbol{v}}_n^{\,\boldsymbol{b}_n}(\boldsymbol{\xi})\coloneqq
	\begin{cases}
		\boldsymbol{v}_n(\boldsymbol{\xi}) & \text{if $\boldsymbol{\xi}\in \img(\boldsymbol{y}_n,\Omega)$,} \\
		\boldsymbol{b}_n & \text{if $\boldsymbol{\xi}\in \imt(\boldsymbol{y}_n,\Omega)\setminus \img(\boldsymbol{y}_n,\Omega)$,}
	\end{cases}
\end{equation*}
satisfies 
\begin{equation}
	\label{eqn:jn}
	\haus\big(J_{\widetilde{\boldsymbol{v}}_n^{\boldsymbol{b}_n}}\cap \imt(\boldsymbol{y}_n,\Omega)\big)=\haus\big(J_{\boldsymbol{v}_n}\cap \imt(\boldsymbol{y}_n,\Omega)\big)+\mathcal{S}(\boldsymbol{y}_n)
\end{equation}
for every $n \in \N$ and, in addition, $\boldsymbol{b}_n \to \boldsymbol{b}$. Now, let $V \subset \subset \imt(\boldsymbol{y},\Omega)$ be open. By Proposition \ref{prop:conv-topim}(i) we find $V \subset \subset \imt(\boldsymbol{y}_n,\Omega)$   for every $n \in \N$ along a not relabeled subsequence.  From \eqref{eqn:cav-compactness2},   \eqref{eqn:cav-compactness4}, and \eqref{eqn:v-imt}--\eqref{eqn:imc-conv},  we have
\begin{equation*}
	\text{$\widetilde{\boldsymbol{v}}_n^{\boldsymbol{b}_n}\to \widetilde{\boldsymbol{v}}^{\boldsymbol{b}}$ a.e.\ in $V$,}\qquad 
	\text{$\nabla \widetilde{\boldsymbol{v}}_n^{\boldsymbol{b}_n} \wk \nabla \widetilde{\boldsymbol{v}}^{\boldsymbol{b}}$ in $L^2(V; \R^{M \times N}  )$.}
\end{equation*}
Applying Theorem \ref{thm:ambrosio-compactness}, we obtain
\begin{equation*}
	\haus\big(J_{\widetilde{\boldsymbol{v}}^{\boldsymbol{b}}}\cap V\big)\leq \liminf_{n \to \infty} 	\haus\big(J_{\widetilde{\boldsymbol{v}}_n^{\boldsymbol{b}_n}}\cap V\big)\leq \liminf_{n \to \infty} 	\haus\big(J_{\widetilde{\boldsymbol{v}}_n^{\boldsymbol{b}_n}}\cap \imt(\boldsymbol{y}_n,\Omega)\big),
\end{equation*}
from which, taking the supremum among all open sets $V \subset \subset \imt(\boldsymbol{y},\Omega)$, we deduce
\begin{equation*}
	\haus\big(J_{\widetilde{\boldsymbol{v}}^{\boldsymbol{b}}}\cap \imt(\boldsymbol{y},\Omega)\big)\leq \liminf_{n \to \infty} 	\haus\big(J_{\widetilde{\boldsymbol{v}}_n^{\boldsymbol{b}_n}}\cap \imt(\boldsymbol{y}_n,\Omega)\big).
\end{equation*}
Eventually, recalling \eqref{eqn:j}--\eqref{eqn:jn}, we get
\begin{equation}
	\label{eqn:cav-surface}
	\haus\big(J_{\boldsymbol{v}}\cap \imt(\boldsymbol{y},\Omega)\big)+\mathcal{S}(\boldsymbol{y})\leq \liminf_{n \to \infty}  \haus\big(J_{\boldsymbol{v}_n}\cap \imt(\boldsymbol{y}_n,\Omega)\big) + \mathcal{S}(\boldsymbol{y}_n).
\end{equation}
The combination of \eqref{eqn:cav-exchange}--\eqref{eqn:cav-elastic},  \eqref{eqn:cav-surface},    Theorem \ref{thm:determinant-surface-energy},  and the assumption $\lambda_1\geq \lambda_3>0$  gives \eqref{eqn:cav-lsc}.
\end{proof}

The next example motivates our choice to penalize only the measure of  $J_{\boldsymbol{v}}\cap \imt(\boldsymbol{y},\Omega)$ rather than that of the whole set $J_{\boldsymbol{v}}$. Also, it shows that, in general, the functional $(\boldsymbol{y},\boldsymbol{v})\mapsto \haus(J_{\boldsymbol{v}})$ is not lower semicontinuous with respect to the relevant topology.

\begin{figure}
	\begin{tikzpicture}[baseline, remember picture, scale=1.2]
		\usetikzlibrary{calc}
	\draw[fill=black!20]  (-4.5,2) rectangle (-3.5,-2);
	\node[above right] at (-3.5,1.7) {\large $\Omega$};
	\node[above right] at (1,1.7) {\large $\img(\boldsymbol{y},\Omega)$};
	\draw[fill=black!20] (180:1) coordinate (beta) arc (180:-180:1) coordinate (alpha) -- (-180:2) arc (-180:180:2) -- cycle;
	\coordinate  (x) at (2.5,0);
	\coordinate (o) at (0,0);
	\foreach \a in {35,60,90,118}
	{
		\draw[blue,->,thick] ({1.33*cos(\a)},{1.33*sin(\a)}) -- ({(1.33*cos(\a)-(\a/145)*sin(\a))},{(1.33*sin(\a)+(\a/145)*cos(\a))}); 
		\draw[blue,->,thick] ({1.66*cos(\a)},{1.66*sin(\a)}) -- ({(1.66*cos(\a)-(\a/145)*sin(\a))},{(1.66*sin(\a)+(\a/145)*cos(\a))}); 
		\draw[blue,->,thick] ({1.33*cos(\a)},{-1.33*sin(\a)}) -- ({(1.33*cos(\a)-(\a/145)*sin(\a))},{(-1.33*sin(\a)-(\a/145)*cos(\a))}); 
		\draw[blue,->,thick] ({1.66*cos(\a)},{-1.66*sin(\a)}) -- ({(1.66*cos(\a)-(\a/145)*sin(\a))},{(-1.66*sin(\a)-(\a/145)*cos(\a))}); 
	}
	\foreach \a in {15}
	{
		\draw[blue,->,thick] ({1.33*cos(\a)},{1.33*sin(\a)}) -- ({(1.33*cos(\a)-(.15)*sin(\a))},{(1.33*sin(\a)+(.15)*cos(\a))}); 
		\draw[blue,->,thick] ({1.66*cos(\a)},{1.66*sin(\a)}) -- ({(1.66*cos(\a)-(.15)*sin(\a))},{(1.66*sin(\a)+(.15)*cos(\a))}); 
		\draw[blue,->,thick] ({1.33*cos(\a)},{-1.33*sin(\a)}) -- ({(1.33*cos(\a)-(.15)*sin(\a))},{(-1.33*sin(\a)-(.15)*cos(\a))}); 
		\draw[blue,->,thick] ({1.66*cos(\a)},{-1.66*sin(\a)}) -- ({(1.66*cos(\a)-(.15)*sin(\a))},{(-1.66*sin(\a)-(.15)*cos(\a))}); 
	}
	\foreach \a in {145}
	{
		\draw[blue,->,thick] ({1.33*cos(\a)},{1.33*sin(\a)}) -- ({(1.33*cos(\a)-(\a/170)*sin(\a))},{(1.33*sin(\a)+(\a/170)*cos(\a))}); 
		\draw[blue,->,thick] ({1.66*cos(\a)},{1.66*sin(\a)}) -- ({(1.66*cos(\a)-(\a/170)*sin(\a))},{(1.66*sin(\a)+(\a/170)*cos(\a))}); 
		\draw[blue,->,thick] ({1.33*cos(\a)},{-1.33*sin(\a)}) -- ({(1.33*cos(\a)-(\a/170)*sin(\a))},{(-1.33*sin(\a)-(\a/170)*cos(\a))}); 
		\draw[blue,->,thick] ({1.66*cos(\a)},{-1.66*sin(\a)}) -- ({(1.66*cos(\a)-(\a/170)*sin(\a))},{(-1.66*sin(\a)-(\a/170)*cos(\a))}); 
	}
	\node[blue, ultra thick] at (1.33,0) {.}; 	\node[blue, ultra thick] at (1.66,0) {.};
	\node[red] at (-.8,0) {$\Sigma$};
	\node[blue] at (.25,1.6) {\large $\boldsymbol{v}$};
	\draw[red,thick] (-2,0)--(-1,0);
	\end{tikzpicture}
	\hspace*{30pt} 
	\begin{tikzpicture}[baseline, remember picture, scale=1.2]
	\draw[fill=black!20] (160:1) coordinate (beta) arc (160:-160:1) coordinate (alpha) -- (-160:2) arc (-160:160:2) -- cycle;
	\node[above right] at (1,1.7) {\large $\img(\boldsymbol{y}_n,\Omega)$};
	\coordinate  (x) at (2.5,0);
	\coordinate (o) at (0,0);
	\foreach \a in {35,60,90,118}
	{
		\draw[blue,->,thick] ({1.33*cos(\a)},{1.33*sin(\a)}) -- ({(1.33*cos(\a)-(\a/145)*sin(\a))},{(1.33*sin(\a)+(\a/145)*cos(\a))}); 
		\draw[blue,->,thick] ({1.66*cos(\a)},{1.66*sin(\a)}) -- ({(1.66*cos(\a)-(\a/145)*sin(\a))},{(1.66*sin(\a)+(\a/145)*cos(\a))}); 
		\draw[blue,->,thick] ({1.33*cos(\a)},{-1.33*sin(\a)}) -- ({(1.33*cos(\a)-(\a/145)*sin(\a))},{(-1.33*sin(\a)-(\a/145)*cos(\a))}); 
		\draw[blue,->,thick] ({1.66*cos(\a)},{-1.66*sin(\a)}) -- ({(1.66*cos(\a)-(\a/145)*sin(\a))},{(-1.66*sin(\a)-(\a/145)*cos(\a))}); 
	}
	\foreach \a in {15}
	{
		\draw[blue,->,thick] ({1.33*cos(\a)},{1.33*sin(\a)}) -- ({(1.33*cos(\a)-(.15)*sin(\a))},{(1.33*sin(\a)+(.15)*cos(\a))}); 
		\draw[blue,->,thick] ({1.66*cos(\a)},{1.66*sin(\a)}) -- ({(1.66*cos(\a)-(.15)*sin(\a))},{(1.66*sin(\a)+(.15)*cos(\a))}); 
		\draw[blue,->,thick] ({1.33*cos(\a)},{-1.33*sin(\a)}) -- ({(1.33*cos(\a)-(.15)*sin(\a))},{(-1.33*sin(\a)-(.15)*cos(\a))}); 
		\draw[blue,->,thick] ({1.66*cos(\a)},{-1.66*sin(\a)}) -- ({(1.66*cos(\a)-(.15)*sin(\a))},{(-1.66*sin(\a)-(.15)*cos(\a))}); 
	}
	\foreach \a in {145}
	{
		\draw[blue,->,thick] ({1.33*cos(\a)},{1.33*sin(\a)}) -- ({(1.33*cos(\a)-(\a/170)*sin(\a))},{(1.33*sin(\a)+(\a/170)*cos(\a))}); 
		\draw[blue,->,thick] ({1.66*cos(\a)},{1.66*sin(\a)}) -- ({(1.66*cos(\a)-(\a/170)*sin(\a))},{(1.66*sin(\a)+(\a/170)*cos(\a))}); 
		\draw[blue,->,thick] ({1.33*cos(\a)},{-1.33*sin(\a)}) -- ({(1.33*cos(\a)-(\a/170)*sin(\a))},{(-1.33*sin(\a)-(\a/170)*cos(\a))}); 
		\draw[blue,->,thick] ({1.66*cos(\a)},{-1.66*sin(\a)}) -- ({(1.66*cos(\a)-(\a/170)*sin(\a))},{(-1.66*sin(\a)-(\a/170)*cos(\a))}); 
	}
	\node[blue, ultra thick] at (1.33,0) {.}; 	\node[blue, ultra thick] at (1.66,0) {.};
	\node[blue] at (.25,1.6) {\large $\boldsymbol{v}_n$};
	\end{tikzpicture}
	\caption{The deformations in Example \ref{ex:jump-not-contained-imt} with the corresponding Eulerian fields. }
	\label{fig:jump-not-contained-imt}
\end{figure}
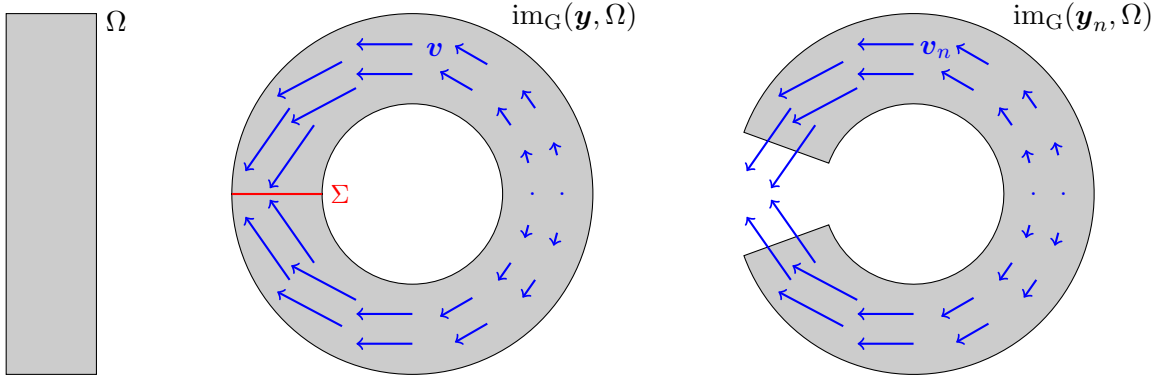

\begin{example}[Jump points of Eulerian maps] \label{ex:jump-not-contained-imt}
Let $N=M=2$. We employ the notation in \eqref{eqn:q-norm}--\eqref{eqn:q-ball-annulus}. Also, we write \begin{equation*}
\vartheta(\boldsymbol{\xi})\coloneqq \mathrm{sgn}(\xi_2)\arccos\left(\frac{\xi_1}{|\boldsymbol{\xi}|}  \right) \quad \text{for all $\boldsymbol{\xi}\in  \R^2 \setminus \lbrace 0 \rbrace $.}
\end{equation*}
In this way, $-\pi \leq \vartheta(\boldsymbol{\xi})\leq \pi$ for all $\boldsymbol{\xi}\in   \R^2 \setminus \lbrace 0 \rbrace$.
Set $\Omega\coloneqq (1/2,1)\times (-1,1)$. The deformation $\boldsymbol{y}\colon \Omega \to \R^2$	defined by $\boldsymbol{y}(\boldsymbol{x})\coloneqq x_1(\cos (\pi x_2),\sin (\pi x_2))$ is a diffeomorphism onto its image and transforms $\Omega$ into an annulus (Figure \ref{fig:jump-not-contained-imt}). Precisely, we have 
\begin{equation*}
\imt(\boldsymbol{y},\Omega)=\img(\boldsymbol{y},\Omega)=\left \{ \boldsymbol{\xi}\in A(1/2,1): \hspace{2pt}  -\pi <\vartheta(\boldsymbol{\xi})<\pi \right\}=A(1/2,1)\setminus \Sigma,
\end{equation*}
where $\Sigma\coloneqq (-1,-1/2)\times \{0\}$. Observe that $\img(\boldsymbol{y},\Omega)^{(1)}=A(1/2,1)$, so that $\imt(\boldsymbol{y},\Omega)\subset \img(\boldsymbol{y},\Omega)^{(1)}$  with strict inclusion. For every $n\in \N$, let $\boldsymbol{y}_n\colon \Omega \to \R^2$ be defined as $\boldsymbol{y}_n(\boldsymbol{x})\coloneqq x_1\left (\cos \left( \alpha_n x_2 \right) ,\sin \left(\alpha_n x_2 \right) \right )$, where $\alpha_n\coloneqq \frac{\pi n}{n+1}$. This map is also a diffeomorphism and transforms $\Omega$ into an annular sector (Figure \ref{fig:jump-not-contained-imt}).  We have
\begin{equation*}
	\imt(\boldsymbol{y}_n,\Omega)=\img(\boldsymbol{y}_n,\Omega)=\left \{ \boldsymbol{\xi}\in A(1/2,1):  \hspace{2pt} -\alpha_n<\vartheta(\boldsymbol{\xi}) <\alpha_n   \right \}.
\end{equation*}
One checks that \eqref{eqn:weak}--\eqref{eqn:weak-det} hold.
Now, define  $\boldsymbol{v}\colon \img(\boldsymbol{y},\Omega)\to\R^2$ as  
$\boldsymbol{v}(\boldsymbol{\xi})\coloneqq (\vartheta(\boldsymbol{\xi})/\pi)(-\xi_2,\xi_1)^\top/|\boldsymbol{\xi}|$ 
Then, $J_{\boldsymbol{v}}=\Sigma$. Indeed, for every $-1<\xi_1<-1/2$, we have  $\lim_{\xi_2 \to 0^\pm} \boldsymbol{v}(\boldsymbol{\xi})=  \mp  \boldsymbol{e}_2$.    In particular, $J_{\boldsymbol{v}}\subset \R^2 \setminus \imt(\boldsymbol{y},\Omega)$. Define also $\boldsymbol{v}_n\colon \img(\boldsymbol{y}_n,\Omega)\to \R^2$ as $\boldsymbol{v}_n\coloneqq \boldsymbol{v}\restr{\img(\boldsymbol{y}_n,\Omega)}$ for all $n\in \N$. Then, the convergences in \eqref{eqn:cav-compactness1}--\eqref{eqn:cav-compactness5} hold and  $J_{\boldsymbol{v}_n}=\emptyset$ for every $n\in\N$. Thus,
\begin{equation*}
\mathscr{H}^1(J_{\boldsymbol{v}})=\mathscr{H}^1(\Sigma)=\frac{1}{2}>0=\liminf_{n\to \infty} \mathscr{H}^1(J_{\boldsymbol{v}_n}).
\end{equation*}
\end{example}

In the next example, we illustrate the observation made in Remark \ref{rem:cav}(c). Again, we look at cavities squashing onto lower-dimensional sets along converging sequences of deformations.

\begin{figure}
	\centering
	\begin{tikzpicture}[scale=3.8]
		\draw[fill=black!20] (0,1) -- (1,0) -- (0,-1) -- (-1,0) -- (0,1);
		\draw[fill=white] (0,.2) -- (.5,0) -- (0,-.2) -- (-.5,0) -- (0,.2);
		\foreach \a in {90,120,150,180,210,240,270}
		{
			\draw[blue,->,thick] ({-.5+.05*cos(\a)},{.05*sin(\a)}) -- ({-.5+.35*cos(\a)},{.35*sin(\a)});
		}
		\foreach \a in {90,60,30,0,-30,-60,-90}
		{
			\draw[blue,->,thick] ({.5+.05*cos(\a)},{.05*sin(\a)}) -- ({.5+.35*cos(\a)},{.35*sin(\a)});
		}	
		\draw[blue,->,thick] (-.33,.2) -- (-.33,.5); 
		\draw[blue,->,thick] (-.16,.4) -- (-.16,.7); 
		\draw[blue,->,thick] (.33,.2) -- (.33,.5); 
		\draw[blue,->,thick] (.16,.4) -- (.16,.7); 
		\draw[blue,->,thick] (0,.55) -- (0,.85); 
		\draw[blue,->,thick] (-.33,-.2) -- (-.33,-.5); 
		\draw[blue,->,thick] (-.16,-.4) -- (-.16,-.7); 
		\draw[blue,->,thick] (.33,-.2) -- (.33,-.5); 
		\draw[blue,->,thick] (.16,-.4) -- (.16,-.7); 
		\draw[blue,->,thick] (0,-.55) -- (0,-.85);
		\node[right] at (.15,1) {\large $\img(\boldsymbol{y}_n,\Omega)$};
		\node[blue] at (.43,.45) {\large $\boldsymbol{v}_n$};
	\end{tikzpicture}
	\hspace{30pt}
	\begin{tikzpicture}[scale=3.8]
		\usetikzlibrary{calc}
		\draw[fill=black!20] (0,1) -- (1,0) -- (0,-1) -- (-1,0) -- (0,1);
		\draw[red,thick] (-.5,0)--(.5,0);
		\foreach \a in {90,120,150,180,210,240,270}
		{
			\draw[blue,->,thick] ({-.5+.05*cos(\a)},{.05*sin(\a)}) -- ({-.5+.35*cos(\a)},{.35*sin(\a)});
		}
		\foreach \a in {90,60,30,0,-30,-60,-90}
		{
			\draw[blue,->,thick] ({.5+.05*cos(\a)},{.05*sin(\a)}) -- ({.5+.35*cos(\a)},{.35*sin(\a)});
		}	
		\draw[blue,->,thick] (-.33,.2) -- (-.33,.5); 
		\draw[blue,->,thick] (-.16,.4) -- (-.16,.7); 
		\draw[blue,->,thick] (.33,.2) -- (.33,.5); 
		\draw[blue,->,thick] (.16,.4) -- (.16,.7); 
		\draw[blue,->,thick] (0,.55) -- (0,.85); 
		\draw[blue,->,thick] (-.33,-.2) -- (-.33,-.5); 
		\draw[blue,->,thick] (-.16,-.4) -- (-.16,-.7); 
		\draw[blue,->,thick] (.33,-.2) -- (.33,-.5); 
		\draw[blue,->,thick] (.16,-.4) -- (.16,-.7); 
		\draw[blue,->,thick] (0,-.55) -- (0,-.85);
		\node[right] at (.15,1) {\large $\img(\boldsymbol{y},\Omega)$};
		\node[blue] at (.43,.45) {\large $\boldsymbol{v}$};
		\node[red] at (.25,.1) {\large $\Sigma$};
	\end{tikzpicture}
	\caption{The Eulerian fields in Example \ref{ex:jump-not-lsc}.}
\end{figure}
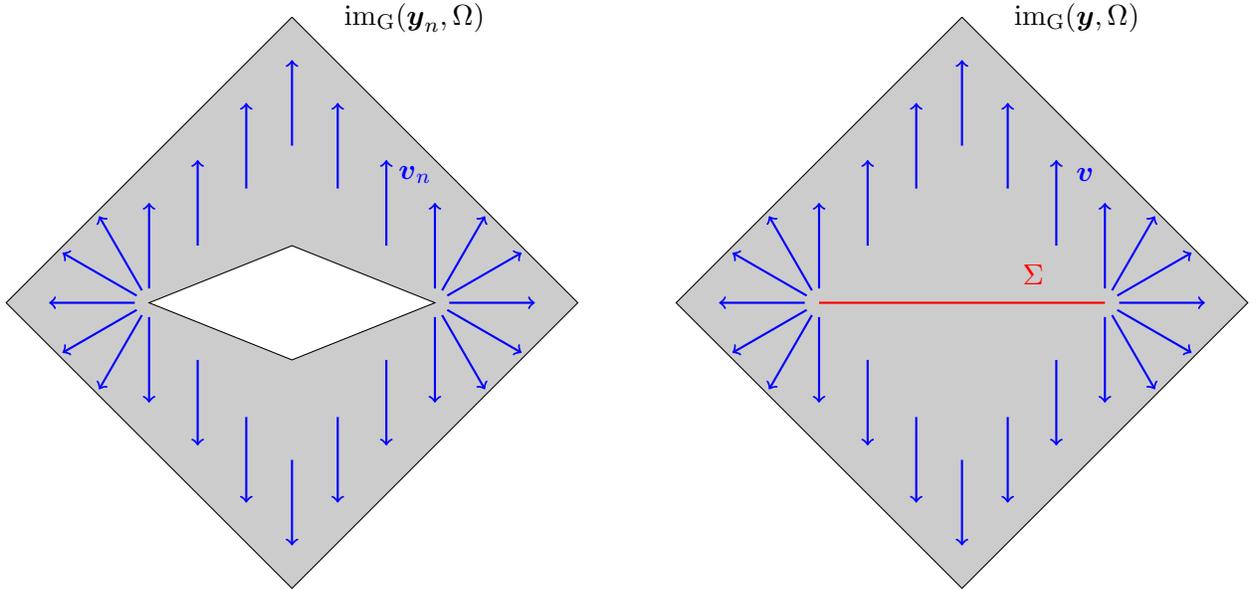

\begin{example}[Jump points of Eulerian maps inside the topological image]\label{ex:jump-not-lsc}  
Let $N=M=2$. We employ the notation in \eqref{eqn:q-ball-annulus}.  Let $\Omega$, $Q_n$, $\Sigma$, $\boldsymbol{y}_n$, and $\boldsymbol{y}$ be as in Example \ref{ex:cavity-segment}. Recall that
\begin{equation*}
	\text{ $\mathcal{S}(\boldsymbol{y}_n)=4\sqrt{ 1/4 +2^{-2n}}$ \quad for all $n\in \N$, \qquad  $\mathcal{S}(\boldsymbol{y})=0$. }
\end{equation*}
Define the Eulerian field
 $\boldsymbol{v}\colon \img(\boldsymbol{y},\Omega)\to \R^2$ by setting
\begin{equation*}
	\boldsymbol{v}(\boldsymbol{\xi})\coloneqq \begin{cases}
	\sgn(\xi_2)\boldsymbol{e}_2 & \text{if $|\xi_1|\leq 1/2$,} \\
	\displaystyle \frac{\boldsymbol{\xi}-\sgn(\xi_1)(\boldsymbol{e}_1/2)}{|\boldsymbol{\xi}-\sgn(\xi_1)(\boldsymbol{e}_1/2)|} & \text{if $1/2<|\xi_1|<1$.}
	\end{cases}
\end{equation*} 
This map is smooth within its domain $\img(\boldsymbol{y},\Omega)= B_1  \setminus \Sigma$ and  $J_{\boldsymbol{v}}=\Sigma \subset \imt(\boldsymbol{y},\Omega)$.
Setting $\boldsymbol{v}_n\coloneqq \boldsymbol{v}\restr{\img(\boldsymbol{y}_n,\Omega)}\colon \img(\boldsymbol{y}_n,\Omega)\to \R^2$ for every $n\in \N$,  where $\img(\boldsymbol{y}_n,\Omega)=B_1 \setminus Q_n$,   we observe that $J_{\boldsymbol{v}_n}=\emptyset$ and  the convergences in \eqref{eqn:cav-compactness1}--\eqref{eqn:cav-compactness5} hold. However,
\begin{equation*}
	\mathscr{H}^1(J_{\boldsymbol{v}}\cap \imt(\boldsymbol{y},\Omega))=\mathscr{H}^1(\Sigma)=1>0=\liminf_{n\to\infty} \mathscr{H}^1(J_{\boldsymbol{v}_n}\cap \imt(\boldsymbol{y}_n,\Omega)),
\end{equation*}
while
\begin{equation*}
	\mathscr{H}^1(J_{\boldsymbol{v}}\cap \imt(\boldsymbol{y},\Omega)) + \mathcal{S(\boldsymbol{y})}=1< 2  =\liminf_{n \to \infty} \left \{ \mathscr{H}^1(J_{\boldsymbol{v}_n}\cap \imt(\boldsymbol{y}_n,\Omega)) + \mathcal{S}(\boldsymbol{y}_n)\right  \}.
\end{equation*}
\end{example}

In Theorem \ref{thm:cav}  and Remark \ref{rem:cav}(f),  the assumption $\widetilde{\boldsymbol{v}}\in W^{1,2}(\imt(\boldsymbol{y},\Omega);\R^M)$  formally corresponds to the case in which  $\boldsymbol{v}$ satisfies   $\boldsymbol{v}=\boldsymbol{0}$ on $\partial\, \img(\boldsymbol{y},\Omega)\cap \imt(\boldsymbol{y},\Omega)$.  However, such boundary conditions may or may not be acceptable  depending on the physical meaning of the variable $\boldsymbol{v}$ in specific models, see Remark \ref{rem:nem}(a), Remark \ref{rem:ph}(a), and Remark \ref{rem:mag}(a). 

For this reason, we propose below an alternative approach. Unfortunately,  this poses a restriction on the class of admissible deformations.
 Given $\kappa>0$, we consider the class $\mathcal{Y}^{\rm cav}_{p,\kappa}(\Omega)$ in \eqref{eqn:Y-cav-kappa} and we recall  
Definition~\ref{def:material-image}. For every $\boldsymbol{y} \in \mathcal{Y}^{\rm cav}_{p,\kappa}(\Omega)$   with   $\mathcal{S}(\boldsymbol{y})<+\infty$, Remark \ref{rem:INV-top-im}(a) implies $\mathscr{H}^0({ C_{\boldsymbol{y}}}) < +\infty$, and thus   the material image $\imm(\boldsymbol{y},\Omega)$ is an open set with $\imm(\boldsymbol{y},\Omega)\cong \img(\boldsymbol{y},\Omega)$ by  items  (a) and (b) in Remark \ref{rem:matim}. Thus, we can consider $\boldsymbol{v}\in W^{1,2}(\imm(\boldsymbol{y},\Omega); Z)$. In this setting, the functional in \eqref{eqn:F-EL}--\eqref{eqn:F-EL2} can be rewritten as
\begin{equation}
	\label{eqn:F-EL-cav2}
	 \mathcal{F}_2^{\rm cav}(\boldsymbol{y},\boldsymbol{v})= \int_{\Omega} W(D \boldsymbol{y},\boldsymbol{v}\circ \boldsymbol{y})\,\d\boldsymbol{x}+ \lambda_1\mathcal{S}(\boldsymbol{y}) + \int_{\imm(\boldsymbol{y},\Omega)} |D \boldsymbol{v}|^2\,\d\boldsymbol{\xi}. 
\end{equation}

For this functional, we have the following result.

\begin{theorem}[Coercivity and lower semicontinuity: Sobolev deformations  with $\boldsymbol{\kappa}>0$] \label{thm:cavk}
Let  $\kappa>0$ and $ \mathcal{F}_2^{\rm cav}$  be defined as in \eqref{eqn:F-EL-cav2} with $W$ satisfying  {\rm (W1)}--{\rm (W3)} and $\lambda_1>0$.   Let $(\boldsymbol{y}{}_n)_n\subset \mathcal{Y}^{\rm cav}_{p,\kappa}(\Omega)$ and $(\boldsymbol{v}_n)_n$ be a sequence of maps $\boldsymbol{v}_n \in W^{1,2}( \imm(\boldsymbol{y}_n,  \Omega); Z)$. Suppose that 
\begin{equation}\label{eqn:cav2-bound}
	\sup \left \{ \|\boldsymbol{y}_n\|_{L^p(\Omega;\R^N)}+ \|\boldsymbol{v}_n\|_{L^2(\imm(\boldsymbol{y}_n,\Omega);\R^M)}+ \mathcal{F}_2^{\rm cav}(\boldsymbol{y}_n,\boldsymbol{v}_n) \right \}<+\infty.
\end{equation}
 Eventually,  assume that the sequence $(\boldsymbol{v}_n \circ \boldsymbol{y}_n)_n$ is equi-integrable. Then, there exist  $\boldsymbol{y}\in  \mathcal{Y}_{p,\kappa}^{\rm cav}(\Omega) $ and $\boldsymbol{v}\in W^{1,2}(\imm(\boldsymbol{y},\Omega);\R^M)$ satisfying $\boldsymbol{v}\circ \boldsymbol{y}\in L^1(\Omega;\R^M)$  such that, up to subsequences, we have:
\begin{align}
\label{eqn:cav2-compactness1}
\text{$\boldsymbol{y}_n \wk \boldsymbol{y}$ in $W^{1,p}(\Omega;\R^N)$}, \qquad &\text{$\det D\boldsymbol{y}_n \wk \det D\boldsymbol{y}$ in $L^1(\Omega)$},\\
\label{eqn:cav2-compactness2}
\text{$\chi_{\imm(\boldsymbol{y}_n,\Omega)}\boldsymbol{v}_n \to \chi_{\imm(\boldsymbol{y},\Omega)}\boldsymbol{v}$ } &\text{in $L^1(\R^N;\R^M)$,}\\ 
\label{eqn:cav2-compactness3}
\text{$\chi_{\imm(\boldsymbol{y}_n,\Omega)}\boldsymbol{v}_n \wk \chi_{\imm(\boldsymbol{y},\Omega)}\boldsymbol{v}$ } &\text{in $L^2(\R^N;\R^M)$,}  \\
\label{eqn:cav2-compactness4}
\text{$\chi_{\imm(\boldsymbol{y}_n,\Omega)}D\boldsymbol{v}_n \wk \chi_{\imm(\boldsymbol{y},\Omega)}D\boldsymbol{v}$ } &\text{in $L^2(\R^N;\R^{M \times N})$,}\\
\label{eqn:cav2-compactness5}
\text{$\boldsymbol{v}_n \circ \boldsymbol{y}_n \to \boldsymbol{v}\circ \boldsymbol{y}$ } &\text{in $L^1(\Omega;\R^M)$.}
\end{align}
 Moreover, if $\boldsymbol{v}(\boldsymbol{\xi})\in Z$ for almost every $\boldsymbol{\xi}\in\imm(\boldsymbol{y},\Omega)$, then we have
\begin{equation}\label{eqn:F-EL-cav2-lsc}
	\mathcal{F}^{\rm cav}_2(\boldsymbol{y},\boldsymbol{v})\leq \liminf_{n\to \infty}\mathcal{F}^{\rm cav}_2(\boldsymbol{y}_n,\boldsymbol{v}_n).
\end{equation}  
\end{theorem}

\begin{remark}\label{rem:lsc-SobolevII}
	\begin{enumerate}[(a)]
		\item As mentioned in  Remark \ref{rem:Y-cav-kappa}(c),  the a priori lower bound   in \eqref{eqn:Y-cav-kappa} excludes the possibility for cavities to close along a  converging sequence of deformations. The  latter phenomenon  may  lead  to a relaxation effect for which a sequence of  converging  Eulerian maps  exhibits  jumps in the limit. An explicit case is discussed in Example \ref{ex:cavity-segment-relaxation} below. 
		\item  As before,   the assumption $\boldsymbol{v}(\boldsymbol{\xi})\in  Z $ for almost every $\boldsymbol{\xi}\in \imm(\boldsymbol{y},\Omega)$ is superfluous if $ Z$ is closed.
	\end{enumerate}
\end{remark}

\begin{proof}
The proof is divided into three steps.

\textbf{Step 1 (Compactness of deformations).} Arguing as in Step 1 of the proof of Theorem \ref{thm:cav}, we find $\boldsymbol{y}\in \mathcal{Y}^{\rm cav}_{p}(\Omega)$ satisfying \eqref{eqn:cav2-compactness1} and also
\begin{equation}
	\label{eqn:kcav-imm}
	\text{$\chi_{\imm(\boldsymbol{y}_n,\Omega)} \to \chi_{\imm(\boldsymbol{y},\Omega)}$ in $L^1(\R^N)$}
\end{equation}
thanks to  Remark \ref{rem:matim}(a).  By claim (ii) of Lemma \ref{lem:cavities-Det2},  we actually have $\boldsymbol{y}\in \mathcal{Y}^{\rm cav}_{p,\kappa}(\Omega)$.

\textbf{Step 2 (Compactness of  Eulerian maps).} The strategy in this step follows the proof of \cite[Proposition 7.1]{barchiesi.henao.moracorral} and \cite[Theorem 3.2]{bresciani}.  Given \eqref{eqn:cav2-bound},  there exist ${\boldsymbol{w}}  \in L^2(\R^N;\R^M)$ and ${\boldsymbol{W}} \in L^2(\R^N;\R^{M \times N})$ such that, up to subsequences, we have
\begin{equation}
	\label{eqn:kcav-w}
	\text{$\chi_{\imm(\boldsymbol{y}_n,\Omega)}\boldsymbol{v}_n \wk {\boldsymbol{w}}$ in $L^2(\R^N;\R^M)$, \qquad $\chi_{\imm(\boldsymbol{y}_n,\Omega)}D\boldsymbol{v}_n \wk {\boldsymbol{W}}$ in $L^2(\R^N;\R^{M \times N})$. }
\end{equation}
Let $V \subset \subset \imm(\boldsymbol{y},\Omega)$ be open and smooth. By Proposition \ref{prop:conv-matim},  we have $V \subset \subset \imm(\boldsymbol{y}_n,\Omega)$ for every $n \in \N$ along a not relabeled subsequence. Observe that
\begin{equation*}
	\sup_{n \in \N} \left\{ \int_V |\boldsymbol{v}_n|^2\,\d\boldsymbol{\xi}+ \int_V |D\boldsymbol{v}_n|^2\,\d\boldsymbol{\xi}  \right\} \leq \sup_{n \in \N} \left\{ \|\boldsymbol{v}_n\|_{L^2(\imm(\boldsymbol{y},\Omega);\R^M)} +  \mathcal{F}_2^{\rm cav} (\boldsymbol{y}_n,\boldsymbol{v}_n) \right\},
\end{equation*}
where the right-hand side is finite because of \eqref{eqn:cav2-bound}. 
Hence, there exist $\widehat{ \boldsymbol{v}}\in W^{1,2}(V;\R^M)$ and a not relabeled subsequence for which
\begin{equation}
	\label{eqn:kcav-loc}
	\text{$\boldsymbol{v}_n \to \widehat{ \boldsymbol{v}}$ a.e.\ in $V$, \qquad $\boldsymbol{v}_n \wk \widehat{ \boldsymbol{v}}$ in $L^2(V;\R^M)$, \qquad $D\boldsymbol{v}_n \wk D \widehat{ \boldsymbol{v}}$ in $L^2(V;\R^{M \times N})$.}
\end{equation}
Combining \eqref{eqn:kcav-w} and \eqref{eqn:kcav-loc}, we deduce that $\widehat{ \boldsymbol{v}}\cong  {\boldsymbol{w}}$ and $D\widehat{ \boldsymbol{v}}\cong {\boldsymbol{W}}$ in $V$.  Thus, the map $\boldsymbol{v}\coloneqq {\boldsymbol{w}}\restr{\imm(\boldsymbol{y},\Omega)}$ belongs to the space $W^{1,2}(\imm(\boldsymbol{y},\Omega);\R^M)$.

By means of a diagonal argument, we select a not relabeled subsequence satisfying
\begin{equation*}
	\text{$\boldsymbol{v}_n \to {\boldsymbol{v}}$ a.e.\ in $V$, \qquad $\boldsymbol{v}_n \wk {\boldsymbol{v}}$ in $L^2(V;\R^M)$, \qquad $D\boldsymbol{v}_n \wk D \boldsymbol{v}$ in $L^2(V;\R^{M \times N})$}
\end{equation*}
for every open set $V \subset \subset \imm(\boldsymbol{y},\Omega)$. From this and \eqref{eqn:kcav-imm}, we deduce \eqref{eqn:cav2-compactness2}--\eqref{eqn:cav2-compactness4} by arguing as in Step 2 of the proof of Theorem~\ref{thm:cav}. To establish \eqref{eqn:cav2-compactness5},  we   also  argue   as in Step 2 of the proof of Theorem~\ref{thm:cav}.      

\textbf{Step 3 (Lower semicontinuity).} From \eqref{eqn:cav2-compactness4}, we immediately obtain
\begin{equation*}
	\int_{\imm(\boldsymbol{y},\Omega)}|D\boldsymbol{v}|^2\,\d\boldsymbol{\xi}\leq \liminf_{n \to \infty} \int_{\imm(\boldsymbol{y}_n,\Omega)} |D\boldsymbol{v}_n|^2\,\d\boldsymbol{\xi}.
\end{equation*}
The lower semicontinuity of the elastic energy is established by  exploiting the condition $p>N-1$ and weak continuity of minors \cite[Theorem 8.20, Part 4]{dacorogna}, and  by  applying Theorem \ref{thm:lscT}  as in Step 3 of the proof of Theorem \ref{thm:cav}. Eventually, the lower semicontinuity of the surface energy follows from \eqref{eqn:cav2-compactness1} by Theorem~\ref{thm:determinant-surface-energy},  so that \eqref{eqn:F-EL-cav2-lsc} is proved. 
\end{proof}

   The next example shows that  the a  priori lower bound on the cavities volume is necessary in order to exclude jumps of the limiting Eulerian map. 

\begin{example}[Cavities squashing onto a segment and relaxation]\label{ex:cavity-segment-relaxation}
Let $N=M=2$. We employ the notation in \eqref{eqn:q-ball-annulus}. Consider $\Omega$, $Q_n$, $\Sigma$, $\boldsymbol{y}_n$, $\boldsymbol{y}$, $\boldsymbol{v}_n$, and $\boldsymbol{v}$ as in Examples \ref{ex:cavity-segment} and \ref{ex:jump-not-lsc}. Then, we have
 $\imm(\boldsymbol{y}_n,\Omega)=B_1\setminus Q_n$ for all $n \in \N$, and $\imm(\boldsymbol{y},\Omega)=B_1$. Indeed,  $C_{\boldsymbol{y}}=\emptyset$ because $\imt(\boldsymbol{y},\boldsymbol{0})=\Sigma$ and $\mathscr{L}^2(\Sigma)=0$. Note that $\boldsymbol{v}_n\in W^{1,2}(\imm(\boldsymbol{y}_n,\Omega);\R^2)$ for all $n\in \N$, while  $\boldsymbol{v}\in SBV^2(\imm(\boldsymbol{y},\Omega);\R^2)$. In particular, \eqref{eqn:cav2-compactness1}--\eqref{eqn:cav2-compactness5} hold.  However, $\boldsymbol{v}\notin W^{1,2}(\imm(\boldsymbol{y},\Omega);\R^2)$ because $J_{\boldsymbol{v}}=\Sigma$.  
\end{example}

\subsection{Eulerian-Lagrangian energies for  deformations   with bounded variation}
In this subsection,  admissible deformations belong to the  class $\mathcal{Y}_{p}^{\rm frac}(\Omega)$ defined in \eqref{eqn:Y-frac}. Additionally, we impose an a priori $L^\infty$-bound on deformations. 
 Given $b>0$, we set
\begin{equation}
	\label{eqn:Y-frac-M}
	\mathcal{Y}_{p,b}^{\rm frac}(\Omega)\coloneqq \big\{ \boldsymbol{y}\in \mathcal{Y}_p^{\rm frac}(\Omega):  \hspace{4pt} \|\boldsymbol{y}\|_{L^\infty(\Omega;\R^N)}\leq b\big\}.
\end{equation}
Note that this restriction is standard practice when  working with $SBV$-functions in fracture mechanics \cite[Section 4]{ambrosio.fusco.pallara}. 
Given $\boldsymbol{y}\in \mathcal{Y}_{p}^{\rm frac}(\Omega)$ and $\boldsymbol{v}\in L^2(\img(\boldsymbol{y},\Omega); Z)$ almost everywhere approximately differentiable,  the functional in \eqref{eqn:F-EL} and \eqref{eqn:F-EL3} takes the form
\begin{equation}
	\label{eqn:F-EL-frac}
		\begin{split}
		 \mathcal{F}^{\rm frac}  (\boldsymbol{y},\boldsymbol{v})\coloneqq & \int_{\Omega} W(\nabla \boldsymbol{y},\boldsymbol{v}\circ \boldsymbol{y})\,\d\boldsymbol{x}+ \lambda_1   \mathcal{S}(\boldsymbol{y})+ \lambda_2 \per\big(\img(\boldsymbol{y},\Omega)\big)+\haus(J_{\boldsymbol{y}}) \\
		&+  \int_{\img(\boldsymbol{y},\Omega)} |\nabla \boldsymbol{v}|^2\,\d\boldsymbol{\xi}+ \lambda_3 \haus(J_{\boldsymbol{v}}).
		\end{split} 
\end{equation}

We study this functional  under the assumption that the extension of $\boldsymbol{v}$ to the whole space by zero belongs to the space $GSBV^2(\R^N;\R^M)$.

\begin{theorem}[Coercivity and lower semicontinuity: $\boldsymbol{SBV}$-deformations]
\label{thm:frac}
Let $b>0$, and $ \mathcal{F}^{\rm frac} $ be defined as in	\eqref{eqn:F-EL-frac} with $W$ satisfying {\rm (W1)}--{\rm (W3)},  $\lambda_1>0$, and $\lambda_2\geq\lambda_3>0$.   
Let $(\boldsymbol{y}_n)_n\subset \mathcal{Y}_{p,b}^{\rm frac}(\Omega)$ and let $(\boldsymbol{v}_n)_n$ be a sequence of maps $\boldsymbol{v}_n \in L^2(\img(\boldsymbol{y}_n,\Omega); Z  )$  with   $\widebar{\boldsymbol{v}}_n \in GSBV^2(\R^N;\R^M)$, where $\widebar{\boldsymbol{v}}_n$ denotes the extension of $\boldsymbol{v}_n$ to the whole space by zero. Suppose that 
\begin{equation}
	\label{eqn:frac-bound}
	\sup_{n \in \N}  \left \{  \|\boldsymbol{v}_n\|_{L^2(\img(\boldsymbol{y}_n,\Omega);\R^M)}+  \mathcal{F}^{\rm frac}(\boldsymbol{y}_n,\boldsymbol{v}_n) \right \}<+\infty.
\end{equation}
Eventually, assume that  the sequence $(\boldsymbol{v}_n \circ \boldsymbol{y}_n)_n$ is equi-integrable.   Then, there exist $\boldsymbol{y}\in \mathcal{Y}_{p,b}^{\rm frac}(\Omega)$ and $\boldsymbol{v}\in L^2(\img(\boldsymbol{y},\Omega);\R^M)$ with  $\boldsymbol{v}\circ \boldsymbol{y}\in L^1(\Omega;  \R^M )$ and $\widebar{\boldsymbol{v}}\in GSBV^2(\R^N;\R^M)$, where $\widebar{\boldsymbol{v}}$ denotes the extension  of $\boldsymbol{v}$ to the whole space by zero,  such that, up to subsequences, we have:
\begin{align}
	\label{eqn:frac-compactness0}
	\text{$\boldsymbol{y}_n \to \boldsymbol{y}$ in $L^1(\Omega;\R^N)$},  \quad 
	\text{$\nabla \boldsymbol{y}_n \wk \nabla \boldsymbol{y}$ }& \text{in $L^p(\Omega;\rnn)$}, \quad \text{$\det \nabla\boldsymbol{y}_n \wk \det \nabla\boldsymbol{y}$ }\text{in $L^1(\Omega)$}, \\
	\label{eqn:frac-compactness1}
	\text{$\chi_{\img(\boldsymbol{y}_n,\Omega)}\boldsymbol{v}_n \to \chi_{\img(\boldsymbol{y},\Omega)}\boldsymbol{v}$ }&\text{in $L^1(\R^N;\R^M)$,}	\\
	\label{eqn:frac-compactness2}
	\text{$\chi_{\img(\boldsymbol{y}_n,\Omega)}\boldsymbol{v}_n \wk \chi_{\img(\boldsymbol{y},\Omega)}\boldsymbol{v}$ }&\text{in $L^2(\R^N;\R^M)$,} \\
	\label{eqn:frac-compactness3}
	\text{$\chi_{\img(\boldsymbol{y}_n,\Omega)}\nabla\boldsymbol{v}_n \wk \chi_{\img(\boldsymbol{y},\Omega)}\nabla\boldsymbol{v}$ }&\text{in $L^2(\R^N;\R^{M \times N})$}\\
	\label{eqn:frac-compactness4}
	\text{$\boldsymbol{v}_n \circ \boldsymbol{y}_n \to \boldsymbol{v}\circ \boldsymbol{y}$ }&\text{in $L^1(\Omega;\R^M)$.}
\end{align}
Moreover, if $\boldsymbol{v}(\boldsymbol{\xi})\in Z$ for almost every $\boldsymbol{\xi}\in \img(\boldsymbol{y},\Omega)$, then we have 
	\begin{equation}
		\label{eqn:frac-lsc}
		\mathcal{F}^{\rm frac}(\boldsymbol{y},\boldsymbol{v})\leq \liminf_{n \to \infty} \mathcal{F}^{\rm frac}(\boldsymbol{y}_n,\boldsymbol{v}_n). 
	\end{equation}
\end{theorem}

\begin{remark}
	\label{rem:frac}
	\begin{enumerate}[(a)]
		\item The assumption $\widebar{\boldsymbol{v}}\in GSBV^2(\R^N;\R^M)$ is an additional regularity requirement, but the energy functional in \eqref{eqn:F-EL-frac} depends only of the values of $\boldsymbol{v}$.   Loosely speaking, this assumption means that $\haus(J_{\boldsymbol{v}})<+\infty$  and  $\boldsymbol{v}$ has finite Dirichlet energy on $\img(\boldsymbol{y},\Omega)$. 
		\item Since $SBV$-deformations   do not enjoy any kind of continuity, the topological degree is not available for  such   maps. Thus, we do not have any candidate open set for imposing the Sobolev regularity of Eulerian maps. Also, even if the Eulerian maps  do not  jump  inside the  geometric image along the sequence, jumps of the limiting Eulerian map  inside the geometric image of the limiting deformation cannot be excluded,  see  Example \ref{ex:jump-frac} below. 
		\item   Both conditions $\lambda_1>0$ and $\lambda_2\geq \lambda_3>0$ are crucial. As already noted, the first condition is needed to deduce the weak continuity of the Jacobian since this is not ensured by the boundedness of the perimeter of geometric images, see \cite[Section 1]{henao.moracorral.invertibility}. The second condition refers to an observation analogous to one made in Remark \ref{rem:cav}(c).   In the proof of Theorem \ref{thm:frac} we show that
		\begin{align}\label{eqn:GG}
			\haus(J_{\boldsymbol{v}})+\per\big(\img(\boldsymbol{y},\Omega)\big)&\leq \liminf_{n \to \infty} \left\{ \haus(J_{\boldsymbol{v}_n}) + \per\big(\img(\boldsymbol{y}_n,\Omega) \big) \right\},\\ \label{eqn:GG2}
			\per \big(\img(\boldsymbol{y},\Omega)\big)&\leq \liminf_{n \to \infty} \per\big(\img(\boldsymbol{y}_n,\Omega) \big).
		\end{align}
		While the lower semicontinuity of the perimeter of the geometric image is  standard,   the one of the functional $\boldsymbol{v}\mapsto \haus(J_{\boldsymbol{v}})$ does not  generally  hold,   see   Example \ref{ex:jump-not-contained-imt} and    Example \ref{ex:jump-frac} below. Multiplying the terms in \eqref{eqn:GG}--\eqref{eqn:GG2} by $\lambda_3>0$ and $\lambda_2-\lambda_3\geq 0$, respectively, the lower semicontinuity of the two surface terms in \eqref{eqn:F-EL-frac} follows. 
		\item Also here, as in  Remark \ref{rem:cav}(d), the assumption $\boldsymbol{v}(\boldsymbol{\xi})\in  Z $ for almost every $\boldsymbol{\xi}\in \img(\boldsymbol{y},\Omega)$ is trivially satisfied when $ Z  $ is closed.
		\item   If $\sup_{n \in \N}\|\boldsymbol{v}_n\|_{L^\infty(\img(\boldsymbol{y}_n,\Omega);\R^M)}<+\infty$, then $\widebar{\boldsymbol{v}}_n\in SBV^2(\R^N;\R^M)$ for every $n\in \N$ and also $\widebar{\boldsymbol{v}}\in SBV^2(\R^N;\R^M)$. Additionally, if each $\widebar{\boldsymbol{v}}_n$ is piecewise-constant, then so is $\widebar{\boldsymbol{v}}$. The claim is proved as in Remark \ref{rem:cav}(e)  by exploiting \eqref{eqn:PC}--\eqref{eqn:SBVcapLinfty}, \eqref{eqn:frac-compactness1}, and \eqref{eqn:frac-compactness3}. 
		\item  In contrast to Theorem \ref{thm:cav} and Theorem \ref{thm:cavk}, here   it is possible to relax the condition on the exponent $p$ from $p>N-1$ to $p \geq N-1$ at the price of additionally imposing in \eqref{eqn:growth} a superlinear growth in the cofactor.  To this end, in the proof given below  one has to apply \cite[Lemma 2]{henao.moracorral.invertibility} instead of \cite[Corollary~5.31]{ambrosio.fusco.pallara} to recover the weak continuity of the Jacobian cofactor.
		\item In principle, it would be possible to remove the a priori $L^\infty$-bound in \eqref{eqn:Y-frac-M} by considering deformations  $\boldsymbol{y}\in GSBV^p(\Omega;\R^N)$  and resorting   to the compactness result  in \cite{friedrich}. To do so, one should check that, starting from a sequence of injective deformations, the sequence of modifications in \cite[Theorem 3.1]{friedrich} can be constructed in such a way that the global invertibility is preserved.  
	\end{enumerate}
\end{remark}

\begin{proof}[Proof of Theorem \ref{thm:frac}]
The proof is divided into three steps.

\textbf{Step 1 (Compactness of deformations).}
By \eqref{eqn:growth} and \eqref{eqn:frac-bound}, we have
\begin{equation*}
	\sup_{n \in \N} \left\{ \| \nabla  \boldsymbol{y}_n\|_{L^p(\Omega;\rnn)}+\|\gamma(\det  \nabla   \boldsymbol{y}_n)\|_{L^1(\Omega)}+\haus(J_{\boldsymbol{y}_n}) \right\}<+\infty. 
\end{equation*}
Given the  $L^\infty$-bound in \eqref{eqn:Y-frac-M}, thanks to Theorem \ref{thm:ambrosio-compactness}, the first two convergences in \eqref{eqn:frac-compactness0} hold for some $\boldsymbol{y}\in SBV^p(\Omega;\R^N)$ along a not relabeled subsequence. In view of the second condition in \eqref{eqn:growth-gamma} and the De la Vallée Poussin's criterion, there exists $h \in L^1(\Omega)$ such that, up to subsequences, $\det \nabla \boldsymbol{y}_n \wk h$ in $L^1(\Omega)$. As each $\boldsymbol{y}_n$ satisfies $\det D \boldsymbol{y}_n>0$ almost everywhere, there holds $h \geq 0$. With the same contradiction argument based on \eqref{eqn:growth-gamma} employed in Step 1 of the proof of Theorem \ref{thm:cav}, we conclude that $h>0$ almost everywhere.   As \eqref{eqn:frac-bound} gives a uniform bound on $(\mathcal{S}(\boldsymbol{y}_n))_n$, by applying  \cite[Corollary 5.31]{ambrosio.fusco.pallara} and  Theorem \ref{thm:determinant-surface-energy}  we deduce that $\boldsymbol{y}$ is almost everywhere injective and  $h\cong \det \nabla \boldsymbol{y}$.   This completes the proof of  \eqref{eqn:frac-compactness0} and shows that $\boldsymbol{y}\in \mathcal{Y}_{p,b}^{\rm frac}(\Omega)$,  thanks to the lower semicontinuity of the $L^\infty$-norm.   Also,  claim (i) of Proposition \ref{prop:approx-diff} yields  
\begin{equation}
	\label{eqn:SBV-img}
	\text{$\chi_{\img(\boldsymbol{y}_n,\Omega)}\to \chi_{\img(\boldsymbol{y},\Omega)}$ in $L^1(\R^N)$.}	
\end{equation}

\textbf{Step 2 (Compactness of Eulerian maps).} From  \eqref{eqn:frac-bound}, thanks to Lemma \ref{lem:approx-diff-extension}, we  get that 
\begin{equation*}
	\begin{split}
 \sup_{n \in \N} \left \{ \|\widebar{\boldsymbol{v}}_n\|_{L^2(\R^N;\R^M)}+\|\nabla \widebar{\boldsymbol{v}}_n\|_{L^2(\R^N;\R^{M \times N})} \right \}= \sup_{n \in \N} \left  \{ \|{\boldsymbol{v}}_n\|_{L^2(\img(\boldsymbol{y}_n,\Omega);\R^M)}+\|\nabla {\boldsymbol{v}}_n\|_{L^2(\img(\boldsymbol{y}_n,\Omega);\R^{M \times N})} \right \}
	\end{split}
\end{equation*}
 is finite,   while Lemma \ref{lem:jump-extension} yields
\begin{equation*}
	\begin{split}
		\sup_{n \in \N} \haus(J_{\overline{\boldsymbol{v}}_n})&= \sup_{n \in \N} \left\{ \haus\big(J_{\boldsymbol{v}_n}) + \haus(J_{\overline{\boldsymbol{v}}_n}\cap \partial^* \img(\boldsymbol{y}_n,\Omega)\big) \right\} 	\\
		&	\leq \sup_{n \in \N} \left\{ \haus(J_{\boldsymbol{v}_n}) + \per\big(\img(\boldsymbol{y}_n,\Omega)\big) \right\}<+\infty,
	\end{split}  
\end{equation*}
 again  thanks to \eqref{eqn:frac-bound}.  
By Theorem \ref{thm:ambrosio-compactness}, we find ${\boldsymbol{w}} \in GSBV^2(\R^N; \R^M )$ and a not relabeled subsequence for which
\begin{equation}
	\label{eqn:vbar}
	\text{$\widebar{\boldsymbol{v}}_n \to  {\boldsymbol{w}}$ a.e.\ in $\R^N$,} \qquad \text{$\widebar{\boldsymbol{v}}_n \wk  {\boldsymbol{w}}$ in $L^2(\R^N;  \R^M  )$,} \qquad \text{$\nabla \widebar{\boldsymbol{v}}_n \wk \nabla  {\boldsymbol{w}}$ in $L^2(\R^N; \R^{M \times N} )$.}
\end{equation}
From \eqref{eqn:SBV-img}, arguing as in Step 2 of the proof of Theorem~\ref{thm:cav}, we deduce that $ {\boldsymbol{w}} \cong \boldsymbol{0}$ and $\nabla  {\boldsymbol{w}} \cong \boldsymbol{O}$ in $\R^N \setminus \img(\boldsymbol{y},\Omega)$. Thus, letting  $\boldsymbol{v}\coloneqq  {\boldsymbol{w}} \restr{\img(\boldsymbol{y},\Omega)}  \in  L^2(\img(\boldsymbol{y},\Omega);\R^M)$ we have $ {\boldsymbol{w}} \cong \chi_{\img(\boldsymbol{y},\Omega)}\boldsymbol{v}$.   Thus, the last two convergences in \eqref{eqn:vbar} prove  \eqref{eqn:frac-compactness2}--\eqref{eqn:frac-compactness3}, while the first one yields \eqref{eqn:frac-compactness1} by Vitali's convergence theorem. 
Eventually, proceeding once again as in Step 2 of the proof of Theorem \ref{thm:cav}, we establish \eqref{eqn:frac-compactness4}. 

\textbf{Step 3 (Lower semicontinuity).} Applying \cite[Corollary 5.31]{ambrosio.fusco.pallara} and recalling the third convergence in \eqref{eqn:frac-compactness0}, we obtain
\begin{equation*}
 \text{$\adj_r \nabla \boldsymbol{y}_n \wk \adj_r \nabla \boldsymbol{y}$ in $L^1\left(\Omega;\R^{\binom{N}{r}\times \binom{N}{r}}  \right)$ \quad for all 	  $r=1,\dots,N$.}
\end{equation*}
As in  Step 3 of  the proof of Theorem \ref{thm:cav}, the assumption $\boldsymbol{v}(\boldsymbol{\xi})\in  Z$ for almost every $\boldsymbol{\xi}\in \img(\boldsymbol{y},\Omega)$ yields $\boldsymbol{v}(\boldsymbol{y}(\boldsymbol{x}))\in  Z$ for almost every $\boldsymbol{x}\in \Omega$ thanks to Lusin's condition (N${}^{-1}$),  see Remark \ref{rem:federer}(b). Hence, assumptions \eqref{eqn:polyconvex1}--\eqref{eqn:polyconvex2} and \eqref{eqn:frac-compactness4} yield the lower semicontinuity of the elastic energy by  Theorem \ref{thm:lscT}.

 By Theorem \ref{thm:ambrosio-compactness},  we have
\begin{equation*}
	\haus(J_{\boldsymbol{y}}) \leq \liminf_{n \to \infty} \haus(J_{\boldsymbol{y}_n}).
\end{equation*}
The lower semicontinuity of the surface term  $\mathcal{S}$  follows from Theorem \ref{thm:determinant-surface-energy}, while the one of the perimeter term is easily deduced from \eqref{eqn:SBV-img}.
Also,  \eqref{eqn:frac-compactness3} immediately yields
\begin{equation*}
	\int_{\img(\boldsymbol{y},\Omega)}|\nabla \boldsymbol{v}|^2\,\d\boldsymbol{\xi}\leq \liminf_{n \to \infty} \int_{\img(\boldsymbol{y}_n,\Omega)}|\nabla \boldsymbol{v}_n|^2\,\d\boldsymbol{\xi}.
\end{equation*}
To establish the lower semicontinuity of the remaining energy terms, we proceed similarly to Step 3 of the proof of Theorem \ref{thm:cav}. Let $\boldsymbol{\nu}_{\img(\boldsymbol{y},\Omega)}\colon \partial^*\,\img(\boldsymbol{y},\Omega)\to S$ be the outer unit normal   which is well defined since $\img(\boldsymbol{y},\Omega)$ has finite perimeter.  The lateral traces of $\widebar{\boldsymbol{v}}$ at $\boldsymbol{\xi}_0\in \partial^*\,\img(\boldsymbol{y},\Omega)$ are defined  similarly as in \eqref{XXX}, and we observe  that  $\widebar{\boldsymbol{v}}^-(\boldsymbol{\xi}_0)=\boldsymbol{0}$, so that the jump of $\widebar{\boldsymbol{v}}$ across $\partial^*\,\img(\boldsymbol{y},\Omega)$ equals $\left[ \widebar{\boldsymbol{v}} \right]=\widebar{\boldsymbol{v}}^+$. As the map $\widebar{\boldsymbol{v}}^+\colon \partial^*\,\img(\boldsymbol{y},\Omega)\to \R^M$ is  Borel, the sets  $ \left\{ \boldsymbol{\xi}\in \partial^*\,\img(\boldsymbol{y},\Omega): \hspace{4pt} \widebar{\boldsymbol{v}}^+(\boldsymbol{\xi})=\boldsymbol{b} \right\}$ are Borel for all $\boldsymbol{b}\in \R^M$. Since $\haus(\partial^* \img(\boldsymbol{y},\Omega))<+\infty$, we have that for all $\boldsymbol{b}\in \R^M$, except  for   an at most countable number,  
\begin{equation*}
	\haus \left( \left\{ \boldsymbol{\xi}\in \partial^*\,\img(\boldsymbol{y},\Omega): \hspace{4pt} \widebar{\boldsymbol{v}}^+(\boldsymbol{\xi})=\boldsymbol{b} \right\} \right)=0.
\end{equation*}
Choose $\boldsymbol{b}\in \R^M$ satisfying the previous condition and define $\widebar{\boldsymbol{v}}^{\boldsymbol{b}}\colon \R^N \to \R^M$ by setting
\begin{equation}\label{YYY}
	\widebar{\boldsymbol{v}}^{\boldsymbol{b}}(\boldsymbol{\xi})\coloneqq \begin{cases}
		\boldsymbol{v}(\boldsymbol{\xi}) & \text{if $\boldsymbol{\xi}\in \img(\boldsymbol{y},\Omega)$,}\\
		\boldsymbol{b} & \text{if $\boldsymbol{\xi}\in \R^N \setminus \img(\boldsymbol{y},\Omega)$.}
	\end{cases}
\end{equation}
As $\widebar{\boldsymbol{v}}\in GSBV^2(\R^N;\R^M)$ and $\img(\boldsymbol{y},\Omega)$ has finite perimeter, we have $\widebar{\boldsymbol{v}}^{\boldsymbol{b}}\in GSBV^2(\R^N;\R^M)$. Also, by Lemma \ref{lem:jump-extension}, there holds
\begin{equation}
	\label{eqn:frac-Jv}
	\haus(J_{\overline{\boldsymbol{v}}^{\boldsymbol{b}}})= \haus(J_{\boldsymbol{v}})+\per\big(\img(\boldsymbol{y},\Omega)\big).
\end{equation}
Analogously, we find a sequence $(\boldsymbol{b}_n)_n\subset \R^M$  with $\boldsymbol{b}_n \to \boldsymbol{b}$  such that the maps $ (\overline{\boldsymbol{v}}_n^{\boldsymbol{b}_n})_n  \subset GSBV^2(\R^N;\R^M)$ defined   analogously to \eqref{YYY} (with $\boldsymbol{v}_n$, $\boldsymbol{b}_n$, and $\img(\boldsymbol{y}_n,\Omega)$  in place of $\boldsymbol{v}$, $\boldsymbol{b}$, and $\img(\boldsymbol{y},\Omega)$)   
 satisfy
 \begin{equation}
 	\label{eqn:frac-Jvn}
 	\haus(J_{\overline{\boldsymbol{v}}_n^{\boldsymbol{b}_n}})= \haus(J_{\boldsymbol{v}_n})+\per\big(\img(\boldsymbol{y}_n,\Omega)\big)
 \end{equation}
for every $n \in \N$. In this way, from \eqref{eqn:frac-compactness1} and \eqref{eqn:SBV-img}, we obtain 
\begin{equation*}
	\text{$\widebar{\boldsymbol{v}}_n^{\boldsymbol{b}_n}\to \widebar{\boldsymbol{v}}^{\boldsymbol{b}}$ a.e.\ in $\R^N$.}
\end{equation*}
Thus, in view of \eqref{eqn:frac-compactness3} and Lemma \ref{lem:approx-diff-extension},   Theorem \ref{thm:ambrosio-compactness}  yields
\begin{equation*}
	\haus(J_{\overline{\boldsymbol{v}}^{\boldsymbol{b}}})\leq \liminf_{n \to \infty} \haus(J_{\overline{\boldsymbol{v}}_n^{\boldsymbol{b}_n}}).
\end{equation*}
Given  \eqref{eqn:frac-Jv}--\eqref{eqn:frac-Jvn}, this entails 
\begin{equation*}
	\haus(J_{\boldsymbol{v}})+\per\big(\img(\boldsymbol{y},\Omega) \big)\leq \liminf_{n \to \infty} \left\{ \haus(J_{\boldsymbol{v}_n})+\per\big(\img(\boldsymbol{y}_n,\Omega)\big) \right\}
\end{equation*}
which concludes the proof  of \eqref{eqn:frac-lsc} in view of the condition $\lambda_2\geq \lambda_3$ and \eqref{eqn:GG2}.  
\end{proof}

The next example illustrates the observations made in Remark \ref{rem:frac}, items (b) and (c).  

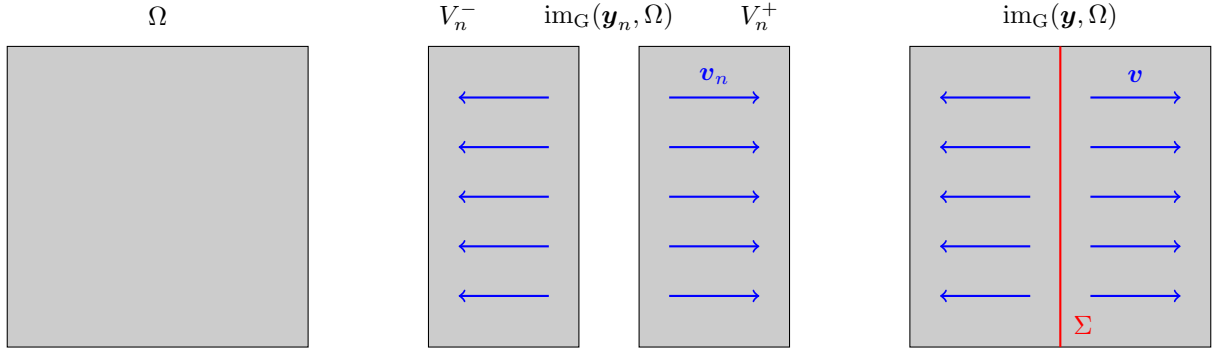
\begin{figure}
	\begin{tikzpicture}[scale=2]
		\node at (-3,1.2) {$\Omega$};
		\draw[fill=black!20] (-4,-1) rectangle (-2,1);
		\node at (0,1.2) {$\img(\boldsymbol{y}_n,\Omega)$};
		\node at (-1,1.2) {$V_n^-$}; \node at (1,1.2) {$V_n^+$};
		\draw[fill=black!20] (-1.2,1) rectangle (-.2,-1);
		\draw[fill=black!20] (1.2,1) rectangle (.2,-1);
		\node[blue] at (.7,.8) {$\boldsymbol{v}_n$};
		\node at (3,1.2) {$\img(\boldsymbol{y},\Omega)$};
		\draw[fill=black!20] (4,-1) rectangle (2,1);
		\foreach \a in {0,.33,.66,-.33,-.66}{
		\draw[blue,thick,->] (.4,\a)--(1,\a);
		\draw[blue,thick,->] (-.4,\a)--(-1,\a);
		\draw[blue,thick,->] (3.2,\a)--(3.8,\a);
		\draw[blue,thick,->] (2.8,\a)--(2.2,\a);
		};
		\draw[red,thick] (3,1)--(3,-1);
		\node[blue] at (3.5,.8) {$\boldsymbol{v}$};
		\node[red] at (3.15,-.85) {$\Sigma$};
	\end{tikzpicture}
	\caption{The deformations in Example \ref{ex:jump-frac} and the corresponding Eulerian fields.}
	\label{fig:jump-frac}
\end{figure}

\begin{example}[Jump of Eulerian maps and fractures] \label{ex:jump-frac}
Let $N=M=2$.  We employ the notation in \eqref{eqn:q-ball-annulus}.   Set $\Omega \coloneqq   B_\infty$ and for all $n\in \N$ define $\boldsymbol{y}_n\in SBV(\Omega;\R^2)$ by setting 
\begin{equation*}
	\boldsymbol{y}_n(\boldsymbol{x})\coloneqq \begin{cases} 
		\boldsymbol{x}+\frac{1}{n}\boldsymbol{e}_1 & \text{if $x_1>0$,}\\ \boldsymbol{x} &  \text{if $x_1=0$,}\\ 
		\boldsymbol{x}-\frac{1}{n}\boldsymbol{e}_1 & \text{if $x_1<0$.}
	\end{cases}
\end{equation*}	
Then, $\boldsymbol{y}_n \in \mathcal{Y}_{p,b}^{\rm frac}(\Omega)$ with $b=2$ and $\img(\boldsymbol{y}_n,\Omega)= V^-_n \cup V^+_n$, where  (see Figure~\ref{fig:jump-frac}):
\begin{equation*}
	 V^+_n\coloneqq \left (\frac{1}{n},1+\frac{1}{n}\right )\times (-1,1), \qquad V^-_n\coloneqq \left (-1-\frac{1}{n},-\frac{1}{n}\right) \times (-1,1).
\end{equation*}
Also, $\per(\img(\boldsymbol{y}_n,\Omega))=12$ for all $n\in \N$.
Clearly, \eqref{eqn:frac-compactness0} holds with  $\boldsymbol{y}\coloneqq\boldsymbol{id}\restr{\Omega}$. In particular, $\per(\img(\boldsymbol{y},\Omega))=8$.
For all $n \in \N$,  define $\boldsymbol{v}_n \in L^2(\img(\boldsymbol{y}_n,\Omega);\R^N)$ by setting
\begin{equation*}
	\boldsymbol{v}_n(\boldsymbol{\xi})\coloneqq \begin{cases}
		\boldsymbol{e}_1 & \text{if $\boldsymbol{\xi}\in V^+_n$,}\\ -\boldsymbol{e}_1 & \text{if $\boldsymbol{\xi}\in V^-_n$.}
	\end{cases}
\end{equation*}
Similarly, let $\boldsymbol{v}\in L^2(\img(\boldsymbol{y},\Omega);\R^2)$ be given by
\begin{equation*}
	\boldsymbol{v}(\boldsymbol{\xi})\coloneqq \begin{cases}
		\boldsymbol{e}_1 & \text{if $\xi_1>0$,}\\ \boldsymbol{\xi} & \text{if $\xi_1=0$,} \\ -\boldsymbol{e}_1 & \text{if $\xi_1<0$.}
	\end{cases}
\end{equation*} 
 Clearly, \eqref{eqn:frac-compactness1}--\eqref{eqn:frac-compactness4} hold.  
Also,  $J_{\boldsymbol{v}_n}=\emptyset$ for all $n \in \N$, while $J_{\boldsymbol{v}}=\Sigma$, where $\Sigma\coloneqq \{0\}\times (-1,1)$. Therefore,
\begin{equation*}
	\haus(J_{\boldsymbol{v}})=2>0=\liminf_{n \to \infty} \haus(J_{\boldsymbol{v}_n}),
\end{equation*}
but
\begin{equation*}
	\haus(J_{\boldsymbol{v}})+\per\big(\img(\boldsymbol{y},\Omega)\big)=10<12=\liminf_{n \to \infty} \left\{ \mathscr{H}^{ N-1}(J_{\boldsymbol{v}_n})+\per\big(\img(\boldsymbol{y}_n,\Omega) \big)\right\}.
\end{equation*}
\end{example}
\section{Applications} \label{sec:appl}

In this section, we apply the results of Section \ref{sec:EL} to establish the existence of minimizers for variational problems featuring mixed Eulerian-Lagrangian formulations.

Henceforth,  $\Omega \subset \widetilde{\Omega} \subset \R^N$ are  bounded Lipschitz domains  with $\haus ( \partial \Omega \cap  \widetilde{\Omega})>0$   and $p>N-1$. We consider $\boldsymbol{d}\in L^1(\widetilde{\Omega};\R^N)$  as boundary datum,  which will be imposed on $\widetilde{\Omega} \setminus \overline{\Omega}$ for $SBV$-deformations.  To cover also the case of Sobolev deformations, for which a weakly continuous  trace operator is defined, we additionally assume that $\boldsymbol{d}\restr{\Gamma}\in L^1(\Gamma;\R^N)$, where  $\Gamma\coloneqq \partial \Omega \cap \widetilde{\Omega}$.    The trace of a deformation $\boldsymbol{y}\in \mathcal{Y}^{\rm cav}_p(\Omega)$ on the boundary of $\Omega$ will be denoted by $\mathrm{tr}_{\partial \Omega}(\boldsymbol{y})$.  In the rest of the section, whenever imposing Dirichlet boundary conditions on the deformations, we tacitly assume $\boldsymbol{d}$ to be sufficiently regular so that the resulting class of admissible deformations is nonempty.  

The elastic energy densities that we consider involve a  function $\Phi \colon \rnn_+ \to [0,+\infty]$ satisfying properties analogous to the ones in Section \ref{sec:EL}, namely:
\begin{enumerate}[($\Phi$1)]
	\item \textbf{Continuity:} The function $\Phi \colon \rnn_+ \to [0,+\infty]$ is continuous;
	\item \textbf{Coercivity:} There exists a constant $C>0$ and a Borel function $\gamma \colon (0,+\infty)\to [0,+\infty]$ satisfying
	\begin{equation}
		\label{eqn:growth-gamma-Phi}
		\lim_{h \to 0^+} \gamma(h)=\lim_{h \to +\infty} \frac{\gamma(h)}{h}=+\infty
	\end{equation}
	such that 
	\begin{equation}
		\label{eqn:growth-Phi}
	\Phi(\boldsymbol{G})\geq C |\boldsymbol{G}|^p+\gamma(\det \boldsymbol{G}) \quad \text{ for all } 	  \boldsymbol{G}\in \rnn_+;
	\end{equation}
	\item \textbf{Polyconvexity:} The function $\Phi$ is polyconvex, that is, there exists a   convex function $\widehat{\Phi}\colon \prod_{r=1}^{N-1} \R^{\binom{N}{r}\times \binom{N}{r}}\times (0,+\infty)\to [0,+\infty]$   such that
	\begin{equation}
		\label{eqn:polyconvex-Phi}
\Phi(\boldsymbol{G})=\widehat{\Phi}(\adj_1\boldsymbol{G}, \dots,\adj_{N-1}\boldsymbol{G},\adj_N\boldsymbol{G})  \quad \text{ for all } 	  \boldsymbol{G}\in \rnn_+.
	\end{equation}
\end{enumerate}

In the following subsections, we discuss three applications concerning nematic elastomers, phase transitions, and magnetoelasticity, respectively.
 We will refer to   the classes $\mathcal{Y}(\Omega)$, $\mathcal{Y}_p^{\rm cav}(\Omega)$,  $\mathcal{Y}_{p,\kappa}^{\rm cav}(\Omega)$, and $\mathcal{Y}_{p,b}^{\rm frac}({\Omega})$ defined in \eqref{eqn:Y}, \eqref{eqn:Y-cav}, \eqref{eqn:Y-cav-kappa}, and \eqref{eqn:Y-frac-M},  respectively.   

\subsection{Nematic elastomers} \label{subsec:nem}
For the first application, we focus on the Oseen-Frank energy model for nematic elastomers   \cite{agostiniani.desimone,barchiesi.desimone,barchiesi.henao.moracorral,desimone.teresi,henao.stroffolini,warner.tarentjev}. In the theory of liquid crystals, this corresponds to the so-called one-constant approximation \cite{ball.mathlc}.  The energy functional depends on two variables:   the deformation $\boldsymbol{y}\in \mathcal{Y}(\Omega)$  and the nematic director $\boldsymbol{n}\colon \img(\boldsymbol{y},\Omega)\to  S$,   where we recall \eqref{eqn:B-S}.  The latter describes the local average orientation of the constituent molecules of the liquid crystal and, hence, complies with the constraint $|\boldsymbol{n}|\cong1$ in $\img(\boldsymbol{y},\Omega)$. Also, the energy should not depend on the orientation of $\boldsymbol{n}$ but only on its direction, that is, the energy should not change by replacing $\boldsymbol{n}$ with $-\boldsymbol{n}$.

As in \cite{barchiesi.desimone,barchiesi.henao.moracorral,henao.stroffolini}, we take an elastic density $W^{\rm nem}\colon \rnn_+ \times  S  \to [0,+\infty]$ of the form
\begin{equation}
	\label{eqn:Wnem}
	W^{\rm nem}(\boldsymbol{F},\boldsymbol{z})\coloneqq \Phi \left( \boldsymbol{A}^{-1}\left(\boldsymbol{z}\right)\boldsymbol{F} \right), \qquad \boldsymbol{A}(\boldsymbol{z})\coloneqq \alpha \boldsymbol{z}\otimes \boldsymbol{z} + \alpha^{-\frac{1}{N-1}} (\boldsymbol{I}-\boldsymbol{z}\otimes \boldsymbol{z}),
\end{equation}
where $\alpha>0$ and $\Phi\colon \rnn_+ \to [0,+\infty]$ satisfies ($\Phi$1)--($\Phi$3). 

The existence theories established in \cite{barchiesi.desimone,barchiesi.henao.moracorral,henao.stroffolini} account for purely elastic deformations and nematic variables that are not allowed to jump. For rigid bodies, models with nematic directors in $SBV$  have also been investigated in the literature \cite{ball.mathlc,ball.bedford,bedford}. The next proposition extends all these results to elastomers possibly undergoing material failure. Both cases of nematic directors with Sobolev or $SBV$-regularity are discussed.

\begin{proposition}[Existence for nematic elastomers] \label{prop:nem}
Let $W^{\rm nem}$ be  as in \eqref{eqn:Wnem} with $\alpha>0$ and $\Phi$  satisfying  {\rm ($ \Phi$1)}--{\rm ($\Phi$3)}.  
\begin{enumerate}[(i)]
\item  For $\lambda_1\geq \lambda_3>0$,  the functional
\begin{equation*}
	(\boldsymbol{y},\boldsymbol{n})\mapsto \int_{\Omega} W^{\rm nem}(D\boldsymbol{y},\boldsymbol{n}\circ \boldsymbol{y})\,\d\boldsymbol{x}+ \lambda_1  \mathcal{S}(\boldsymbol{y})+\int_{\img(\boldsymbol{y},\Omega)}|\nabla \boldsymbol{n}|^2\,\d\boldsymbol{\xi}+ \lambda_3  \haus(J_{\boldsymbol{n}}\cap \imt(\boldsymbol{y},\Omega))
\end{equation*}
admits minimizers within the class
\begin{equation*}
	\left\{ (\boldsymbol{y},\boldsymbol{n}): \hspace{4pt}\boldsymbol{y}\in \mathcal{Y}_p^{\rm cav}(\Omega), \hspace{4pt} \text{$\tr_{\partial \Omega}(\boldsymbol{y})\simeq\boldsymbol{d}$ on $\Gamma$}, \hspace{4pt} \boldsymbol{n}\in L^2(\img(\boldsymbol{y},\Omega); S ), \hspace{4pt} \widetilde{\boldsymbol{n}}\in SBV^2(\imt(\boldsymbol{y},\Omega);\R^N)  \right\},
\end{equation*}
where  $\widetilde{\boldsymbol{n}}$ denotes the extension of $\boldsymbol{n}$ to $\imt(\boldsymbol{y},\Omega)$ by zero.
\item For every $\kappa>0$  and $\lambda_1>0$,  the functional
\begin{equation*}
	(\boldsymbol{y},\boldsymbol{n})\mapsto \int_{\Omega} W^{\rm nem}(D\boldsymbol{y},\boldsymbol{n}\circ \boldsymbol{y})\,\d\boldsymbol{x}+ \lambda_1 \mathcal{S}(\boldsymbol{y})+\int_{\imm(\boldsymbol{y},\Omega)}|D \boldsymbol{n}|^2\,\d\boldsymbol{\xi}
\end{equation*}
admits minimizers within the class
\begin{equation*}
	\left\{ (\boldsymbol{y},\boldsymbol{n}): \hspace{4pt}\boldsymbol{y}\in \mathcal{Y}_{p,\kappa}^{\rm cav}(\Omega), \hspace{4pt} \text{$\tr_{\partial \Omega}(\boldsymbol{y})\simeq\boldsymbol{d}$ on $\Gamma$}, \hspace{4pt} \boldsymbol{n}\in W^{1,2}(\imm(\boldsymbol{y},\Omega); S )  \right\}.
\end{equation*}
\item For every $b>0$,  $\lambda_1>0$, and $\lambda_2\geq \lambda_3>0$,  the functional
\begin{equation*}
	\begin{split}
		(\boldsymbol{y},\boldsymbol{n})\mapsto &\int_{\Omega} W^{\rm nem}(\nabla\boldsymbol{y},\boldsymbol{n}\circ \boldsymbol{y})\,\d\boldsymbol{x}+  \lambda_1 \mathcal{S}(\boldsymbol{y}_n)+\lambda_2  \per\big(\img(\boldsymbol{y},\Omega)\big)+\haus(J_{\boldsymbol{y}})\\
		&+\int_{\img(\boldsymbol{y},\Omega)}|\nabla \boldsymbol{n}|^2\,\d\boldsymbol{\xi}+\lambda_3  \haus(J_{\boldsymbol{n}})
	\end{split}
\end{equation*}
admits minimizers within the class
\begin{equation*}
	\left\{ (\boldsymbol{y},\boldsymbol{n}): \hspace{4pt}\boldsymbol{y}\in \mathcal{Y}_{p,b}^{\rm frac}(\widetilde{\Omega}), \hspace{4pt} \text{$\boldsymbol{y}\cong\boldsymbol{d}$ on $\widetilde{\Omega}\setminus \closure{\Omega}$}, \hspace{4pt} \boldsymbol{n}\in L^2(\img(\boldsymbol{y},\Omega); S ), \hspace{4pt} \widebar{\boldsymbol{n}}\in SBV^2(\R^N;\R^N)  \right\},
\end{equation*}
      where $\widebar{\boldsymbol{n}}$ denotes the extension of $\boldsymbol{n}$ to the whole space by zero. 
\end{enumerate}
In  (i) and (iii), one can also restrict to  $\widetilde{\boldsymbol{n}}\in PC(\imt(\boldsymbol{y},\Omega);\R^N)$ and $\widebar{\boldsymbol{n}}\in PC(\R^N;\R^N)$, respectively.
\end{proposition}

\begin{remark}\label{rem:nem}
\begin{enumerate}[(a)]
	\item Note that,  in presence of cavitations,   $\widetilde{\boldsymbol{n}}\in W^{1,2}(\imt(\boldsymbol{y},\Omega);\R^N)$ and $\boldsymbol{n}\in L^2(\img(\boldsymbol{y},\Omega); S )$ are incompatible assumptions  as $|\widetilde{\boldsymbol{n}}| = 1$ on $\img(\boldsymbol{y},\Omega)$ and $\widetilde{\boldsymbol{n}} = \boldsymbol{0}$ else.  Thus,  the restriction to deformations in $\mathcal{Y}^{\rm cav}_{p,\kappa}(\Omega)$ as in (ii) seems necessary if one wants to exclude jumps of $\boldsymbol{n}$.
	\item  Applied loads  can also be treated.  Precisely, let 
	$\boldsymbol{f}\in L^{q}(\Omega;\R^N)$   and  $\boldsymbol{g}\in L^1(\R^N;\R^N)$ represent applied  forces and  external electric fields, respectively. Their work  is accounted by the functional
	\begin{equation*}
		\mathcal{L}(\boldsymbol{y},\boldsymbol{n})\coloneqq \int_\Omega \boldsymbol{f}\cdot \boldsymbol{y}\,\d\boldsymbol{x}+\int_{\img(\boldsymbol{y},\Omega)} \boldsymbol{g}\cdot \boldsymbol{n}\,\d\boldsymbol{\xi}.
	\end{equation*}
	In (i)--(ii), we take   $q  =  (Np/(N-p))'$ for   $p<N$,  $q>1$ for $p=N$, and $q=1$ for $p>N$   owing to the Sobolev embedding, while in (iii) we  assume $q=1$ given the $L^\infty$-bound in \eqref{eqn:Y-frac-M}.  Under such integrability assumptions, the functional $\mathcal{L}$  is  continuous with respect to the relevant topology.  
	The total energy equals $\mathcal{F}^{\rm nem}-\mathcal{L}$, where  $\mathcal{F}^{\rm nem}$ denotes one of the functionals in Proposition \ref{prop:nem}. This  result can be easily adapted to establish the existence of minimizers  of  the total energy. Indeed,  its coercivity  can be easily checked by means of a standard application of Young and H\"{o}lder inequalities,   see e.g., \cite[Lemma~3.4]{bresciani.thesis}.
\end{enumerate}
\end{remark}

\begin{proof}
All energies  in the statement   are given by the restriction of the functional defined  in \eqref{eqn:F-EL}--\eqref{eqn:F-EL2}  (with  $W^{\rm nem}$ in place of $W$) to the different classes of admissible states under consideration. 
Claims (i), (ii), and (iii) are proved by applying Theorem \ref{thm:cav}, Theorem \ref{thm:cavk}, and  Theorem \ref{thm:frac}, respectively, to minimizing sequences  in  the various classes.  Therefore, we only have to show that the common assumptions of the three theorems are fulfilled. 

First, we check that $W^{\rm nem}$ satisfies (W1)--(W3). The first property being clear, we look at the second one.  
By  direct computation,  using $(\boldsymbol{z}\otimes\boldsymbol{z})\boldsymbol{z} = \boldsymbol{z}$, we obtain  
\begin{equation*}
 \boldsymbol{A}^{-1}(\boldsymbol{z})=\alpha^{-1}\boldsymbol{z}\otimes\boldsymbol{z}+\alpha^{\frac{1}{N-1}}(\boldsymbol{I}-\boldsymbol{z}\otimes\boldsymbol{z}) \quad  \text{for all $\boldsymbol{z}\in S$.} 
\end{equation*}
Considering an orthonormal basis  of $\R^N$ containing $\boldsymbol{z}$ and using the identities $(\boldsymbol{z}\otimes\boldsymbol{z})\boldsymbol{z} = \boldsymbol{z}$ and $(\boldsymbol{z}\otimes\boldsymbol{z})\boldsymbol{\hat{z}} = \boldsymbol{0}$ for all $\boldsymbol{\hat{z}}  \in S$ orthogonal to $\boldsymbol{z}$,  we get
\begin{equation*}
	\det \boldsymbol{A}^{-1}(\boldsymbol{z})=1 \quad \text{for all $\boldsymbol{z}\in S$.}
\end{equation*}
Thanks to \eqref{eqn:growth-Phi}, observing that $|\boldsymbol{A}(\boldsymbol{z})|\leq   C(N,\alpha)$ for all $\boldsymbol{z}\in S$,  we estimate 
\begin{equation}
	\label{eqn:Wnem-coe}
	\begin{split}
 W^{\rm nem}(\boldsymbol{F},\boldsymbol{z})&  =  \Phi(\boldsymbol{A}^{-1}(\boldsymbol{z})\boldsymbol{F})\geq C|\boldsymbol{A}^{-1}(\boldsymbol{z})\boldsymbol{F}|^p+\gamma(\det \boldsymbol{A}^{-1}(\boldsymbol{z})\det \boldsymbol{F}) \geq C(N,p,\alpha) |\boldsymbol{F}|^p+\gamma(\det \boldsymbol{F})
	\end{split}
\end{equation}
 for all 	$ \boldsymbol{F}\in \rnn_+$ and for all $\boldsymbol{z}\in S$.
Thus, $W^{\rm nem}$ satisfies (W2). Assumption (W3) is checked as in \cite[Lemma 8.1]{barchiesi.henao.moracorral},  we include the proof here for the sake of completeness.   Applying   \cite[Proposition~5.66]{dacorogna}, we have
\begin{equation*}
\adj_r(\boldsymbol{A}^{-1}(\boldsymbol{z})\boldsymbol{F})=\adj_r(\boldsymbol{A}^{-1}(\boldsymbol{z}))\,\adj_r\boldsymbol{F} \quad \text{for all $ r=1,\dots,N$},
\end{equation*} 
thus  \eqref{eqn:polyconvex-Phi} yields
\begin{equation*}
			W^{\rm nem}(\boldsymbol{F},\boldsymbol{z})=\widehat{\Phi}\Big(\adj_1(\boldsymbol{A}^{-1}(\boldsymbol{z}))\,\adj_1\boldsymbol{F},\dots,\adj_{N-1}(\boldsymbol{A}^{-1}(\boldsymbol{z}))\,\adj_{N-1}\boldsymbol{F},\adj_N(\boldsymbol{A}^{-1}(\boldsymbol{z}))\,\adj_N\boldsymbol{F}\Big).
\end{equation*} 
Since the map from $\prod_{r=1}^{N-1} \R^{\binom{N}{r}\times \binom{N}{r}}\times (0,+\infty)$ to $\R$ given by
\begin{equation*}
	\begin{split}
		(\boldsymbol{G}_1,\dots,\boldsymbol{G}_{N-1},G_N)&\mapsto \widehat{\Phi}(\adj_1(\boldsymbol{A}^{-1}(\boldsymbol{z}))\,\boldsymbol{G}_1,\dots,\adj_{N-1}(\boldsymbol{A}^{-1}(\boldsymbol{z}))\,\boldsymbol{G}_{N-1},  \adj_N(\boldsymbol{A}^{-1}(\boldsymbol{z}))  G_N) 
	\end{split}
\end{equation*}
 is clearly  convex for every fixed $\boldsymbol{z}\in  S $, this shows that  $W^{\rm nem}$ satisfies (W3).

Now, if $((\boldsymbol{y}_n,\boldsymbol{n}_n))_n$ is a minimizing sequence in   any of the three cases for the respective  class  of admissible states, then $(\boldsymbol{y}_n)_n$ is bounded in $L^p(\Omega;\R^N)$. 
In  (i)--(ii), this boundedness is deduced by applying the Poincaré inequality with trace term given that $\haus(\Gamma)>0$,  while, in (iii), it is immediate  from \eqref{eqn:Y-frac-M}.  
Also, by claim (i) of Corollary \ref{cor:change-of-variable}, we compute
\begin{equation*}
	\int_{\img(\boldsymbol{y}_n,\Omega)} |\boldsymbol{n}_n|^2\,\d\boldsymbol{\xi}=\leb(\img(\boldsymbol{y},\Omega))=\int_{\Omega} \det \nabla \boldsymbol{y}_n\,\d \boldsymbol{x},
\end{equation*}
 where the right-hand side is uniformly bounded because of the growth in \eqref{eqn:Wnem-coe} with respect to the determinant. Eventually, 
the sequence $(\boldsymbol{n}_n \circ \boldsymbol{y}_n)_n$  is also uniformly bounded   on $\Omega$ and,  in turn, equi-integrable.   

 Due to  the compactness results given by the three theorems above, the sequence $((\boldsymbol{y}_n,\boldsymbol{n}_n))_n$ converges to some admissible state $(\boldsymbol{y},\boldsymbol{n})$ with respect to the relevant topology. The boundary conditions for $\boldsymbol{y}$ are satisfied by standard trace theory. Since
the maps $(\boldsymbol{n}_n)_n$ are constrained to take values in the closed set $ Z=S$, their limit $\boldsymbol{n}$ fulfills the same constraint. Hence, by  lower semicontinuity, we deduce that $(\boldsymbol{y},\boldsymbol{n})$ is a minimizer of  the energy  in the suitable class.
\end{proof}

\subsection{Phase transitions}  \label{subsec:ph}
In this subsection, we discuss the variational model for  phase transitions  with interfacial energies proposed by \v{S}ilhav\'{y}  \cite{grandi.etal,grandi.etal2,silhavy.proc,silhavy}. Our modeling approach resembles the one adopted in \cite{grandi.etal2} for multiphase solids.

The energy depends on the deformation $\boldsymbol{y}\in \mathcal{Y}(\Omega)$   and the phase indicator $\boldsymbol{\phi}\colon \img(\boldsymbol{y},\Omega)\to \Lambda$. Here, $M\in \N$ denotes the number of pure phases and 
\begin{align}\label{lambda}
\Lambda\coloneqq \left\{\boldsymbol{z}\in \R^M: \hspace{4pt} 0\leq z_1,\dots,z_M\leq 1, \hspace{4pt} \sum_{m=1}^{M}z_m=1\right\}.
\end{align}
For every $m=1,\dots,M$, the value of $\phi^m\in [0,1]$, where $\boldsymbol{\phi}=(\phi^1,\dots,\phi^M)^\top$, represents the local  concentration of the $m$-th species. This interpretation   explains  why  $\boldsymbol{\phi}$ is constrained to take values in $\Lambda$. 

Inspired by \cite{grandi.etal2}, we introduce the elastic density $W^{\rm ph}\colon \rnn_+ \times \Lambda \to [0,+\infty]$  defined as 
\begin{equation}
	\label{eqn:Wph}
	W^{\rm ph}(\boldsymbol{F},\boldsymbol{z})\coloneqq \sum_{m=1}^{M} z_m\,\Phi_m(\boldsymbol{F}),
\end{equation}
where the functions $\Phi_1,\dots,\Phi_M\colon \rnn_+ \to [0,+\infty]$  satisfy ($\Phi$1)--($\Phi$3).  Also, we consider a continuous function $g\colon  \Lambda \to [0,+\infty)$ such that $g^{-1}(0)=\{ \boldsymbol{p}_m: \hspace{2pt} m=1,\dots, M\}$, where $(\boldsymbol{p}_1, \dots, \boldsymbol{p}_M)$ denotes the canonical basis in $\R^M$. Thus, $g$ is a multi-well potential with wells corresponding to the pure phases.  Setting
\begin{equation}
	\label{eqn:G}
	\mathcal{G}(\boldsymbol{y},\boldsymbol{\phi})\coloneqq \int_{\img(\boldsymbol{y},\Omega)} g\circ \boldsymbol{\phi}\,\d\boldsymbol{\xi},
\end{equation}   
we immediately observe that the value of this functional does not depend on the representatives of $\boldsymbol{y}$ and $\boldsymbol{\phi}$ by Remark~\ref{rem:geom-dom-im}(c).   

In \cite{grandi.etal2,grandi.etal},  the authors consider purely elastic deformations, while the phase indicators are modeled as Sobolev or piecewise-constant maps. These two instances correspond to diffuse and sharp-interface models, respectively. The case of phase indicators in $SBV$ could represent an intermediate model between the two.
The next result  encompasses  all these situations and accounts for possible  cavitation and fracture phenomena.

\begin{proposition}[Existence for phase transitions] \label{prop:ph}
	Let $W^{\rm ph}$ be  as in \eqref{eqn:Wph} with  $\Phi_1,\dots,\Phi_M$ satisfying  {\rm ($\Phi$1)}--{\rm ($\Phi$3)}  and  let  $g\in C^0(\Lambda;[0,+\infty))$ with $g^{-1}(0)=\{ \boldsymbol{p}_m: \hspace{2pt} m=1,\dots, M\}$.    
	\begin{enumerate}[(i)]
		\item  For $\lambda_1\geq \lambda_3>0$,   the functional
		\begin{equation*}
			\begin{split}
				(\boldsymbol{y},\boldsymbol{\phi})\mapsto &\int_{\Omega} W^{\rm ph}(D\boldsymbol{y},\boldsymbol{\phi}\circ \boldsymbol{y})\,\d\boldsymbol{x}+ \lambda_1  \mathcal{S}(\boldsymbol{y})\\
				&+\int_{\img(\boldsymbol{y},\Omega)}|\nabla \boldsymbol{\phi}|^2\,\d\boldsymbol{\xi}+ \lambda_3  \haus(J_{\boldsymbol{\phi}}\cap \imt(\boldsymbol{y},\Omega))  + \int_{\img(\boldsymbol{y},\Omega)} g \circ \boldsymbol{\phi} \,\d\boldsymbol{\xi }
			\end{split}
		\end{equation*}
		admits minimizers within the class
		\begin{equation*}
			\left\{ (\boldsymbol{y},\boldsymbol{\phi}): \hspace{4pt}\boldsymbol{y}\in \mathcal{Y}_p^{\rm cav}(\Omega), \hspace{4pt} \text{$\tr_{\partial \Omega}(\boldsymbol{y})\simeq\boldsymbol{d}$ on $\Gamma$}, \hspace{4pt} \boldsymbol{\phi}\in L^2(\img(\boldsymbol{y},\Omega);\Lambda), \hspace{4pt} \widetilde{\boldsymbol{\phi}}\in SBV^2(\imt(\boldsymbol{y},\Omega);\R^M)  \right\},
		\end{equation*}
		where $\widetilde{\boldsymbol{\phi}}$ denotes the extension of $\boldsymbol{\phi}$ to $\imt(\boldsymbol{y},\Omega)$ by zero.
		\item For every $\kappa>0$  and $\lambda_1>0$,  the functional
		\begin{equation*}
			(\boldsymbol{y},\boldsymbol{\phi})\mapsto \int_{\Omega} W^{\rm ph}(D\boldsymbol{y},\boldsymbol{\phi}\circ \boldsymbol{y})\,\d\boldsymbol{x}+ \lambda_1  \mathcal{S}(\boldsymbol{y})+\int_{\imm(\boldsymbol{y},\Omega)}|D \boldsymbol{\phi}|^2\,\d\boldsymbol{\xi}  + \int_{\img(\boldsymbol{y},\Omega)} g\circ \boldsymbol{\phi} \,\d\boldsymbol{\xi } 
		\end{equation*}
		admits minimizers within the class
		\begin{equation*}
			\left\{ (\boldsymbol{y},\boldsymbol{\phi}): \hspace{4pt}\boldsymbol{y}\in \mathcal{Y}_{p,\kappa}^{\rm cav}(\Omega), \hspace{4pt} \text{$\tr_{\partial \Omega}(\boldsymbol{y})\simeq\boldsymbol{d}$ on $\Gamma$}, \hspace{4pt} \boldsymbol{\phi}\in W^{1,2}(\imm(\boldsymbol{y},\Omega);\Lambda)  \right\}.
		\end{equation*}
		\item For every $b>0$, $\lambda_1>0$, and $\lambda_2\geq \lambda_3>0$,   the functional
		\begin{equation*}
			\begin{split}
				(\boldsymbol{y},\boldsymbol{\phi})\mapsto &\int_{\Omega} W^{\rm ph}(\nabla\boldsymbol{y},\boldsymbol{\phi}\circ \boldsymbol{y})\,\d\boldsymbol{x}+  \lambda_1 \mathcal{S}(\boldsymbol{y}) +\lambda_2  \per\big(\img(\boldsymbol{y},\Omega)\big) + \haus(J_{\boldsymbol{y}})  \\
				&+\int_{\img(\boldsymbol{y},\Omega)}|\nabla \boldsymbol{\phi}|^2\,\d\boldsymbol{\xi}+ \lambda_3  \haus(J_{\boldsymbol{\phi}})  + \int_{\img(\boldsymbol{y},\Omega)} g\circ \boldsymbol{\phi} \,\d\boldsymbol{\xi}
			\end{split}
		\end{equation*}
		admits minimizers within the class
		\begin{equation*}
			\left\{ (\boldsymbol{y},\boldsymbol{\phi}): \hspace{4pt}\boldsymbol{y}\in \mathcal{Y}_{p,b}^{\rm frac}(\widetilde{\Omega}), \hspace{4pt} \text{$\boldsymbol{y}\cong\boldsymbol{d}$ on $\widetilde{\Omega}\setminus \closure{\Omega}$}, \hspace{4pt} \boldsymbol{\phi}\in L^2(\img(\boldsymbol{y},\Omega);\Lambda), \hspace{4pt} \widebar{\boldsymbol{\phi}}\in SBV^2(\R^N;\R^M)  \right\},
		\end{equation*}
		 where   $\widebar{\boldsymbol{\phi}}$ denotes the extension of $\boldsymbol{\phi}$ to the whole space by zero. 
	\end{enumerate}
	In (i) and (iii), one can also restrict to $\widetilde{\boldsymbol{\phi}}\in PC(\imt(\boldsymbol{y},\Omega);\R^M)$ and $\widebar{\boldsymbol{\phi}}\in PC(\R^N;\R^M)$, respectively.
\end{proposition}

\begin{remark}\label{rem:ph}
	\begin{enumerate}[(a)]
		\item  Similarly to Remark \ref{rem:nem}(a),  the assumption $\widetilde{\boldsymbol{\phi}}\in W^{1,2}(\imt(\boldsymbol{y},\Omega);\R^M)$ is incompatible with $\boldsymbol{\phi}\in L^2(\img(\boldsymbol{y},\Omega);\Lambda)$  as $|\widetilde{\boldsymbol{\phi}}| \ge N^{-1/2}$ on $\img(\boldsymbol{y},\Omega)$ and $\widetilde{\boldsymbol{\phi}} = \boldsymbol{0}$ else.  Hence, one should resort to the setting in (ii) when jumps of the phase indicator are excluded.
		\item  Applied loads can be considered. We refer to  Remark \ref{rem:nem}(b) for more details.  
		\item In our model, the components of $\boldsymbol{\phi}$ represent concentrations of the $M$ species in the solid. Avoiding the normalization given by $\Lambda$, one could equivalently consider phase indicators $\boldsymbol{\phi}$ taking values in a more general compact subset of $\R^M$. In this situation, one generally imposes a mass constraint of the form
		\begin{equation*}
			\int_{\img(\boldsymbol{y},\Omega)} \boldsymbol{\phi}\,\d\boldsymbol{\xi}=\boldsymbol{\mu},
		\end{equation*} 
		where $\boldsymbol{\mu}\in \R^M$ is the prescribed distribution of masses. The existence results in Proposition \ref{prop:ph} can be easily adapted to this setting.
	\end{enumerate}
\end{remark}

\begin{proof}
 The proof works analogously to the one of Proposition \ref{prop:nem}. We consider the functional $\widetilde{\mathcal{F}}^{\rm ph}$ defined  as in \eqref{eqn:F-EL}--\eqref{eqn:F-EL2}, but  with $W^{\rm ph}$ in place of $W$. All energies in the statement are given by the restriction of $\mathcal{F}^{\rm ph}\coloneqq \widetilde{\mathcal{F}}^{\rm ph}+\mathcal{G}$, where $\mathcal{G}$ is defined as in \eqref{eqn:G}, to the different classes of admissible states. Since $\mathcal{G}$ is nonnegative, if $\mathcal{F}^{\rm ph}$ is uniformly bounded along a sequence, then $\widetilde{\mathcal{F}}^{\rm ph}$ is also uniformly bounded along the same sequence.  

We check the common assumptions of Theorems \ref{thm:cav}, \ref{thm:cavk}, and \ref{thm:frac}.  Note that phase indicators take values in the closed set $Z =\Lambda$.
Also here, (W1) is clear. By \eqref{eqn:growth-Phi}  and  \eqref{lambda},   we have
\begin{equation*}
W^{\rm ph}(\boldsymbol{F},\boldsymbol{z})\geq \sum_{m=1}^{M} z_m \left( {C}|\boldsymbol{F}|^p+{\gamma}(\det \boldsymbol{F}) \right)={C}|\boldsymbol{F}|^p+{\gamma}(\det \boldsymbol{F}) \quad \text{for all  $\boldsymbol{F}\in \rnn_+$ and  $\boldsymbol{z}\in \Lambda$},
\end{equation*}
so that (W2)  holds  true. The validity of (W3) follows immediately from ($\Phi$3).
The rest of the proof works as the one of  Proposition \ref{prop:nem},  simply replacing the compact set $S$ by $\Lambda$  in the arguments. 
\end{proof}

\subsection{Magnetoelasticity} \label{subsec:mag}
The last application concerns Brown's variational model of magnetoelasticity \cite{brown} which has been recently explored also from the mathematical perspective \cite{barchiesi.henao.moracorral,bresciani,bresciani.davoli.kruzik,desimone.dolzmann,kruzik.stefanelli.zeman,rybka.luskin}.

 The governing energy depends on the deformation $\boldsymbol{y}\in \mathcal{Y}(\Omega)$ and the magnetization  $\boldsymbol{m}\colon \img(\boldsymbol{y},\Omega) \to \R^N$.  As in \cite{bresciani,desimone.dolzmann,james.kinderlehrer,rybka.luskin},  we interpret the latter variable as the local average of magnetic dipoles per unit volume. Therefore, this is subject to the saturation constraint  (see  \cite[p.\ 73]{brown})  which reads   
\begin{equation}
	\label{eqn:sat}
	\text{$|\boldsymbol{m}\circ \boldsymbol{y}|\det \nabla \boldsymbol{y}\cong 1$ in $\Omega$,}
\end{equation}
or equivalently  
\begin{equation}
	\label{eqn:sat2}
	\text{$|\boldsymbol{m}|\cong \det \nabla \boldsymbol{y}^{-1}$ in $\img(\boldsymbol{y},\Omega)$,}
\end{equation}
 thanks to claim (ii) of Lemma \ref{lem:inverse-differentiable}. 
Note that, according to  \eqref{eqn:sat2} and Lemma \ref{lem:inverse-differentiable}(ii),   magnetizations take values in $\R^N_\times  := \R^N \setminus \lbrace \boldsymbol{0} \rbrace$.   In particular, for incompressible materials,  magnetizations  are sphere valued \cite{kruzik.stefanelli.zeman}. We refer to \cite[Section 1]{bresciani}  for a thorough discussion on the  constraint in \eqref{eqn:sat}--\eqref{eqn:sat2}. 

Following the modeling approach in \cite{bresciani,james.kinderlehrer}, we define the elastic density $W^{\rm mag}\colon \rnn_+ \times \R^N_\times \to [0,+\infty]$ as
\begin{equation}
	\label{eqn:Wmag}
	W^{\rm mag}(\boldsymbol{F},\boldsymbol{z})\coloneqq 
		\Phi \left( \boldsymbol{B}\left(\boldsymbol{z}\right)\boldsymbol{F} \right), \qquad \boldsymbol{B}(\boldsymbol{z})\coloneqq \beta_0 \frac{\boldsymbol{z}}{|\boldsymbol{z}|}\otimes  \frac{\boldsymbol{z}}{|\boldsymbol{z}|} + \beta_1 \left(\boldsymbol{I}-  \frac{\boldsymbol{z}}{|\boldsymbol{z}|} \otimes  \frac{\boldsymbol{z}}{|\boldsymbol{z}|} \right),
\end{equation} 
where $\beta_0,\beta_1>0$ and $\Phi \colon \rnn_+ \to [0,+\infty]$  satisfies ($\Phi$1)--($\Phi$3). Additionally, as in \cite{bresciani,bresciani.thesis},  for technical reasons (see \eqref{technical} below),  we strengthen the first condition in \eqref{eqn:growth-gamma-Phi} by requiring   
\begin{equation}
	\label{eqn:gamma-st}
	\lim_{h \to 0^+} h\gamma(h)=+\infty.
\end{equation}

The magnetoelastic energy also comprises a contribution depending on the stray field $\boldsymbol{h}_{\boldsymbol{m}}\coloneqq -D u_{\boldsymbol{m}}$,  where the potential $u_{\boldsymbol{m}}$ is given by a weak solution to the Maxwell equation
\begin{equation}
	\label{eqn:maxwell}
	\text{$\div \left( -D u_{\boldsymbol{m}}  + \chi_{\img(\boldsymbol{y},\Omega)}\boldsymbol{m} \right)=0$ in $\R^N$.}
\end{equation}
Namely, $u_{\boldsymbol{m}} \in V^{1,2}(\R^N)$ satisfies  
\begin{equation}
	\label{eqn:maxwell-weak}
\int_{\R^N} Du_{\boldsymbol{m}} \cdot D u\,\d\boldsymbol{\xi}= \int_{\img(\boldsymbol{y},\Omega)} \boldsymbol{m}\cdot D u\,\d\boldsymbol{\xi} \quad \text{	for all $u \in V^{1,2}(\R^N)$},
\end{equation}
 where we recall  \eqref{beppo}.  The stray-field energy is given by
\begin{equation}
	\label{eqn:H}
	\mathcal{H}(\boldsymbol{y},\boldsymbol{m})\coloneqq \int_{\R^N}|\boldsymbol{h}_{\boldsymbol{m}}|^2\,\d\boldsymbol{\xi}.
\end{equation}

We collect some observations on the stray-field energy in the next lemma.  Analogous results have been proved in \cite[Proposition 8.8]{barchiesi.henao.moracorral}, \cite[Proposition 3.6]{bresciani}, \cite[Theorem 3.2]{bresciani.davoli.kruzik}, and \cite[Lemma 2.3]{kruzik.stefanelli.zeman}.  

\begin{lemma}[Stray field] \label{lem:stray-field}
Let $\mathcal{H}$ be defined as in \eqref{eqn:H}. Then:
\begin{enumerate}[(i)]
	\item For every $\boldsymbol{y}\in \mathcal{Y}(\Omega)$   and $\boldsymbol{m}\in L^2(\img(\boldsymbol{y},\Omega);\R^N_\times)$, the Maxwell equation \eqref{eqn:maxwell} admits a weak solution $u_{\boldsymbol{m}}\in V^{1,2}(\R^N)$ which is unique up to additive constants and satisfies
	\begin{equation*}
		\|\boldsymbol{h}_{\boldsymbol{m}}\|_{L^2( \R^N;  \R^N)}=\|D u_{\boldsymbol{m}}\|_{L^2(\R^N;\R^N)} \leq \|\chi_{\img(\boldsymbol{y},\Omega)}\boldsymbol{m}\|_{L^2(\R^N;\R^N)}.
	\end{equation*}
	In particular,   $\mathcal{H}(\boldsymbol{y},\boldsymbol{m})  < + \infty $.  
	\item Let $(\boldsymbol{y}_n)_n\subset  \mathcal{Y}(\Omega)$  and $(\boldsymbol{m}_n)_n$ be a sequence of maps $\boldsymbol{m}_n \in L^2(\img(\boldsymbol{y}_n,\Omega);\R^N_\times)$. Also, let $\boldsymbol{y}\in \mathcal{Y}(\Omega)$ and $\boldsymbol{m}\in L^2(\img(\boldsymbol{y},\Omega);\R^N_\times)$. If 
	\begin{equation}\label{MMM}
		\text{$\chi_{\img(\boldsymbol{y}_n,\Omega)}\boldsymbol{m}_n\wk \chi_{\img(\boldsymbol{y},\Omega)}\boldsymbol{m}$ in $L^2(\R^N;\R^N)$,}
	\end{equation}
	then   $\boldsymbol{h}_{\boldsymbol{m}_n}  \wk \boldsymbol{h}_{\boldsymbol{m}}$ in $L^2(\R^N; \R^N)$ and thus 
	\begin{equation*}
		\mathcal{H}(\boldsymbol{y},\boldsymbol{m})\leq \liminf_{n \to \infty} \mathcal{H}(\boldsymbol{y}_n,\boldsymbol{m}_n).
	\end{equation*}
\end{enumerate}
\end{lemma}  
\begin{remark}
	\begin{enumerate}[(a)]
		\item If $(\boldsymbol{y}_1,\boldsymbol{m}_1)$ and $(\boldsymbol{y}_2,\boldsymbol{m}_2)$ satisfy  $\boldsymbol{y}_1 \cong \boldsymbol{y}_2$ and $\boldsymbol{m}_1 \cong \boldsymbol{m}_2$ in $\img(\boldsymbol{y}_1,\Omega)\cap \img(\boldsymbol{y}_2,\Omega)$, then $\chi_{\img(\boldsymbol{y}_1,\Omega)}\boldsymbol{m}_1\cong \chi_{\img(\boldsymbol{y}_2,\Omega)}\boldsymbol{m}_2$. In this case, $\boldsymbol{h}_{\boldsymbol{m}_1}\cong \boldsymbol{h}_{\boldsymbol{m}_2}$ as a consequence of claim (i).   In particular,
		the value of $\mathcal{H}(\boldsymbol{y},\boldsymbol{m})$ does not depend on the representatives of $\boldsymbol{y}$ and $\boldsymbol{m}$.
		\item In (ii), if the convergence of $(\chi_{\img(\boldsymbol{y}_n,\Omega)}\boldsymbol{m}_n)$ is strong in $L^2(\R^N;\R^N)$, then $\boldsymbol{h}_{\boldsymbol{m}_n}\to \boldsymbol{h}_{\boldsymbol{m}}$ in $L^2(\R^N;\R^N)$ and, in turn, $\mathcal{H}(\boldsymbol{y}_n,\boldsymbol{m}_n)\to\mathcal{H}(\boldsymbol{y},\boldsymbol{m})$, see \cite[Remark 3.7]{bresciani} for a proof.
	\end{enumerate}
\end{remark}

\begin{proof}
(i) The claim is proved in \cite[Proposition 8.8]{barchiesi.henao.moracorral}.

(ii)  From \eqref{MMM}, by claim (i), we see that $ \boldsymbol{h}_{\boldsymbol{m}_n}=-D u_{{\boldsymbol{m}}_n}  \wk \boldsymbol{h}$ in $L^2(\R^N;\R^N)$ for some $\boldsymbol{h}\in L^2(\R^N;\R^N)$ up to subsequences.  Passing to the limit $n\to \infty$ on both sides of \eqref{eqn:maxwell-weak},  {(with $\boldsymbol{m}_n$ and $\boldsymbol{y}_n$ in place of  $\boldsymbol{m}$ and $\boldsymbol{y}$, respectively)} with the aid of \eqref{MMM}, we deduce that $\boldsymbol{h}= - Du_{\boldsymbol{m}} =\boldsymbol{h}_{\boldsymbol{m}}$.  Thus, the lower semicontinuity of $\mathcal{H}$ follows from the one of the $L^2$-norm.
\end{proof}

All existence results for minimizers of the magnetoelastic energy available in the literature concern purely elastic deformations.  Most of these are framed in the setting of micromagnetics \cite{brown}, so that magnetizations are not allowed to jump \cite{barchiesi.henao.moracorral,bresciani,bresciani.davoli.kruzik,kruzik.stefanelli.zeman,rybka.luskin}. Piecewise-constant magnetizations can be employed to  describe  the formation of magnetic-domain structures as done in \cite{anzellotti.baldo.visintin,grandi.etal2} for the high-anisotropy limit and \cite{james.kinderlehrer} for the large-body limit.   Models with magnetizations in $SBV$ have also been investigated \cite{acerbi.fonseca.mingione,santos}. The next theorem   extends all these results  to the setting of material failure.

\begin{proposition}[Existence for magnetoelasticity] \label{prop:mag}
Let  $W^{\rm mag}$ be  as in \eqref{eqn:Wmag} with $\beta_0,\beta_1>0$ and $\Phi$ satisfying {\rm ($\Phi$1)}--{\rm ($\Phi$3)} together with \eqref{eqn:gamma-st}. 
	\begin{enumerate}[(i)]
		\item  For $\lambda_1>0$,  the functional
		\begin{equation*}
			(\boldsymbol{y},\boldsymbol{m})\mapsto \int_{\Omega} W^{\rm mag}(D\boldsymbol{y},\boldsymbol{m}\circ \boldsymbol{y})\,\d\boldsymbol{x}+ \lambda_1  \mathcal{S}(\boldsymbol{y})+\int_{\img(\boldsymbol{y},\Omega)}|D \boldsymbol{m}|^2\,\d\boldsymbol{\xi}+\int_{\R^N} |\boldsymbol{h}_{\boldsymbol{m}}|^2\,\d\boldsymbol{\xi}
		\end{equation*}
		admits minimizers within the class
		\begin{equation*}
			\begin{split}
				\big\{ (\boldsymbol{y},\boldsymbol{m}): \hspace{8pt}\boldsymbol{y}\in \mathcal{Y}_p^{\rm cav}(\Omega), \hspace{8pt} \text{$\tr_{\partial \Omega}(\boldsymbol{y})\simeq\boldsymbol{d}$ on $\Gamma$}, \hspace{8pt} \boldsymbol{m}\in L^2(\img(\boldsymbol{y},\Omega);\R^N_\times),&\\
				\text{$ |\boldsymbol{m}\circ \boldsymbol{y}|   \det D \boldsymbol{y} \cong1$ in $\Omega$},\hspace{8pt} \widetilde{\boldsymbol{m}}\in W^{1,2}(\imt(\boldsymbol{y},\Omega);\R^N)  \big\},&
			\end{split}
		\end{equation*}
		 where   $\widetilde{\boldsymbol{m}}$ denotes the extension of $\boldsymbol{m}$ to $\imt(\boldsymbol{y},\Omega)$ by zero.
		\item  For $\lambda_1 \geq \lambda_3>0$,  the functional
		\begin{equation*}
			\begin{split}
				(\boldsymbol{y},\boldsymbol{m})\mapsto &\int_{\Omega} W^{\rm mag}(D\boldsymbol{y},\boldsymbol{m}\circ \boldsymbol{y})\,\d\boldsymbol{x}+ \lambda_1  \mathcal{S}(\boldsymbol{y})\\
				&+\int_{\img(\boldsymbol{y},\Omega)} |\nabla \boldsymbol{m}|^2\,\d\boldsymbol{\xi}+\int_{\R^N} |\boldsymbol{h}_{\boldsymbol{m}}|^2\,\d\boldsymbol{\xi}+ \lambda_3  \haus(J_{\boldsymbol{m}}\cap \imt(\boldsymbol{y},\Omega))
			\end{split}
		\end{equation*}
		admits minimizers within the class
		\begin{equation*}
			\begin{split}
				\big\{ (\boldsymbol{y},\boldsymbol{m}): \hspace{8pt}\boldsymbol{y}\in \mathcal{Y}_p^{\rm cav}(\Omega), \hspace{4pt} \text{$\tr_{\partial \Omega}(\boldsymbol{y})\simeq\boldsymbol{d}$ on $\Gamma$}, \hspace{8pt} \boldsymbol{m}\in L^2(\img(\boldsymbol{y},\Omega);\R^N_\times),&\\
				 \text{$|  \boldsymbol{m}\circ \boldsymbol{y}|  \det D \boldsymbol{y}  \cong1$ in $\Omega$},\hspace{8pt} \widetilde{\boldsymbol{m}}\in GSBV^2(\imt(\boldsymbol{y},\Omega);\R^N)  \big\},&
			\end{split}
		\end{equation*}
where	   $\widetilde{\boldsymbol{m}}$ denotes the extension of $\boldsymbol{m}$ to $\imt(\boldsymbol{y},\Omega)$ by zero.
		\item For every $\kappa>0$  and $\lambda_1>0$,  the functional
		\begin{equation*}
			(\boldsymbol{y},\boldsymbol{m})\mapsto \int_{\Omega} W^{\rm mag}(D\boldsymbol{y},\boldsymbol{m}\circ \boldsymbol{y})\,\d\boldsymbol{x}+ \lambda_1  \mathcal{S}(\boldsymbol{y})+\int_{\imm(\boldsymbol{y},\Omega)}|D \boldsymbol{m}|^2\,\d\boldsymbol{\xi}+\int_{\R^N} |\boldsymbol{h}_{\boldsymbol{m}}|^2\,\d\boldsymbol{\xi}
		\end{equation*}
		admits minimizers within the class
		\begin{equation*}
			\left\{ (\boldsymbol{y},\boldsymbol{m}): \hspace{4pt}\boldsymbol{y}\in \mathcal{Y}_{p,\kappa}^{\rm cav}(\Omega), \hspace{4pt} \text{$\tr_{\partial \Omega}(\boldsymbol{y})\simeq\boldsymbol{d}$ on $\Gamma$}, \hspace{4pt} \boldsymbol{m}\in W^{1,2}(\imm(\boldsymbol{y},\Omega);\R^N_\times), \hspace{4pt} \text{$|\boldsymbol{m}\circ \boldsymbol{y}|\det D \boldsymbol{y}\cong 1$ in $\Omega$}  \right\}.
		\end{equation*}
		\item For every $b>0$,  $\lambda_1>0$, and $\lambda_2\geq \lambda_3>0$,   the functional
		\begin{equation*}
			\begin{split}
				(\boldsymbol{y},\boldsymbol{m})\mapsto &\int_{\Omega} W^{\rm mag}(\nabla\boldsymbol{y},\boldsymbol{m}\circ \boldsymbol{y})\,\d\boldsymbol{x}+ \lambda_1 \mathcal{S}(\boldsymbol{y})+ \lambda_2  \per\big(\img(\boldsymbol{y},\Omega)\big)+ \haus(J_{\boldsymbol{y}})\\
				&+\int_{\imm(\boldsymbol{y},\Omega)}|\nabla \boldsymbol{m}|^2\,\d\boldsymbol{\xi}+\int_{\R^N} |\boldsymbol{h}_{\boldsymbol{m}}|^2\,\d\boldsymbol{\xi}+
				+ \lambda_3  \haus(J_{\boldsymbol{m}})
			\end{split}
		\end{equation*}
		admits minimizers within the class
		\begin{equation*}
			\begin{split}
				\big \{ (\boldsymbol{y},\boldsymbol{m}): \hspace{4pt}\boldsymbol{y}\in \mathcal{Y}_{p,b}^{\rm frac}(\widetilde{\Omega}), \hspace{4pt} \text{$\boldsymbol{y}\cong\boldsymbol{d}$ on $\widetilde{\Omega}\setminus \closure{\Omega}$}, \hspace{4pt} \boldsymbol{m}\in L^2(\img(\boldsymbol{y},\Omega);\R^N_\times),&\\ \hspace{4pt}  \text{$|\boldsymbol{m}\circ \boldsymbol{y}|\det \nabla \boldsymbol{y}\cong 1$ in $\Omega$},  \hspace{4pt} \widebar{\boldsymbol{m}}\in GSBV^2(\R^N;\R^N)  \big \},&
			\end{split}
		\end{equation*}
		where $\widebar{\boldsymbol{m}}$ denotes the extension of $\boldsymbol{m}$ to the whole space by zero. 
	\end{enumerate}
	In  (ii) and (iv), one can also restrict to $\widetilde{\boldsymbol{m}}\in PC(\imt(\boldsymbol{y},\Omega);\R^N)$ and $\widebar{\boldsymbol{m}}\in PC(\R^N;\R^N)$, respectively.
\end{proposition}
\begin{remark}\label{rem:mag}
	\begin{enumerate}[(a)]
		\item   The setting in (i) is compatible with the saturation constraint  but poses a restriction on the class of admissible deformations. Precisely, the assumption $\widetilde{\boldsymbol{m}}\in W^{1,2}(\imt(\boldsymbol{y},\Omega);\R^N)$ and \eqref{eqn:sat2} force the Jacobian determinant of $\boldsymbol{y}$ to explode in correspondence of cavitation points to avoid jumps of $\widetilde{\boldsymbol{m}}$ on the boundary of  cavities. To see this, let $\boldsymbol{a}\in C_{\boldsymbol{y}}$ and $\boldsymbol{\nu}_{\imt(\boldsymbol{y},\boldsymbol{a})}\colon \partial^* \imt(\boldsymbol{y},\boldsymbol{a})\to S$ be the outer unit normal, where we recall Theorem \ref{thm:INV-top-im}(i). 
		Fix $\boldsymbol{\xi}_0\in \partial^*\imt(\boldsymbol{y},\boldsymbol{a})$ and set $\boldsymbol{\nu}_0\coloneqq \boldsymbol{\nu}_{\imt(\boldsymbol{y},\boldsymbol{a})}(\boldsymbol{\xi}_0)$. 
		 Since $\widetilde{\boldsymbol{m}}\cong \boldsymbol{0}$ in $\imt(\boldsymbol{y},\boldsymbol{a})$,  we have  $		\widetilde{\boldsymbol{m}}^-(\boldsymbol{\xi}_0)= \boldsymbol{0}$,  while \eqref{eqn:sat2} yields
		\begin{equation*}
			|\widetilde{\boldsymbol{m}}^+(\boldsymbol{\xi}_0)|=\aplim_{\substack{\boldsymbol{\xi}\to \boldsymbol{\xi}_0 \\ \boldsymbol{\xi}\in H^+(\boldsymbol{\xi}_0,\boldsymbol{\nu}_0)}} |\widetilde{\boldsymbol{m}}(\boldsymbol{\xi})|=\aplim_{\substack{\boldsymbol{\xi}\to \boldsymbol{\xi}_0 \\ \boldsymbol{\xi}\in H^+(\boldsymbol{\xi}_0,\boldsymbol{\nu}_0)}} |\boldsymbol{m}(\boldsymbol{\xi})|=\aplim_{\substack{\boldsymbol{\xi}\to \boldsymbol{\xi}_0 \\ \boldsymbol{\xi}\in H^+(\boldsymbol{\xi}_0,\boldsymbol{\nu}_0)}} \det \nabla \boldsymbol{y}^{-1}(\boldsymbol{\xi}).
		\end{equation*}
		Therefore, the regularity of $\widetilde{\boldsymbol{m}}$ gives the identity
		\begin{equation*}
			\aplim_{\substack{\boldsymbol{\xi}\to \boldsymbol{\xi}_0 \\ \boldsymbol{\xi}\in H^+(\boldsymbol{\xi}_0,\boldsymbol{\nu}_0)}} \det \nabla \boldsymbol{y}^{-1}(\boldsymbol{\xi})=0
		\end{equation*}
		for $\haus$-almost all $\boldsymbol{\xi}_0\in \partial^*\imt(\boldsymbol{y},\boldsymbol{a})$. 
		From this, recalling Lemma \ref{lem:inverse-differentiable}(ii) and setting $\boldsymbol{x}_0\coloneqq \boldsymbol{y}^{-1}(\boldsymbol{\xi}_0)$, we   deduce 
		\begin{equation*}
			\aplim_{\boldsymbol{x}\to \boldsymbol{x}_0} \det  D  \boldsymbol{y}(\boldsymbol{x})=+\infty.
		\end{equation*}
		This behavior is discussed in  Example \ref{ex:mag1} below. In particular, given an incompressible  deformation   $\boldsymbol{y}$, i.e., satisfying $\det D \boldsymbol{y}=1$ almost everwhere, with $C_{\boldsymbol{y}}\neq \emptyset$ there is no $\boldsymbol{m}\in L^2(\img(\boldsymbol{y},\Omega);\R^N_\times)$ satisfying both \eqref{eqn:sat}--\eqref{eqn:sat2} and $\widetilde{\boldsymbol{m}}\in W^{1,2}(\imt(\boldsymbol{y},\Omega);\R^N)$.  
		\item In contrast  to  Proposition \ref{prop:nem} and Proposition \ref{prop:ph}, here it is possible to relax the assumption ($\Phi3$) by  requiring the convexity of the map $(\boldsymbol{G}_1,\dots,\boldsymbol{G}_{N-1})\mapsto \widehat{\Phi}  (\boldsymbol{G}_1,\dots,\boldsymbol{G}_{N-1},d)$ for every fixed $d>0$. Indeed, if
		$((\boldsymbol{y}_n,\boldsymbol{m}_n))_n$ and  $(\boldsymbol{y},\boldsymbol{m})$ are admissible states for which
		  the composition of magnetizations with deformations converge strongly $L^1(\Omega;\R^N)$, then, by exploiting the constraint in \eqref{eqn:sat}, one can prove  that $\det \nabla \boldsymbol{y}_n \to \det \nabla \boldsymbol{y}$ in $L^1(\Omega)$, see  \cite[Theorem 3.2]{bresciani} for details.
		\item  Applied loads as in Remark \ref{rem:nem}(b) can be treated.     In this case,  $\boldsymbol{g}$ represents an external magnetic field which needs to belong to $ L^2(\R^N;\R^N)$. 
	\end{enumerate}
\end{remark}
\begin{proof}
 The proof has the same structure of the one of Proposition \ref{prop:ph}.   We consider the functional $\widetilde{\mathcal{F}}^{\rm mag}$ defined as in \eqref{eqn:F-EL}--\eqref{eqn:F-EL2}, but with $W^{\rm mag}$ in place of $W$.  All the energies in the statement  are given by the restriction of $\mathcal{F}^{\rm mag}\coloneqq \widetilde{\mathcal{F}}^{\rm mag}+\mathcal{H}$, where $\mathcal{H}$ is as in \eqref{eqn:H}, to  the  different classes of admissible states. Note that, if $\mathcal{F}^{\rm mag}$ is uniformly bounded along a sequence of admissible states, then the same holds for $\widetilde{\mathcal{F}}^{\rm mag}$ along the same sequence since $\mathcal{H}$ is  nonnegative. 

We check that  the common assumptions of Theorem \ref{thm:cav}, Theorem \ref{thm:cavk}, and Theorem \ref{thm:frac} are fulfilled. 
 The conditions  (W1)--(W3) are verified as in the proof of Proposition \ref{prop:nem},  given that the structure of $W^{\rm mag}$ and $	W^{\rm nem}$ are essentially the same.  
  Let $  ((\boldsymbol{y}_n,\boldsymbol{m}_n))_n$ be a minimizing sequence for $\mathcal{F}^{\rm mag}$ with respect to any of the classes of admissible states. The boundedness of $(\boldsymbol{y}_n)_n$ in $L^p(\Omega;\R^N)$  follows as in Proposition \ref{prop:nem}, either by a  Poincaré inequality or the definition  in \eqref{eqn:Y-frac-M}.  

We show that $(\boldsymbol{m}_n \circ \boldsymbol{y}_n)_n$ is equi-integrable. Recall that
\begin{equation}
	\label{eqn:satn}
\text{$|\boldsymbol{m}_n \circ \boldsymbol{y}_n|\det \nabla \boldsymbol{y}_n \cong 1$ in $\Omega$ \ \  for all $n \in \N$.}
\end{equation}
As a consequence, magnetizations are constrained to take values in the set $ Z  =\R^N_\times$. 
 We  define $\widehat{\gamma}\colon (0,+\infty)\to [0,+\infty]$ by setting $\widehat{\gamma}(t)\coloneqq \gamma(1/t)$. Then, \eqref{eqn:gamma-st} yields 
\begin{align}\label{technical}
\lim_{t \to +\infty} \widehat{\gamma}(t)/t=+\infty.
\end{align}
Exploiting \eqref{eqn:satn}, we compute
\begin{equation*}
	\int_\Omega \widehat{\gamma}(|\boldsymbol{m}_n \circ \boldsymbol{y}_n|)\,\d\boldsymbol{x}=\int_\Omega \widehat{\gamma}\left( \frac{1}{\det \nabla \boldsymbol{y}_n}\right)\,\d\boldsymbol{x}=\int_{\Omega} \gamma(\det \nabla \boldsymbol{y}_n)\,\d\boldsymbol{x},
\end{equation*}
where the right-hand side is uniformly bounded due to \eqref{eqn:growth} and the boundedness of $\widetilde{\mathcal{F}}^{\rm mag}$. Thus, $(\boldsymbol{m}_n \circ \boldsymbol{y}_n)_n$ is equi-integrable by the De la Vallée {Poussin} criterion and, in particular, it is bounded in $L^1(\Omega;\R^N)$. From this, we easily deduce the boundedness of $(\chi_{\img(\boldsymbol{y}_n,\Omega)}\boldsymbol{m}_n)$  in $L^2(\R^N;\R^N)$. Indeed, exploiting again \eqref{eqn:satn} and applying Corollary \ref{cor:change-of-variable}(i), we obtain
\begin{equation*}
	\int_{\img(\boldsymbol{y}_n,\Omega)} |\boldsymbol{m}_n|^2\,\d\boldsymbol{\xi}=\int_{\Omega} |\boldsymbol{m}_n \circ \boldsymbol{y}_n|^2\,\det \nabla \boldsymbol{y}_n\,\d\boldsymbol{x}=\int_\Omega |\boldsymbol{m}_n \circ \boldsymbol{y}_n|\,\d\boldsymbol{x}.
\end{equation*}
 This allows us to apply the compactness part of our main theorems.   We find  a pair $(\boldsymbol{y},\boldsymbol{m})$ with  $\boldsymbol{y}\in \mathcal{Y}(\Omega)$  in the respective class of deformations  and $\boldsymbol{m}\in L^2(\img(\boldsymbol{y},\Omega);\R^N)$  such that, up to subsequences,    $((\boldsymbol{y}_n,\boldsymbol{m}_n) )_n$ converges to $(\boldsymbol{y},\boldsymbol{m})$ with respect to the relevant topology.  
We check that $(\boldsymbol{y},\boldsymbol{m})$ is admissible, namely, it satisfies \eqref{eqn:sat}  and    $\boldsymbol{m}$ takes values in $\R^N_\times$.   We have
\begin{equation} \label{eqn:yy}
	\text{$\boldsymbol{y}_n \to \boldsymbol{y}$ a.e.\ in $\Omega$, \qquad $\det \nabla \boldsymbol{y}_n \wk \det \nabla \boldsymbol{y}$ in $L^1(\Omega)$},
\end{equation}
\begin{equation}\label{eqn:mm}
	\text{$\chi_{\img(\boldsymbol{y}_n,\Omega)}\boldsymbol{m}_n \to \chi_{\img(\boldsymbol{y},\Omega)}\boldsymbol{m}$ a.e.\ and in $L^1(\R^N;\R^N)$,}
\end{equation}
\begin{equation}\label{eqn:mm2}
	\text{$\chi_{\img(\boldsymbol{y}_n,\Omega)}\boldsymbol{m}_n \wk \chi_{\img(\boldsymbol{y},\Omega)}\boldsymbol{m}$ in $L^2(\R^N;\R^N)$.}
\end{equation}
Let $\varphi \in C^\infty_{\rm c}(\Omega)$. Applying \eqref{eqn:satn} and   Corollary \ref{cor:change-of-variable}(i), we compute
\begin{equation}\label{eqn:cm}
	\int_\Omega \varphi\,\d\boldsymbol{x}=\int_\Omega \varphi |\boldsymbol{m}_n \circ \boldsymbol{y}_n|\det \nabla \boldsymbol{y}_n\,\d\boldsymbol{x}=\int_{\img(\boldsymbol{y}_n,\Omega)} \varphi \circ \boldsymbol{y}_n^{-1}\,|\boldsymbol{m}_n|\,\d\boldsymbol{\xi}.
\end{equation}
From \eqref{eqn:yy}--\eqref{eqn:mm}   and  Proposition \ref{prop:approx-diff}(i), we see that
\begin{equation*}
	\text{$\chi_{\img(\boldsymbol{y}_n,\Omega)}  \, \varphi \circ  \boldsymbol{y}_n^{-1}  \,|\boldsymbol{m}_n|\to \chi_{\img(\boldsymbol{y},\Omega)} \, \varphi \circ   \boldsymbol{y}^{-1}  \,|\boldsymbol{m}|$ a.e.\ in $\R^N$}
\end{equation*}
for a not relabeled subsequence, thanks to the continuity of $\varphi$. 
As
 $${\chi_{\img(\boldsymbol{y}_n,\Omega)} \, |\varphi \circ  \boldsymbol{y}_n^{-1}  |\,|\boldsymbol{m}_n|\leq \chi_{\img(\boldsymbol{y}_n,\Omega)}\|\varphi\|_{C^0_{\rm b}( \Omega  )}\,|\boldsymbol{m}_n|,}$$
  where the sequence of the right-hand side converges strongly in $L^1(\R^N)$ as a consequence of \eqref{eqn:mm}, we can pass to the limit, as $n \to \infty$, in \eqref{eqn:cm} by applying the dominated convergence theorem. We obtain
\begin{equation*}
	\int_{\Omega} \varphi\,\d \boldsymbol{x}=\int_{\img(\boldsymbol{y},\Omega)} \varphi \circ \boldsymbol{y}^{-1}\,|\boldsymbol{m}|\,\d\boldsymbol{\xi}=\int_\Omega \varphi\, |\boldsymbol{m}\circ \boldsymbol{y}|\,\det \nabla \boldsymbol{y}\,\d\boldsymbol{x},
\end{equation*}
where we applied once again  Corollary \ref{cor:change-of-variable}(i). Given the arbitrariness of $\varphi \in C^\infty_{\rm c}(\Omega)$, the previous identity yields $|\boldsymbol{m}\circ \boldsymbol{y}|\,\det \nabla \boldsymbol{y}\cong 1$ in $\Omega$, as desired.  Note that this also directly  shows that  $\boldsymbol{m}$ takes values in $\R^N_\times$.  

At this point, the lower semicontinuity claim of Theorems \ref{thm:cav}, \ref{thm:cavk}, and \ref{thm:frac} yields
\begin{equation*}
	\widetilde{\mathcal{F}}^{\rm mag}(\boldsymbol{y},\boldsymbol{m})\leq \liminf_{n \to \infty} \widetilde{\mathcal{F}}^{\rm mag}(\boldsymbol{y}_n,\boldsymbol{m}_n). 
\end{equation*}
Then, recalling \eqref{eqn:mm2}  and applying claim (ii) of Lemma \ref{lem:stray-field}, we obtain
\begin{equation*}
	{\mathcal{F}}^{\rm mag}(\boldsymbol{y},\boldsymbol{m})\leq \liminf_{n \to \infty} {\mathcal{F}}^{\rm mag}(\boldsymbol{y}_n,\boldsymbol{m}_n),
\end{equation*}
which shows that $(\boldsymbol{y},\boldsymbol{m})$ is a minimizer of $\mathcal{F}^{\rm mag}$.
\end{proof}

The next example  concerns  the observation made in Remark \ref{rem:mag}(a).

\begin{example}\label{ex:mag1}
Let $N=2$ and $\Omega\coloneqq B$, where we employ the notation in \eqref{eqn:q-ball-annulus}--\eqref{eqn:B-S}. Consider  $\boldsymbol{y}_1\colon \Omega \to \R^2$ defined as
\begin{equation*}
	\text{$\boldsymbol{y}_1(\boldsymbol{x})\coloneqq \frac{1}{2} \left( |\boldsymbol{x}|+1 \right) \frac{\boldsymbol{x}}{|\boldsymbol{x}|}$ \quad for all $\boldsymbol{x}\in \Omega\setminus \{ \boldsymbol{0} \}$.}
\end{equation*}   
By  Lemma \ref{lem:radial-degree},   $\boldsymbol{y}_1\in \mathcal{Y}_p^{\rm cav}(\Omega)$. 
This deformation fixes the boundary $S$ of $B$ and opens a concentric cavity of radius $1/2$ at the origin (Figure \ref{fig:cavity-point}).  Precisely,  $\img(\boldsymbol{y}_1,\Omega)=A(1/2,1)$ and $\imt(\boldsymbol{y}_1,\Omega)=B$, so that $\partial\, \img(\boldsymbol{y}_1,\Omega) \cap \imt(\boldsymbol{y}_1,\Omega)=S(1/2)$. We observe that $\boldsymbol{y}\restr{\Omega \setminus \{ \boldsymbol{0} \}}$ is a smooth diffeomorphism whose inverse $\boldsymbol{y}^{-1}_1\colon \img(\boldsymbol{y}_1,\Omega)\to \R^N$ is given by
\begin{equation}\label{eqn:mag-inv}
	\text{$\boldsymbol{y}^{-1}_1(\boldsymbol{\xi})\coloneqq (2|\boldsymbol{\xi}|-1)\,\frac{\boldsymbol{\xi}}{|\boldsymbol{\xi}|}$ \quad for all $\boldsymbol{\xi}\in \img(\boldsymbol{y}_1,\Omega)$.}
\end{equation}
Using  Corollary \ref{cor:radial-jacobian}, we compute
\begin{equation*}
	\text{$\det D  \boldsymbol{y}_1(\boldsymbol{x})=\frac{1}{4}\left (\frac{|\boldsymbol{x}|+1}{|\boldsymbol{x}|}\right )$ \quad  for all $\boldsymbol{x}\in \Omega \setminus \{ \boldsymbol{0} \}$.}
\end{equation*}
Then, applying Lemma \ref{lem:inverse-differentiable}(ii) and substituting \eqref{eqn:mag-inv}, we obtain 
\begin{equation*}
	\text{$\det \nabla \boldsymbol{y}^{-1}_1(\boldsymbol{\xi})=4 \left( \frac{|\boldsymbol{y}^{-1}_1(\boldsymbol{\xi})|}{|\boldsymbol{y}^{-1}_1(\boldsymbol{\xi})|+1} \right)=2 \left(\frac{2|\boldsymbol{\xi}|-1}{|\boldsymbol{\xi}|} \right)$ \quad for all $\boldsymbol{\xi}\in \img(\boldsymbol{y},\Omega)$. }
\end{equation*}
Define $\boldsymbol{m}_1\colon \img(\boldsymbol{y}_1,\Omega)\to \R^2_\times$ by setting
\begin{equation*}
	\text{$\boldsymbol{m}_1(\boldsymbol{\xi})\coloneqq 2 \left(\frac{2|\boldsymbol{\xi}|-1}{|\boldsymbol{\xi}|} \right) \,\frac{\boldsymbol{\xi}}{|\boldsymbol{\xi}|}$ \quad for all $\boldsymbol{\xi}\in \img(\boldsymbol{y}_1,\Omega)$.}
\end{equation*}
Clearly,  $\boldsymbol{m}_1\in W^{1,2}(\img(\boldsymbol{y}_1,\Omega);\R^2_\times)$ and \eqref{eqn:sat2} is satisfied. Also,  $\lim_{|\boldsymbol{\xi}|\to 1/2} \boldsymbol{m}_1(\boldsymbol{\xi})=\boldsymbol{0}$.  Therefore,    the extension of $\boldsymbol{m}_1$ to $\imt(\boldsymbol{y}_1,\Omega)$ by zero belongs to $W^{1,2}(\imt(\boldsymbol{y}_1,\Omega);\R^2)$, so that $(\boldsymbol{y}_1,\boldsymbol{m}_1)$ is an admissible competitor for Proposition \ref{prop:mag}(i).

Instead, define $\boldsymbol{y}_2 \colon \Omega \to \R^2$ by setting
\begin{equation*}
	\text{$\boldsymbol{y}_2(\boldsymbol{x})\coloneqq \sqrt{  |\boldsymbol{x}|^2+\frac{1}{4}}\,\frac{\boldsymbol{x}}{|\boldsymbol{x}|}$ \quad for all $\boldsymbol{x}\in \Omega \setminus \{ \boldsymbol{0}  \}$.}
\end{equation*}
 This deformation maps $\Omega$ into $B(\sqrt{5}/2)$ while opening a concentric cavity of radius $1/2$. As before, $\boldsymbol{y}_2\in \mathcal{Y}_p^{\rm cav}(\Omega)$ by   Lemma \ref{lem:radial-degree}   with  $\img(\boldsymbol{y}_2,\Omega)=A(1/2,\sqrt{5}/2)$, and $\imt(\boldsymbol{y}_2,\Omega)=B(\sqrt{5}/2)$.  By Corollary \ref{cor:radial-jacobian}, $\det D \boldsymbol{y}_2=1$ in $\Omega\setminus \{ \boldsymbol{0} \}$ and, hence, $\det \nabla \boldsymbol{y}^{-1}_2=1$ in $\img(\boldsymbol{y}_2,\Omega)$. Now, if  $\boldsymbol{m}_2\in L^2(\img(\boldsymbol{y}_2,\Omega);\R^2_\times)$ satisfies the saturation constraint \eqref{eqn:sat2},  then $|\boldsymbol{m}_2|\cong1$ in $\img(\boldsymbol{y}_2,\Omega)$. Thus, $\lim_{|\boldsymbol{\xi}|\to 1/2} |\boldsymbol{m}_2(\boldsymbol{\xi})|=1$, so that the extension of $\boldsymbol{m}_2$ to $\imt(\boldsymbol{y}_2,\Omega)$ by zero jumps on $\partial\, \img(\boldsymbol{y}_2,\Omega) \cap \imt(\boldsymbol{y}_2,\Omega)=S(1/2)$ and, in turn, cannot belong to $W^{1,2}(\imt(\boldsymbol{y}_2,\Omega);\R^2)$. In conclusion, there exists no $\boldsymbol{m}_2\in L^2(\img(\boldsymbol{y}_2,\Omega);\R^2_\times)$ satisfying the saturation constraint which is an admissible competitor for Proposition \ref{prop:mag}(i). 
\end{example}

\section*{Appendix: radial deformations}

\setcounter{section}{0}
\setcounter{theorem}{0}

\setcounter{equation}{0}

\renewcommand{\thetheorem}{A.\arabic{theorem}}

\renewcommand{\theequation}{A.\arabic{equation}}

In this appendix, we collect some results on radial deformations with respect to the $\nu$-norm. For $\nu=1,2,\infty$, maps of this kind  have been considered in the examples presented in  Sections  \ref{sec:conv}, \ref{sec:EL}, and \ref{sec:appl}.  

Let $1\leq \mu,\nu \leq \infty$ and recall the notation in \eqref{eqn:q-norm}. By the equivalence of all norms on finite-dimensional spaces,  we have
\begin{equation}\label{eqn:equivalence}
	C_1(\mu,\nu,N)|\boldsymbol{x}|_\nu \leq |\boldsymbol{x}|_\mu \leq C_2(\mu,\nu,N)|\boldsymbol{x}|_\nu \quad \text{for all $\boldsymbol{x}  \in \R^N$,}
\end{equation} 
 for constants   $C_1,C_2>0$.  
Concerning the measure of balls and spheres, employing the notation in \eqref{eqn:q-ball-annulus}--\eqref{eqn:q-sphere}, we have
\begin{equation}\label{eqn:volume-area}
	\leb(B_\nu(r))=C_3(N,\nu) r^N, \quad \haus(S_\nu(r))=C_4(N,\nu)r^{N-1} \quad \text{for all $r>0$,}
\end{equation}
for some  constants $C_3,C_4>0$ with $C_4= C_3N $.

We are concerned with  maps $\boldsymbol{u}\colon B_\nu \to \R^N$ of the form
\begin{equation}\label{eqn:radial-def}
	\boldsymbol{u}(\boldsymbol{x})\coloneqq \eta(|\boldsymbol{x}|_\nu) \frac{\boldsymbol{x}}{|\boldsymbol{x}|_\nu} \quad \text{for all $\boldsymbol{x}\in B_\nu \setminus \{ \boldsymbol{0} \}$}
\end{equation}
for some  Borel  function $\eta \colon (0,1)\to \R$.   For definiteness, we  may set $\boldsymbol{u}(\boldsymbol{0})\coloneqq \boldsymbol{0}$.  

The next result provides necessary and sufficient conditions for the Sobolev  regularity of such maps. The case $\nu=2$ is classical  \cite[Lemma 4.1]{ball.cavitation}, but we have not found the general case $1\leq \nu\leq \infty$ in the literature, so we include details here.   Partial results for $\nu=\infty$ are given in   \cite[Counterexample 7.4]{ball.murat}  and \cite[Section 11]{mueller.spector}.    

\begin{proposition}[Sobolev regularity of radial maps]\label{prop:radial-regularity}
Let $\eta\colon (0,1)\to \R$ be a  Borel  function. Also, let $1\leq \nu \leq \infty$ and $\boldsymbol{u}\colon B_\nu \to \R^N$ be defined as in \eqref{eqn:radial-def}.
Then,  $\boldsymbol{u}\in W^{1,p}(B_\nu;\R^N)$ for some $1\leq p< \infty$ if and only if $\eta$ is locally absolutely continuous and satisfies 
\begin{equation}\label{eqn:radial-integrability}
	\int_0^1 \left (   \frac{|\eta(r)|^p}{r^p}  +
	|\eta'(r)|^p \right ) r^{N-1}\,\d r<+\infty.
\end{equation}
In that case, for all $i,j=1,\dots,N$ and for almost every $\boldsymbol{x}\in B_\nu$, we have
\begin{equation}\label{eqn:weak-der-radial}
	\partial_j u^i(\boldsymbol{x})=  
	\begin{cases} \displaystyle
		 \delta_{ij} \frac{\eta(|\boldsymbol{x}|_\nu)}{|\boldsymbol{x}|_\nu}+  \frac{|\boldsymbol{x}|_\nu  \eta'(|\boldsymbol{x}|_\nu) - \eta(|\boldsymbol{x}|_\nu)}{|\boldsymbol{x}|_\nu^{\nu+1}}\, x_i\,\sgn(x_j)|x_j|^{\nu-1} 	& \text{if $1\leq \nu <\infty$,}\\[10pt] \displaystyle
		\delta_{ij} \frac{\eta(|\boldsymbol{x}|_\infty)}{|\boldsymbol{x}|_\infty}+  \frac{|\boldsymbol{x}|_\infty  \eta'(|\boldsymbol{x}|_\infty) - \eta(|\boldsymbol{x}|_\infty)}{|\boldsymbol{x}|_\infty^2} \, x_i\,\sgn(x_j)\chi_{\Delta_j}(\boldsymbol{x})  & \text{if $\nu=\infty$,}
	\end{cases}
\end{equation}
where we set $\Delta_j\coloneqq \{\boldsymbol{w}\in \R^N: \hspace{2pt}|\boldsymbol{w}|_\infty=|w_j|\}$.
\end{proposition} 
\begin{remark}\label{rem:radial-regularity}
\begin{enumerate}[(a)]
\item If $\eta$ is locally absolutely continuous, then $\boldsymbol{u}$ is continuous on $B_\nu \setminus \{\boldsymbol{0} \}$ and, for $1<\nu<\infty$, it is   
differentiable on the union of the spheres   $S_\nu(r)$ among all $0<r<1$ for which $\eta$ is differentiable at $r$. For $\nu=1$ and $\nu=\infty$, one has to additionally exclude the points where the $\nu$-norm is not differentiable. Thus,  $\boldsymbol{u}$ is almost everywhere differentiable for all $1\leq \nu \leq \infty$.
\item  For $\nu=2$, we recover  	\cite[Equation (4.3)]{ball.cavitation}, namely
\begin{equation*}
	D\boldsymbol{u}(\boldsymbol{x})=\frac{\eta(|\boldsymbol{x}|)}{|\boldsymbol{x}|} \boldsymbol{I}+ \frac{|\boldsymbol{x}|\eta'(|\boldsymbol{x}|)-\eta(|\boldsymbol{x}|)}{|\boldsymbol{x}|^3}\, \boldsymbol{x}\otimes \boldsymbol{x}\quad \text{for almost all $\boldsymbol{x}\in B$},
\end{equation*}	
where we recall \eqref{eqn:B-S}.
\end{enumerate}
\end{remark}
\begin{proof}
First of all,  for every $j=1,\dots,N$, we have
\begin{equation*}
	\partial_j(|\boldsymbol{x}|_\nu)=\begin{cases} \displaystyle 
		\frac{\sgn(x_j)|x_j|^{\nu-1}}{|\boldsymbol{x}|_\nu^{\nu-1}} & \text{if $1\leq \nu<\infty$,}\\[10pt]
		\sgn(x_j)\chi_{\Delta_j}(\boldsymbol{x}) & \text{if $\nu=\infty$,}
	\end{cases} \quad \text{for all $\boldsymbol{x}\neq \boldsymbol{0}$.}
\end{equation*}	
The map $\boldsymbol{x}\mapsto \frac{\boldsymbol{x}}{|\boldsymbol{x}|_\nu}$ belongs to $W^{1,q}_\loc(\R^N;\R^N)$ for all $1\leq q<N$. In particular, for every $i,j=1,\dots,N$, we have
\begin{equation*}
	\partial_j \left(\frac{x_i}{|\boldsymbol{x}|_\nu}\right)=\begin{cases}
		\displaystyle \frac{\delta_{ij}}{|\boldsymbol{x}|_\nu}-\frac{x_i \sgn(x_j)|x_j|^{\nu-1}}{|\boldsymbol{x}|_{\nu}^{\nu+1}} & \text{if $1\leq \nu <\infty$,}\\[10pt]
		 \displaystyle \frac{\delta_{ij}}{|\boldsymbol{x}|_\nu}-\frac{x_i \sgn(x_j)\chi_{\Delta_j}(\boldsymbol{x})}{|\boldsymbol{x}|_{\infty}^2} & \text{if $ \nu =\infty$,}
	\end{cases}\quad \text{for all $\boldsymbol{x} \neq \boldsymbol{0}$.}
\end{equation*}

The proof is divided into two steps.

\textbf{Step 1 (Necessity).}	
	Assume that $\boldsymbol{u}\in W^{1,p}(B_\nu;\R^N)$.  First, we show that $\eta \in W^{1,p}_{\loc}(0,1)$.
	Consider the map  $\boldsymbol{\Upsilon}\colon (0,+\infty)\times (-\pi/2,\pi/2)^{N-1} \to \R^N$ defined as
	\begin{align*}
		\Upsilon^i(r,\boldsymbol{\vartheta})\coloneqq r \prod_{k=1}^{i-1} \sin \vartheta_k\,\cos \vartheta_i \quad \text{for all $i= 1,  \dots,N-1$,} \quad \quad 
		\Upsilon^N(r,\boldsymbol{\vartheta})\coloneqq r \prod_{k=1}^{N-1} \sin \vartheta_k,
	\end{align*} 
	where $\boldsymbol{\vartheta}=(\vartheta_1,\dots,\vartheta_{N-1})$.  The map $\boldsymbol{\Upsilon}$ is a parametrization of the upper {hemisphere} of $S$.   Note that $\boldsymbol{\Upsilon}(r,\boldsymbol{\vartheta})=r\boldsymbol{e}_{\boldsymbol{\vartheta}}$ with $\boldsymbol{e}_{\boldsymbol{\vartheta}}\coloneqq \boldsymbol{\Upsilon}(1,\boldsymbol{\vartheta})\in S$.
	The map $\boldsymbol{\Upsilon}$ is injective and Lipschitz, and so is its inverse. 
	 Let $0<R_1<R_2<1$. For  $0<\delta<\pi/2$ sufficiently small, the set $\boldsymbol{\Upsilon} \big( (R_1,R_2)\times (-\delta,\delta)^{N-1} \big)$ is a portion of a spherical cap which is contained in $B_\nu$. In this case, $\boldsymbol{u} \circ \boldsymbol{\Upsilon} \in W^{1,p}((R_1,R_2)\times (-\delta,\delta)^{N-1} ; \R^N)$ by \cite[Theorem~11.53]{leoni}. Then, by \cite[Theorem~11.45]{leoni},  the map $v_{\boldsymbol{\vartheta}}\colon (R_1,R_2) \to \R$ defined as $v_{\boldsymbol{\vartheta}}(r)\coloneqq \boldsymbol{u}(\boldsymbol{\Upsilon}(r,\boldsymbol{\vartheta}))\cdot \boldsymbol{e}_{\boldsymbol{\vartheta}}$ satisfies $v_{\boldsymbol{\vartheta}}\in W^{1,p}(R_1,R_2)$ for $\mathscr{L}^{N-1}$-almost every $\boldsymbol{\vartheta}\in (-\delta,\delta)^{N-1}$. Choosing any such $\boldsymbol{\vartheta}$, we observe that
	\begin{equation*}
		v_{\boldsymbol{\vartheta}}(r)=\boldsymbol{u}(r\boldsymbol{e}_{\boldsymbol{\vartheta}})  \cdot \boldsymbol{e}_{\boldsymbol{\vartheta}}  =\frac{1}{|\boldsymbol{e}_{\boldsymbol{\vartheta}}|_\nu}\eta(r |\boldsymbol{e}_{\boldsymbol{\vartheta}}|_\nu) \quad \text{for all $r\in (R_1,R_2)$.}
	\end{equation*}
	Therefore,  we deduce that $\eta \in W^{1,p}(R_1 |\boldsymbol{e}_{\boldsymbol{\vartheta}}|_\nu, R_2|\boldsymbol{e}_{\boldsymbol{\vartheta}}|_\nu)$.  Since $||\boldsymbol{e}_{\boldsymbol{\vartheta}}|_\nu - 1| \le c_\delta$ with $c_\delta \to 0$ as $\delta \to 0^+$,  this entails $\eta \in W^{1,p}(R_1(1+c_\delta), R_2(1-c_\delta))$.  As $R_1,R_2 \in (0,1)$   were arbitrary,  by sending $\delta \to 0^+$,  this proves that $\eta \in W^{1,p}_{\loc}(0,1)$ and, in turn,   that  $\eta$ is locally absolutely continuous.
	
	At this point, we know that $\boldsymbol{u}$ is almost everywhere differentiable in $B_\nu$, see Remark~\ref{rem:radial-regularity}(a). Since we assume $\boldsymbol{u}\in W^{1,p}(B_\nu;\R^N)$, its pointwise  and weak derivatives coincide. Thus, by  a  direct computation, we determine the expression in \eqref{eqn:weak-der-radial} for the weak derivatives.  In particular,
	 \begin{equation}\label{eqn:radial-div}
		\div \,  \boldsymbol{u}(\boldsymbol{x})=(N-1)\frac{\eta(|\boldsymbol{x}|_\nu)}{|\boldsymbol{x}|_\nu}+\eta'(|\boldsymbol{x}|_\nu) \quad \text{for almost all $\boldsymbol{x}\in B_\nu$.}
	\end{equation}
	Using \eqref{eqn:volume-area} and the coarea formula, we obtain
	\begin{equation}\label{eqn:g-lp-weight}
		\int_{0}^{1} \left |(N-1)\frac{\eta(r)}{r}+\eta'(r)\right |^p r^{N-1}\,\d r=C_4^{-1} \int_{B_\nu} |\div \,  \boldsymbol{u}|^p\d\boldsymbol{x}<+\infty.
	\end{equation}
	However, the previous estimate does not  yet  yield \eqref{eqn:radial-integrability} since both functions inside the modulus may have a sign.

	Define $a\colon B_\nu \to (0,+\infty)$ as $a(\boldsymbol{x})\coloneqq |\boldsymbol{u}(\boldsymbol{x})|_\nu$. By approximating $\boldsymbol{u}$ with smooth functions, we check that $a\in W^{1,p}(B_\nu)$ and the chain rule holds. Thus, for every $j=1,\dots, N$ and for almost all $\boldsymbol{x}\in B_\nu$, we have 
	\begin{equation*}
		|\partial_j a(\boldsymbol{x} )|=
		\begin{cases} \displaystyle
			\frac{|x_j|^{\nu-1}}{|\boldsymbol{x}|_\nu^{\nu-1}}|\eta'(|\boldsymbol{x}|_\nu)| & \text{if $1\leq \nu <\infty$,}\\[10pt] \displaystyle
			\chi_{\Delta_j}(\boldsymbol{x})|\eta'(|\boldsymbol{x}|_\infty)| & \text{if $\nu=\infty$.}
		\end{cases}
	\end{equation*}
	In particular,
	\begin{equation*}
		|D a(\boldsymbol{x})|_1=\sum_{j=1}^{N} |\partial_j a(\boldsymbol{x})|= \begin{cases} \displaystyle 
			\frac{\sum_{j=1}^{N} |x_j|^{\nu -1}}{|\boldsymbol{x}|_\nu^{\nu -1}}|\eta'(|\boldsymbol{x}|_\nu)| & \text{if $1\leq \nu <\infty$,}\\[10pt] 
			\sum_{j=1}^{N} \chi_{\Delta_j}(\boldsymbol{x})|\eta'(|\boldsymbol{x}|_\infty)|	& \text{if $\nu=\infty$}.
		\end{cases} 
	\end{equation*}
	Since the map  $\mu \mapsto ( \sum_{j=1}^{N} |x_j|^\mu)^{1/\mu}$ is nonincreasing and we have  the identity $\sum_{j=1}^{N}\chi_{\Delta_j}(\boldsymbol{x})=1$ for almost all $\boldsymbol{x}\in \R^N$, we deduce $|Da(\boldsymbol{x})|_1\geq|\eta'(|\boldsymbol{x}|_\nu)|$ for almost all $\boldsymbol{x}\in B_\nu$. Using the coarea formula together with \eqref{eqn:equivalence}--\eqref{eqn:volume-area}, this inequality yields
	\begin{equation}
		\label{eqn:f'-lp-weight}
		\begin{split}
			\int_0^1 |\eta'(r)|^p\,r^{N-1}\,\d r=C_4^{-1}\int_{ B_\nu} |\eta'(|\boldsymbol{x}|_\nu)|^p\,\d \boldsymbol{x} \leq C_1^{-p} C_4^{-1} \int_{B_\nu} |D a(\boldsymbol{x})|^p\,\d\boldsymbol{x}<+\infty.
		\end{split}
	\end{equation}  
	Then, combining \eqref{eqn:g-lp-weight}--\eqref{eqn:f'-lp-weight}, we obtain
	\begin{equation*}
		\begin{split}
			\int_0^1\left | \frac{\eta(r)}{r} \right |^p r^{N-1}\,\d r<+\infty,
		\end{split}
	\end{equation*}
	so that \eqref{eqn:radial-integrability} is proved.

	 \textbf{Step 2 (Sufficiency).} Assume that $\eta$ is  locally  absolutely continuous and satisfies \eqref{eqn:radial-integrability}. 
	 Using \eqref{eqn:equivalence}--\eqref{eqn:volume-area} and the coarea formula, we estimate
	 \begin{equation*}
	 	\int_{B_\nu} |\boldsymbol{u}(\boldsymbol{x})|^p\,\d\boldsymbol{x}\leq C_2^p\int_{B_\nu} |\eta(|\boldsymbol{x}|_\nu)|^p\,\d\boldsymbol{x}=C_2^pC_4 \int_0^1 |\eta(r)|^p r^{N-1}\,\d r<+\infty,
	 \end{equation*} 
	 where the integral on the right-hand side is finite in view of \eqref{eqn:radial-integrability}. Thus, $\boldsymbol{u}\in L^p(B_\nu;\R^N)$.
	 
	 We extend $\eta$  by zero outside of $(0,1)$  without renaming it and we consider the family $(\rho_\varepsilon)_{\varepsilon>0}$ of standard radial mollifiers. By assumption, $\eta\in W^{1,p}_\loc(0,1)$. Setting $\eta_\varepsilon \coloneqq \eta \ast \rho_\varepsilon$, we have $\eta_\varepsilon \to \eta$  in $W^{1,p}_\loc(0,1)$. Define $\boldsymbol{u}_\varepsilon \colon B_\nu \setminus \{\boldsymbol{0}\} \to \R^N$ by setting
	 \begin{equation*}
	 	\boldsymbol{u}_\varepsilon(\boldsymbol{x})\coloneqq \eta_\varepsilon(|\boldsymbol{x}|_\nu)\frac{\boldsymbol{x}}{|\boldsymbol{x}|_\nu} \quad \text{for all $\boldsymbol{x}\in B_\nu\setminus \{ \boldsymbol{0} \}$.}
	 \end{equation*}  
	 Then, $\boldsymbol{u}_\varepsilon \in W^{1,\infty}(B_\nu \setminus \{\boldsymbol{0}\};\R^N)$. 
	 Let $0<R_1<R_2<1$. By \eqref{eqn:equivalence} and the coarea formula
	 \begin{equation*}
	 	\begin{split}
	 		\int_{A_\nu(R_1,R_2)} |\boldsymbol{u}_\varepsilon(\boldsymbol{x})-\boldsymbol{u}(\boldsymbol{x}) |^p\,\d\boldsymbol{x}&\leq    C_1^{-p}   \int_{A_\nu(R_1,R_2)} |\boldsymbol{u}_\varepsilon(\boldsymbol{x})-\boldsymbol{u}(\boldsymbol{x}) |^p_\nu\,\d\boldsymbol{x} =C_1^{-p} \int_{R_1}^{R_2} |\eta_\varepsilon(r)-\eta(r)|^p r^{N-1}\,\d r\\
	 		&\leq  C_1^{-p} R_2^{N-1} \int_{R_1}^{R_2} |\eta_\varepsilon(r)-\eta(r)|^p\,\d r.
	 	\end{split}
	 \end{equation*}
	 For every $i,j=1,\dots,N$ and $\boldsymbol{x}\in A_\nu(R_1,R_2)$, we have
	 \begin{equation*} \label{eqn:weak-der-radial-eps}
	 	\partial_j u_\varepsilon^i(\boldsymbol{x})=  
	 	\begin{cases} \displaystyle
	 		\delta_{ij}\frac{\eta_\varepsilon(|\boldsymbol{x}|_\nu)}{|\boldsymbol{x}|_\nu}  + \frac{ |\boldsymbol{x}|_\nu \eta'_\varepsilon(|\boldsymbol{x}|_\nu) - \eta_\varepsilon(|\boldsymbol{x}|_\nu) }{|\boldsymbol{x}|_\nu^{\nu+1}} \, x_i \sgn(x_j)|x_j|^{\nu-1} 	& \text{if $1\leq \nu <\infty$,}\\[10pt] \displaystyle
	  \delta_{ij} 		\frac{\eta_\varepsilon(|\boldsymbol{x}|_\infty)}{|\boldsymbol{x}|_\infty}  + \frac{ |\boldsymbol{x}|_\infty \eta'_\varepsilon(|\boldsymbol{x}|_\infty)  -\eta_\varepsilon(|\boldsymbol{x}|_\infty) }{|\boldsymbol{x}|_\infty^2}\, x_i \sgn(x_j)\chi_{\Delta_j}(\boldsymbol{x})& \text{if $\nu=\infty$.}
	 	\end{cases} 
	 \end{equation*} 
		Define $w^i_j\colon B_\nu \to \R$ with $w^i_j(\boldsymbol{x})$ being given by the right-hand side of \eqref{eqn:weak-der-radial}. Employing again \eqref{eqn:volume-area} and the coarea formula, we estimate 
	 \begin{equation*}
	 	\begin{split}
	 		\int_{A_\nu(R_1,R_2)} |\partial_j u^i_\varepsilon(\boldsymbol{x})-w^i_j(\boldsymbol{x})|^p\,\d\boldsymbol{x}
	 		&\leq C(N,p)  \int_{A_\nu(R_1,R_2)} \left( 1 + \frac{|x_i|\, |x_j|^{\nu-1}}{|\boldsymbol{x}|_\nu^\nu}\right)^p \frac{|\eta_\varepsilon(|\boldsymbol{x}|_\nu)-\eta(|\boldsymbol{x}|_\nu)|^p}{|\boldsymbol{x}|_\nu^p}\,\d\boldsymbol{x}\\
	 		& \ \ \ +C(N,p) \int_{A_\nu(R_1,R_2)} \left( \frac{|x_i|\, |x_j|^{\nu-1}}{|\boldsymbol{x}|_\nu^\nu} \right)^p |\eta_\varepsilon'(|\boldsymbol{x}|_\nu)-\eta'(|\boldsymbol{x}|_\nu)|^p\,\d\boldsymbol{x}\\
	 		&\leq C(N,p,\nu,R_1,R_2) \int_{A_\nu(R_1,R_2)} |\eta_\varepsilon(|\boldsymbol{x}|_\nu)-\eta(|\boldsymbol{x}|_\nu)|^p \,\d\boldsymbol{x} \\
	 		& \ \ \ + C(N,p,\nu,R_1,R_2) \int_{A_\nu(R_1,R_2)} |\eta'_\varepsilon(|\boldsymbol{x}|_\nu) -\eta'(|\boldsymbol{x}|_\nu)|^p \,\d\boldsymbol{x}\\
	 		&\leq C(N,p,\nu,R_1,R_2) \, C_4 \int_{R_1}^{R_2}  \Big (|\eta_\varepsilon(r)-\eta(r)|^p+|\eta'_\varepsilon(r)-\eta'(r)|^p \Big ) r^{N-1}\,\d r\\
	 		&\leq C(N,p,\nu,R_1,R_2) \, C_4 R_2^{N-1} \int_{R_1}^{R_2}  \Big (|\eta_\varepsilon(r)-\eta(r)|^p+|\eta'_\varepsilon(r)-\eta'(r)|^p \Big ) \,\d r.
	 	\end{split}
	 \end{equation*}
	 The previous two estimates show that $\boldsymbol{u}_\varepsilon \to \boldsymbol{u}$ in $W^{1,p}_\loc(B_\nu\setminus \{\boldsymbol{0}\};\R^N)$. Thus,  $\boldsymbol{u} \in W^{1,p}_\loc(B_\nu \setminus \{\boldsymbol{0}\};\R^N)$ with weak derivatives as in \eqref{eqn:weak-der-radial}.  Applying  \eqref{eqn:volume-area},  \eqref{eqn:weak-der-radial},  and the coarea formula, we obtain
	 \begin{equation*}
	 	\begin{split}
	 		\int_{B_\nu} |D \boldsymbol{u}(\boldsymbol{x})|^p\,\d \boldsymbol{x}&\leq C(N,p) \int_{B_\nu} \left( \frac{|\eta(|\boldsymbol{x}|_\nu)|^p}{|\boldsymbol{x}|_\nu^p}+|\eta'(|\boldsymbol{x}|_\nu)|^p \right)\,\d\boldsymbol{x}\\
	 		&=C(N,p) C_4 \int_0^1 \left( \frac{|\eta(r)|^p}{r^p} + |\eta'(r)|^p \right)r^{N-1}\,\d r,
	 	\end{split}
	 \end{equation*}
	 where the integral on the right-hand side is finite by \eqref{eqn:radial-integrability}. This shows that $D\boldsymbol{u}\in L^p(B_\nu;\rnn)$.      Finally, we need to check that $\boldsymbol{u}\in W^{1,p}(B_\nu;\R^N)$. This can be done in a classical way via integration by parts, we include a short alternative proof here.  We consider the sequence  $(\boldsymbol{u}_n)_n\subset GSBV^p(B_\nu;\R^N)$ defined by $\boldsymbol{u}_n \coloneqq \boldsymbol{u} \chi_{A(1/n,1)}$. We can apply Theorem \ref{thm:ambrosio-compactness} to see that $\boldsymbol{u} \in GSBV^p(B_\nu;\R^N)$. Since $$\haus(J_{\boldsymbol{u}})\leq \liminf_{n \to \infty} \haus(J_{\boldsymbol{u}_n}) \le \liminf_{n \to \infty} C_4 \left (\frac{1}{n} \right )^{N-1} = 0,$$
	 and $\boldsymbol{u}\in L^\infty(\Omega;\R^N)$, from \eqref{eqn:SBVcapLinfty} we immediately find $\boldsymbol{u} \in W^{1,p}(B_\nu;\R^N)$. 
\end{proof}

From Proposition~\ref{prop:radial-regularity}, we deduce the following formula for the Jacobian determinant of radial deformations.

\begin{corollary}[Jacobian determinant of radial  maps]\label{cor:radial-jacobian}
	Let $\eta\colon (0,1)\to \R$ be a locally absolutely continuous function  satisfying \eqref{eqn:radial-integrability} for some $1\leq p<\infty$. Also, let $1\leq \nu\leq \infty$ and $\boldsymbol{u}\colon B_\nu\to \R^N$ be defined as in \eqref{eqn:radial-def}.
	Then
	\begin{equation}\label{eqn:radial-jacobian}
		\det D \boldsymbol{u}(\boldsymbol{x})= \left(\frac{\eta(|\boldsymbol{x}|_\nu)}{|\boldsymbol{x}|_\nu}\right)^{N-1}\eta'(|\boldsymbol{x}|_\nu) \quad \text{for almost all $\boldsymbol{x}\in  B_\nu$.}
	\end{equation} 
	In particular, $\det D \boldsymbol{u}\in L^1(B_\nu)$ if and only if  $\eta$  satisfies
	\begin{equation}\label{eqn:radial-jacobian-coarea}
		\int_0^1 \left( {\eta(r)}\right)^{N-1}\eta'(r)\,\d r<+\infty.
	\end{equation}
\end{corollary}
\begin{remark}\label{rem:radial-jacobian}
\begin{enumerate}[(a)]
	\item From \eqref{eqn:radial-jacobian}, we see that $\det D \boldsymbol{u}(\boldsymbol{x})\neq 0$  if and only if both $\eta(|\boldsymbol{x}|_\nu)$ and $\eta'(|\boldsymbol{x}|_\nu)$ are nonzero. If we assume that $\eta$ is positive, then $\det D \boldsymbol{u}>0$ almost everywhere in $B_\nu$ if and only if $\eta'$ is  positive almost everywhere or, equivalently, if  $\eta$ is strictly increasing.
	\item For $\nu=2$, we recover the formula in  \cite[p.~568]{ball.cavitation}. 
\end{enumerate}	
\end{remark}
\begin{proof}
Observe that $\boldsymbol{u}\in W^{1,p}(B_\nu;\R^N)$ by Proposition \ref{prop:radial-regularity}. Once  \eqref{eqn:radial-jacobian} is proven, the equivalence between $\det D \boldsymbol{u}\in L^1(B_\nu)$ and  \eqref{eqn:radial-jacobian-coarea} follows by the coarea formula.
 For convenience, we set
\begin{equation*}
	\boldsymbol{p}_\nu(\boldsymbol{x})\coloneqq \displaystyle 
	\begin{cases} 
		\frac{1}{|\boldsymbol{x}|_\nu^{\nu-1}}\sum_{k=1}^{N}\sgn(x_k)|x_k|^{\nu-1}\boldsymbol{e}_k & \text{if $1\leq \nu <\infty$,}\\[10pt]
		\sum_{k=1}^{N}\sgn(x_k)\chi_{\Delta_k}(\boldsymbol{x})\boldsymbol{e}_k & \text{if $\nu=\infty$,} 
	\end{cases} \quad \text{for all $\boldsymbol{x}\in\R^N$.}
\end{equation*}
With this notation,  \eqref{eqn:weak-der-radial} gives
\begin{equation}
	\label{eqn:weak-der-radial-p}
	D\boldsymbol{u}(\boldsymbol{x})=
		\displaystyle \frac{\eta(|\boldsymbol{x}|_\nu)}{|\boldsymbol{x}|_\nu}\,\boldsymbol{I}+\frac{|\boldsymbol{x}|_\nu \eta'(|\boldsymbol{x}|_\nu)-\eta(|\boldsymbol{x}|_\nu)}{|\boldsymbol{x}|_\nu^{2}}\,  \boldsymbol{x}\otimes \boldsymbol{p}_{\nu}(\boldsymbol{x})
	 \quad \text{for almost all $\boldsymbol{x}\in B_\nu$.}
\end{equation} 
Let $0<r<1$ be such that $\eta$ is differentiable at $r$. At almost all $\boldsymbol{x}\in S_\nu(r)$, the map $\boldsymbol{u}$ is differentiable  by Remark \ref{rem:radial-regularity}(a) and the tangent space of $S_\nu(r)$ at that point is well defined as the hyperplane orthogonal  to $\boldsymbol{p}_{\nu}(\boldsymbol{x})$. Let $\boldsymbol{\tau}_1,\dots,\boldsymbol{\tau}_{N-1} \in \R^N$ form an orthonormal basis of this tangent space. For $j=1,\dots,N-1$, using \eqref{eqn:weak-der-radial-p}, we compute the tangential derivatives
\begin{equation*}
	\frac{\partial \boldsymbol{u}}{\partial \boldsymbol{\tau}_j}(\boldsymbol{x})=\big (D\boldsymbol{u}(\boldsymbol{x})\big )\boldsymbol{\tau}_j=\frac{\eta(r)}{r}\,\boldsymbol{\tau}_j.
\end{equation*}
This equation shows that each $\boldsymbol{\tau}_j$ is an eigenvector of $D\boldsymbol{u}(\boldsymbol{x})$ with corresponding eigenvalue $\lambda_j\coloneqq \frac{\eta(r)}{r}$. Let $\lambda_N$ be the last eigenvalue of $D\boldsymbol{u}(\boldsymbol{x})$. By \eqref{eqn:radial-div}, we have
\begin{equation*}
	\lambda_N=   {\rm tr} \big( D\boldsymbol{u}(\boldsymbol{x}) \big) -\sum_{j=1}^{N-1} \lambda_j =  \div \,  \boldsymbol{u}(\boldsymbol{x})-\sum_{j=1}^{N-1} \lambda_j=\eta'(r).
\end{equation*}
Therefore,
\begin{equation*}
	\det D \boldsymbol{u}(\boldsymbol{x})=\prod_{j=1}^N \lambda_j=\left( \frac{\eta(r)}{r}\right)^{N-1} \eta'(r).
\end{equation*}
 This concludes the proof. 
\end{proof}

Henceforth, we will assume  $\eta$ to be positive.  Observe that deformations violating this  condition correspond to  eversions    and, in turn,   exhibit interpenetration of matter   in the sense that they  do  not fulfill condition (INV). See \cite[Remark~4, p.17]{mueller.spector} for a related example.  

We make the following simple observation about the injectivity of radial maps.

\begin{lemma}[Injectivity of radial maps]\label{lem:radial-injective}
Let  $\eta\colon (0,1)\to (0,+\infty)$ be a continuous function. Also, let $1\leq \nu\leq \infty$ and $\boldsymbol{u}\colon B_\nu \to \R^N$ be defined as in \eqref{eqn:radial-def}. Then,   $\boldsymbol{u}\restr{B_\nu \setminus \{ \boldsymbol{0} \}}$  is injective if and only if $\eta$ is strictly monotone.    
\end{lemma}
\begin{proof}
First, assume that  $\eta$   is strictly monotone. Let $\boldsymbol{x}_1,\boldsymbol{x}_2\in B_\nu \setminus \{\boldsymbol{0}\}$. If $\boldsymbol{u}(\boldsymbol{x}_1)=\boldsymbol{u}(\boldsymbol{x}_2)$, then, taking the $\nu$-norm at both sides, we obtain $\eta(|\boldsymbol{x}_1|_\nu)=\eta(|\boldsymbol{x}_2|_\nu)$. Since $\eta$ is injective, this yields $|\boldsymbol{x}_1|_\nu=|\boldsymbol{x}_2|_\nu$. Hence, from the identity  $\boldsymbol{u}(\boldsymbol{x}_1)=\boldsymbol{u}(\boldsymbol{x}_2)$, we conclude $\boldsymbol{x}_1=\boldsymbol{x}_2$.
Conversely, assume that $\boldsymbol{u}\restr{B_\nu \setminus \{\boldsymbol{0}\}}$ is injective. Let $0<r_1,r_2<1$ be such that $\eta(r_1)=\eta(r_2)$. Set $\boldsymbol{x}_i\coloneqq r_i \frac{\boldsymbol{e}_1}{|\boldsymbol{e}_1|_\nu}$ for $i=1,2$. Then, $\boldsymbol{u}(\boldsymbol{x}_1)=\boldsymbol{u}(\boldsymbol{x}_2)$, which yields $\boldsymbol{x}_1=\boldsymbol{x}_2$ by assumption.  In particular,  $r_1=|\boldsymbol{x}_1|_\nu=|\boldsymbol{x}_2|_\nu=r_2$. Therefore, $\eta$ is injective and, in turn, strictly monotone.
\end{proof}

For ease of exposition, we register a preliminary result on the topological degree of Sobolev deformations in the supercritical regime. 

\begin{lemma}[Degree of Sobolev maps with supercritical integrability]
	\label{lem:degree-supercritical}
Let $A \subset \R^N$ be a bounded Lipschitz domain, let $\boldsymbol{h}\in W^{1,q}(A;\R^N) \cap C^0(\closure{A};\R^N)$ with $N<q\leq \infty$ and $\det D \boldsymbol{h}\geq 0$ almost everywhere. Setting $\Lambda_{\boldsymbol{h}}\coloneqq \{ \det D \boldsymbol{h}=0  \}$, suppose that $\boldsymbol{h}$ is almost everywhere injective on $A \setminus  \Lambda_{\boldsymbol{h}}$. 
Then,
\begin{equation*}
	\deg(\boldsymbol{h},A,\cdot)=\chi_V \quad \text{in $\R^N \setminus \boldsymbol{h}(\partial A)$,}
\end{equation*}
where $V\coloneqq \boldsymbol{h}(A) \setminus \boldsymbol{h}(\partial A)$. In particular, $\imt(\boldsymbol{h},A)=V$.  
\end{lemma}
\begin{proof}
First, observe that $\deg(\boldsymbol{h},A,\cdot)=0$ on $\R^N \setminus \closure{V}$ because of Remark~\ref{rem:top-deg}(d). 
For $\boldsymbol{\xi}_0\in V$,  consider a nonnegative function $\psi \in C^0_{\rm c}(\R^N)$ supported in the connected component of $V$ containing $\boldsymbol{\xi}_0$ and having unit integral. Then, using the integral formula for the degree \cite[Theorem~5.38]{fonseca.gangbo}  together with Proposition~\ref{prop:change-of-variable}, we compute
\begin{equation*}
	\begin{split}
		\deg(\boldsymbol{h},A,\boldsymbol{\xi}_0)&=\int_A \psi \circ \boldsymbol{h}\,\det D \boldsymbol{h}\,\d \boldsymbol{x}=\int_{A \setminus  \Lambda_{\boldsymbol{h}}  } \psi \circ \boldsymbol{h}\,\det D \boldsymbol{h}\,\d \boldsymbol{x}  
		=\int_{V \setminus \boldsymbol{h}(\Lambda_{\boldsymbol{h}}) } \psi \,\d\boldsymbol{\xi}=\int_{V  } \psi \,\d\boldsymbol{\xi}=1,
	\end{split}
\end{equation*}
where  the second last equality  follows as $\leb(\boldsymbol{h}(\Lambda_{\boldsymbol{h}}))=0$ by \cite[Corollary~3]{marcus.mizel}.  
\end{proof}

In the next result, we compute the topological degree of deformations as in \eqref{eqn:radial-def}.  Note that, at this stage, we cannot resort to Lemma~\ref{lem:orientation} as we do not know yet whether  maps as in \eqref{eqn:radial-def} satisfy condition (INV).

\begin{lemma}[Degree of radial  maps]
	\label{lem:radial-degree}
Let $\eta \colon (0,1) \to (0,+\infty)$ be a locally {absolutely} continuous  and strictly increasing function satisfying \eqref{eqn:radial-integrability} for some $p>N-1$ and \eqref{eqn:radial-jacobian-coarea}. Also, let $1\leq \nu\leq \infty$ and $\boldsymbol{u}\colon B_\nu \to \R^N$ be defined as in \eqref{eqn:radial-def}. Then, $\boldsymbol{u}\in W^{1,p}(B_\nu;\R^N)$ with $\det D \boldsymbol{u}\in L^1_+(B_\nu)$. Moreover, 
 for every domain $U \subset \subset \Omega$  with  $\partial U \subset \R^N \setminus \{  \boldsymbol{0} \}$ such that $\R^N \setminus \partial U$ has exactly two connected components, there holds
\begin{equation}
	\label{eqn:radial-deg}
	\deg(\boldsymbol{u},U,\cdot)=\chi_V \quad \text{in $\R^N \setminus \boldsymbol{u}(\partial U)$,}
\end{equation}
where $V$ denotes the unique bounded connected component of $\R^N \setminus \boldsymbol{u}(\partial U)$, and also
\begin{equation}
	\label{eqn:radial-inclusion}
	\imt(\boldsymbol{u},U)=V, \qquad \boldsymbol{u}(U\setminus \{ \boldsymbol{0} \}) \subset V, \qquad \boldsymbol{u}(B_\nu \setminus (U  \cup \{\boldsymbol{0}\}))\subset \R^N \setminus V.
\end{equation} 
In particular,   $\boldsymbol{u}$ satisfies condition {\rm (INV)}.	
\end{lemma}
\begin{remark}
Suppose that $\closure{U} \subset \R^N \setminus \{ \boldsymbol{0} \}$. By Remark~\ref{rem:radial-regularity}(a), 	 $\boldsymbol{u}\restr{U}\in C^0(\closure{U};\R^N)$ and this map is injective by Lemma~\ref{lem:radial-injective}. Thus, $\boldsymbol{u}\restr{U}$ is a homeomorphism by the invariance of domain  theorem \cite[Theorem~3.30]{fonseca.gangbo}. In this case,  combining \eqref{eqn:radial-deg} with \cite[Theorem~3.35]{fonseca.gangbo}, we deduce that $\boldsymbol{u}(U)=V$.
\end{remark}
\begin{proof}
Note that $\boldsymbol{u}\in W^{1,p}(B_\nu;\R^N)$ with $\det D \boldsymbol{u}\in L^1_+(B_\nu)$ by Proposition~\ref{prop:radial-regularity} and Corollary ~\ref{cor:radial-jacobian}. Consider $U$ as in the statement.  
Since $\boldsymbol{u}\restr{\partial U}$ is injective by Lemma~\ref{lem:radial-injective}, by Jordan's separation theorem \cite[Theorem~3.29]{fonseca.gangbo}, the set $\R^N \setminus \boldsymbol{u}(\partial U)$ has exactly two   connected components: $V$ and $V_0$, where the  former  is bounded and the latter is unbounded.

Let $0<\delta \ll 1$ be such that $\partial U \subset A_\nu(\delta,1-\delta)$. For  $\varepsilon \in (0,\delta] $, define $\eta_\varepsilon \colon [0,1-\varepsilon] \to [0,+\infty)$ and $\boldsymbol{u}_\varepsilon\colon {B}_\nu(1-\varepsilon) \to \R^N$ as
\begin{equation}
	\label{eqn:ueps}
	\eta_\varepsilon(r)\coloneqq \begin{cases} \frac{\eta(\varepsilon)}{\varepsilon}r & \text{if $0\leq r\leq \varepsilon$,}\\
		\eta(r) & \text{if $\varepsilon \leq r \leq 1-\varepsilon$,}
		
	\end{cases}  \qquad \boldsymbol{u}_\varepsilon (\boldsymbol{x})\coloneqq \eta_\varepsilon(|\boldsymbol{x}|_\nu) \frac{\boldsymbol{x}}{|\boldsymbol{x}|_\nu}.
\end{equation}
We have  $\boldsymbol{u}_\varepsilon \in C^0(\closure{B}_\nu(1-\varepsilon);\R^N)$  and $\boldsymbol{u}_\varepsilon \restr{\partial {U}}=\boldsymbol{u}\restr{\partial {U}}$. Observe that  $\boldsymbol{u}_\varepsilon$ is a homeomorphism as a consequence of Lemma~\ref{lem:radial-injective} and the invariance of domain theorem \cite[Theorem~3.30]{fonseca.gangbo}. Then, thanks to Remark~\ref{rem:top-deg}(d) and \cite[Theorem~3.35]{fonseca.gangbo}, we realize that  $\boldsymbol{u}_\varepsilon(U)=V$ as well as $V_0=\R^N \setminus \boldsymbol{u}_\varepsilon(\closure{U})$.  
At this point,  \eqref{eqn:radial-deg} is proved if we show the that
\begin{equation}
	\label{eqn:deg1}
	\deg(\boldsymbol{u},U,\cdot)=\chi_{\boldsymbol{u}_\varepsilon(U)} \quad \text{in $\R^N \setminus \boldsymbol{u}(\partial U)$.}
\end{equation}
By Definition~\ref{def:top-deg}, we have
\begin{equation}
	\label{eqn:deg-deg}
	\deg(\boldsymbol{u},U,\cdot)=\deg(\boldsymbol{u}_\varepsilon,U,\cdot) \quad \text{in $\R^N \setminus \boldsymbol{u}(\partial U)$.}
\end{equation}
Observe that $\deg(\boldsymbol{u}_\varepsilon,U,\cdot)=0$ on $V_0$  by \eqref{eqn:deg1}--\eqref{eqn:deg-deg}, which gives $\deg(\boldsymbol{u},U,\cdot)=0$ on $\R^N \setminus \boldsymbol{u}_\varepsilon(\closure{U})$ by \eqref{eqn:deg-deg}.
Let $\boldsymbol{\xi} \in \boldsymbol{u}_\varepsilon(U)$.
By \cite[Theorem~2.7(1)]{fonseca.gangbo}, we have
\begin{equation}
	\label{eqn:additivity}
	\deg({\boldsymbol{u}}_\varepsilon,B_\nu(1-\varepsilon),\boldsymbol{\xi})=\deg({\boldsymbol{u}}_\varepsilon,U,\boldsymbol{\xi})+ \deg({\boldsymbol{u}}_\varepsilon,B_\nu(1-\varepsilon)\setminus \closure{U},\boldsymbol{\xi}).
\end{equation}
We claim that 
\begin{equation}
	\label{eqn:deg-a}
	\deg({\boldsymbol{u}}_\varepsilon,B_\nu(1-\varepsilon),\boldsymbol{\xi})=1
\end{equation}
and
\begin{equation}
	\label{eqn:deg-b}
	\deg({\boldsymbol{u}}_\varepsilon,B_\nu(1-\varepsilon)\setminus \closure{U},\boldsymbol{\xi})=0.
\end{equation}
Plugging in \eqref{eqn:deg-a}--\eqref{eqn:deg-b} into \eqref{eqn:additivity},  we obtain
$\deg(\boldsymbol{u}_\varepsilon,U,\boldsymbol{\xi})=1$, so that  \eqref{eqn:deg-deg} yields $\deg(\boldsymbol{u},U,\boldsymbol{\xi})=1$, which concludes the proof of \eqref{eqn:deg1}.
To check \eqref{eqn:deg-a}, define $\widetilde{\boldsymbol{u}}_\varepsilon \in C^1(\closure{B}_\nu(1-\varepsilon);\R^N)$ as $\widetilde{\boldsymbol{u}}_\varepsilon(\boldsymbol{x})\coloneqq \frac{\eta(1-\varepsilon)}{1-\varepsilon}\boldsymbol{x}$. 
Observe that $\widetilde{\boldsymbol{u}}_\varepsilon\restr{S_\nu(1-\varepsilon)}={\boldsymbol{u}}_\varepsilon\restr{S_\nu(1-\varepsilon)}$ and $\widetilde{\boldsymbol{u}}_\varepsilon(B_\nu(1-\varepsilon))={\boldsymbol{u}}_\varepsilon(B_\nu(1-\varepsilon))=B_\nu(\eta(1-\varepsilon))$, where we use that $\eta_\varepsilon$ is increasing and $\eta_\varepsilon((0,1-\varepsilon)) = (0, \eta(1-\varepsilon))$. 
Given  that $\widetilde{\boldsymbol{u}}_\varepsilon$ is linear, the classical formula for the degree \cite[Definition~2.1]{fonseca.gangbo} yields $\deg(\widetilde{\boldsymbol{u}}_\varepsilon,B_\nu(1-\varepsilon),\cdot)=\chi_{B_\nu(\eta(1-\varepsilon))}$ in $\R^N \setminus S_\nu(\eta(1-\varepsilon))$. Thus, \eqref{eqn:deg-a} follows by \cite[Theorem~2.4]{fonseca.gangbo}. 
Claim \eqref{eqn:deg-b} holds as a consequence of \cite[Theorem~2.1]{fonseca.gangbo}
considering that $\boldsymbol{\xi}\in \boldsymbol{u}_\varepsilon(U) \subset \R^N \setminus \boldsymbol{u}_\varepsilon(B_\nu(1-\varepsilon)\setminus \overline{U})$ because of the injectivity of $\boldsymbol{u}_\varepsilon$.

We move to the proof of \eqref{eqn:radial-inclusion}.  First,   we note that the identity in \eqref{eqn:radial-inclusion} is immediate from \eqref{eqn:radial-deg}.
For the two inclusions in \eqref{eqn:radial-inclusion},  let $\boldsymbol{x}\in B_\nu \setminus \{ \boldsymbol{0} \}$ and  choose  $\varepsilon \in (0, \delta]$  such that $\varepsilon < |\boldsymbol{x}|_\nu < 1-\varepsilon$. In this case, $\boldsymbol{u}(\boldsymbol{x})=\boldsymbol{u}_\varepsilon(\boldsymbol{x})$.  If $\boldsymbol{x}\in U$, then 
$\boldsymbol{u}(\boldsymbol{x})=\boldsymbol{u}_\varepsilon(\boldsymbol{x})\in \boldsymbol{u}_\varepsilon(U)=V$. Instead, if $\boldsymbol{x}\notin U$, then we have 
\begin{equation*}
	\boldsymbol{u}(\boldsymbol{x})=\boldsymbol{u}_\varepsilon(\boldsymbol{x})\in \boldsymbol{u}_\varepsilon(B_\nu \setminus (U \cup \{ \boldsymbol{0} \})) \subset \R^N \setminus \boldsymbol{u}_\varepsilon(U)=\R^N \setminus V,
\end{equation*} 
owing to the injectivity of $\boldsymbol{u}_\varepsilon$. Thus, \eqref{eqn:radial-inclusion} is proved. 
\end{proof}

All the examples discussed in Sections \ref{sec:conv},  \ref{sec:EL}, and  \ref{sec:appl} concern maps defined as the composition of a radial deformation with a Lipschitz transformation. 
First, we discuss right compositions corresponding to changes of  the  reference configuration. See \cite[Theorem 9.1]{mueller.spector} for a related result.

\begin{proposition}[Changes of reference configuration and radial maps]\label{prop:radial-coref}
Let $\eta\colon (0,1)\to (0,+\infty)$ be a locally absolutely continuous and strictly increasing function  satisfying \eqref{eqn:radial-integrability} for some $p>N-1$ and \eqref{eqn:radial-jacobian-coarea}. Also, let $1\leq \nu\leq \infty$ and    $\boldsymbol{u}\colon B_\nu \to \R^N$ be defined as in \eqref{eqn:radial-def}. 
Eventually, let $\Omega \subset \R^N$ be a bounded domain  and let   $\boldsymbol{f}\colon {\Omega} \to \boldsymbol{f}({\Omega})$ be an injective Lipschitz map with $\boldsymbol{f}({\Omega}) \subset B_\nu$ and $\det D \boldsymbol{f}>0$ almost everywhere whose inverse is also  Lipschitz. Then, setting $\boldsymbol{y}\coloneqq \boldsymbol{u}\circ \boldsymbol{f}$, we have $\boldsymbol{y}\in W^{1,p}(\Omega;\R^N)$ and $\det D \boldsymbol{y}\in L^1_+(\Omega)$. Moreover, $\boldsymbol{y}$ satisfies condition {\rm (INV)}. 	
\end{proposition}
\begin{proof}
 From Proposition~\ref{prop:radial-regularity} and Corollary~\ref{cor:radial-jacobian}, we have $\boldsymbol{u}\in W^{1,p}(\Omega;\R^N)$ with $\det D \boldsymbol{u}\in L^1_+(\Omega)$.   
By \cite[Theorem 11.53]{leoni}, we deduce that  $\boldsymbol{y}\in W^{1,p}(\Omega;\R^N)$ with $D\boldsymbol{y}=(D\boldsymbol{u}\circ \boldsymbol{f})  (D\boldsymbol{f})$ almost everywhere. In particular,  $\det D \boldsymbol{y}=  ((\det D \boldsymbol{u})\circ \boldsymbol{f})  \det D \boldsymbol{f}>0$ almost everywhere  thanks to  Lusin's condition (N${}^{-1}$) satisfied by $\boldsymbol{f}$.   Moreover,  using Corollary \ref{cor:change-of-variable} and  Lemma \ref{lem:radial-injective}, we estimate  
	\begin{equation*}
		\int_\Omega \det D \boldsymbol{y}\,\d\boldsymbol{x}=\int_{\boldsymbol{f}(\Omega)} \det D \boldsymbol{u}\,\d\boldsymbol{\xi}\leq \int_{B_\nu} \det D \boldsymbol{u}\,\d\boldsymbol{\xi},
	\end{equation*}
	where  the integral on the right-hand side is finite.     Hence, 
	$\det D \boldsymbol{y}\in L^1_+(\Omega)$. 

  We show that, for every domain $U\subset \subset \Omega$ with $\boldsymbol{f}(\partial U) \subset \R^N \setminus \{ \boldsymbol{0} \}$ such that $\R^N \setminus \partial U$ has exactly two connected components, we have
\begin{equation}
	\label{eqn:radial-deg-cor}
	\deg(\boldsymbol{y},U,\cdot)=\chi_V \quad \text{in $\R^N \setminus \boldsymbol{y}(\partial U)$,}
\end{equation}
where $V$ denotes the unique bounded connected component of $\R^N \setminus \boldsymbol{y}(\partial U)$, and it holds
\begin{equation}
	\label{eqn:radial-inclusion-cor}
	\imt(\boldsymbol{y},U)=V,\qquad \boldsymbol{y}(U \setminus \{ \boldsymbol{f}^{-1}(\boldsymbol{0})  \}) \subset V, \qquad \boldsymbol{y}(\Omega \setminus (U \cup \{ \boldsymbol{f}^{-1}(\boldsymbol{0})  \}) )\subset \R^N \setminus V.
\end{equation}
  Once this is proved, condition (INV) for $\boldsymbol{y}$  follows  by choosing $U$  in \eqref{eqn:radial-inclusion-cor} as a ball.

Let us prove \eqref{eqn:radial-deg-cor}--\eqref{eqn:radial-inclusion-cor}. Since $\boldsymbol{f}\restr{\partial U}$ is injective, by Jordan's separation theorem \cite[Theorem~3.29]{fonseca.gangbo}, the set $\R^N \setminus \boldsymbol{f}(\partial U)$ has exactly two connected components: $W$ and $W_0$, the former being bounded and the latter being unbounded. As  $\boldsymbol{f}$ is  Lipschitz  with $\det D \boldsymbol{f}>0$ almost everywhere, from Lemma~\ref{lem:degree-supercritical} we see that $W=\boldsymbol{f}(U)$ and we compute 
 \begin{equation}
 	\label{eqn:deg-f}
 	\deg(\boldsymbol{f},U,\cdot)=\chi_W \quad \text{in $\R^N \setminus \boldsymbol{f}(\partial U)$.} 
 \end{equation}
 In particular, $\partial W=\boldsymbol{f}(\partial U)$. Since $\boldsymbol{u}\restr{\partial W}$ is injective, again by Jordan's separation theorem \cite[Theorem~3.29]{fonseca.gangbo}, we can write $\R^N \setminus \boldsymbol{u}(\partial W)=V \cup V_0$, where $V$ and $V_0$ are the bounded and unbounded connected component, respectively.  
 
 Let  $0<\varepsilon \ll 1$ be such that $\partial W \subset A_\nu(\varepsilon,1-\varepsilon)$, and define $\eta_\varepsilon \in C^0([0,1-\varepsilon])$ and $\boldsymbol{u}_\varepsilon\in C^0(\closure{B}_\nu(1-\varepsilon);\R^N)$ as in \eqref{eqn:ueps}. Thus, $\boldsymbol{u}_\varepsilon\restr{\partial W}=\boldsymbol{u}\restr{\partial W}$ and, in turn, $\boldsymbol{y}_\varepsilon\restr{\partial U}=\boldsymbol{y}\restr{\partial U}$, where we set $\boldsymbol{y}_\varepsilon\coloneqq \boldsymbol{u}_\varepsilon \circ \boldsymbol{f}\in C^0(\closure{U};\R^N)$.  In particular, Definition~\ref{def:top-deg} yields
 \begin{equation}
 	\label{eqn:deg-deg-y}
 	\deg(\boldsymbol{y},U,\cdot)=\deg(\boldsymbol{y}_\varepsilon,U,\cdot) \quad \text{in $\R^N \setminus \boldsymbol{y}(\partial U)$.}
 \end{equation}
 Now, from Lemma~\ref{lem:radial-degree} and Definition~\ref{def:top-deg}, we have 
 \begin{equation}
 	\label{eqn:deg-u}
 	\deg(\boldsymbol{u}_\varepsilon,W,\cdot)=\deg(\boldsymbol{u},W,\cdot)=\chi_V \quad \text{in $\R^N \setminus \boldsymbol{y}(\partial U)$}
 \end{equation}
 and
 \begin{equation}
 	\label{eqn:radial-inclusion-u}
 	\boldsymbol{u}(W \setminus \{ \boldsymbol{0} \}) \subset V, \qquad \boldsymbol{u}(B_\nu \setminus (W \cup \{ \boldsymbol{0}  \})) \subset \R^N \setminus V.
 \end{equation}
 Using the multiplication formula for the degree \cite[Theorem~2.10]{fonseca.gangbo} together with \eqref{eqn:deg-f} and \eqref{eqn:deg-u}, we get
 \begin{equation*}
 	\begin{split}
 		\deg(\boldsymbol{y}_\varepsilon,U,\cdot)&=\deg(\boldsymbol{u}_\varepsilon,W,\cdot)\,\deg(\boldsymbol{f},U,W)+\deg(\boldsymbol{u}_\varepsilon,W_0,\cdot)\,\deg(\boldsymbol{f},U,W_0)
 		=\chi_V \quad \text{in $\R^N \setminus \boldsymbol{y}(\partial U)$,}
 	\end{split}
 \end{equation*}
 where $\deg(\boldsymbol{f},U,W)$ and $\deg(\boldsymbol{f},U,W_0)$ denote the value of $\deg(\boldsymbol{f},U,\boldsymbol{\xi})$ for any $\boldsymbol{\xi}\in W$ and $\boldsymbol{\xi}\in W_0$, respectively, see \eqref{eqn:deg-f}.    Combining the previous equation with \eqref{eqn:deg-deg-y}, we obtain \eqref{eqn:radial-deg-cor}. From this, the equality in \eqref{eqn:radial-inclusion-cor} immediately follows. The two inclusions in \eqref{eqn:radial-inclusion-cor} are deduced from \eqref{eqn:radial-inclusion-u} recalling that $W=\boldsymbol{f}(U)$.

\end{proof}

Next, we look at left compositions corresponding to the superposition of a Lipschitz transformation to a radial deformation. 

\begin{proposition}[Superposition and radial maps]\label{prop:radial-superposition}
Let $\eta\colon (0,1)\to (0,+\infty)$ be a locally absolutely continuous and strictly increasing function satisfying \eqref{eqn:radial-integrability} for  $p>N-1$ and \eqref{eqn:radial-jacobian-coarea}. Also, let $1\leq \nu \leq \infty$  and  $\boldsymbol{u}\colon B_\nu \to \R^N$ be defined as in \eqref{eqn:radial-def}.  Eventually, let $\boldsymbol{g}\colon   \R^N  \to \R^N$ be   a Lipschitz map with $\det D \boldsymbol{g}\geq 0$ almost everywhere and, for $\Lambda_{\boldsymbol{g}}\coloneqq \{ \det D \boldsymbol{g}=0 \}$, suppose that $\boldsymbol{g}$ is injective in $\R^N \setminus \closure{\Lambda}_{\boldsymbol{g}}$, $\leb(\partial \Lambda_{\boldsymbol{g}})=0$,   and $\boldsymbol{u}(B_\nu \setminus \{ \boldsymbol{0}  \}) \cap \closure{\Lambda}_{\boldsymbol{g}}=\emptyset$.   
Then, setting $\boldsymbol{y}\coloneqq \boldsymbol{g}\circ \boldsymbol{u}$, we have $\boldsymbol{y} \in W^{1,p}(B_\nu;\R^N)$ and $\det D \boldsymbol{y}\in L^1_+(B_\nu)$.  Moreover,  $\boldsymbol{y}$ satisfies condition {\rm (INV)}.  
\end{proposition}

For the proof of Proposition \ref{prop:radial-superposition}, we need a preliminary result on the chain rule for the composition of Sobolev and Lipschitz maps.  For scalar-valued compositions, the validity of the chain rule is a classical result, see, e.g., \cite[Theorem~12.69]{leoni}. Instead, the vector-valued case  is much more  delicate  \cite{leoni.morini}.   For our purposes, the following sufficient conditions are satisfactory. The next result is probably known to experts, but we include it with proof as we did not find any reference for it into the literature.

\begin{lemma}[Superposition]\label{lem:superposition}
	Let $A\subset \R^N$ be a bounded domain	and $\boldsymbol{u}\in W^{1,q}(A;\R^N)$ with $1\leq q<\infty$ satisfying $\det D \boldsymbol{u}>0$ almost everywhere in $A$. Then, for every Lipschitz map $\boldsymbol{g}\colon \R^N \to \R^M$, we have $\boldsymbol{g} \circ \boldsymbol{u}\in W^{1,q}(A;\R^M)$.   Moreover, the chain rule holds, namely
	\begin{equation*}
		\text{$D(\boldsymbol{g}\circ \boldsymbol{u})=  (D\boldsymbol{g}\circ \boldsymbol{u})  \,(D\boldsymbol{u})$ \quad almost everywhere in $A$.}
	\end{equation*} 
\end{lemma}
\begin{remark}
	\begin{enumerate}[(a)]
		\item It is sufficient to assume that $\boldsymbol{u}$ satisfies Lusin's condition (N$^{-1}$)  instead of $\det D \boldsymbol{u}>0$ almost everywhere.   	
		\item If $A$ is not bounded, then  $\boldsymbol{g}\circ \boldsymbol{u}\in W^{1,q}_{\rm loc}(A;\R^M)$ with $D(\boldsymbol{g}\circ \boldsymbol{u})\in L^q(A;\R^{M\times N})$. Additionally, if $\boldsymbol{g}(\boldsymbol{0})=\boldsymbol{0}$, we also have $\boldsymbol{g}\circ \boldsymbol{u}\in L^{q}(A;\R^M)$ and, in turn, $\boldsymbol{g}\circ \boldsymbol{u}\in W^{1,q}(A;\R^M)$.  
	\end{enumerate}	
\end{remark}
\begin{proof}
	First, note that $\boldsymbol{u}$ satisfies Lusin's condition (N$^{-1}$) thanks to Remark \ref{rem:federer}(b). By \cite[Theorem 11.45]{leoni}, up to the choice of a representative, the map $\boldsymbol{u}$ is absolutely continuous on almost every line in $\Omega$ parallel to the coordinate axes. Precisely, for every $i=1,\dots, N$ and for almost every $\boldsymbol{s}\in \R^{N-1}$, the map $\boldsymbol{u}^{(i)}_{\boldsymbol{s}}\colon A^{(i)}_{\boldsymbol{s}} \to \R^N$ defined as $\boldsymbol{u}^{(i)}_{\boldsymbol{s}}(t)\coloneqq \boldsymbol{u}(\boldsymbol{s},t)$ is absolutely continuous. Here,  $A^{(i)}_{\boldsymbol{s}}\coloneqq \left \{t\in \R: \hspace{2pt} (\boldsymbol{s},t)\in A\right \}$, where, with a slight abuse of notation, we write $(\boldsymbol{s},t)=(s_1,\dots,s_{i-1},t,s_{i+1}, \dots,s_{N-1})^\top$ with obvious modifications for the cases $i=1,N$. In particular, $A^{(i)}_{\boldsymbol{s}}$ is an open set and  $\boldsymbol{u}^{(i)}_{\boldsymbol{s}}$ is differentiable on $A^{(i)}_{\boldsymbol{s}}\setminus Z^{(i)}_{\boldsymbol{s}}$ for some set $Z^{(i)}_{\boldsymbol{s}} \subset A^{(i)}_{\boldsymbol{s}}$ with $\mathscr{L}^1(Z^{(i)}_{\boldsymbol{s}})=0$. As a consequence, $\boldsymbol{u}$ admits  partial  derivatives  in the direction $\boldsymbol{e}_i$  with  $\partial_i \boldsymbol{u}(\boldsymbol{s},t)=(\boldsymbol{u}^{(i)}_{\boldsymbol{s}})'(t)$ at  $(\boldsymbol{s},t)$ for all $t\in A^{(i)}_{\boldsymbol{s}}\setminus Z^{(i)}_{\boldsymbol{s}}$.
	It is enough to consider $M=1$. Let $g\colon \R^N \to \R$ be Lipschitz. For all $i$ and $\boldsymbol{s}$ as above, the map $g \circ \boldsymbol{u}^{(i)}_{\boldsymbol{s}}$ is absolutely continuous on $A^{(i)}_{\boldsymbol{s}}$. Thus, $g \circ \boldsymbol{u}^{(i)}_{\boldsymbol{s}}$ is almost everywhere differentiable on $A^{(i)}_{\boldsymbol{s}}$, but we do not necessarily have the chain rule. 
	Define $\Sigma_ {g}\subset \R^N$ as the set of points where $g$ is not differentiable. By Rademacher's theorem,  $\leb(\Sigma_g)=0$, so that $\leb(\boldsymbol{u}^{-1}(\Sigma_g))=0$ by Lusin's condition  (N$^{-1}$). Using the coarea formula and standard properties of the Hausdorff measure  \cite[Proposition~2.49(iv)]{ambrosio.fusco.pallara}, we see that $\mathscr{L}^1(P^{(i)}_{\boldsymbol{s}})=0$, where $P^{(i)}_{\boldsymbol{s}}\coloneqq \{ t\in A^{(i)}_{\boldsymbol{s}}: \hspace{2pt} (t,\boldsymbol{s})\in \boldsymbol{u}^{-1}(\Sigma_g)  \}$, for all $i=1,\dots,N$ and for almost every $\boldsymbol{s}\in \R^{N-1}$. Therefore, for all $i=1,\dots,N$ and for almost all $\boldsymbol{s}\in \R^{N-1}$, we have that $\boldsymbol{u}^{(i)}_{\boldsymbol{s}}$ if differentiable at $t$ and  $g$  is differentiable at $\boldsymbol{u}^{(i)}_{\boldsymbol{s}}(t)$ for all $t\in A^{(i)}_{\boldsymbol{s}}\setminus (Z^{(i)}_{\boldsymbol{s}} \cup P^{(i)}_{\boldsymbol{s}})$. In that case,  $g \circ \boldsymbol{u}^{(i)}_{\boldsymbol{s}}$ is differentiable at $t$ with $(g\circ \boldsymbol{u}^{(i)}_{\boldsymbol{s}})'(t)=Dg(\boldsymbol{u}^{(i)}_{\boldsymbol{s}}(t))\cdot (\boldsymbol{u}^{(i)}_{\boldsymbol{s}})'(t)$. At this point, we deduce that $g \circ \boldsymbol{u}$ admits all partial derivatives at almost every point of $A$ with $D(g \circ \boldsymbol{u})= (Dg \circ \boldsymbol{u})  \cdot (D\boldsymbol{u})$. From   
	this formula, we see that $D(g \circ \boldsymbol{u})\in L^q(A;\R^N)$ thanks to the uniform boundedness of $Dg$. Also, exploiting the Lipschitz continuity of $g$ and the boundedness of $A$, we easily check that $g \circ \boldsymbol{u}\in L^q(A)$. Therefore, by \cite[Theorem 11.45]{leoni}, we conclude that $D(g \circ \boldsymbol{u})$  coincides  with the weak gradient of $g \circ \boldsymbol{u}$, so that $g \circ \boldsymbol{u}\in W^{1,q}(A)$. 
\end{proof}

\begin{proof}[Proof of Proposition \ref{prop:radial-superposition}]
	By Lemma \ref{lem:superposition}, we have 
	$\boldsymbol{y}\in W^{1,p}(B_\nu;\R^N)$ with $D\boldsymbol{y}=(D\boldsymbol{g}\circ \boldsymbol{u})(D\boldsymbol{u})$ almost everywhere. Thus, $\det D \boldsymbol{y}=(\det D \boldsymbol{g} \circ \boldsymbol{u})\det D \boldsymbol{u}>0$ almost everywhere thanks to  $\boldsymbol{u}(B_\nu \setminus \{ \boldsymbol{0} \}) \cap \closure{\Lambda}_{\boldsymbol{g}}=\emptyset$ and Remark \ref{rem:radial-jacobian}(a). Using Corollary \ref{cor:change-of-variable}(i) and Lemma \ref{lem:radial-injective}, we estimate   
	\begin{equation*}
	\begin{split}
	\int_{B_\nu}\det D \boldsymbol{y}\,\d\boldsymbol{x}&=\int_{\boldsymbol{u}(B_\nu)} \det D \boldsymbol{g}\, \d \boldsymbol{\xi} \leq C\,\|D\boldsymbol{g}\|_{L^\infty(B_\nu;\rnn)}^N \, \leb(\boldsymbol{u}(B_\nu\setminus \{\boldsymbol{0}\}))\leq C(\boldsymbol{g})\int_{B_\nu} \det D \boldsymbol{u}\,\d\boldsymbol{x},
	\end{split}
	\end{equation*}
	where the right-hand side is finite by Corollary \ref{cor:radial-jacobian}. This shows that $\det D \boldsymbol{y}\in L^1_+(B_\nu)$.
	
	 We show that for every domain $U \subset \subset \Omega$ with $\partial U \subset \R^N \setminus \{  \boldsymbol{0} \}$ such that $\R^N \setminus \partial U$ has exactly two connected  components, we have 
	\begin{equation}
		\label{eqn:radial-deg-superpos}
		\deg(\boldsymbol{y},U,\cdot )=\chi_V \quad \text{in $\R^N \setminus \boldsymbol{y}(\partial U)$,}
	\end{equation}
	where $V$ denotes the unique bounded connected component of $\R^N \setminus \boldsymbol{y}(\partial U)$, and it holds 
	\begin{equation}
		\label{eqn:radial-inclusion-superpos}
		\imt(\boldsymbol{y},U)=V, \qquad \boldsymbol{y}(U \setminus \{ \boldsymbol{0} \})\subset V, \qquad \boldsymbol{y}(B_\nu \setminus (U \cup \{ \boldsymbol{0} \}))\subset \R^N \setminus V.
	\end{equation}
Once this is proved,  $\boldsymbol{y}$ satisfies condition (INV)   by choosing $U$ in \eqref{eqn:radial-inclusion-superpos} to be a ball.

We now show \eqref{eqn:radial-deg-superpos}--\eqref{eqn:radial-inclusion-superpos}. 	As $\boldsymbol{u}\restr{\partial U}$ is injective by Lemma~\ref{lem:radial-injective},  thanks to  Jordan's separation theorem \cite[Theorem~3.29]{fonseca.gangbo}, we can write  $\R^N \setminus \boldsymbol{u}(\partial U)=W \cup W_0$, where $W$ and $W_0$ are the bounded and unbounded connected component, respectively. Also, 	Lemma~\ref{lem:radial-degree}   yields 
	\begin{equation}
		\label{eqn:deg-u-W}
		\deg(\boldsymbol{u},U,\cdot)=\chi_W \quad \text{in $\R^N \setminus \boldsymbol{u}(\partial U)$}
	\end{equation}
	and
	\begin{equation}
		\label{eqn:uu}
		\boldsymbol{u}(B_\nu \setminus \{ \boldsymbol{0}\}) \subset W, \qquad \boldsymbol{u}(B_\nu \setminus (U \cup \{ \boldsymbol{0} \}))\subset \R^N \setminus  W. 
	\end{equation}
 Given that $\boldsymbol{u}(\partial U) \cap \closure{\Lambda}_{\boldsymbol{g}}=\emptyset$, by assumption,    $\boldsymbol{g}\restr{\boldsymbol{u}(\partial U)}$ is also injective and again by Jordan's separation theorem \cite[Theorem~3.29]{fonseca.gangbo}, we  have  $\R^N \setminus \boldsymbol{y}(\partial U)=V \cup V_0$, where $V$ and $V_0$ are connected components with $V$ bounded and $V_0$ unbounded. Note that $\boldsymbol{g}$ is almost everywhere injective in $\R^N \setminus \Lambda_{\boldsymbol{g}}$ as $\leb(\partial \Lambda_{\boldsymbol{g}})=0$. 
 Thus,  we can use Lemma~\ref{lem:degree-supercritical} to  see that $\boldsymbol{g}(W)=V$  and there holds
	\begin{equation}
		\label{eqn:deg-g}
		\deg(\boldsymbol{g},W,\cdot)=\chi_V \quad \text{in $\R^N \setminus \boldsymbol{y}(\partial U)$.}
	\end{equation} 
	Now, let $0<\varepsilon  \ll 1$ be such that $\partial U \subset A_\nu(\varepsilon,1-\varepsilon)$, and define $\eta_\varepsilon \in C^0([0,1-\varepsilon])$ and $\boldsymbol{u}_\varepsilon \in C^0(\closure{B}_\nu(1-\varepsilon);\R^N)$ as in \eqref{eqn:ueps}.   As $\boldsymbol{u}_\varepsilon\restr{\partial U}=\boldsymbol{u}\restr{\partial U}$, by Definition~\ref{def:top-deg}, we have \eqref{eqn:deg-deg}. Then,  \eqref{eqn:deg-u-W} gives
	\begin{equation}
		\label{eqn:deg-ueps}
		\deg(\boldsymbol{u}_\varepsilon,U,\cdot)=\chi_W \quad \text{in $\R^N \setminus \boldsymbol{u}(\partial U)$.}
	\end{equation} 
	Setting $\boldsymbol{y}_\varepsilon\coloneqq \boldsymbol{g}\circ \boldsymbol{u}_\varepsilon \in C^0(\closure{B}_\nu;\R^N)$, we have $\boldsymbol{y}_\varepsilon\restr{\partial U}=\boldsymbol{y}\restr{\partial U}$, so that   	$\deg(\boldsymbol{y},U,\cdot)=\deg(\boldsymbol{y}_\varepsilon,U,\cdot)$ in $\R^N \setminus \boldsymbol{y}(\partial U)$ by   Definition~\ref{def:top-deg}.  Applying the multiplication formula for the degree \cite[Theorem~2.10]{fonseca.gangbo} taking into account \eqref{eqn:deg-g}--\eqref{eqn:deg-ueps}, we obtain
	\begin{equation*}
		\deg(\boldsymbol{y}_\varepsilon,U,\cdot)=\deg(\boldsymbol{g},W,\cdot)\,\deg(\boldsymbol{u}_\varepsilon,U,W)+\deg(\boldsymbol{g},W_0,\cdot)\,\deg(\boldsymbol{u}_\varepsilon,U,W_0)=\chi_V \quad \text{in $\R^N \setminus  \boldsymbol{y}(\partial U)$,  }
	\end{equation*} 
	where $\deg(\boldsymbol{u}_\varepsilon,U,W)$ and $\deg(\boldsymbol{u}_\varepsilon,U,W_0)$ denote the value of $\deg(\boldsymbol{u}_\varepsilon,U,\boldsymbol{\xi})$ for $\boldsymbol{\xi}\in W$ and $\boldsymbol{\xi}\in W_0$, respectively, see \eqref{eqn:deg-u-W} and \eqref{eqn:deg-ueps}. 	Thus,  $\deg(\boldsymbol{y},U,\cdot)=\deg(\boldsymbol{y}_\varepsilon,U,\cdot)$ in $\R^N \setminus \boldsymbol{y}(\partial U)$  and the previous equation give \eqref{eqn:radial-deg-superpos}. 
	At this point, the equality in \eqref{eqn:radial-inclusion-superpos} is immediate, while the two inclusion follow from \eqref{eqn:uu} by applying $\boldsymbol{g}$ at both sides.
\end{proof}

\section*{Acknowledgements}
M.~Bresciani and M.~Friedrich have been  supported by the DFG project FR 4083/5-1 and by the Deutsche Forschungsgemeinschaft (DFG, German Research Foundation) under Germany's Excellence Strategy EXC 2044 -390685587, Mathematics M\"unster: Dynamics--Geometry--Structure.  

C.~Mora-Corral has been supported by the Agencia Estatal de Investigaci\'{o}n of the Spanish Ministry of
Research and Innovation, through project PID2021-124195NB-C32 and the Severo Ochoa Programme for
Centres of Excellence in R\&D CEX2019-000904-S, by the Madrid Government (Comunidad de Madrid,
Spain) under the multiannual Agreement with UAM in the line for the Excellence of the University Research
Staff in the context of the V PRICIT   (Regional Programme of Research and Technological Innovation), and
by the ERC Advanced Grant 834728.  

The authors  are thankful to the Erwin Schr\"odinger Institute in Vienna, where part of this work has been developed during the workshop ``Between Regularity and Defects: Variational and Geometrical Methods in Materials Science".

\end{document}